\input xy
\xyoption{all}
\xyoption{2cell}
\UseAllTwocells

\documentclass{memo-l}

\usepackage{amsmath}
\usepackage{amscd}
\usepackage{rotating}

\def\pr{{\operatorname{pr}}}
\def\id{{\operatorname{id}}}
\def\Fr{{\operatorname{Fr}}}

\newtheorem{theorem}{Theorem}[chapter]
\newtheorem{proposition}[theorem]{Proposition}
\newtheorem{corollary}[theorem]{Corollary}
\newtheorem{lemma}[theorem]{Lemma}

\theoremstyle{definition}
\newtheorem{definition}[theorem]{Definition}
\newtheorem{example}[theorem]{Example}

\newtheorem*{main}{Theorem}
\newtheorem*{newmain}{Definition}

\theoremstyle{remark}

\numberwithin{section}{chapter} \numberwithin{equation}{chapter}

\makeindex

\begin{document}
\frontmatter

\title{Points on Quantum Projectivizations}

\author{Adam Nyman}
\address{Department of Mathematics, University of Montanta, Missoula, MT 59812-0864}
\email{nymana@mso.umt.edu}

\date{September 5, 2002}
\thanks{2000 {\it Mathematical Subject Classification. } Primary 14A22; Secondary 18D30}
\keywords{Non-commutative algebraic geometry, Point module,
Non-commutative projective bundle}

\begin{abstract}
The use of geometric invariants has recently played an important
role in the solution of classification problems in non-commutative
ring theory.  We construct geometric invariants of non-commutative
projectivizataions, a significant class of examples in
non-commutative algebraic geometry.  More precisely, if $S$ is an
affine, noetherian scheme, $X$ is a separated, noetherian
$S$-scheme, $\mathcal{E}$ is a coherent
${\mathcal{O}}_{X}$-bimodule  and $\mathcal{I} \subset
T(\mathcal{E})$ is a graded ideal then we develop a compatibility
theory on adjoint squares in order to construct the functor
$\Gamma_{n}$ of flat families of truncated
$T(\mathcal{E})/\mathcal{I}$-point modules of length $n+1$.  For
$n \geq 1$ we represent $\Gamma_{n}$ as a closed subscheme of
${\mathbb{P}}_{X^{2}}({\mathcal{E}}^{\otimes n})$.  The
representing scheme is defined in terms of both
${\mathcal{I}}_{n}$ and the bimodule Segre embedding, which we
construct.

Truncating a truncated family of point modules of length $i+1$ by taking its first $i$ components defines a morphism $\Gamma_{i} \rightarrow \Gamma_{i-1}$ which makes the set $\{\Gamma_{n}\}$ an inverse system.  In order for the point modules of $T(\mathcal{E})/\mathcal{I}$ to be parameterizable by a scheme, this system must be eventually constant.  In \cite{pointsonqrs}, we give sufficient conditions for this system to be constant and show that these conditions are satisfied when ${\sf{Proj}} T(\mathcal{E})/\mathcal{I}$ is a quantum ruled surface.  In this case, we show the point modules over $T(\mathcal{E})/\mathcal{I}$ are parameterized by the closed points of ${\mathbb{P}}_{X^{2}}(\mathcal{E})$.
\end{abstract}

\maketitle

\tableofcontents

\mainmatter
\chapter{Introduction} Many solutions to classification problems in
non-commutative ring theory have recently been attained with the
aid of geometric invariants.  By associating to a non-commutative
ring a scheme whose geometric properties reflect the properties of
the ring, we may use the geometry of schemes to study the ring.
The discovery of geometric invariants associated to
non-commutative rings has proved enormously helpful in classifying
non-commutative rings of low dimension.  In \cite{stafford,
groth}, this technique is used to classify domains of dimension
two and to classify three dimensional non-commutative rings which
are close to commutative polynomial rings, including
three-dimensional Sklyanin algebras, the enveloping algebra of the
Heisenberg Lie algebra, and a significant class of
three-dimensional iterated Ore extensions.

\section{Geometric invariants in the absolute case} We recount
the construction and use of geometric invariants of graded rings
over a field $k$ given in \cite{groth}.  Suppose $V$ is a finite
dimensional $k$-vector space, $T(V)$ is the tensor algebra of $V$,
$I \subset T(V)$ is a homogeneous ideal and $A = T(V)/I$.

\begin{definition} \cite[Definition 3.8, p. 45]{groth}
A graded right $A$-module $M$ is a {\bf point module} if:
\begin{enumerate}
\item{it is generated in degree zero,}
\item{$M_{0}=k$, and}
\item{$\operatorname{dim }M_{i}=1$ for all $i \geq 0$.}
\end{enumerate}
A {\bf truncated point module of length ${\mathbf{n+1}}$} is a module $M$ satisfying (1) and (2) above and whose Hilbert function is $1$ if $0 \leq i \leq n$ and $0$ otherwise.
\end{definition}
If $A$ is commutative, an $A$-point module is isomorphic to a polynomial ring in one variable, and is thus the coordinate ring of a closed point in $\operatorname{Proj }A$.  Furthermore, if $k$ is algebraically closed, the point modules of $A$ are parameterized by $\operatorname{Proj }A$.
\begin{definition} \cite[Definition 3.8, p. 45]{groth}
Let $R$ be a commutative $k$-algebra.  A {\bf flat family of $\mathbf{A}$-point modules parameterized by $\mathbf{\operatorname{\bf{Spec }}R}$} is a graded $R \otimes_{k}A$ module $M$, generated in degree zero, such that $M_{0} \cong R$ and $M_{i}$ is locally free of rank one for each $i$.  Families of truncated point modules are defined similarly.
\end{definition}
Let $\sf{k}$ denote the category of affine noetherian $k$-schemes.
\begin{definition}
The assignment $\Gamma_{n}:\sf{k} \rightarrow \sf{Sets}$ sending $U$ to
$$
\{ \mbox{isomorphism classes of truncated $U$-families of length n+1} \}
$$
and sending $f:V \rightarrow U$ to the map $\Gamma_{n}(f)$ defined by $\Gamma_{n}(f)[M] = [f^{*}M]$, is the {\bf functor of flat families of truncated $\mathbf{A}$-point modules of length $\mathbf{n+1}$}.
\end{definition}

\begin{proposition} \label{prop.start} \cite[Proposition 3.9, p. 46]{groth}
The functor $\Gamma_{n}$ of flat families of truncated $A$-modules of length $n+1$ is representable.  We abuse notation by letting $\Gamma_{n}$ denote the scheme representing the functor $\Gamma_{n}$.
\end{proposition}
The scheme $\Gamma_{n}$ is a subscheme of ${\mathbb{P}}_{{\operatorname{Spec }}k}(V)^{n}$ which depends on $I_{n}$, the $n$th graded piece of $I$.

Truncating a truncated family of point modules of length $i+1$ by taking its first $i$ components defines a morphism $\Gamma_{i} \rightarrow \Gamma_{i-1}$ which makes the set $\{\Gamma_{n}\}$ an inverse system.  In order for the point modules of $A$ to be parameterizable by a scheme, this system must be eventually constant.

If $A=T(V)$, then it is not hard to show that
$\Gamma_{n}={\mathbb{P}}_{{\operatorname{Spec }}k}(V)^{n}$ for all
$n$.  In this case, the inverse system $\{\Gamma_{n}\}$ is never
constant, and the point modules over $A$ are not parameterizable
by a scheme.  If $A=S(V)$, the symmetric algebra of $V$, then the
truncation morphisms $\Gamma_{i} \rightarrow \Gamma_{i-1}$ are
isomorphisms for all $i>1$.  In particular, $\Gamma_{n} \cong
{\mathbb{P}}_{{\operatorname{Spec }}k}(V)$ for all $n \geq 1$.
Thus, the inverse limit, $\Gamma$, of the inverse system
$\{\Gamma_{n}\}$ is just ${\mathbb{P}}_{{\operatorname{Spec
}}k}(V)=\operatorname{Proj }A$.  This is not surprising in light
of the fact that when $A$ is commutative and $k$ is algebraically
closed, the closed points of $\Gamma$ correspond to the closed
points of $\operatorname{Proj }A$.  Although $\Gamma_{n}$ is
relatively easy to compute in case $A=T(V)$ or $A=S(V)$, the
computation of $\Gamma_{n}$ is usually more complicated.
Nevertheless, Artin, Tate and Van den Bergh describe sufficient
conditions for the inverse system $\{\Gamma_{n}\}$ to be
eventually constant \cite[Propositions 3.5, 3.6, and 3.7, pp.
44-45]{groth}.  They then show that these conditions are satisfied
when $A$ is a regular algebra of global dimension two or three
generated in degree one.  If $A$ is a regular algebra of global
dimension two generated in degree one, the truncation map
$\Gamma_{n} \rightarrow \Gamma_{2}$ is an isomorphism, and
$\Gamma_{2}$ is the graph of an automorphism of
${\mathbb{P}}^{1}$.  Furthermore,  two such algebras are
isomorphic if and only if the conjugacy classes of the
corresponding automorphisms of ${\mathbb{P}}^{1}$ are the same
(\cite[p.47-48]{groth}).  Regular algebras of dimension three
generated in degree one are also classified by their associated
point scheme, but the description of their classification is more
complicated than the dimension two case \cite[Theorem 3, p.
36]{groth}.

While regular algebras of global dimension four generated in degree one are far from classified, important classes of such algebras are characterized by their point schemes and by higher dimensional geometric invariants.  For a brief discussion of work in this direction, see \cite{linemod, linemod2}.

\section{Bimodules and algebras} The aim of this paper is to
obtain a relative version of Proposition \ref{prop.start}.  More
precisely, if $S$ is an affine, noetherian scheme, $X$ is a
separated, noetherian $S$-scheme, $\mathcal{E}$ is a coherent
${\mathcal{O}}_{X}$-bimodule, $T(\mathcal{E})$ is the tensor
algebra of $\mathcal{E}$ (Example \ref{example.dave}) and
$\mathcal{I} \subset T(\mathcal{E})$ is an ideal, then we
associate to $\mathcal{B}=T(\mathcal{E})/\mathcal{I}$ a scheme of
points which specializes to Artin, Tate and Van den Bergh's
construction when $S=X$ is the spectrum of a field.

Important examples in noncommutative algebraic geometry are constructed from $\mathcal{B}$.  For example, suppose $X$ is a smooth curve over $k$ and $\mathcal{E}$ is an ${\mathcal{O}}_{X}$-bimodule which is locally free of rank two.  If $\mathcal{Q} \subset \mathcal{E} \otimes_{{\mathcal{O}}_{X}} \mathcal{E}$ is an invertible bimodule, then the quotient $\mathcal{B} = T(\mathcal{E})/(\mathcal{Q})$ is a bimodule algebra.  A quantum ruled surface\index{quantum ruled surface|textbf} was originally defined by Van den Bergh to be the quotient of ${\sf{Grmod }}\mathcal{B}$ by direct limits of modules which are zero in high degree \cite[p.33]{qrs}, ${\sf{Proj }}\mathcal{B}$\index{Proj B@${\sf Proj }\mathcal{B}$|textbf}.  In order that $\mathcal{B}$ has desired regularity properties, Patrick imposes the condition of {\it admissibility} on $\mathcal{Q}$ \cite[Section 2.3]{qrs} and birationally classifies many instances of such quantum ruled surfaces \cite[Section 3]{2side}.  In \cite{pointsonqrs} we follow a suggestion of Van den Bergh by insisting that $\mathcal{Q}$ be nondegenerate \cite[Definition 2.17]{pointsonqrs}.  Quotients of four-dimensional Sklyanin algebras by central homogeneous elements of degree two provide important examples of coordinate rings of quantum ruled surfaces in this sense \cite[Theorem 7.4.1, p. 29]{translation}.

In order to describe our generalization of Proposition \ref{prop.start}, we first describe our generalization of a flat family of point modules.  Just as when $S=X$ is the spectrum of a field, in order to define such a family, we need to know that if $f:V \rightarrow U$ is a map of affine $S$-schemes, $\mathcal{B}$ is an ${\mathcal{O}}_{U \times X}$-bimodule algebra, and if $\tilde{f}=f \times \id_{X}:V \times X \rightarrow U \times X$, then $(\tilde{f} \times \tilde{f})^{*}\mathcal{B}$ inherits an ${\mathcal{O}}_{V \times X}$-bimodule algebra structure from $\mathcal{B}$.  Because the definitions of the bimodule tensor product (Definition \ref{def.tensorproduct}), the associativity isomorphisms of this product (Proposition \ref{prop.tensor}), and the left and right unit morphisms (\ref{eqn.unitary}) are quite complicated, the fact that a bimodule algebra can be lifted to a different base poses technical difficulties.  However, the compatibilities one needs in order to lift a bimodule algebra can be neatly described in the language of indexed categories.  The fact that a bimodule algebra can be lifted reduces to the fact that associativity of the bimodule tensor product and the left and right unit morphisms are indexed natural transformations.  In order to prove these facts, we define, in Chapter 2, the category of squares.  After reviewing the definitions of indexed category, indexed functor and indexed natural transformation, we give necessary conditions under which families of functors between indexed categories are indexed (Proposition \ref{prop.indexedfunct}) and under which transformations associated to squares of indexed categories are indexed (Proposition \ref{cor.reduce}).  Since both the bimodule tensor product associativity maps and the left and right unit maps are a composition of transformations associated to squares, they can be shown to be indexed (Proposition \ref{prop.assindex} and Proposition \ref{prop.unitindex}).

\section{Geometric invariants in the relative case} Now suppose
$\mathcal{B}$ is an ${\mathcal{O}}_{S \times X}$-bimodule algebra,
$U$ is a noetherian affine $S$ scheme, $f:U \rightarrow S$ is the
structure map, $d:U \rightarrow U \times U$ is the diagonal map
and $\tilde{f}=f \times \id_{X}:U \times X \rightarrow S \times
X$.  Let ${\mathcal{B}}^{U} = (\tilde{f} \times
\tilde{f})^{*}\mathcal{B}$.

\begin{newmain}
A {\bf family of $\mathcal{B}$-point modules parameterized by $\mathit{\mathbf{U}}$} \index{family!of point modules|textbf}is a graded ${\mathcal{B}}^{U}$-module ${\mathcal{M}}_{0} \oplus {\mathcal{M}}_{1} \oplus \cdots $ generated in degree zero such that, for each $i\geq 0$ there exists a map
$$
q_{i}: U \rightarrow X
$$
and an invertible ${\mathcal{O}}_{U}$-module ${\mathcal{L}}_{i}$ with ${\mathcal{L}}_{0} \cong {\mathcal{O}}_{U}$ and
$$
{\mathcal{M}}_{i} \cong (\id_{U} \times q_{i})_{*}d_{*}{\mathcal{L}}_{i}.
$$
A {\bf family of truncated $\mathcal{B}$-point modules of length $\mathbf{n+1}$ parameterized by $\mathbf{U}$} \index{family!of truncated point modules|textbf}is a graded ${\mathcal{B}}^{U}$-module ${\mathcal{M}}_{0} \oplus {\mathcal{M}}_{1} \oplus \cdots $generated in degree zero, such that for each $i\geq 0$ there exists a map
$$
q_{i}: U \rightarrow X
$$
and an invertible ${\mathcal{O}}_{U}$-module ${\mathcal{L}}_{i}$ with ${\mathcal{L}}_{0} \cong {\mathcal{O}}_{U}$,
$$
{\mathcal{M}}_{i} \cong (\id_{U} \times q_{i})_{*}d_{*} {\mathcal{L}}_{i}
$$
for $i \leq n$ and ${\mathcal{M}}_{i} = 0$ for $i > n$.
\end{newmain}

\begin{newmain}
Let $\sf{S}$ be the category of affine, noetherian $S$-schemes.  The assignment $\Gamma_{n}:{\sf{S}} \rightarrow \sf{Sets}$\index{Gamman@$\Gamma_{n}$|textbf} sending $U$ to
$$
\{ \mbox{isomorphism classes of truncated $U$-families of length n+1} \}
$$
and sending $f:V \rightarrow U$ to the map $\Gamma_{n}(f)$ defined by $\Gamma_{n}(f)[\mathcal{M}] = [{\tilde{f}}^{*}\mathcal{M}]$, is the {\bf functor of flat families of truncated $\mathcal{B}$-point modules of length $\mathbf{n+1}$}\index{functor of flat families of truncated $\mathcal{B}$-point modules|textbf}.
\end{newmain}
We show that the functor $\Gamma_{n}$ is representable.  In order to describe the representing scheme, we introduce a generalization of the Segre embedding\index{Segre embedding}, which we now proceed to describe.  First, we need some definitions.  If $W$ is a scheme and $\mathcal{A}$ is a quasi-coherent ${\mathcal{O}}_{W}$-module, then we let $\operatorname{SSupp} \mathcal{A}$ denote the scheme theoretic support of $\mathcal{A}$\index{scheme theoretic support}.  Let $p_{\mathcal{A}}:{\mathbb{P}}_{W}(\mathcal{A}) \rightarrow W$ be the structure map.

\begin{newmain}
Let $X$, $Y$, and $Z$ be schemes, suppose $\mathcal{E}$ is a quasi-coherent ${\mathcal{O}}_{X \times Y}$-module and $\mathcal{F}$ is a quasi-coherent ${\mathcal{O}}_{Y \times Z}$-module.  We say that $\mathcal{E}$ and $\mathcal{F}$ have the {\bf affine direct image property}\index{affine direct image property|textbf} if the restriction of the projection map $\pr_{13}$ to
$$
\operatorname{SSupp}({\pr_{12}}^{*}\mathcal{E} \otimes {\pr_{23}}^{*}\mathcal{F})
$$
is affine.
\end{newmain}
We denote by
$$
{\mathbb{P}}_{X \times Y}(\mathcal{E}) \otimes_{Y} {\mathbb{P}}_{Y \times Z}(\mathcal{F})
$$
the fiber product in the following diagram:
$$
\xymatrix{
& {\mathbb{P}}_{X \times Y}(\mathcal{E}) \otimes_{Y} {\mathbb{P}}_{Y \times Z}(\mathcal{F}) \ar[dr]^{q_{\mathcal{F}}} \ar[dl]_{q_{\mathcal{E}}}  & \\
{\mathbb{P}}_{X \times Y}(\mathcal{E}) \ar[d]_{p_{\mathcal{E}}} & & {\mathbb{P}}_{Y \times Z}(\mathcal{F}) \ar[d]^{p_{\mathcal{F}}} \\
X \times Y \ar[dr]_{\pr_{2}} & & Y \times Z \ar[dl]^{\pr_{1}} \\
& Y &.
}
$$
We prove the following Theorem (Theorem \ref{theorem.bigsegre}):

\begin{main}
Let $\mathcal{E}$ be a quasi-coherent ${{\mathcal{O}}_{X \times Y}}$-module and let $\mathcal{F}$ be a quasi-coherent ${{\mathcal{O}}_{Y \times Z}}$-module such that $\mathcal{E}$ and $\mathcal{F}$ have the affine direct image property.  Then there exists a canonical map
$$
s:{\mathbb{P}}_{X\times Y}(\mathcal{E}) \otimes_{Y} {\mathbb{P}}_{Y \times Z}(\mathcal{F}) \rightarrow {\mathbb{P}}_{X \times Z}({\mathcal{E}} \otimes_{{\mathcal{O}}_{Y}} {\mathcal{F}})
$$
such that $s$ is a closed immersion\index{bimodule Segre embedding!is a closed immersion} which is functorial\index{bimodule Segre embedding!is functorial}, compatible with base change\index{bimodule Segre embedding!is compatible with base change}, and associative\index{bimodule Segre embedding!is associative}.
\end{main}
We call the map $s$ in the previous theorem the {\bf bimodule Segre embedding.}\index{bimodule Segre embedding}  The bimodule Segre embedding is uniquely determined by algebraic data (Theorem \ref{theorem.bigsegre}).

Our generalization of Proposition \ref{prop.start} takes on the following form: if we let $X=Y$ and define the trivial bimodule Segre embedding $s:{\mathbb{P}}_{X^{2}}(\mathcal{E})^{\otimes 1} \rightarrow {\mathbb{P}}_{X^{2}}({\mathcal{E}}^{\otimes 1})$ as the identity map, then we have Theorem \ref{theorem.bigone}:

\begin{main}

For $n \geq 1$, $\Gamma_{n}$ is represented by the pullback of the diagram
$$
\xymatrix{
& {\mathbb{P}}_{X^{2}}(\mathcal{E})^{\otimes n} \ar[d]^{s} \\
{\mathbb{P}}_{X^{2}}({\mathcal{E}}^{\otimes n}/{\mathcal{I}}_{n}) \ar[r] & {\mathbb{P}}_{X^{2}}({\mathcal{E}}^{\otimes n})
}
$$
\end{main}
Now suppose $S= \operatorname{Spec }k$.
\begin{definition}
A {\bf point module}\index{point module|textbf} over $\mathcal{B}$ is an ${\mathbb{N}}$-graded $\mathcal{B}$-module ${\mathcal{M}}_{0} \oplus {\mathcal{M}}_{1} \oplus \cdots $ such that, for each $i\geq 0$, the multiplication map ${\mathcal{M}}_{i} \otimes_{{\mathcal{O}}_{U \times X}} {\mathcal{B}}_{1} \rightarrow {\mathcal{M}}_{i+1}$ is an epimorphism and ${\mathcal{M}}_{i} \cong {\mathcal{O}}_{p_{i}}$ for some closed point $p_{i} \in X$.
\end{definition}

As in the case when $S=X=\operatorname{Spec }k$, truncating a truncated family of point modules of length $i+1$ by taking its first $i$ components defines a morphism $\Gamma_{i} \rightarrow \Gamma_{i-1}$.  This collection of morphisms makes the set $\{\Gamma_{n}\}$ an inverse system.  In order for the point modules of $\mathcal{B}=T(\mathcal{E})/\mathcal{I}$ to be parameterized by a scheme, this system must be eventually constant.  In \cite{pointsonqrs}, we prove, as suggested by Van den Bergh, analogues of \cite[Propositions 3.5, 3.6, and 3.7, pp. 44-45]{groth} in order to give sufficient conditions for the inverse system $\{ \Gamma_{n} \}$ to be eventually constant.  Our main result in this direction is \cite[Proposition 2.22, p. 775]{pointsonqrs}:

\begin{proposition} \label{proposition.newaddition}
\begin{enumerate}
\item{}
For $n \geq m$, Let $p_{m}^{n}:\Gamma_{n} \rightarrow \Gamma_{m}$ denote the map induced by truncation.  Assume that, for some $d$, $p_{d}^{d+1}$ defines a closed immersion from $\Gamma_{d+1}$ to $\Gamma_{d}$, thus identifying $\Gamma_{d+1}$ with a closed subscheme $E \subset \Gamma_{d}$.  Then $\Gamma_{d+1}$ defines a map $\sigma:E \rightarrow \Gamma_{d}$ such that if $(p_{1}, \ldots, p_{d})$ is a point in $E$,
$$
\sigma(p_{1}, \ldots, p_{d}) = (p_{2}, \ldots, p_{d+1}),
$$
where $(p_{1}, \ldots, p_{d+1})$ is the unique point of $\Gamma_{d+1}$ lying over $(p_{1}, \ldots, p_{d}) \in \Gamma_{d}$.
\item{}
If, in addition, $\sigma(E) \subset E$, $\mathcal{I}$ is generated in degree $\leq d$ and $\Gamma$ has linear fibres \index{linear fibres}\cite[Definition 2.17, p. 773]{pointsonqrs}, then $p_{d}^{n}:\Gamma_{n} \rightarrow E$ is an isomorphism for every $n \geq d+1$.
\end{enumerate}
\end{proposition}
In \cite{pointsonqrs}, we use Proposition \ref{proposition.newaddition} to show that, when ${\sf{Proj }}\mathcal{B}$ is a ruled surface in the sense of \cite{pointsonqrs}, the point modules over $\mathcal{B}$ are parameterized by the closed points of ${\mathbb{P}}_{X^{2}}(\mathcal{E})$.

Van den Bergh has developed another definition of quantum ruled surface \cite[Definition 11.4, p.35]{p1bundles}\index{quantum ruled surface}, based on the notion of a non-commutative symmetric algebra \index{non-commutative symmetric algebra}generated by $\mathcal{E}$, which does not depend on $\mathcal{Q}$!  Van den Bergh proves \cite{p1bundles} that the functor of flat families of $\mathcal{B}$-point modules is represented by ${\mathbb{P}}_{X^{2}}(\mathcal{E})$, without using Theorem \ref{theorem.bigone}, then uses this parameterization to show that the category of graded modules over a non-commutative symmetric algebra is noetherian.  While Van den Bergh's proof is specific to non-commutative $\mathbb{P}^{1}$-bundles\index{non-commutative $\mathbb{P}^{1}$-bundles}, our results apply to any non-commutative projectivization.

\section{Organization of the paper}

In Chapter 2 we study the category of squares, and use the resulting theory to provide conditions under which a family of functors between indexed categories is indexed, and conditions under which a family of natural transformations associated to a square of indexed categories is indexed.  In Chapter 3, after reviewing the definitions and basic properties of bimodules and bimodule algebras, we use the results of Chapter 2 to show that a bimodule algebra, and a module over a bimodule algebra, may be lifted to a different base.  This allows us to define $\Gamma_{n}$.  In Chapter 4 we show that $\Gamma_{n}$ is compatible with descent (Definition \ref{def.compwithdescent}).  This fact allows us to deduce that $\Gamma_{n}$ is representable by showing that a suitable subfunctor of $\Gamma_{n}$ imbeds in the scheme of points which ostensibly represents $\Gamma_{n}$.

Now, suppose $(\mathcal{I})_{i}=0$ for $i < n$.  In Chapter 5, we represent $\Gamma_{i}$ for $0 \leq i < n$.  The proof is technical and employs the compatibility result Corollary \ref{cor.adcom2}.  In Chapter 6 we define the bimodule Segre embedding and show that it is functorial, associative and compatible with base change.  These properties allow us to show the map $s:{\mathbb{P}}_{X^{2}}(\mathcal{E})^{\otimes n} \rightarrow {\mathbb{P}}_{X^{2}}({\mathcal{E}}^{\otimes n})$ is well defined.  We then show, in Chapter 7, that $\Gamma_{i}$ is representable for $i \geq n$.

\section{Advice to the reader}

Our attention to compatibilities tends to remove emphasis from the basic idea behind the proof of Theorem \ref{theorem.bigone}.  For the reader who wishes to ignore issues of compatibility, we offer a list of items which can be read independently of the rest of the paper, and which gives the reader the idea behind the proof of Theorem \ref{theorem.bigone}.  For those who want to read the complete proof of Theorem \ref{theorem.bigone}, we recommend first reading this list of items before attempting to read the entire paper.

\begin{itemize}
\item{Chapter 2:  Introduction, Definition \ref{def.bigsquares} through the definition of the ``2-cell dual to $\Phi$", Examples \ref{example.projform} and \ref{example.basechange}.}
\item{Chapter 3:  Introduction, Section 3.1, Section 3.2, Section 3.4 through the statement of Lemma \ref{lem.function}, the statement of Proposition \ref{prop.gamman}, and Definition \ref{def.gamma}.}
\item{Chapter 4:  Through Section 4.3, Lemma \ref{lem.sup1} and Lemma \ref{lem.sup2}.}
\item{Chapter 5:  Through the proof of Proposition \ref{prop.trans}, Step 2, and Proposition \ref{prop.refreport}.}
\item{Chapter 6:  Through Section 6.1 (ignoring Example \ref{example.refrep}), Lemmas \ref{lem.eisenbud} and \ref{lem.diamond}, and Propositions \ref{prop.cruciso} and \ref{prop.segretensor}.}
\item{Chapter 7:  Everything except the proof that $\Sigma^{\Fr}$ is natural in Proposition \ref{prop.injy}.}
\end{itemize}

\section{Notation and conventions}

We will denote indexed categories by $\mathfrak{A}$, $\mathfrak{B}$, $\mathfrak{C}$, \ldots, indexed functors by $\mathfrak{F}$, $\mathfrak{G}$, $\mathfrak{H}$, $\ldots$ and indexed natural transformations (and ordinary natural transformations) by $\Delta$, $\Theta$, $\Omega$, $\ldots$.  We will denote categories and objects in 2-categories by $\sf{A}$, $\sf{B}$, $\sf{C}$, $\ldots$ and functors by $F$, $G$, $H$, $\ldots$.  We will denote objects in an ordinary category by $\mathcal{A}$, $\mathcal{B}$, $\mathcal{C}$, $\ldots$.  We will denote morphisms between objects $\alpha$, $\beta$, $\gamma$, \ldots.  In a 2-category, 2-cells will be denoted by $\Longrightarrow$, while 1-cells will be denoted by $\rightarrow$.  In an ordinary category, morphisms will always be denoted by $\rightarrow$.

There are four exceptions to the above conventions.  If $(F,G)$ is an adjoint pair of functors with $F$ left adjoint to $G$, then we denote its unit 2-cell by $\eta:\operatorname{id} \Longrightarrow GF$ and its counit 2-cell by $\epsilon:FG \Longrightarrow \operatorname{id}$, thus using lowercase Greek letters instead of our usual uppercase.  In addition, $\Gamma_{n}$ will denote a functor, violating the usual convention that upper-case Greek denotes a natural transformation.  If $f:Y \rightarrow X$ is a map of schemes, then $f^{*}$ and $f_{*}$ will denote extension and restriction of scalars, respectively.  Similarly, if $f:V \rightarrow U$ is a morphism in a category $\sf{S}$ and $\mathfrak{A}$ is an indexed category, then we write the induced functor $f^{*}{\sf{A}}^{U} \rightarrow {\sf{A}}^{V}$.  Finally, the map of sheaves of rings associated to $f$, ${\mathcal{O}}_{X} \rightarrow f_{*}{\mathcal{O}}_{Y}$ will sometimes be denoted $f^{\#}$.

If $\sf{C}$ and $\sf{D}$ are categories and $F:\sf{C} \rightarrow \sf{D}$ and $G:\sf{D} \rightarrow \sf{C}$ are functors, then we write $(F,G)$ if $F$ is left adjoint to $G$.  In this case, we will often make use of the fact that the compositions

\begin{equation} \label{eqn.defad}
\xymatrix{
F \ar@{=>}[r]^{F*\eta} & FGF \ar@{=>}[r]^{\epsilon*F} & F & & G \ar@{=>}[r]^{\eta*G} & GFG \ar@{=>}[r]^{G*\eta} & G
}
\end{equation}
are identities (\cite[IV.1 Theorem 1 (ii)]{cat}).

If $X$ is a scheme, then we let ${\sf{Qcoh}}X$\index{qcohx@${\sf{Qcoh}} X$} denote the category of quasi-coherent ${\mathcal{O}}_{X}$-modules and we let ${\sf{Coh}}X$\index{cohx@${\sf{Coh }}X$} denote the category of coherent ${\mathcal{O}}_{X}$-modules.  We say $\mathcal{M}$ is an ${\mathcal{O}}_{X}$-module if $\mathcal{M}$ is a quasi-coherent ${\mathcal{O}}_{X}$-module.

Throughout this paper, $S$ will denote an affine, noetherian
scheme and, from Chapter 3 onwards, $\sf{S}$ will denote the
category of affine, noetherian $S$-schemes.  Unless otherwise
stated, all schemes are $S$ schemes so that products are over $S$.
If $f:V \rightarrow U$ is a morphism in $\sf{S}$, we let
$\tilde{f}=f \times \id_{X}:V \times X \rightarrow U \times X$,
and we let ${\tilde{f}}^{2}:(V\times X) \times (V\times X)
\rightarrow (U \times X) \times (U \times X)$ denote the map
$\tilde{f} \times \tilde{f}$.  We will reintroduce this notation
throughout to avoid confusion.

If $U$ is a scheme and $Y$ is a $U$-scheme, then we denote the $n$-fold product of $Y$ over $U$ by $(Y)_{U}^{n}$.  Thus, for example, $(X \times Y)_{Y}^{2}=(X \times Y) \times_{Y} (X \times Y)$.  This is, roughly speaking, the subscheme of $(X \times Y) \times (X \times Y)$ where the $Y$-coordinates are equal.

We will often abuse notation as follows:  if $q: U \rightarrow X$ and $r:U \rightarrow Y$ are morphisms of schemes and if $d:U \rightarrow U \times U$ is the diagonal morphism, then we sometimes write $q \times r$ to denote the morphism $(q \times r) \circ d:U \rightarrow X \times Y$.  When this abuse is employed, we write $q \times r:U \rightarrow X \times Y$ or we indicate that $(q \times r)_{*}$ is a functor from the category of quasi-coherent ${\mathcal{O}}_{U}$-modules to the category of quasi-coherent ${\mathcal{O}}_{X \times Y}$-modules.

{\bf Acknowledgments.}  We thank M. Artin, C. Ingalls, D. Patrick, P. Perkins, S.P. Smith, M. Van den Bergh and J. Zhang for numerous helpful conversations and feedback on earlier versions of this paper.  We thank the referee for their numerous useful comments.  We are especially indebted to S.P. Smith for suggesting the study of points on quantum ruled surfaces, and for providing invaluable guidance and support throughout the formulation of the ideas in this paper.

\chapter{Compatibilities on Squares}

When working with several functors simultaneously, one sometimes finds different ways of constructing canonical natural transformations between their compositions.  For example, let $F: \sf{A} \rightarrow \sf{B}$ be a functor and suppose $(F,G)$ is an adjoint pair with unit and counit $(\eta,\epsilon)$.  Then there are two obvious canonical morphisms between $F$ and itself.  There is an identity transformation $\id_{F}:F \Longrightarrow F$, and there is also the transformation
$$
\xymatrix{
F \ar@{=>}[r]^{F * \eta} & FGF \ar@{=>}[r]^{\epsilon * F} & F,
}
$$
where the symbol ``$*$'' denotes the horizontal product of natural transformations (Definition \ref{def.horprod}).  This map equals the identity transformation (Equation \ref{eqn.defad}).  We call equalities of this kind {\it compatibilities}\index{compatibilities|textbf}, although our use of this term is not standard.  More generally, we call any commuting diagram of canonical morphisms between objects in a category a compatibility.  Authors faced with potential compatibilities sometimes assume all such diagrams commute since tedious proofs of the commutativity of such ``natural'' diagrams exist on a case by case basis.  For an example of confessions of these assumptions, see \cite[Chap. II.6]{residues}.

The purpose of this chapter is to establish many of the compatibilities we need in order to define and represent the functor $\Gamma_{n}$.  Suppose $X$ is a noetherian separated scheme, $f:V \rightarrow U$ is a morphism of affine, noetherian schemes, $\tilde{f}=f \times \id_{X}$ and $\mathcal{B}$ is an ${\mathcal{O}}_{U \times X}$-bimodule algebra.  In order to define $\Gamma_{n}$, we must show that $\tilde{f}^{2*}\mathcal{B}$ inherits an ${\mathcal{O}}_{V \times X}$-bimodule algebra structure from $\mathcal{B}$.  The proof of this fact reduces to the verification that several large diagrams, one involving the associativity of the bimodule tensor product, the others involving the left and right scalar multiplication maps, commute.  Rather then check these complicated compatibilities directly, we note that both the associativity of the tensor product and the scalar multiplication maps are composed of canonical maps (Proposition \ref{prop.tensor}, Proposition \ref{prop.scalar}), each of which can be described as data associated to a {\it square} (Definition \ref{def.bigsquares}).  Motivated by this observation, we study squares in general.  Since the compatibilities we must verify are neatly stated in the language of indexed categories, we introduce this language (Definitions \ref{def.index}, \ref{def.indexfunct} and \ref{def.intrans}).  Then, using the theory of squares we have developed, we reduce the compatibilities we must verify to more easily checked compatibilities (Propositions \ref{prop.indexedfunct} and \ref{cor.reduce}).

The definition of a square (Definition \ref{def.bigsquares}) is inspired by the following situation in algebraic geometry:  Let $Y, Y', Z$ and $Z'$ be schemes with maps between them as in the following diagram:

\begin{equation} \label{eqn.commute}
\xymatrix{
Y' \ar[r]^{p'} \ar[d]_{f'} & Y \ar[d]^{f} \\
Z' \ar[r]_{p} & Z.
}
\end{equation}
Suppose that diagram (\ref{eqn.commute}) commutes.  We then have the equality
$$
\Phi:  f_{*}{p'}_{*} = p_{*}{f'}_{*}.
$$
Applying $p^{*}$ to the left of this equation and ${p'}^{*}$ to the right, we have the equality

$$
p^{*}f_{*}{p'}_{*}{p'}^{*}=p^{*}p_{*}{f'}_{*}{p'}^{*}.
$$
Composing this with various adjointness maps, we have the {\bf 2-cell dual to $\Phi$ induced by (\ref{eqn.commute})}:\index{2-cell dual to phi@2-cell dual to $\Phi$!induced by a diagram of schemes|textbf}
$$
\Psi:  p^{*}f_{*}\Longrightarrow  p^{*}f_{*}{p'}_{*}{p'}^{*}=p^{*}p_{*}{f'}_{*}{p'}^{*} \Longrightarrow {f'}_{*}{p'}^{*}.
$$
We show (Lemma \ref{prop.duality}) that $\Phi$ and $\Psi$ are dual to each other.  It is this duality which allows us to reduce a property of $\Psi$ to its dual $\Phi$.  In the example above, since $\Phi$ is trivial, we expect $\Phi$ to be an easier map to study then $\Psi$.  This expectation is not disappointed (Corollary \ref{cor.indexcan}).

We employ various results from the theory of squares elsewhere.
Corollary \ref{cor.adcom2} allows us to establish useful
compatibilities in the category of sets, and is crucial in the
proof that the bimodule Segre embedding has desired properties. We
also use the Corollary to show that our representation of
$\Gamma_{n}$ is natural (Proposition \ref{prop.trans}).
Propositions \ref{prop.adcom3} and \ref{prop.horizontal} allow us
to relate the canonical and dual 2-cell of a diagram to the
canonical and dual 2-cells of various subdiagrams.  These results
reflect the duality between $\Phi$ and $\Psi$.

\section{2-Categories}
We review the definition of 2-categories since many of our results involve the 2-category of all categories.  The following definition is from \cite[Section 5]{stacks}.

\begin{definition} \label{def.2cat}
A {\bf 2-category \index{2 category@$2$-\mbox{category}|textbf}$\mathbb{T}$} consists of the following data:
\begin{enumerate}
\item{A class of objects $\operatorname{ob}\mathbb{T}$,}
\item{for each pair $\sf{X}$, ${\sf{Y}} \in \operatorname{ob}\mathbb{T}$, a category $\operatorname{Hom}({\sf{X}},{\sf{Y}})$ whose objects, called 1-morphisms, are represented by $\xymatrix{{\sf{X}} \ar[r]^{F} & {\sf{Y}}}$ and whose morphisms, called 2-morphisms or 2-cells, are represented by
$$
\xymatrix{
{\sf{X}} \rtwocell^F_G{\Delta} & {\sf{Y}}\\
}
$$
subject to the following conditions:}

\begin{enumerate}
\item{}
({\it Composition of 1-morphisms})  Given a diagram
$$
\xymatrix{
{\sf{X}} \ar[r]^{F} & {\sf{Y}} \ar[r]^{G} & {\sf{Z}}
}
$$
there exists a diagram
$$
\xymatrix{
{\sf {X}} \ar[r]^{G \circ F} & {\sf{Z}}
}
$$
and this composition is associative,
\item{}
({\it Identity for 1-morphisms})  For each object ${\sf{X}}$ there is a 1-morphism $\operatorname{id}_{{\sf{X}}}:{\sf{X}} \rightarrow {\sf{X}}$ such that if $F: {\sf{X}} \rightarrow {\sf{Y}}$ is a 1-morphism, then $F \circ \operatorname{id}_{{\sf{X}}}=\operatorname{id}_{{\sf{Y}}} \circ F=F$,
\item{}
({\it Vertical composition of 2-morphisms}) Given a diagram
$$
\xymatrix{
{\sf{X}}\ruppertwocell^F{\Delta} \rlowertwocell_H{\theta} \ar[r]_{\hskip -.2in G} & {\sf{Y}} & \mbox{there exists} & {\sf{X}}\rrtwocell^F_H{\hskip .1in \theta \circ \Delta}  &  & {\sf{Y}}
}
$$
and this composition is associative,
\item{}
({\it Horizontal composition of 2-morphisms}) Given a diagram
$$
\xymatrix{
{\sf{X}} \rtwocell^F_{F'}{\Delta} & {\sf{Y}} \rtwocell^G_{G'}{\Theta} & {\sf{Z}} & \mbox{there exists} & {\sf{X}} \rrtwocell^{G \circ F}_{G' \circ F'}{\hskip .1in \Theta * \Delta} & & {\sf{Z}}
}
$$
and this composition is associative,
\item{}
({\it Identity for 2-morphisms})  For every 1-morphism $F$ there is a 2-morphism $\operatorname{id}_{F}$ (or sometimes just $F$ or $\operatorname{id}$ when the context is clear) such that, given a diagram
$$
\xymatrix{
{\sf{X}} \rtwocell^F_{F'}{\Delta} & {\sf{Y}} \rtwocell^G_{G}{\operatorname{id}} & {\sf{Z}}
}
$$
$\Delta \circ \operatorname{id}_{F}=\operatorname{id}_{F'} \circ \Delta$ and $\operatorname{id}_{G}*\operatorname{id}_{F}=\operatorname{id}_{G \circ F}$,
\item{}
(\it Compatibility between horizontal and vertical composition of 2-cells) Given a diagram
$$
\xymatrix{
{\sf{X}}\ruppertwocell^F{\Delta} \rlowertwocell_{F''}{\Delta'} \ar[r]_{\hskip -.2in F'} & {\sf{Y}} \ruppertwocell^G{\Theta} \rlowertwocell_{G''}{\Theta'} \ar[r]_{\hskip -.2in G'} & {\sf{Z}}
}
$$
we have $(\Theta' \circ \Theta)*(\Delta'\circ \Delta)=(\Theta' * \Delta') \circ (\Theta * \Delta)$.
\end{enumerate}
\end{enumerate}

\end{definition}
Let ${\sf Cat}$\index{cat@${\sf Cat}$|textbf} denote the category of categories, whose objects are (small) categories and whose morphisms are functors between categories.  If $\sf{A}$ and $\sf{B}$ are categories, then the collection of functors between them, $\operatorname{Hom}(\sf{A},\sf{B})$ forms a category which has natural transformations between functors as morphisms.  We describe a horizontal composition of natural transformations, $*$, (2d in Definition \ref{def.2cat}) which is associative and satisfies 4e,f.
\begin{definition} \label{def.horprod}\index{cat@${\sf Cat}$!horizontal composition of 2-morphisms in|textbf}
Suppose we are given a diagram of categories, functors, and 2-cells:
$$
\xymatrix{
\sf{X} \rtwocell^F_{F'}{\Delta} & \sf{Y} \rtwocell^G_{G'}{\Theta} & \sf{Z}.
}
$$
The {\bf horizontal composition of $\Delta$ and $\Theta$}, denoted $\Theta * \Delta$ is defined, for every object $X$ of $\sf{X}$, by the formula

$$(\Theta * \Delta)_{X} = \Theta_{F'X} \circ G(\Delta_{X})=G'(\Delta_{X}) \circ \Theta_{FX}$$
\end{definition}
This composition makes $\sf{Cat}$ a 2-category \cite[Theorem 1, p. 44]{cat}.

\section{The category of squares} \index{category of squares|(}
We define the category of squares and describe some of its
properties.  Applications of the theory will be deferred to later
chapters.
\begin{definition} \label{def.bigsquares}\index{category of squares|textbf}
Let $\operatorname{Sq}$ be the category defined as follows:
\begin{itemize}
\item{}
The objects of $\operatorname{Sq}$ consist of
\begin{enumerate}
\item{}
a diagram of categories and functors
$$
\xymatrix{
\sf{A} \ar@<1ex>[r]^{G'} \ar[d]_{H} & \sf{B} \ar@<1ex>[l]^{F'} \ar[d]^{I}\\
\sf{C} \ar@<1ex>[r]^{G}  & \sf{D} \ar@<1ex>[l]^{F}.
}
$$
such that $(F,G)$ and $(F',G')$ are adjoint pairs with units/counits $(\eta,\epsilon)$ and $(\eta',\epsilon')$ respectively, and
\item{}
a pair of natural transformations
$$
\xymatrix{
\Phi:IG' \ar@{=>}[r] & GH & & \Psi:FI \ar@{=>}[r] & HF'
}
$$
such that $\Psi$ is the composition
\begin{equation} \label{eqn.psi}
\xymatrix{
FI \ar@{=>}[rr]^{FI*\eta} & & FIG'F' \ar@{=>}[rr]^{F * \Phi * F'} & & FGHF' \ar@{=>}[rr]^{\epsilon * HF'} & & HF'.
}
\end{equation}
An object composed of this data will be denoted
\begin{equation} \label{eqn.definee}
\xymatrix{
\sf{A} \ar@<1ex>[r]^{G'} \ar[d]_{H} \ar@{=>}[dr]|{\Phi}\hole & \sf{B} \ar@<1ex>[l]^{F'} \ar[d]^{I}\\
\sf{C} \ar@<1ex>[r]^{G}  & \sf{D} \ar@<1ex>[l]^{F}.
}
\end{equation}
We call $\Phi$ the {\bf canonical 2-cell associated to (\ref{eqn.definee})}\index{canonical 2-cell associated to a square|textbf} and we call $\Psi$ the {\bf 2-cell dual to (\ref{eqn.definee})} or the {\bf 2-cell dual to $\Phi$}\index{2-cell dual to phi@2-cell dual to $\Phi$|textbf}.
\end{enumerate}
\item{}
A morphism $(J,\Theta)$
\begin{equation} \label{eqn.morphine}
\xymatrix{
\sf{A} \ar@<1ex>[r]^{G'} \ar[d]_{H} \ar@{=>}[dr]|{\Phi}\hole & \sf{B} \ar@<1ex>[l]^{F'} \ar[d]^{I} & \ar[r]  & & \overline{\sf{A}} \ar@<1ex>[r]^{\overline{G'}} \ar[d]_{\overline{H}} \ar@{=>}[dr]|{\overline{\Phi}}\hole & \overline{\sf{B}} \ar@<1ex>[l]^{\overline{F'}} \ar[d]^{\overline{I}} \\
\sf{C} \ar@<1ex>[r]^{G}  & \sf{D} \ar@<1ex>[l]^{F} & & & \overline{\sf{C}} \ar@<1ex>[r]^{\overline{G}}  & \overline{\sf{D}} \ar@<1ex>[l]^{\overline{F}}.
}
\end{equation}
in $\operatorname{Sq}$ consists of
\begin{enumerate}
\item{}
functors
$$
J^{\sf{A,B,C,D}}: \sf{A},\sf{B},\sf{C},\sf{D} \rightarrow \overline{A},\overline{B},\overline{C},\overline{D}
$$
and
\item{}
natural transformations
$$
\xymatrix{
\overline{F}J^{\sf{D}} \ar@{=>}[r]^{\Theta_{\overline{F}}} & J^{\sf{C}}F & &\overline{F'}J^{\sf{B}} \ar@{=>}[r]^{\Theta_{\overline{F'}}} & J^{\sf{A}}F' \\
J^{\sf{D}}G \ar@{=>}[r]^{\Theta_{G}} & \overline{G}J^{\sf{C}} & & J^{\sf{B}}G' \ar@{=>}[r]^{\Theta_{G'}} & \overline{G'}J^{\sf{A}} \\
\overline{H}J^{\sf{A}} \ar@{=>}[r]^{\Theta_{\overline{H}}} & J^{\sf{C}}H & &\overline{I}J^{\sf{B}} \ar@{=>}[r]^{\Theta_{\overline{I}}} & J^{\sf{D}}I
}
$$
such that the diagrams
\begin{enumerate}
\item{}
$$
\xymatrix{
J^{\sf{B}} \ar@{=>}[r]^{\overline{\eta'}*J^{\sf{B}}} \ar@{=>}[d]_{J^{\sf{B}* \eta'}} & \overline{G'}\overline{F'}J^{\sf{B}} \ar@{=>}[d]^{\overline{G'}*\Theta_{\overline{F'}}} & & \overline{F}J^{\sf{D}}G \ar@{=>}[r]^{\overline{F}*\Theta_{G}} \ar@{=>}[d]_{\Theta_{\overline{F}}* G} & \overline{F} \overline{G}J^{\sf{C}} \ar@{=>}[d]^{\overline{\epsilon}*J^{\sf{C}}} \\
J^{\sf{B}}G'F' \ar@{=>}[r]_{\Theta_{G'}*F'} &  \overline{G'}J^{\sf{A}}F' & & J^{\sf{C}}FG \ar@{=>}[r]_{J^{\sf{C}}*\epsilon} & J^{\sf{C}}
}
$$
and
\item{}
$$
\xymatrix{
& \overline{I}J^{B}G' \ar@{=>}[dr]^{\Theta_{\overline{I}}*G'} \ar@{=>}[dl]_{\overline{I}*\Theta_{G'}} & \\
\overline{I}\overline{G'}J^{\sf{A}} \ar@{=>}[d]_{\overline{\Phi}*J^{\sf{A}}}   & & J^{\sf{D}}IG' \ar@{=>}[d]^{J^{\sf{D}}*\Phi} \\
\overline{G}\overline{H}J^{\sf{A}} \ar@{=>}[dr]_{\overline{G}*\Theta_{\overline{H}}}   & & J^{\sf{D}}GH \ar@{=>}[dl]^{\Theta_{G}*H}  \\
& \overline{G}J^{\sf{C}}H  &
}
$$
commute.
\end{enumerate}
\end{enumerate}
\end{itemize}
Suppose $(J,\Theta)$ is a morphism (\ref{eqn.morphine}) and $(K,\Delta)$ is a morphism whose domain equals the codomain of $(J, \Theta)$.  Suppose the codomain of $(K,\Delta)$ is the square
$$
\xymatrix{
\sf{\underline{A}} \ar@<1ex>[r]^{\underline{G'}} \ar[d]_{\underline{H}} \ar@{=>}[dr]|{\underline{\Phi}}\hole & \sf{\underline{B}} \ar@<1ex>[l]^{\underline{F'}} \ar[d]^{\underline{I}}\\
\sf{\underline{C}} \ar@<1ex>[r]^{\underline{G}}  & \sf{\underline{D}} \ar@<1ex>[l]^{\underline{F}}.
}
$$
Then the composition $(K,\Delta) \circ (J,\Theta)$ is defined to be the pair $(KJ,\Xi)$ where $\Xi$ is defined as follows:
$$
\xymatrix{
\Xi_{\underline{F}} = (K^{\sf{C}}*\Theta_{\overline{F}}) \circ (\Delta_{\underline{F}}*J^{\sf{D}})
}
$$
and $\Xi_{\underline{F'},\underline{H},\underline{I}}$ are defined similarly, and
$$
\xymatrix{
\Xi_{\underline{G}} =  (\Delta_{\overline{G}}*J^{\sf{C}})\circ (K^{\sf{D}}*\Theta_{G})
}
$$
and $\Xi_{\underline{G'}}$ is defined similarly.  It is an easy but tedious exercise to show that the composition $(KJ,\Xi)$ is actually a morphism.  We call $\operatorname{Sq}$ the {\bf category of squares}.
\end{definition}
{\it Remark.}  An arbitrary square (\ref{eqn.definee}) will have an arbitrary $\Phi$.  Hence, calling $\Phi$ canonical might seem ill advised.  However, the following examples show that there exists a canonical choice for $\Phi$ in important situations.
\begin{example} \label{example.projform}
Let $Y, Y', Z$ and $Z'$ be schemes with maps between them as in the following diagram:

\begin{equation} \label{eqn.commutee}
\xymatrix{
Y' \ar[r]^{p'} \ar[d]_{f'} & Y \ar[d]^{f} \\
Z' \ar[r]_{p} & Z.
}
\end{equation}
Suppose (\ref{eqn.commutee}) commutes.  If
$$
\Phi:  f_{*}{p'}_{*} \Longrightarrow p_{*}{f'}_{*}.
$$
is equality, then
$$
\xymatrix{
{\sf{Qcoh}} Y' \ar@<1ex>[r]^{{p'}_{*}} \ar[d]_{{f'}_{*}} \ar@{=>}[dr]|{\Phi}\hole & {\sf{Qcoh}} Y \ar@<1ex>[l]^{{p'}^{*}} \ar[d]^{f_{*}}\\
{\sf{Qcoh}} Z' \ar@<1ex>[r]^{p_{*}}  & {\sf{Qcoh}} Z \ar@<1ex>[l]^{p^{*}}.
}
$$
is a square.  We sometimes call $\Phi$ the {\bf canonical 2-cell induced by (\ref{eqn.commutee})} or the {\bf canonical 2-cell of the square induced by (\ref{eqn.commutee})}.  We sometimes say $\Psi$ is the {\bf dual 2-cell to (\ref{eqn.commutee})} or the {\bf dual 2-cell to the square induced by (\ref{eqn.commutee})}.
\end{example}

\begin{example} \label{example.basechange}
Let $f:Y \rightarrow X$ be a map of schemes.  For each $\mathcal{P}$ in ${\sf{Qcoh}}Y$, let
$$
\Phi^{\mathcal{P}}:  f_{*}\mathcal{P} \otimes f_{*}- \Longrightarrow f_{*}(\mathcal{P} \otimes -).
$$
be the natural transformation
$$
f_{*}\mathcal{P} \otimes f_{*}- \Longrightarrow f_{*}f^{*}(f_{*}\mathcal{P} \otimes f_{*}- ) \Longrightarrow f_{*}(f^{*}f_{*}\mathcal{P} \otimes f^{*}f_{*} -) \Longrightarrow f_{*}(\mathcal{P} \otimes -).
$$
Then
\begin{equation} \label{eqn.project}
\xymatrix{
{\sf{Qcoh}}Y  \ar@<1ex>[r]^{f_{*}} \ar[d]_{f_{*}(\mathcal{P}\otimes -)} \ar@{=>}[dr]|{\Phi^{\mathcal{P}}}\hole  & {\sf{Qcoh}} X \ar@<1ex>[l]^{f^{*}} \ar[d]^{f_{*}\mathcal{P} \otimes -}\\
{\sf{Qcoh}}X \ar@<1ex>[r]^{\operatorname{id}}  & {\sf{Qcoh}}X \ar@<1ex>[l]^{\operatorname{id}}.
}
\end{equation}
is a square, and the dual to (\ref{eqn.project}),
$$
\Psi^{\mathcal{P}}:f_{*}\mathcal{P} \otimes - \Longrightarrow f_{*}(\mathcal{P} \otimes f^{*}-)
$$
is the {\it projection formula}\index{projection formula|textbf}.
\end{example}

\subsection{Elementary properties of squares}\index{canonical 2-cell associated to a square|(}\index{2-cell dual to phi@2-cell dual to $\Phi$|(}

\begin{lemma} \label{lem.precom}
If
$$
\xymatrix{
\sf{A} \ar@<1ex>[r]^{G'} \ar[d]_{H} \ar@{=>}[dr]|{\Phi}\hole & \sf{B} \ar@<1ex>[l]^{F'} \ar[d]^{I}\\
\sf{C} \ar@<1ex>[r]^{G}  & \sf{D} \ar@<1ex>[l]^{F}.
}
$$
is an object in $\operatorname{Sq}$, then the diagrams
$$
\xymatrix{
I \ar@{=>}[rr]^{I * \eta'} \ar@{=>}[d]_{\eta * I} & & IG'F' \ar@{=>}[d]^{ \eta*(IG'F')}  &  & H & & HF'G' \ar@{=>}[ll]_{H*\epsilon'}   \\
GFI \ar@{=>}[rr]_{(GFI)* \eta'} & & GFIG'F' & & FGH \ar@{=>}[u]^{\epsilon*H} & & FGHF'G' \ar@{=>}[u]_{\epsilon*(HF'G')} \ar@{=>}[ll]^{(FGH)*\epsilon'}
}
$$
commute.
\end{lemma}

\begin{proof}
We first prove the left diagram commutes.  To prove this assertion, we apply 2f of Definition \ref{def.2cat} to the diagram
$$
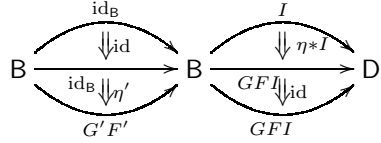
\xymatrix{
{\sf{B}}\rruppertwocell^{\operatorname{id}_{\sf{B}}}{\operatorname{id}} \rrlowertwocell_{G'F'}{\eta'} \ar[rr]_{\hskip -.25in \operatorname{id}_{\sf{B}}}  & & {\sf{B}} \rruppertwocell^I{\hskip .15in \eta*I} \rrlowertwocell_{\hskip -.12in GFI}{\operatorname{id}} \ar[rr]_{\hskip -.25in GFI}  & & {\sf{D}}
}
$$
and conclude that $(\eta * I)*\eta'=(GFI*\eta') \circ (\eta * I)$.  But $(\eta*I)*\eta'=\eta*(I*\eta')$ by 2d of Definition \ref{def.2cat}.  We claim $\eta*(I*\eta')=(\eta*IG'F') \circ (I*\eta')$.  In fact, this follows easily by applying 2f of Definition \ref{def.2cat} to the diagram
$$
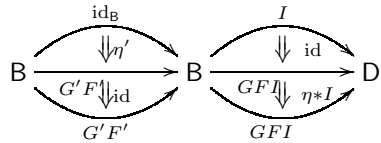
\xymatrix{
{\sf{B}}\rruppertwocell^{\operatorname{id}_{\sf{B}}}{\eta'} \rrlowertwocell_{G'F'}{\operatorname{id}} \ar[rr]_{\hskip -.25in G'F'}  & & {\sf{B}} \rruppertwocell^I{\hskip .15in \operatorname{id}} \rrlowertwocell_{\hskip -.12in GFI}{\hskip .2in \eta*I} \ar[rr]_{\hskip -.25in GFI}  & & {\sf{D}}
}
$$
Thus, the left square does indeed commute.  The proof that the right diagram commutes is similar, and we omit it.
\end{proof}

Objects in the category $\operatorname{Sq}$ enjoy the following compatibilities:

\begin{proposition} \label{prop.adcom}
If
$$
\xymatrix{
\sf{A} \ar@<1ex>[r]^{G'} \ar[d]_{H} \ar@{=>}[dr]|{\Phi}\hole & \sf{B} \ar@<1ex>[l]^{F'} \ar[d]^{I}\\
\sf{C} \ar@<1ex>[r]^{G}  & \sf{D} \ar@<1ex>[l]^{F}.
}
$$
is an object in $\operatorname{Sq}$, then the diagrams
$$
\xymatrix{
I \ar@{=>}[r]^{\eta * I} \ar@{=>}[d]_{I * \eta'} & GFI \ar@{=>}[d]^{G*\Psi} & &FIG' \ar@{=>}[r]^{\Psi*G'} \ar@{=>}[d]_{F*\Phi} & HF'G' \ar@{=>}[d]^{H*\epsilon'}  \\
IG'F' \ar@{=>}[r]_{\Phi * F'} & GHF' & & FGH \ar@{=>}[r]_{\epsilon*H} & H
}
$$
commute.
\end{proposition}

\begin{proof}
We show the left diagram commutes.  The argument that the right diagram commutes is similar, so we omit it.  We claim
$$
\xymatrix{
I \ar@{=>}[rr]^{I * \eta'} \ar@{=>}[d]_{\eta * I} & & IG'F' \ar@{=>}[rr]^{\Phi * F'} \ar@{=>}[d]^{\eta  * (IG'F')}  & & GHF' \ar@{=>}[d]^{\eta*GHF'}   \\
GFI \ar@{=>}[rr]_{(GFI)*\eta} & & GFIG'F' \ar@{=>}[rr]_{GF*\Phi*F'} & & GFGHF'
}
$$
commutes.  For, by Lemma \ref{lem.precom}, the left square commutes, while the right square commutes by naturality of $\eta$.  Since
$$
\xymatrix{
G \ar@{=>}[r]^{\eta * G} & GFG \ar@{=>}[r]^{G * \epsilon} & G
}
$$
is the identity map, the diagram
$$
\xymatrix{
I \ar@{=>}[rr]^{I * \eta'} \ar@{=>}[d]_{\eta * I} & & IG'F' \ar@{=>}[rr]^{\Phi * F'} & & GHF'  \\
GFI \ar@{=>}[rr]_{(GFI)*\eta} & & GFIG'F' \ar@{=>}[rr]_{GF*\Phi*F'} & & GFGHF' \ar@{=>}[u]_{G*\epsilon*HF'}
}
$$
commutes also, which is just what we needed to establish.
\end{proof}

\begin{lemma} \label{prop.duality}
(Duality)\index{duality of squares|textbf}.  Suppose
\begin{equation} \label{eqn.firstsquare}
\xymatrix{
\sf{A} \ar@<1ex>[r]^{G'} \ar[d]_{H} \ar@{=>}[dr]|{\Phi}\hole & \sf{B} \ar@<1ex>[l]^{F'} \ar[d]^{I}\\
\sf{C} \ar@<1ex>[r]^{G}  & \sf{D} \ar@<1ex>[l]^{F}.
}
\end{equation}
is a square.  Then the map
$$
f:{\sf Nat}(IG',GH) \rightarrow {\sf Nat}(FI,HF')
$$
sending $\Phi:  IG' \Longrightarrow GH$ to the natural transformation $\Psi:FI \Longrightarrow HF'$ given by the composition
$$
\xymatrix{
FI \ar@{=>}[r]^{FI*\eta'} & FIG'F' \ar@{=>}[r]^{F * \Phi * F'} & FGHF' \ar@{=>}[r]^{\epsilon * HF'} & HF'
}
$$
is a bijection whose inverse, $g$, sends $\Psi:FI \Longrightarrow HF'$ to the composition
$$
\xymatrix{
IG' \ar@{=>}[r]^{\eta*IG'} & GFIG' \ar@{=>}[r]^{G * \Psi * G'} & GHF'G' \ar@{=>}[r]^{GH*\epsilon'} & GH.
}
$$
\end{lemma}

\begin{proof}
We show that $gf(\Phi)=\Phi$.  We must show that the diagram
\begin{equation} \label{eqn.orig}
\xymatrix{
IG' \ar@{=>}[rr]^{\Phi} \ar@{=>}[d]_{\eta*IG'} & & GH \\
GFIG' \ar@{=>}[d]_{GFI*\eta'*G'} & & GHF'G' \ar@{=>}[u]_{GH*\epsilon'} \\
GFIG'F'G' \ar@{=>}[rr]_{GF*\Phi*F'G'} & & GFGHF'G' \ar@{=>}[u]_{G*\epsilon*HF'G'}
}
\end{equation}
commutes.  We note that
$$
\xymatrix{
IG' \ar@{=>}[rr]^{\Phi} \ar@{=>}[d]_{\eta*IG'} & & GH \ar@{=>}[d]^{\eta*GH} \\
GFIG' \ar@{=>}[rr]_{GF*\Phi} & & GFGH  \\
GFIG'F'G' \ar@{=>}[rr]_{GF*\Phi*F'G'} \ar@{=>}[u]^{GFIG'*\epsilon'}  & & GFGHF'G' \ar@{=>}[u]_{GFGH*\epsilon'}
}
$$
commutes by naturality of $\eta$ and $\epsilon'$, so that
$$
\xymatrix{
IG' \ar@{=>}[rr]^{\Phi} \ar@{=>}[d]_{\eta*IG'} & & GH  \\
GFIG' \ar@{=>}[rr]_{GF*\Phi} \ar@{=>}[d]_{GFI*\eta'*G'} & & GFGH \ar@{=>}[u]_{G*\epsilon*H} \\
GFIG'F'G' \ar@{=>}[rr]_{GF*\Phi*F'G'}  & & GFGHF'G' \ar@{=>}[u]_{GFGH*\epsilon'}
}
$$
commutes also.  Thus, the result will follow provided
$$
\xymatrix{
H & \\
FGH \ar@{=>}[u]^{\epsilon * H} & HF'G' \ar@{=>}[ul]_{H*\epsilon'} \\
FGHF'G' \ar@{=>}[u]^{FGH*\epsilon'} \ar@{=>}[ur]_{\epsilon*HF'G'}
}
$$
commutes, since we may attach $G$ applied to this diagram to the right of the previous one and recover diagram \ref{eqn.orig}.  This last diagram commutes by Lemma \ref{lem.precom}.  Thus $gf=\operatorname{id}$ as desired.
\end{proof}
{\it Remark.}  Compare \cite[2.1, p.167]{2cell}.

\begin{corollary} \label{cor.solu}
If
$$
\xymatrix{
\sf{A} \ar@<1ex>[r]^{G'} \ar[d]_{H} \ar@{=>}[dr]|{\Phi}\hole & \sf{B} \ar@<1ex>[l]^{F'} \ar[d]^{I}\\
\sf{C} \ar@<1ex>[r]^{G}  & \sf{D} \ar@<1ex>[l]^{F}.
}
$$
is a square, then $\Delta:FI \Longrightarrow HF'$ has the property that
\begin{equation} \label{eqn.mockup}
\xymatrix{
FIG' \ar@{=>}[r]^{\Delta*G'} \ar@{=>}[d]_{F*\Phi} & HF'G' \ar@{=>}[d]^{H*\epsilon'} \\
FGH \ar@{=>}[r]_{\epsilon*H} & H
}
\end{equation}
commutes if and only if $\Delta=\Psi$, and $\Theta:IG' \Longrightarrow GH$ has the property that
$$
\xymatrix{
I \ar@{=>}[r]^{\eta * I} \ar@{=>}[d]_{I * \eta'} & GFI \ar@{=>}[d]^{G*\Psi} \\
IG'F' \ar@{=>}[r]_{\Theta * F'} & GHF' \\
}
$$
commutes if and only if $\Theta=\Phi$.
\end{corollary}

\begin{proof}
We establish the first assertion.  The latter fact is proved in a similar fashion.  The diagram
$$
\xymatrix{
FI \ar@{=>}[d] \ar@{=>}[rr]^{\Delta} & & HF' \ar@{=>}[d] \\
FIG'F' \ar@{=>}[d]_{F* \Phi * F'} \ar@{=>}[rr]^{\Delta*G'F'} & & HF'G'F' \ar@{=>}[d]^{H*\epsilon'*F'} \\
FGHF' \ar@{=>}[rr]_{\epsilon*HF'} & & HF'
}
$$
commutes:  the top commutes by the naturality of $\eta'$ while the bottom commutes since (\ref{eqn.mockup}) commutes.  The right vertical is the identity while the left route is $\Psi$, so $\Phi=\Delta$.  The converse holds by Proposition \ref{prop.adcom}.
\end{proof}

\begin{corollary} \label{cor.adcom2}
If
$$
\xymatrix{
\sf{A} \ar@<1ex>[r]^{G'} \ar[d]_{H} \ar@{=>}[dr]|{\Phi}\hole & \sf{B} \ar@<1ex>[l]^{F'} \ar[d]^{I}\\
\sf{C} \ar@<1ex>[r]^{G}  & \sf{D} \ar@<1ex>[l]^{F}.
}
$$
is a square, then the diagram
\begin{equation} \label{eqn.digge}
\xymatrix{
\mbox{Hom}_{\sf{A}}(F'\mathcal{B},\mathcal{A}) \ar[r] \ar[d]_{H(-)} & \mbox{Hom}_{\sf{B}}(\mathcal{B},G'\mathcal{A}) \ar[d]^{I(-)} \\
\mbox{Hom}_{\sf{C}}(HF'\mathcal{B},H{\mathcal{A}}) \ar[d]_{- \circ \Psi_{\mathcal{B}}} & \mbox{Hom}_{\sf{D}}(I\mathcal{B},IG'\mathcal{A}) \ar[d]^{\Phi_{\mathcal{A}} \circ -} \\
\mbox{Hom}_{\sf{C}}(FI\mathcal{B},H\mathcal{A}) \ar[r] & \mbox{Hom}_{\sf{D}}(I\mathcal{B},GH\mathcal{A})
}
\end{equation}
whose horizontal maps are the adjoint isomorphisms, commutes.
\end{corollary}

\begin{proof}
Let
$$
\mu:F'\mathcal{B} \rightarrow \mathcal{A}.
$$
Then $\mu$ goes, via the left route of the diagram, to the map
$$
\xymatrix{
I\mathcal{B} \ar[r]^{(\eta * I)_{\mathcal{B}}} & GFI\mathcal{B} \ar[r]^{(G*\Psi)_{\mathcal{B}}} & GHF'\mathcal{B} \ar[r]^{GH\mu} & GH\mathcal{A}
}
$$
On the other hand, $\mu$ goes, via the right route of the diagram, to the map
$$
\xymatrix{
I\mathcal{B} \ar[r]^{(I*\eta')_{\mathcal{B}}} & IG'F'\mathcal{B} \ar[r]^{IG'\mu} \ar[d]_{(\Phi*F')_{\mathcal{B}}} & IG'\mathcal{A} \ar[d]^{\Phi_{\mathcal{A}}} \\
& GHF'\mathcal{B} \ar[r]_{GH\mu} & GH \mathcal{A}
}
$$
where the right square commutes by naturality of $\Phi$.  To show that these two maps are the same, it suffices to show that the diagram of functors
$$
\xymatrix{
I \ar@{=>}[r]^{I*\eta'} \ar@{=>}[d]_{\eta*I} & IG'F' \ar@{=>}[d]^{\Phi*F'} \\
GFI \ar@{=>}[r]_{G* \Psi} & GHF'
}
$$
commutes.  The commutivity of this diagram follows from Proposition \ref{prop.adcom}.
\end{proof}
The Corollary asserts that the data in the hypothesis of Proposition \ref{prop.adcom} gives an {\it adjoint square} \cite[ex. 4. p.101]{cat}.

\begin{lemma} \label{lem.simplecomp2}
If $J^{\sf{A,B,C,D}}:\sf{A,B,C,D} \rightarrow \sf{A,B,C,D}$ are identity functors, $\Theta_{\overline{H},\overline{I}}$ are identity morphisms, $\Theta_{\overline{F},\overline{F'},G,G'}$ are isomorphisms, and if $(J,\Theta)$ defines a morphism
$$
\xymatrix{
\sf{A} \ar@<1ex>[r]^{G'} \ar[d]_{H} \ar@{=>}[dr]|{\Phi}\hole & \sf{B} \ar@<1ex>[l]^{F'} \ar[d]^{I} & \ar[r]  & & \sf{A} \ar@<1ex>[r]^{\overline{G'}} \ar[d]_{{H}} \ar@{=>}[dr]|{\overline{\Phi}}\hole & \sf{B} \ar@<1ex>[l]^{\overline{F'}} \ar[d]^{{I}} \\
\sf{C} \ar@<1ex>[r]^{G}  & \sf{D} \ar@<1ex>[l]^{F} & & & \sf{C} \ar@<1ex>[r]^{\overline{G}}  & \sf{D} \ar@<1ex>[l]^{\overline{F}}.
}
$$
in $\operatorname{Sq}$ then the diagrams
\begin{equation} \label{eqn.adjointpic}
\xymatrix{
\operatorname{id} \ar@{=>}[r]^{\overline{\eta'}} \ar@{=>}[d]_{\eta'} & \overline{G'}\overline{F'} \ar@{=>}[d]^{\overline{G'}*\Theta_{\overline{F'}}} & & \overline{F}G \ar@{=>}[r]^{\overline{F}*\Theta_{G}} \ar@{=>}[d]_{\Theta_{\overline{F}}* G} & \overline{F} \overline{G} \ar@{=>}[d]^{\overline{\epsilon}} \\
G'F' \ar@{=>}[r]_{\Theta_{G'}*F'} &  \overline{G'}F' & & FG \ar@{=>}[r]_{\epsilon} & \operatorname{id}
}
\end{equation}
and
\begin{equation} \label{eqn.psipic}
\xymatrix{
\overline{F}I \ar@{=>}[d]_{\overline{\Psi}} \ar@{=>}[rr]^{\Theta_{\overline{F}}*I} & & FI \ar@{=>}[d]^{\Psi} \\
H\overline{F'} \ar@{=>}[rr]_{H*\Theta_{\overline{F}}} & & HF'
}
\end{equation}
commute.
\end{lemma}

\begin{proof}
In the proof, all 2-morphisms will be either a unit, a counit or $\Theta$.  Since it will be clear which is which, we will not label the maps.  Suppose $\Theta$ is a morphism.  Then, by definition of a morphism, (\ref{eqn.adjointpic}) and
$$
\xymatrix{
IG' \ar@{=>}[r] \ar@{=>}[d] & I\overline{G'} \ar@{=>}[d] \\
GH \ar@{=>}[r] & \overline{G}H
}
$$
commute.  Thus, the top center of the diagram
$$
\xymatrix{
FI \ar@{=>}[r] & FIG'F' \ar@{=>}[r] \ar@{=>}[d] & FGHF' \ar@{=>}[r] \ar@{=>}[d] & HF' \\
& FI\overline{G'}F' \ar@{=>}[r] & F\overline{G}HF' & \\
\overline{F}I  \ar@{=>}[uu] \ar@{=>}[r] & \overline{F}I\overline{G'}\overline{F'} \ar@{=>}[r] \ar@{=>}[u] & \overline{F}\overline{G}H\overline{F'} \ar@{=>}[r] \ar@{=>}[u] & H \overline{F'} \ar@{=>}[uu]
}
$$
commutes.  The bottom center of this diagram commutes by naturality of $\Theta_{\overline{F'}}$.  The fact that the left and right rectangles of the diagram commute follows readily from (\ref{eqn.adjointpic}) so we omit the details.
\end{proof}

\begin{definition} \label{def.squares}
Fix categories $\sf{A}$, $\sf{B}$, $\sf{C}$ and $\sf{D}$.  Let $\operatorname{Sq}(\sf{A},\sf{B},\sf{C},\sf{D})$ be the subset of $Ob$$\operatorname{Sq}$ consisting of all objects of the form
$$
\xymatrix{
\sf{A} \ar@<1ex>[r]^{G'} \ar[d]_{H} \ar@{=>}[dr]|{\Phi}\hole & \sf{B} \ar@<1ex>[l]^{F'} \ar[d]^{I}\\
\sf{C} \ar@<1ex>[r]^{G}  & \sf{D} \ar@<1ex>[l]^{F}.
}
$$
\end{definition}
We introduce three composition laws between various $\operatorname{Sq}$-categories.  There is a {\it vertical} composition from an appropriate subset of
$$
\operatorname{Sq}(\sf{A},\sf{B},\sf{C},\sf{D}) \times \operatorname{Sq}(\sf{C},\sf{D},\sf{E},\sf{F})
$$
to $\operatorname{Sq}(\sf{A},\sf{B},\sf{E},\sf{F})$ a {\it horizontal} composition from a subset of
$$
\operatorname{Sq}(\sf{A},\sf{B},\sf{C},\sf{D}) \times \operatorname{Sq}(\sf{B},\sf{E},\sf{D},\sf{F})
$$
to $ \operatorname{Sq}(\sf{A},\sf{E},\sf{C},\sf{F})$ and a {\it rotation} from a subset of $\operatorname{Sq}(\sf{A},\sf{B},\sf{C},\sf{D})$ to $\operatorname{Sq}(\sf{B},\sf{D},\sf{A},\sf{C})$.  We now describe these laws.
\begin{proposition} \label{prop.adcom3}
(Vertical composition of squares)\index{vertical composition of squares|textbf}.  Suppose
$$
\xymatrix{
\sf{A} \ar@<1ex>[r]^{G''} \ar[d]_{H} \ar@{=>}[dr]|{\Phi^{t}}\hole & \sf{B} \ar@<1ex>[l]^{F''} \ar[d]^{I}\\
\sf{C} \ar@<1ex>[r]^{G'}  & \sf{D} \ar@<1ex>[l]^{F'}.
}
$$
is an element of $\operatorname{Sq}(\sf{A},\sf{B},\sf{C},\sf{D})$ while
$$
\xymatrix{
\sf{C} \ar@<1ex>[r]^{G'} \ar[d]_{J} \ar@{=>}[dr]|{\Phi^{b}}\hole & \sf{D} \ar@<1ex>[l]^{F'} \ar[d]^{K}\\
\sf{E} \ar@<1ex>[r]^{G}  & \sf{F} \ar@<1ex>[l]^{F}.
}
$$
is an element of $\operatorname{Sq}(\sf{C},\sf{D},\sf{E},\sf{F})$.  If we define
$$
\Phi^{o} = (\Phi^{b}*H) \circ (K * \Phi^{t})
$$
then
$$
\xymatrix{
\sf{A} \ar@<1ex>[r]^{G''} \ar[d]_{JH} \ar@{=>}[dr]|{\Phi^{o}}\hole & \sf{B} \ar@<1ex>[l]^{F''} \ar[d]^{KI} \\
\sf{E} \ar@<1ex>[r]^{G}  & \sf{F} \ar@<1ex>[l]^{F}.
}
$$
is an element of $\operatorname{Sq}(\sf{A},\sf{B},\sf{E},\sf{F})$ which has
$$
\Psi^{o}=(J*\Psi^{t})\circ (\Psi^{b}*I).
$$
\end{proposition}

\begin{proof}
Our strategy is to show three diagrams commute, then arrange them left to right and note that the top circuit of this large diagram equals
$$
(J*\Psi^{t})\circ (\Psi^{b}*I)
$$
while the bottom circuit equals $\Psi^{o}$.  First, the diagram
$$
\xymatrix{
FKI \ar@{=>}[rr]^{FK*\eta'*I} \ar@{=>}[d]_{FKI*\eta''} & & FKG'F'I \ar@{=>}[d]^{FKG'*\Psi^{t}} \\
FKIG''F'' \ar@{=>}[rr]_{FK*\Phi^{t} *F''} & & FKG'HF'' \\
}
$$
commutes by Proposition \ref{prop.adcom}.  Next, the diagram
$$
\xymatrix{
FKG'F'I \ar@{=>}[rr]^{F*\Phi^{b}*F'I} \ar@{=>}[d]_{FKG'*\Psi^{t}} & & FGJF'I \ar@{=>}[dd]^{FGJ*\Psi^{t}} \\
FKG'HF'' \ar@{=>}[drr]^{\hskip .35in F*\Phi^{b}*HF''} & & \\
FKIG''F'' \ar@{=>}[u]^{FK*\Phi^{t}*F''} \ar@{=>}[rr]_{F*\Phi^{o}*F''} & & FGJHF''
}
$$
commutes, as follows:  The top commutes by 2ef of Definition \ref{def.2cat}, while the bottom commutes by hypothesis.  Finally, the diagram
$$
\xymatrix{
FGJF'I \ar@{=>}[d]_{FGJ*\Psi^{t}} \ar@{=>}[rr]^{\epsilon*JF'I}& & JF'I \ar@{=>}[d]^{J*\Psi^{t}}  \\
FGJHF'' \ar@{=>}[rr]_{\epsilon*JHF''} & & JHF''
}
$$
commutes, by naturality of $\epsilon$.  Combining these three diagrams left to right, we note that
$$
\Psi^{b}*I = (\epsilon * JF'I) \circ (F * \Phi^{b} *F'I) \circ (FK*\eta'*I)
$$
is just the top row of morphisms of this large diagram, while
$$
\Psi^{o} = (\epsilon * JHF'') \circ (F * \Phi^{o} * F'') \circ (FKI*\eta''),
$$
is the bottom route.  The assertion follows.
\end{proof}
The proof of the following result is similar, so we omit it.
\begin{proposition} \label{prop.horizontal}
(Horizontal composition of squares)\index{horizontal composition of squares|textbf}.    Suppose
$$
\xymatrix{
\sf{A} \ar@<1ex>[r]^{G'} \ar[d]_{H} \ar@{=>}[dr]|{\Phi^{l}}\hole & \sf{B} \ar@<1ex>[l]^{F'} \ar[d]^{I} & \mbox{and}  &  \sf{B} \ar@<1ex>[r]^{G'''} \ar[d]_{I} \ar@{=>}[dr]|{\Phi^{r}}\hole & \sf{E} \ar@<1ex>[l]^{F'''} \ar[d]^{J} \\
\sf{C} \ar@<1ex>[r]^{G}  & \sf{D} \ar@<1ex>[l]^{F} & &  \sf{D} \ar@<1ex>[r]^{G''}  & \sf{F} \ar@<1ex>[l]^{F''}.
}
$$
are elements of $\operatorname{Sq}(\sf{A},\sf{B},\sf{C},\sf{D})$ and $\operatorname{Sq}(\sf{B},\sf{E},\sf{D},\sf{F})$ respectively.  If we let
$$
\Phi^{o} = (G'''*\Phi^{l}) \circ (\Phi^{r} * G')
$$
then
$$
\xymatrix{
\sf{B} \ar@<3ex>[r]^{G'''G'} \ar[d]_{H} \ar@{=>}[dr]|{\Phi^{o}}\hole & \sf{E} \ar@<-1ex>[l]^{F'F'''} \ar[d]^{J} \\
\sf{D} \ar@<-1ex>[r]^{G''G}  & \sf{F} \ar@<3ex>[l]^{FF''}
}
$$
is an element of $\operatorname{Sq}(\sf{B},\sf{E},\sf{D},\sf{F})$ such that
$$
\Psi^{o} = (\Psi^{l}*F''')\circ (F * \Psi^{r}).
$$
\end{proposition}
{\it Remark.}  We have assumed that the unit and counit of $(F'F''',G'''G')$ are the canonical unit and counit one gets by composing adjoint pairs \cite[Theorem 1, p.101]{cat}.

\begin{proposition} \label{prop.rotation}
(Rotation)\index{rotation of squares|textbf}.  Let
$$
\xymatrix{
\sf{A} \ar@<1ex>[r]^{G'} \ar[d]_{H} \ar@{=>}[dr]|{\Phi}\hole & \sf{B} \ar@<1ex>[l]^{F'} \ar[d]^{I}\\
\sf{C} \ar@<1ex>[r]^{G}  & \sf{D} \ar@<1ex>[l]^{F}
}
$$
be an object in $\operatorname{Sq}(\sf{A},\sf{B},\sf{C},\sf{D})$ such that $I$ has left adjoint $\overline{I}$ and unit/counit $(\eta_{I},\epsilon_{I})$ while $H$ has left adjoint $\overline{H}$ and unit/counit $(\eta_{H},\epsilon_{H})$.  If $\Phi^{R} = \Psi$, then
$$
\xymatrix{
\sf{B} \ar@<1ex>[r]^{I} \ar[d]_{F'} \ar@{=>}[dr]|{\Phi^{R}}\hole & \sf{D} \ar@<1ex>[l]^{\overline{I}} \ar[d]^{F}\\
\sf{A} \ar@<1ex>[r]^{H}  & \sf{C} \ar@<1ex>[l]^{\overline{H}}
}
$$
is an object of $\operatorname{Sq}(\sf{B},\sf{D},\sf{A}, \sf{C})$ such that

\begin{equation} \label{eqn.rotate}
\xymatrix{
\overline{H}FIG'F'\overline{I} \ar@{=>}[rrr]^{\overline{H}F*\Phi*F'\overline{I}} & & & \overline{H}FGHF'\overline{I} \ar@{=>}[d]^{(\epsilon_{H} \circ (\overline{H}*\epsilon*H))*F'\overline{I}} \\
\overline{H}F \ar@{=>}[rrr]_{\Psi^{R}} \ar@{=>}[u]^{\overline{H}F*((I*\eta'*\overline{I}) \circ \eta_{I})} & & & F'\overline{I}
}
\end{equation}
commutes.
\end{proposition}

\begin{proof}
By Proposition \ref{prop.adcom} applied to the square
\begin{equation} \label{eqn.first}
\xymatrix{
\sf{A} \ar@<1ex>[r]^{G'} \ar[d]_{H} \ar@{=>}[dr]|{\Phi}\hole & \sf{B} \ar@<1ex>[l]^{F'} \ar[d]^{I}\\
\sf{C} \ar@<1ex>[r]^{G}  & \sf{D}  \ar@<1ex>[l]^{F}
}
\end{equation}
we know the diagram
\begin{equation} \label{eqn.big}
\xymatrix{
\overline{H}F \ar@{=>}[r]^{\overline{H}F*\eta_{I}} & \overline{H}FI\overline{I} \ar@{=>}[rr]^{\overline{H}F*\eta*I\overline{I}} \ar@{=>}[d]_{\overline{H}FI*\eta'*\overline{I}} & & \overline{H}FGFI\overline{I} \ar@{=>}[d]^{\overline{H}FG*\Psi*\overline{I}} & & \\
& \overline{H}FIG'F'\overline{I} \ar@{=>}[rr]_{\overline{H}F*\Phi*F'\overline{I}} & & \overline{H}FGHF'\overline{I} \ar@{=>}[d]^{\overline{H}*\epsilon*HF'\overline{I}} & & \\
& & &\overline{H}HF'\overline{I} \ar@{=>}[rr]_{\epsilon_{H}*F'\overline{I}} & & F' \overline{I} \
}
\end{equation}
commutes.  Since the bottom route of this diagram is exactly the top route of (\ref{eqn.rotate}), it suffices to show that the top route of (\ref{eqn.big}) equals $\Psi^{R}$.  If we can show that the second, third and fourth maps of the top route of (\ref{eqn.big}) equals $\overline{H}*\Psi*\overline{I}$, then the assertion will follow from the definition of $\Psi^{R}$.  To show the second, third and fourth map of the top route of (\ref{eqn.big}) equals $\overline{H}*\Psi*\overline{I}$, it suffices to show that the composition
$$
\xymatrix{
& FI \ar@{=>}[d]^{F*\eta*I} \\
& FGFI \ar@{=>}[d]^{FG*\Psi} \\
HF' & FGHF' \ar@{=>}[l]^{\epsilon*HF'}
}
$$
equals $\Psi$.  This follows from the commutivity of
$$
\xymatrix{
& FI \ar@{=>}[d]^{F*\eta*I} \\
FI \ar@{=>}[d]_{\Psi} & FGFI \ar@{=>}[d]^{FG*\Psi} \ar@{=>}[l]_{\epsilon*FI}\\
HF' & FGHF' \ar@{=>}[l]^{\epsilon*HF'}
}
$$
\end{proof}\index{category of squares|)}\index{canonical 2-cell associated to a square|)}\index{2-cell dual to phi@2-cell dual to $\Phi$|)}

\section{Indexed categories}
We define the terms {\it indexed category}, {\it indexed functor},
and {\it indexed natural transformation}.  We give relevant
examples of the above terms.  In the next section we will use some
of the consequences of Proposition \ref{prop.adcom} to give
criteria as to when a family of functors or natural
transformations is indexed (Propositions \ref{prop.indexedfunct}
and \ref{cor.reduce}).

\begin{definition} \cite[p.63]{Coherence} \label{def.index}\index{indexed!category|textbf}
Let $\sf{S}$ be a category with finite limits.  An {\bf{ $\sf{S}$-indexed category, $\mathfrak{A}$}}, consists of the following data:
\begin{enumerate}
\item{for each $W \in \sf{S}$, a category ${\sf{A}}^{W}$,}
\item{for each $f:W \rightarrow V$ in $\sf{S}$ a functor $f^{*}:{\sf{A}}^{V} \rightarrow {\sf{A}}^{W}$, called an {\bf{induced functor}}}
\item{for each composable pair
$$
\xymatrix{
W \ar[r]^{g} & V \ar[r]^{f} & U}
$$
in $\sf{S}$, a natural isomorphism $\Xi_{f,g}:g^{*}f^{*} \Longrightarrow (fg)^{*}$,}
\item{for each $W$ in $\sf{S}$, a natural isomorphism $\Upsilon_{W}:(1_{W})^{*} \Longrightarrow 1_{{\sf{A}}^{W}}$,}
\end{enumerate}
subject to the following compatibilities:
\newenvironment{labelitemize}{%
    \renewcommand{\labelitemi}{(C1)}%
    \begin{itemize}}{\end{itemize}}
\begin{labelitemize}
\item{for each composable triple $\xymatrix{X \ar[r]^{h} & W \ar[r]^{g} & V \ar[r]^{f} & U}$ in $\sf{S}$, the diagram
$$
\xymatrix{
h^{*}g^{*}f^{*} \ar@{=>}[r]^{h^{*}\Xi_{f,g}} \ar@{=>}[d]_{\Xi_{g,h}f^{*}} & h^{*}(fg)^{*} \ar@{=>}[d]^{\Xi_{fg,h}} \\
(gh)^{*}f^{*} \ar@{=>}[r]_{\Xi_{f,gh}} & (fgh)^{*}
}
$$
commutes, and}
\end{labelitemize}
\newenvironment{labelitemize2}{%
    \renewcommand{\labelitemi}{(C2)}%
    \begin{itemize}}{\end{itemize}}
\begin{labelitemize2}
\item{for each $f:V \rightarrow U$ in $\sf{S}$,
$$
\xymatrix{
\Xi_{1_{U},f}=f^{*}\Upsilon_{U}.
}
$$
}
\end{labelitemize2}
\end{definition}
Note that we have not specified a category structure on $\mathfrak{A}$, so that an indexed category is {\it not} a category.

\begin{lemma}  \cite[(3.6.2), p.104]{derived} \label{lem.adjointnew}
Suppose $X$, $Y$ and $Z$ are ringed spaces with maps between them
$$
\xymatrix{
X \ar[r]^{g} & Y \ar[r]^{f} & Z.
}
$$
Then there is a unique isomorphism $g^{*}f^{*} \Longrightarrow (fg)^{*}$ making the diagrams
$$
\xymatrix{
(fg)^{*}(fg)_{*} \ar@{=>}[d] & g^{*}f^{*}(fg)_{*} \ar@{=>}[l]_{\cong} \ar@{=>}[r]^{\cong} & g^{*}f^{*}f_{*}g_{*} \ar@{=>}[d]^{=} \\
\operatorname{id}_{{\sf{Qcoh}}X} & g^{*}g_{*} \ar@{=>}[l] & g^{*}f^{*}f_{*}g_{*} \ar@{=>}[l]
}
$$
and
$$
\xymatrix{
(fg)_{*}(fg)^{*}  & (fg)_{*}g^{*}f^{*} \ar@{=>}[l]_{\cong} \ar@{=>}[r]^{\cong} & f_{*}g_{*}g^{*}f^{*} \ar@{=>}[d]^{=} \\
\operatorname{id}_{{\sf{Qcoh}}Z} \ar@{=>}[r] \ar@{=>}[u] & f_{*}f^{*}  \ar@{=>}[r] & f_{*}g_{*}g^{*}f^{*}
}
$$
commute.
\end{lemma}

\begin{proof}
Since $(fg)_{*} \cong f_{*}g_{*}$, the unique isomorphism $g^{*}f^{*} \Longrightarrow (fg)^{*}$ is due to the uniqueness of the left adjoint of a given right adjoint.
\end{proof}
The following example will be used in the next chapter to define
$\Gamma_{n}$.
\begin{example} \label{example.spaces}
Suppose $I=\{1,\ldots, n\}$, ${\{X_{i} \}}_{i \in I}$ are schemes and $U$ is an affine noetherian scheme.  We remind the reader that $\times_{i}(U \times X_{i})_{U}$ is the product of the $n$ spaces $U \times_{S} X_{i}$ over $U$.  Let $\sf{S}$ denote the category of affine noetherian schemes, and let ${\mathfrak{\times_{i} X_{i}}}$ denote the indexed category defined by the following data: for each $U \in \sf{S}$, let
$$
{\sf{\times_{i} X_{i}}}^{U}={\sf{Qcoh}}\times_{i}(U \times_{S} X_{i})_{U}.
$$
For each $f:V \rightarrow U$ in $\sf{S}$, let
$$
\tilde{f} = \times_{i}(f \times 1_{X_{i}}):\times_{i}(V \times X_{i})_{V} \rightarrow \times_{i}(U \times X_{i}),
$$
and let
$$
f^{*} = {\tilde{f}}^{*}:({\sf{\times_{i}X_{i}}})^{U} \rightarrow ({\sf{\times_{i} X_{i}}})^{V}.
$$
For each composable pair
$$
\xymatrix{
W \ar[r]^{g} & V \ar[r]^{f} & U
}
$$
in $\sf{S}$ and each object $\mathcal{M}$ in $({\sf{\times_{i}X_{i}}})^{U}$,
let
$$
\Xi_{f,g,\mathcal{M}}:g^{*}f^{*}\mathcal{M} \rightarrow (fg)^{*}\mathcal{M}
$$
be the isomorphism constructed in Lemma \ref{lem.adjointnew}.  Finally, for each $U$ in $\sf{S}$ and each object $\mathcal{M}$ in $({\sf{\times_{i}X_{i}}})^{U}$, let
$$
\Upsilon_{U,\mathcal{N}}:(1_{U})^{*}\mathcal{M} \rightarrow 1_{{\sf{\times_{i}X_{i}}}^{U}}\mathcal{M}
$$
be the canonical isomorphism of ${\mathcal{O}}_{\times_{i}(U \times X_{i})_{U}}$-modules
$$
{\mathcal{O}}_{\times_{i}(U\times X_{i})_{U}}\otimes_{{\mathcal{O}}_{\times_{i}(U\times X_{i})_{U}}} {\mathcal{M}} \rightarrow  {\mathcal{M}}.
$$
To verify that $\mathfrak{\times_{i}X_{i}}$ is an $\sf{S}$-indexed category, we must verify that our morphisms $\Xi$ and $\Upsilon$ satisfy conditions $(C1)$ and $(C2)$.  These conditions may be checked locally, and $\mathfrak{\times_{i}X_{i}}$ is indeed an $\sf{S}$-indexed category.
\end{example}

\begin{definition}\cite[p.63-64]{Coherence} \label{def.indexfunct}
Suppose $\mathfrak{A}$ and $\mathfrak{B}$ are indexed categories.  A {\bf weakly indexed functor} \index{weakly indexed functor|textbf}$\mathfrak{F}:\mathfrak{A} \rightarrow \mathfrak{B}$ consists of the following data:
\begin{enumerate}
\item{for each $W \in \sf{S}$, a functor ${F}^{W}$,}
\item{for each $f:W \rightarrow V$ in $\sf{S}$ a natural transformation
$$
\Lambda_{f}:f^{*}F^{V} \Longrightarrow F^{W}f^{*}
$$
subject to the condition}
\end{enumerate}
\newenvironment{labelitemize}{%
    \renewcommand{\labelitemi}{(F)}%
    \begin{itemize}}{\end{itemize}}
\begin{labelitemize}
\item{
for each composable pair
$$
\xymatrix{
W \ar[r]^{g} & V \ar[r]^{f} & U}
$$
in $\sf{S}$, the diagram
$$
\xymatrix{
g^{*}f^{*}F^{U} \ar@{=>}[r]^{\Xi_{f,g}F^{U}} \ar@{=>}[d]_{g^{*}\Lambda_{f}} & (fg)^{*}F^{U} \ar@{=>}[dd]^{\Lambda_{fg}} \\
g^{*}F^{V}f^{*} \ar@{=>}[d]_{\Lambda_{g}f^{*}} & \\
F^{W}g^{*}f^{*} \ar@{=>}[r]_{F^{W}\Xi_{f,g}} & F^{W}(fg)^{*}
}
$$
commutes.  If $\Lambda_{f}$ is an isomorphism for all $f$ in $\sf{S}$, then we say $\mathfrak{F}$ is {\bf indexed}\index{indexed!functor|textbf}.
}
\end{labelitemize}
\end{definition}

\begin{example} \label{example.functors}
We retain the notation from Example \ref{example.spaces}.  Suppose $X$ and $Y$ are noetherian schemes.  We define an ${\sf{S}}$-indexed functor, ${\mathfrak{F}}:{\mathfrak{X}} \rightarrow {\mathfrak{Y}}$.  Let $\mathcal{F}$ be an ${\mathcal{O}}_{X}-{\mathcal{O}}_{Y}$-bimodule (\ref{def.bimodule}).  For $U, V \in {\sf{S}}$, suppose $p:(U \times X) \times_{U} (U \times Y)\rightarrow X \times Y$ is the projection map and, for each $U \in \sf{S}$, let ${\pr_{i}}^{U}:(U \times X)\times_{U}(U \times Y) \rightarrow U \times X, U \times Y$ be the standard projection map.  If
$$
F^{U} = {{\pr_{2}}^{U}}_{*}({{\pr_{1}}^{U}}^{*}- \otimes {p}^{*}\mathcal{F}):{\sf{Qcoh}} U \times X \rightarrow {\sf{Qcoh}} U \times Y,
$$
then we check that $\mathfrak{F}$ satisfies Definition $\ref{def.indexfunct}$, 2.  Suppose $f:V \rightarrow U$ is a morphism in $\sf{S}$, $f_{X}=f \times \id_{X}$, $f_{Y}=f \times \id_{Y}$, $f_{XY}=f_{X} \times f_{Y}$ and suppose $\mathcal{M}$ is an ${\mathcal{O}}_{U \times X}$-module.  Since
\begin{equation} \label{eqn.indexfunct}
\xymatrix{
(V \times X) \times_{V}(V \times Y) \ar[r]^{f_{XY}} \ar[d]_{{\pr_{2}}^{V}} & (U \times X)\times_{U}(U \times Y) \ar[d]^{{\pr_{2}^{U}}} \\
V \times Y \ar[r]_{{f_{Y}}} & U \times Y
}
\end{equation}
is a pull-back diagram, we will later show (Proposition \ref{prop.canon}) the dual 2-cell of the square induced by (\ref{eqn.indexfunct}), applied to the functor ${{\pr_{1}}^{U}}^{*}- \otimes p^{*}\mathcal{F}$,
$$
\Psi:  {{\pr_{2}}^{V}}_{*}f_{XY}^{*}({{\pr_{1}}^{U}}^{*}-\otimes {p}^{*}\mathcal{F}) \Longrightarrow
$$
$$
{f_{Y}}^{*}{{\pr_{2}}^{U}}_{*}({{\pr_{1}}^{U}}^{*}- \otimes {p}^{*}\mathcal{F}) = {{f_{Y}}}^{*}F^{U},
$$
is an isomorphism.  On the other hand, there is a natural isomorphism
$$
 {{\pr_{2}}^{V}}_{*}f_{XY}^{*} ({{\pr_{1}}^{U}}^{*}- \otimes {p}^{*}\mathcal{F}) \Longrightarrow F^{V} {{f_{X}}}^{*}.
$$
Composing these maps, we have a natural isomorphism
$$
\Lambda_{f}:{f_{Y}}^{*}F^{U} \Longrightarrow F^{V}{f_{X}}^{*}
$$
as desired.  We will show in Lemma \ref{lem.isommonoid} that these maps satisfy (F).\end{example}

\begin{definition} \label{def.intrans} \cite[p.64]{Coherence}
Suppose $\mathfrak{F}, \mathfrak{F}':\mathfrak{A} \rightarrow \mathfrak{B}$ are two weakly indexed functors.  An {\bf indexed natural transformation $\Delta:\mathfrak{F} \Longrightarrow \mathfrak{F}'$}\index{indexed!natural transformation|textbf} consists of a natural transformation $\Delta^{U}:F^{U} \Longrightarrow {F'}^{U}$ for each $U \in \sf{S}$ such that for each $f:V \rightarrow U$ in $\sf{S}$,
\begin{equation} \label{eqn.lippy}
\xymatrix{
f^{*}F^{U} \ar@{=>}[d]_{\Lambda_{f}} \ar@{=>}[r]^{f^{*}\Delta^{U}} & f^{*}{F'}^{U} \ar@{=>}[d]^{\Lambda_{f}} \\
F^{V}f^{*} \ar@{=>}[r]_{\Delta^{V}f^{*}} & {F'}^{V}f^{*}
}
\end{equation}
commutes.
\end{definition}
We mention a convention we will employ.  Suppose we are given, for each $U \in \sf{S}$, a category ${\sf{A}}^{U}$, and for each morphism $f:V \rightarrow U$ in $\sf{S}$ we are given a functor $f^{*}:{\sf{A}}^{U} \rightarrow {\sf{A}}^{V}$.  If we show that this assignment satisfies the axioms of an indexed category, then we will write $\mathfrak{A}$ for the resulting indexed category.

$\sf{S}$-indexed categories, functors, and natural transformations form a 2-category in an obvious way \cite[p.45]{2cat}.  For example, if $\mathfrak{A}$, $\mathfrak{B}$ and $\mathfrak{C}$ are $\sf{S}$-indexed categories, $f:V \rightarrow U$ is a morphism in $\sf{S}$, $\mathfrak{F}:\mathfrak{A} \rightarrow \mathfrak{B}$ is a weakly indexed functor with indexed structure $\Lambda^{F}_{f}:f^{*}F^{U} \Longrightarrow F^{V}f^{*}$ and $\mathfrak{G}:\mathfrak{B} \rightarrow \mathfrak{C}$ is a weakly indexed functor with indexed structure $\Lambda^{G}_{f}:f^{*}G^{U} \Longrightarrow G^{V}f^{*}$, then the indexed structure of $\mathfrak{GF}$ is given by the composition
\begin{equation} \label{eqn.indexcomp}
\xymatrix{
f^{*}G^{U}F^{U} \ar@{=>}[r]^{\Lambda^{G}_{f}*F^{U}} & G^{V}f^{*}F^{U} \ar@{=>}[r]^{G^{V}*\Lambda^{F}_{f}} & G^{V}F^{V}f^{*}.
}
\end{equation}
We make use of this fact routinely.

\begin{example} \label{example.indexcheck}
Retain the notation of Example \ref{example.spaces}.  Suppose $J \subset I=\{1,. . . ,n \}$ and, for each $U \in \sf{S}$, $\pr_{J}^{U}:\times_{i \in I}(U \times X_{i})_{U} \rightarrow \times_{i \in J}(U \times X_{i})_{U}$ denotes projection.  Then the $\sf{S}$-indexed set $\{\pr_{J}^{U*}\}$ induces, given the obvious $\Lambda$ in Definition \ref{def.indexfunct}, an indexed functor $\pr_{J}^{*}:\times_{i \in J}\mathfrak{{X_{i}}} \rightarrow \times_{i \in I}{\mathfrak{X_{i}}}$.  It is easy to prove this by localizing, since the localization of a sheaf under an inverse image functor is easy to describe in terms of the localization of the sheaf.  Similarly, if $\mathcal{P}$ is in ${\sf{Qcoh}}\times_{i}X_{i}$, and if $p_{i}^{U}:U \times X_{i} \rightarrow X_{i}$ is projection, then $-\otimes_{{\mathcal{O}}_{\times_{i}(U\times X_{i})_{U}}}(\times_{i}p_{i}^{U})^{*}\mathcal{P}$ induces an $\sf{S}$-indexed functor, as above.  Thus, compositions of these functors are indexed.  Similarly, if there is some family of natural transformations (indexed by $\sf{S}$) between two such indexed functors, then one can check whether or not such a family of transformations is indexed by localizing.  We will routinely suppress these sorts of computations.
\end{example}

Although we will not make use of this fact, we mention that indexed categories have a coherence theorem (\cite[Theorem, p.61]{Coherence}) for diagrams involving $\Xi$'s (Definition \ref{def.index}, (3)), $\Upsilon$'s (Definition \ref{def.index}, (4)) and $\Theta$'s (Definition \ref{def.indexfunct}, (2)).

\section{Squares of indexed categories}

In the next chapter, we will have many examples of families of functors whose ostensible indexing structure, $\Lambda$ (Definition \ref{def.indexfunct} (2)), happens to be the dual 2-cell of some square of indexed categories (Definition \ref{def.squareof}).  We show that the family of canonical 2-cells associated to a square of indexed categories is indexed if and only if the family of dual 2-cells is indexed (Corollary \ref{cor.reduce}).  We also give sufficient conditions for a family of functors between indexed categories to be weakly indexed (Proposition \ref{prop.indexedfunct}).

\begin{definition}
An indexed category is said to be {\bf adjointed} \index{indexed!adjointed category|textbf}if for each $f:V \rightarrow U$ in $\sf{S}$, there is an adjoint pair $(f^{*},f_{*})$ with unit and counit $(\eta_{f},\epsilon_{f})$ where $f^{*}$ is the functor induced by $f$, such that there is a natural isomorphism $\Omega_{f,g}:(fg)_{*} \Longrightarrow f_{*}g_{*}$ which induces $\Xi_{f,g}:g^{*}f^{*} \Longrightarrow (fg)^{*}$ and such that if
$$
\xymatrix{
W \ar[r]^{g} & V \ar[r]^{f} & U
}
$$
are morphisms in $\sf{S}$, then the pair associated to $fg$ is $((fg)^{*}(fg)_{*})$ with unit and counit
\begin{equation} \label{eqn.unitcounit}
((\Omega_{f,g}^{-1}*\Xi_{f,g}) \circ (f_{*}*\eta_{g}*f^{*}) \circ \eta_{f}, \epsilon_{g} \circ (g^{*}*\epsilon_{f}*g_{*}) \circ \Xi^{-1}_{f,g}*\Omega_{f,g}).
\end{equation}
\end{definition}

\begin{example} \label{example.schemes}
By Lemma \ref{lem.adjointnew}, the indexed category ${\mathfrak{X}}^{n}$ is adjointed.  More generally, if ${\sf{S}}$ is the category of schemes,  it is easy to show that the assignment $W \in {\sf{S}} \mapsto {\sf{Qcoh}}W$ makes ${\sf{S}}$ an index for an indexed category, $\mathfrak{S}$ which, in light of Lemma \ref{lem.adjointnew}, is naturally adjointed.
\end{example}

\begin{proposition} \label{prop.indexedfunct}
Let $\mathfrak{A}$ and $\mathfrak{B}$ be $\sf{S}$-indexed, adjointed categories.  Suppose, for each $U \in \sf{S}$, there is a functor $H^{U}:{\sf{A}}^{U} \rightarrow {\sf{B}}^{U}$.  In addition, suppose for each morphism $f:V \rightarrow U$ in $\sf{S}$, there is a natural transformation
$$
\Phi_{f}:H^{U}f_{*} \Longrightarrow f_{*}H^{V}
$$
such that if
$$
\xymatrix{
W \ar[r]^{g} & V \ar[r]^{f} & U
}
$$
is a diagram of objects and morphisms in $\sf{S}$, then
$$
\xymatrix{
H^{U}(fg)_{*} \ar@{=>}[dd]_{H^{U}*\Omega} \ar@{=>}[r]^{\Phi_{fg}} & (fg)_{*}H^{W} \ar@{=>}[d]^{\Omega *H^{W}}  \\
& f_{*}g_{*}H^{W} \\
H^{U}f_{*}g_{*} \ar@{=>}[r]_{\Phi_{f}*g_{*}} & f_{*}H^{V}g_{*} \ar@{=>}[u]_{f_{*}*\Phi_{g}}
}
$$
commutes.  Let
$$
\Psi_{f}:f^{*}H^{U} \Longrightarrow H^{V}f^{*}
$$
be the dual 2-cell of the square
$$
\xymatrix{
{\sf{A}}^{V} \ar@<1ex>[r]^{f_{*}} \ar[d]_{H^{V}} \ar@{=>}[dr]|{\Phi_{f}}\hole   & {\sf{A}}^{U} \ar@<1ex>[l]^{f^{*}} \ar[d]^{H^{U}} \\
{\sf{B}}^{V} \ar@<1ex>[r]^{f_{*}}  & {\sf{B}}^{U} \ar@<1ex>[l]^{f^{*}}.
}
$$
Then $\Psi_{f}$ gives the set $\{ H^{U} \}$ the structure of a weakly indexed functor\index{weakly indexed functor}.
\end{proposition}

\begin{proof}
Let
$$
\xymatrix{
W \ar[r]^{g} & V \ar[r]^{f} & U
}
$$
be a diagram of objects and morphisms in $\sf{S}$ and consider the horizontal diagram
$$
\xymatrix{
{\sf{A}}^{W} \ar@<1ex>[r]^{g_{*}} \ar[d]_{H^{W}} \ar@{=>}[dr]|{\Phi_{g}}\hole & {\sf{A}}^{V} \ar@<1ex>[r]^{f_{*}} \ar@<1ex>[l]^{g^{*}} \ar[d]^{H^{V}}  \ar@{=>}[dr]|{\Phi_{f}}\hole & {\sf{A}}^{U} \ar@<1ex>[l]^{f^{*}} \ar[d]^{H^{U}} \\
{\sf{B}}^{W} \ar@<1ex>[r]^{g_{*}}  & {\sf{B}}^{V} \ar@<1ex>[r]^{f_{*}} \ar@<1ex>[l]^{g^{*}}  & {\sf{B}}^{U} \ar@<1ex>[l]^{f^{*}}
}
$$
The horizontal composition of these squares is the square
$$
\xymatrix{
{\sf{A}}^{W} \ar@<3ex>[r]^{f_{*}g_{*}} \ar[d]_{H^{W}} \ar@{=>}[dr]|{\Phi}\hole & {\sf{A}}^{U} \ar@<-1ex>[l]^{g^{*}f^{*}} \ar[d]^{H^{U}} \\
{\sf{B}}^{W} \ar@<-1ex>[r]^{f_{*}g_{*}} & {\sf{B}}^{U} \ar@<3ex>[l]^{g^{*}f^{*}}}
$$
with unit/counits $(\overline{\eta},\overline{\epsilon})$, where
$$
\Phi=(f_{*}* \Phi_{g}) \circ (\Phi_{f}*g_{*}).
$$
Suppose $\Psi$ is dual to $\Phi$.  By Proposition \ref{prop.horizontal}, the diagram
$$
\xymatrix{
g^{*}f^{*}H^{U} \ar@{=>}[rr]^{\Psi} \ar@{=>}[dr]_{g^{*} * \Psi_{f}} & & H^{W}g^{*}f^{*} \\
& g^{*}H^{V}f^{*} \ar@{=>}[ur]_{\Psi_{g} * f^{*}} &
}
$$
commutes.  Thus, to complete the proof, we must show that the diagram
$$
\xymatrix{
g^{*}f^{*}H^{U} \ar@{=>}[rr]^{\Psi} \ar@{=>}[d]_{\Xi_{f,g}*H^{U}} & & H^{W}g^{*}f^{*} \ar@{=>}[d]^{H^{W}\Xi_{f,g}} \\
(fg)^{*}H^{U} \ar@{=>}[rr]_{\Psi_{fg}} & & H^{W}(fg)^{*}
}
$$
commutes.  By Lemma \ref{lem.simplecomp2}, in order to show this diagram commutes, it suffices to show that the data $J^{{\sf{A}}^{W},{\sf{A}}^{U},{\sf{B}}^{W},{\sf{B}}^{U}}=\operatorname{Id}:{{\sf{A}}^{W},{\sf{A}}^{U},{\sf{B}}^{W},{\sf{B}}^{U}} \rightarrow {{\sf{A}}^{W},{\sf{A}}^{U},{\sf{B}}^{W},{\sf{B}}^{U}}$, $\Theta_{(fg)^{*}}=\Xi_{f,g}^{-1}$, $\Theta_{f_{*}g_{*}}=\Omega_{f,g}^{-1}$, $\Theta_{H^{W},H^{U}}$ equal identity transformations, defines a morphism $(J,\Theta)$
$$
\xymatrix{
{\sf{A}}^{W} \ar@<3ex>[r]^{f_{*}g_{*}} \ar[d]_{H^{W}} \ar@{=>}[dr]|{\Phi}\hole & {\sf{A}}^{U} \ar@<-1ex>[l]^{g^{*}f^{*}} \ar[d]^{H^{U}} & \ar[r] & &  {\sf{A}}^{W} \ar@<3ex>[r]^{(fg)_{*}} \ar[d]_{H^{W}} \ar@{=>}[dr]|{\Phi_{fg}}\hole & {\sf{A}}^{U} \ar@<-1ex>[l]^{(fg)^{*}} \ar[d]^{H^{U}} \\
{\sf{B}}^{W} \ar@<-1ex>[r]^{f_{*}g_{*}} & {\sf{B}}^{U} \ar@<3ex>[l]^{g^{*}f^{*}} & & & {\sf{B}}^{W} \ar@<-1ex>[r]^{(fg)_{*}} & {\sf{B}}^{U} \ar@<3ex>[l]^{(fg)^{*}}
}
$$
We proceed to show $(J,\Theta)$ is a morphism.  Since $\mathfrak{A}$ and $\mathfrak{B}$ are adjointed, the diagrams
$$
\xymatrix{
\id \ar@{=>}[rr]^{\eta} \ar@{=>}[d]_{\overline{\eta}} & & (fg)_{*}(fg)^{*} & & g^{*}f^{*}f_{*}g_{*} \ar@{=>}[rr]^{\Xi_{f,g}*\Omega_{f,g}^{-1}} \ar@{=>}[d]_{\overline{\epsilon}} & & (fg)^{*}(fg)_{*} \ar@{=>}[d]^{\epsilon} \\
f_{*}g_{*}g^{*}f^{*} \ar@{=>}[rr]_{\Omega_{f,g}^{-1}*g^{*}f^{*}} & &  (fg)_{*}g^{*}f^{*} \ar@{=>}[u]_{(fg)_{*}*\Xi_{f,g}} & & \id \ar@{=>}[rr]_{=} & & \id
}
$$
commute.  Thus, we need only show that the other data defining a morphism,
$$
\xymatrix{
H^{U}f_{*}g_{*} \ar@{=>}[r] \ar@{=>}[d] & H^{U}(fg)_{*} \ar@{=>}[dd] \\
f_{*}H^{V}g_{*} \ar@{=>}[d] & \\
f_{*}g_{*}H^{W} \ar@{=>}[r] & (fg)_{*}H^{W}
}
$$
commutes.  This is true by hypothesis.  Thus, the hypothesis of Lemma \ref{lem.simplecomp2} are satisfied by $(J,\Theta)$ so $(J,\Theta)$ is a morphism and the result follows.
\end{proof}

\begin{definition} \label{def.squareof}
Suppose $\mathfrak{A}$, $\mathfrak{B}$, $\mathfrak{C}$ and $\mathfrak{D}$ are $\sf{S}$-indexed categories, $\mathfrak{F}$, $\mathfrak{F}'$, $\mathfrak{H}$ and $\mathfrak{I}$ are $\sf{S}$-indexed, $\mathfrak{G}$ and $\mathfrak{G'}$ are weakly $\sf{S}$-indexed, and for each $U \in {\sf{S}}$,
$$
\xymatrix{
{\sf{A}}^{U} \ar@<1ex>[r]^{{G'}^{U}} \ar[d]_{H^{U}} \ar@{=>}[dr]|{\Phi^{U}}\hole & {\sf{B}}^{U} \ar@<1ex>[l]^{{F'}^{U}} \ar[d]^{I^{U}}\\
{\sf{C}}^{U} \ar@<1ex>[r]^{G^{U}}  & {\sf{D}}^{U} \ar@<1ex>[l]^{F^{U}}.
}
$$
is a square.  A {\bf square of indexed categories} \label{square of indexed categories|textbf}consists of the above data, denoted
$$
\xymatrix{
\mathfrak{A} \ar@<1ex>[r]^{\mathfrak{G'}} \ar[d]_{\mathfrak{H}} \ar@{=>}[dr]|{\mathfrak{\Phi}}\hole & \mathfrak{B} \ar@<1ex>[l]^{\mathfrak{F'}} \ar[d]^{\mathfrak{I}}\\
\mathfrak{C} \ar@<1ex>[r]^{\mathfrak{G}}  & \mathfrak{D} \ar@<1ex>[l]^{\mathfrak{F}}
}
$$
such that for each morphism $f:V \rightarrow U$ in $\sf{S}$, the diagrams
\begin{equation} \label{eqn.adjointcomp}
\xymatrix{
f^{*} \ar@{=>}[rr]^{{\eta'}^{V}*f^{*}} \ar@{=>}[d]_{f^{*} * {\eta'}^{U}} & & {G'}^{V}{F'}^{V}f^{*} \ar@{=>}[d]^{{G'}^{V}*{\Lambda_{f}^{F'}}^{-1}} & & F^{V}f^{*}G^{U} \ar@{=>}[rr]^{F^{V}*\Lambda_{f}^{G}} \ar@{=>}[d]_{\Lambda_{f}^{{F}^{-1}}* G^{U}} & &  F^{V}G^{V}f^{*} \ar@{=>}[d]^{{\epsilon}^{V}*f^{*}} \\
f^{*}{G'}^{U}{F'}^{U} \ar@{=>}[rr]_{\Lambda_{f}^{G'}*{F'}^{U}} & & {G'}^{V}f^{*}{F'}^{U} & & f^{*}F^{U}G^{U} \ar@{=>}[rr]_{f^{*}*\epsilon^{U}} & & f^{*}
}
\end{equation}
where $\Lambda$ gives the indexed structure of the various (weakly) indexed functors, commute.
\end{definition}
{\it Remark.}  If
$$
\xymatrix{
\mathfrak{A} \ar@<1ex>[r]^{\mathfrak{G'}} \ar[d]_{\mathfrak{H}} \ar@{=>}[dr]|{\mathfrak{\Phi}}\hole & \mathfrak{B} \ar@<1ex>[l]^{\mathfrak{F'}} \ar[d]^{\mathfrak{I}}\\
\mathfrak{C} \ar@<1ex>[r]^{\mathfrak{G}}  & \mathfrak{D} \ar@<1ex>[l]^{\mathfrak{F}}
}
$$
is a square of indexed functors, then $\epsilon$ and $\eta'$ are indexed natural transformations by (\ref{eqn.adjointcomp}).

\begin{proposition} \label{cor.reduce}
Suppose
\begin{equation} \label{eqn.megasquare}
\xymatrix{
\mathfrak{A} \ar@<1ex>[r]^{\mathfrak{G'}} \ar[d]_{\mathfrak{H}} \ar@{=>}[dr]|{\mathfrak{\Phi}}\hole & \mathfrak{B} \ar@<1ex>[l]^{\mathfrak{F'}} \ar[d]^{\mathfrak{I}}\\
\mathfrak{C} \ar@<1ex>[r]^{\mathfrak{G}}  & \mathfrak{D} \ar@<1ex>[l]^{\mathfrak{F}}
}
\end{equation}
is a square of indexed categories.  Then $\Phi:\mathfrak{IG'} \Longrightarrow \mathfrak{GH}$ is an indexed natural transformation if and only if $\Psi:\mathfrak{FI} \Longrightarrow \mathfrak{HF'}$ is an indexed natural transformation.\index{indexed!natural transformation}
\end{proposition}

\begin{proof}
Since (\ref{eqn.megasquare}) is a square of indexed categories, $\epsilon$ and $\eta'$ are indexed natural transformations.  Thus, by definition of $\Psi$, the fact that $\Phi$ is indexed implies $\Psi$ is indexed.

Conversely, since, for every $f:V \rightarrow U$ in $\sf{S}$, the diagrams (\ref{eqn.adjointcomp}) commute, then by Corollary \ref{cor.solu},
$$
\xymatrix{
{\sf{A}}^{U} \ar@<1ex>[r]^{{G'}^{U}} \ar[d]_{f^{*}} \ar@{=>}[dr]|{\Lambda^{G'}}\hole & {\sf{B}}^{U} \ar@<1ex>[l]^{{F'}^{U}} \ar[d]^{f^{*}} & \mbox{and} & {\sf{C}}^{U} \ar@<1ex>[r]^{{G}^{U}} \ar[d]_{f^{*}} \ar@{=>}[dr]|{\Lambda^{G}}\hole & {\sf{D}}^{U} \ar@<1ex>[l]^{{F}^{U}} \ar[d]^{f^{*}} \\
{\sf{A}}^{V} \ar@<1ex>[r]^{{G'}^{V}}  & {\sf{B}}^{V} \ar@<1ex>[l]^{{F'}^{V}} & & {\sf{C}}^{V} \ar@<1ex>[r]^{G^{V}}  & {\sf{D}}^{V} \ar@<1ex>[l]^{F^{V}}.
}
$$
are squares with duals $(\Lambda^{F'})^{-1}$ and $(\Lambda^{F})^{-1}$ respectively.  Thus, by Proposition \ref{prop.adcom}, the diagrams
$$
\xymatrix{
f^{*} \ar@{=>}[rr]^{{\eta}^{V}*f^{*}} \ar@{=>}[d]_{f^{*} * {\eta}^{U}} & & {G}^{V}{F}^{V}f^{*} \ar@{=>}[d]^{{G}^{V}*{\Lambda_{f}^{{F}}}^{-1}} & & {F'}^{V}f^{*}{G'}^{U} \ar@{=>}[rr]^{{F'}^{V}*\Lambda_{f}^{G'}} \ar@{=>}[d]_{{\Lambda_{f}^{{F'}}}^{-1} * {G'}^{U}} & &  {F'}^{V}{G'}^{V}f^{*} \ar@{=>}[d]^{{\epsilon'}^{V}*f^{*}} \\
f^{*}{G}^{U}{F}^{U} \ar@{=>}[rr]_{\Lambda_{f}^{G}*{F}^{U}} & & {G}^{V}f^{*}{F}^{U} & & f^{*}{F'}^{U}{G'}^{U} \ar@{=>}[rr]_{f^{*}*{\epsilon'}^{U}} & & f^{*}
}
$$
commute also.  Thus, $\eta$ and $\epsilon'$ are indexed so that by definition of $\Psi$, the fact that $\Psi$ is indexed implies $\Phi$ is indexed.
\end{proof}
The following easy result gives criteria for a family of squares to be a square of indexed categories.
\begin{lemma} \label{lem.applicat}
Suppose $\mathfrak{A}$, $\mathfrak{B}$, $\mathfrak{C}$ and $\mathfrak{D}$ are $\sf{S}$-indexed categories, $\mathfrak{F}$, $\mathfrak{F}'$, $\mathfrak{G}$, $\mathfrak{G'}$, $\mathfrak{H}$ and $\mathfrak{I}$ are $\sf{S}$-indexed functors, and, for each $U \in {\sf{S}}$,
\begin{equation} \label{eqn.thisisit}
\xymatrix{
{\sf{A}}^{U} \ar@<1ex>[r]^{{G'}^{U}} \ar[d]_{H^{U}} \ar@{=>}[dr]|{\Phi^{U}}\hole & {\sf{B}}^{U} \ar@<1ex>[l]^{{F'}^{U}} \ar[d]^{I^{U}}\\
{\sf{C}}^{U} \ar@<1ex>[r]^{G^{U}}  & {\sf{D}}^{U} \ar@<1ex>[l]^{F^{U}}.
}
\end{equation}
is a square.  Then this data defines a square of indexed categories if and only if, for all $f:V \rightarrow U$
\begin{equation} \label{eqn.lastsquaress}
\xymatrix{
{\sf{A}}^{U} \ar@<1ex>[r]^{{G'}^{U}} \ar[d]_{f^{*}} \ar@{=>}[dr]|{\Lambda^{G'}}\hole & {\sf{B}}^{U} \ar@<1ex>[l]^{{F'}^{U}} \ar[d]^{f^{*}} & \mbox{and} & {\sf{C}}^{U} \ar@<1ex>[r]^{{G}^{U}} \ar[d]_{f^{*}} \ar@{=>}[dr]|{\Lambda^{G}}\hole & {\sf{D}}^{U} \ar@<1ex>[l]^{{F}^{U}} \ar[d]^{f^{*}} \\
{\sf{A}}^{V} \ar@<1ex>[r]^{{G'}^{V}}  & {\sf{B}}^{V} \ar@<1ex>[l]^{{F'}^{V}} & & {\sf{C}}^{V} \ar@<1ex>[r]^{G^{V}}  & {\sf{D}}^{V} \ar@<1ex>[l]^{F^{V}}.
}
\end{equation}
are squares with duals $(\Lambda^{F'})^{-1}$ and $(\Lambda^{F})^{-1}$ respectively.\index{2-cell dual to phi@2-cell dual to $\Phi$}
\end{lemma}

\begin{proof}
Suppose (\ref{eqn.lastsquaress}) are squares with duals $(\Lambda^{F'})^{-1}$ and $(\Lambda^{F})^{-1}$ respectively.  Then the diagrams
\begin{equation} \label{eqn.thisisit2}
\xymatrix{
f^{*} \ar@{=>}[rr]^{{\eta'}^{V}*f^{*}} \ar@{=>}[d]_{f^{*} * {\eta'}^{U}} & & {G'}^{V}{F'}^{V}f^{*} \ar@{=>}[d]^{{G'}^{V}*{\Lambda^{F'}}^{-1}} & & F^{V}f^{*}G^{U} \ar@{=>}[rr]^{F^{V}*\Lambda^{G}} \ar@{=>}[d]_{{\Lambda^{F}}^{-1}* G^{U}} & &  F^{V}G^{V}f^{*} \ar@{=>}[d]^{{\epsilon}^{V}*f^{*}} \\
f^{*}{G'}^{U}{F'}^{U} \ar@{=>}[rr]_{\Lambda^{G'}*{F'}^{U}} & & {G'}^{V}f^{*}{F'}^{U} & & f^{*}F^{U}G^{U} \ar@{=>}[rr]_{f^{*}*\epsilon^{U}} & & f^{*}
}
\end{equation}
must commute by Proposition \ref{prop.adcom} so that (\ref{eqn.thisisit}) defines a square of indexed categories.

Conversely, if (\ref{eqn.thisisit}) defines a square of indexed categories, (\ref{eqn.thisisit2}) commutes by definition.
\end{proof}\index{square of indexed categories}

\chapter{Construction of the Functor $\Gamma_{n}$} All schemes
introduced in this chapter are separated, noetherian schemes.  As
before, we let ${\sf{S}}$ denote the category of affine noetherian
schemes. We use the results in the previous chapter to define the
functor $\Gamma_{n}$.  After we review the definition of bimodule
and bimodule algebra and give basic properties of these objects,
we show that if $f:V \rightarrow U$ is a morphism in $\sf{S}$,
$\tilde{f}=f \times \id_{X}:V \times X \rightarrow U \times X$,
$\mathcal{B}$ is an ${\mathcal{O}}_{U \times X}$-bimodule algebra
(Definition \ref{def.bimodalg}) and $\mathcal{M}$ is a
$\mathcal{B}$-module (Definition \ref{def.bimodalg}), then
${\tilde{f}}^{2*}\mathcal{B}$ inherits a bimodule algebra
structure from $\mathcal{B}$ and $f^{*}\mathcal{M}$ inherits a
${\tilde{f}}^{2*}\mathcal{B}$-module structure from $\mathcal{M}$
(Theorem \ref{theorem.lift}).  These observations are required in
order to define $\Gamma_{n}$ (Definition \ref{def.gamma}) and this
definition concludes the chapter.

\section{Bimodules}
We define the term {\it bimodule}.  We use a slightly different notion than that appearing in \cite{translation}.

\begin{definition}
Let $Y$ be a scheme, and suppose $\mathcal{M}'$ is a coherent ${\mathcal{O}}_{Y}$-module.  Then the set Supp $\mathcal{M}'$ is closed (\cite[Proposition 1.22 i]{alggeo1}).  From now on, we abuse notation by calling this closed set with the induced reduced scheme structure {\bf Supp $\mathcal{M}'$}\index{support}.
\end{definition}

\begin{definition}
Let $f:Y \rightarrow X$ be a morphism of schemes and let $\mathcal{M}$ be a quasi-coherent ${\mathcal{O}}_{Y}$-module.  We say $\mathcal{M}$ is {\bf relatively locally affine (rla) for $f$}\index{relatively locally affine|textbf}, if for all coherent $\mathcal{M}'\subset \mathcal{M}$ one has that $f|_{\mbox{Supp }\mathcal{M}'}$ is affine.  We call the full subcategory of ${\sf{Qcoh}}Y$ consisting of all quasi-coherent ${\mathcal{O}}_{Y}$-modules which are rla with respect to $f$ ${\sf{rla}}_{f}$.\index{rlaf@${\sf{rla}}_{f}$|textbf}
\end{definition}
Compare \cite[Definition 2.1, p.439]{translation}.

The next result is a variation on \cite[Proposition 2.2]{translation}.  The only difference is that we work with maps that are relatively locally affine, while Van den Bergh works with relatively locally finite maps.  The proof of the following result is the same as that appearing in \cite{translation}.

\begin{proposition} \label{prop.canon}
Let $Z$, $Z'$, and $Y$ be schemes, and let $f:Y \rightarrow Z$, $p:Z' \rightarrow Z$ be morphisms.
\begin{enumerate}
\item{(Projection formula)\index{projection formula}  If $\mathcal{M}$ is an ${\mathcal{O}}_{Y}$-module, rla for $f$, and $\mathcal{N}$ is a quasi-coherent ${\mathcal{O}}_{Z}$-module, then the map
\begin{equation} \label{eqn.projform}
\Psi_{\mathcal{N}}:f_{*}\mathcal{M}\otimes_{{\mathcal{O}}_{Z}}\mathcal{N} \rightarrow f_{*}(\mathcal{M} \otimes_{{\mathcal{O}}_{Y}}f^{*}\mathcal{N})
\end{equation}
appearing in Example \ref{example.basechange} is an isomorphism.}
\item{(Base change)  If there is a pullback diagram
$$
\xymatrix{
Y' \ar[r]^{p'} \ar[d]_{f'} & Y \ar[d]^{f}\\
Z' \ar[r]_{p} & Z
}
$$
and if $\mathcal{M}$ is a quasi-coherent ${\mathcal{O}}_{Y}$-module, rla for $f$, then ${p'}^{*}\mathcal{M}$ is rla for $f'$, and the map
\begin{equation} \label{eqn.pull}
\Psi_{\mathcal{M}}:p^{*}f_{*}\mathcal{M} \rightarrow {f'}_{*}{p'}^{*}\mathcal{M}
\end{equation}
dual to $\Phi$ in Example \ref{example.projform} is an isomorphism\index{base change|textbf}.}
\end{enumerate}
\end{proposition}
{\it Remark.}  In the previous proposition, it is important that the maps (\ref{eqn.projform}) and (\ref{eqn.pull}) are those constructed in Examples \ref{example.projform} and \ref{example.basechange} respectively.  With these choices of maps, one can show that the tensor product of bimodules is associative in such a way that the associativity isomorphisms satisfy MacLane's pentagon axiom\index{MacLane's pentagon axiom} (Proposition \ref{prop.tensor}).  Thus, the category of ${\mathcal{O}}_{X}$-bimodules is a monoidal category\index{monoidal category} with its tensor product.  This allows us to define an ${\mathcal{O}}_{X}$-bimodule algebra as an algebra object in this monoidal category.

We now define {\it bimodules}.  Our definition is essentially the same as \cite[Definition 2.3, p.440]{translation}.
\begin{definition} \label{def.bimodule}
Let $X$ and $Y$ be schemes.  An {\bf ${\mathcal{O}}_{S}$-central ${\mathcal{O}}_{X}-{\mathcal{O}}_{Y}$-bimodule}, or just an {\bf ${\mathcal{O}}_{X}-{\mathcal{O}}_{Y}$-bimodule} when the base $S$ is understood, is a quasi-coherent ${\mathcal{O}}_{X \times_{S} Y}$-module, rla for the projections $\pr_{1,2}:X \times_{S} Y \rightarrow X,Y$\index{bimodule|textbf}.
\end{definition}
For completeness, we define the tensor product of bimodules and prove some basic facts about this product.  The following Definition and Proposition are from \cite{translation}.
For the remainder of Section 3.1, let $\pr_{ij}:X \times Y \times Z \rightarrow X \times Y, Y \times Z, X \times Z$ be the standard projections.
\begin{definition} \label{def.tensorproduct} \cite[p.442]{translation}
If $\mathcal{M}$ is an ${\mathcal{O}}_{X}-{\mathcal{O}}_{Y}$-bimodule and $\mathcal{N}$ is an ${\mathcal{O}}_{Y}-{\mathcal{O}}_{Z}$-bimodule then the {\bf tensor product of $\mathcal{M}$ and $\mathcal{N}$} is an ${\mathcal{O}}_{X}-{\mathcal{O}}_{Z}$-bimodule, defined by
$$
\mathcal{M}\otimes_{{\mathcal{O}}_{Y}} \mathcal{N} = {\pr_{13}}_{*}({\pr_{12}}^{*}\mathcal{M} \otimes_{{\mathcal{O}}_{X \times Y \times Z}} {\pr_{23}}^{*}\mathcal{N}).
$$
\end{definition}\index{tensor product!of bimodules|textbf}
Suppose $\mathcal{M}$ is an ${\mathcal{O}}_{Y}$-module and $\mathcal{N}$ is an ${\mathcal{O}}_{Y}-{\mathcal{O}}_{Z}$-bimodule.  If $\pr_{i}:Y \times Z \rightarrow Y,Z$ is a projection, then we define
$$
\mathcal{M}\otimes_{{\mathcal{O}}_{Y}}\mathcal{N} = \pr_{2*}(\pr_{1}^{*}\mathcal{M}\otimes_{{\mathcal{O}}_{Y \times Z}} \mathcal{N}).
$$
\index{tensor product!of a module and a bimodule|textbf}Associativity of the tensor product is established using only the canonical isomorphisms in Proposition \ref{prop.canon}.  For the readers convenience, we include the proof of this fact appearing in \cite{translation}.
\begin{proposition} \label{prop.tensor} \cite[Proposition 2.5, p.442]{translation}
The tensor product of bimodules is associative\index{tensor product!of bimodules!is associative} and satisfies MacLane's pentagon axiom\index{MacLane's pentagon axiom}.  The tensor product of a coherent module to the left of two coherent bimodules is associative\index{tensor product!of a module and a bimodule!is associative}.
\end{proposition}

\begin{proof}
Let $X_{1}$, $X_{2}$, $X_{3}$ and $X_{4}$ be schemes.  For $I \subset \{1, 2, 3, 4 \}$ denote by $X_{I}$ the product $\times_{i \in I} X_{i}$ and for $I \subset J \subset \{1,2,3,4\}$ let $\pr_{I}^{J}$ be the projection $X_{J} \rightarrow X_{I}$.  Let ${\mathcal{M}}_{i,i+1}$ be an ${\mathcal{O}}_{S}$-central ${\mathcal{O}}_{X_{i}}-{\mathcal{O}}_{X_{i+1}}$-bimodule.  Then
$$
({\mathcal{M}}_{12} \otimes_{{\mathcal{O}}_{X_{2}}} {\mathcal{M}}_{23}) \otimes_{{\mathcal{O}}_{X_{3}}} {\mathcal{M}}_{34}
$$
is by definition,
\begin{equation} \label{eqn.bigten}
\pr_{14*}^{134}(\pr_{13}^{134*}\pr_{13*}^{123}(\pr_{12}^{123*}{\mathcal{M}}_{12} \otimes_{{\mathcal{O}}_{X_{123}}}\pr_{23}^{123*}{\mathcal{M}}_{23}) \otimes_{{\mathcal{O}}_{X_{134}}}\pr_{34}^{134*}{\mathcal{M}}_{34}).
\end{equation}
By Proposition \ref{prop.canon} (2) applied to the diagram
$$
\xymatrix{
X_{1234} \ar[r] \ar[d] & X_{123} \ar[d] \\
X_{134} \ar[r] & X_{13}
}
$$
it follows that
$$
\pr_{13}^{134*}\pr_{13*}^{123} \cong \pr_{134*}^{1234}\pr_{123}^{1234*}
$$
hence (\ref{eqn.bigten}) is isomorphic to
$$
\pr_{14*}^{134}(\pr_{134*}^{1234}\pr_{123}^{1234*}(\pr_{12}^{123*}{\mathcal{M}}_{12} \otimes_{{\mathcal{O}}_{X_{123}}}\pr_{23}^{123*}{\mathcal{M}}_{23}) \otimes_{{\mathcal{O}}_{X_{134}}}\pr_{34}^{134*}{\mathcal{M}}_{34}) \cong
$$
\begin{equation} \label{eqn.bigten2}
\pr_{14*}^{134}(\pr_{134*}^{1234}(\pr_{12}^{1234*}{\mathcal{M}}_{12} \otimes_{{\mathcal{O}}_{X_{1234}}}\pr_{23}^{1234*}{\mathcal{M}}_{23}) \otimes_{{\mathcal{O}}_{X_{134}}}\pr_{34}^{134*}{\mathcal{M}}_{34}).
\end{equation}
Applying the projection formula (Proposition \ref{prop.canon} (1)) for $\pr_{134}^{1234}$ shows that (\ref{eqn.bigten2}) is isomorphic to
$$
\pr_{14*}^{134}(\pr_{134*}^{1234}(\pr_{12}^{1234*}{\mathcal{M}}_{12} \otimes_{{\mathcal{O}}_{X_{1234}}}\pr_{23}^{1234*}{\mathcal{M}}_{23} \otimes_{{\mathcal{O}}_{X_{1234}}}\pr_{134}^{1234*}\pr_{34}^{134*}{\mathcal{M}}_{34})) \cong
$$
$$
\pr_{14*}^{1234}(\pr_{12}^{1234*}{\mathcal{M}}_{12} \otimes_{{\mathcal{O}}_{X_{1234}}}\pr_{23}^{1234*}{\mathcal{M}}_{23} \otimes_{{\mathcal{O}}_{X_{1234}}}\pr_{34}^{1234*}{\mathcal{M}}_{34}).
$$
A similar computation shows that we get the same result if we start from
$$
{\mathcal{M}}_{12} \otimes_{{\mathcal{O}}_{X_{2}}} ({\mathcal{M}}_{23} \otimes_{{\mathcal{O}}_{X_{3}}} {\mathcal{M}}_{34}).
$$
This yields a natural isomorphism between
$$
({\mathcal{M}}_{12} \otimes_{{\mathcal{O}}_{X_{2}}} {\mathcal{M}}_{23}) \otimes_{{\mathcal{O}}_{X_{3}}} {\mathcal{M}}_{34}
$$
and
$$
{\mathcal{M}}_{12} \otimes_{{\mathcal{O}}_{X_{2}}} ({\mathcal{M}}_{23} \otimes_{{\mathcal{O}}_{X_{3}}} {\mathcal{M}}_{34}).
$$
An easy, but very tedious verification shows that this isomorphism satisfies the pentagon axiom.

The final assertion is proven in a similar way.
\end{proof}
${\mathcal{O}}_{X}-{\mathcal{O}}_{Y}$-bimodules also have left and right multiplication maps which we now describe.
\begin{proposition} \cite[Definition 3.1, p.449]{translation} \label{prop.scalar}
Suppose $X$ and $Y$ are schemes, $d:Y \rightarrow Y \times Y$ is the diagonal, $\mathcal{B}$ is any ${\mathcal{O}}_{X \times Y}$-module and $\mathcal{C}$ is any ${\mathcal{O}}_{Y \times X}$-module.  There exist natural isomorphisms ${\mathcal{O}}_{\mu}:d_{*}{\mathcal{O}}_{Y} \otimes_{{\mathcal{O}}_{Y}}\mathcal{C} \rightarrow \mathcal{C}$ and ${}_{\mu}\mathcal{O}:\mathcal{B} \otimes_{{\mathcal{O}}_{Y}} d_{*}{\mathcal{O}}_{Y} \rightarrow \mathcal{B}$ called the left and right ${\mathcal{O}}_{X}$-linear multiplication maps or the left and right scalar multiplication maps.  Similarly, if $\mathcal{M}$ is any ${\mathcal{O}}_{Y}$-module, there exists a natural isomorphism ${}_{\mu}\mathcal{O}:\mathcal{M} \otimes_{{\mathcal{O}}_{Y}} d_{*}{\mathcal{O}}_{Y} \rightarrow \mathcal{M}$.\index{linear multiplication map|textbf}

\end{proposition}

\begin{proof}
We describe the left multiplication map.  The construction of the right multiplication map is similar, so we omit it.  Let $e = \id_{X} \times d$.  We note that there is a pullback diagram
$$
\xymatrix{
X \times Y \ar[r]^{\pr_{2}} \ar[d]_{e} & Y \ar[d]^{d} \\
X \times Y^{2} \ar[r]_{\pr_{23}} & Y^{2}.
}
$$
Thus,
\begin{align} \label{eqn.unitary}
\pr_{13*}(\pr_{12}^{*}\mathcal{B} \otimes \pr_{23}^{*}d_{*}{\mathcal{O}}_{X}) & \cong \pr_{13*}(\pr_{12}^{*}\mathcal{B} \otimes e_{*}\pr_{2}^{*}{\mathcal{O}}_{X}) \\
& \cong \pr_{13*}e_{*}(e^{*}\pr_{12}^{*}\mathcal{B} \otimes \pr_{2}^{*}{\mathcal{O}}_{X}) \nonumber \\
& = e^{*}\pr_{23}^{*}\mathcal{B} \otimes \pr_{2}^{*}{\mathcal{O}}_{X} \nonumber.
\end{align}
where the first isomorphism is induced by the canonical isomorphism (Proposition \ref{prop.canon} (2)) while the second isomorphism is the projection formula (Proposition \ref{prop.canon} (1)).  Now, $\pr_{2}^{*}{\mathcal{O}}_{X} \cong {\mathcal{O}}_{X^{2}}$ and $\pr_{12}e=\operatorname{id}$ so $e^{*}\pr_{23}^{*}\mathcal{B}\otimes \pr_{2}^{*}{\mathcal{O}}_{X} \cong \mathcal{B}$.

The assertion involving $\mathcal{M}$ is proved in a similar fashion.
\end{proof}
We show that bimodule tensor products share various properties with ordinary tensor products of modules.   To this end, we introduce the {\it scheme theoretic support} of a bimodule.
\begin{proposition} \label{prop.ssupp} \cite[Proposition 1.27, p. 59]{alggeo1}
Let $\mathcal{F}$ be a coherent ${\mathcal{O}}_{X}$-module.  Let $\mathcal{I}$ be the annihilator of $\mathcal{F}$.  Give Supp $\mathcal{F}$ the topology it inherits from $X$, and let $i:\operatorname{Supp }\mathcal{F} \rightarrow X$ be inclusion.  Then $(\operatorname{Supp }\mathcal{F}, i^{-1}({\mathcal{O}}_{X}/\mathcal{I}))$ is a closed subscheme of $X$.
\end{proposition}

\begin{definition}
Keep notation as in Proposition \ref{prop.ssupp}.  The {\bf scheme theoretic support of ${\mathcal{F}}$}\index{scheme theoretic support|textbf}, denoted {\bf $\operatorname{SSupp } \mathcal{F}$}, is the scheme $(\operatorname{Supp }\mathcal{F}, i^{-1}({\mathcal{O}}_{X}/\mathcal{I}))$.
\end{definition}
\begin{lemma}
Let $i:Z \rightarrow X$ be a closed immersion of schemes which is a homeomorphism on spaces.  If $U \subset Z$ is affine, and if $i(U)$ is not affine in $X$, then there exists a point $p \in U$ such that ${\mathcal{O}}_{Z,p}=0$.  In particular, if $Y$ is noetherian, $\mathcal{A}$ is a coherent ${\mathcal{O}}_{Y}$-module,
$$
i:\operatorname{Supp }\mathcal{A} \rightarrow \operatorname{SSupp }\mathcal{A}
$$
\index{support}is inclusion, and $U \subset \operatorname{Supp }\mathcal{A}$ is affine, then $i(U) \subset \operatorname{SSupp }\mathcal{A}$ is affine.
\end{lemma}

\begin{proof}
We begin by showing that the second claim follows from the first.  If $Z = \mbox{Supp }\mathcal{A}$, and if $p \in Z$, then since ${\mathcal{A}}_{p} \neq 0$, ${\mathcal{O}}_{Z,p} \neq 0$.

To prove the first claim, we note that since $i$ is a closed immersion, there is a surjection ${\mathcal{O}}_{X} \rightarrow i_{*}{\mathcal{O}}_{Z}$, so that we have a short exact sequence
$$
0 \rightarrow \mathcal{J} \rightarrow {\mathcal{O}}_{X} \rightarrow i_{*}{\mathcal{O}}_{Z} \rightarrow 0.
$$

Thus, ${\mathcal{O}}_{Z}(U) \cong {\mathcal{O}}_{X}(i(U))/{\mathcal{J}}(i(U))$.  By hypothesis, Spec ${\mathcal{O}}_{X}(i(U))/{\mathcal{J}}(i(U)) \cong i(U)$.  If $i(U)$ is not affine in $X$ then Spec ${\mathcal{O}}_{X}(i(U)) \neq $ Spec ${\mathcal{O}}_{X}(i(U))/{\mathcal{J}}(i(U)) \cong i(U)$.  Thus, there exists a prime of ${\mathcal{O}}_{X}(i(U))$, $p$, which does not contain ${\mathcal{J}}(i(U))$.  We may then deduce that ${\mathcal{J}}_{p}={\mathcal{O}}_{X,p}$, so that ${\mathcal{O}}_{Z,p}=0$.
\end{proof}

\begin{corollary} \label{cor.sexend}
Suppose $f:Y \rightarrow X$ is a morphism of schemes and suppose $\mathcal{M}$ is a coherent ${\mathcal{O}}_{Y}$-module which is rla with respect to $f$.  Suppose $j: \operatorname{SSupp }\mathcal{M} \rightarrow Y$ is inclusion.  Then $fj$ is affine.
\end{corollary}

\begin{proof}
If $i:\operatorname{Supp }\mathcal{M} \rightarrow \operatorname{SSupp }\mathcal{M}$ and $k:\operatorname{Supp }\mathcal{B} \rightarrow Y$ are inclusions, then the diagram of schemes
$$
\xymatrix{
\operatorname{Supp }\mathcal{M} \ar[r]^{k} \ar[dr]_{i} & Y \ar[r]^{f} & X \\
& \operatorname{SSupp }\mathcal{M} \ar[u]_{j}
}
$$
commutes, so $ji=k$.  Since $fk$ is affine by hypothesis, $fji$ is affine.  Thus, if $V \subset X$ is affine open, then $i^{-1}j^{-1}f^{-1}(V) \subset \operatorname{Supp }\mathcal{M}$ is open affine.  By the previous lemma, $j^{-1}f^{-1}(V)$ is affine.  Thus $fj$ is affine, as desired.
\end{proof}
\begin{lemma} \label{lem.exact}
Let $f:Y \rightarrow X$ be a map of schemes.  If $\mathcal{M}$ and $\mathcal{N}$ are ${\mathcal{O}}_{Y}$-modules which are rla with \index{relatively locally affine}respect to $f$, and if $\phi:\mathcal{M} \rightarrow \mathcal{N}$ is an epimorphism, then $f_{*}\phi$ is an epimorphism.
\end{lemma}

\begin{proof}
Since the image of a coherent module is coherent, $\phi$ determines a direct system of epimorphisms $\phi_{i}:{\mathcal{M}}_{i} \rightarrow {\mathcal{N}}_{i}$ between coherent modules whose direct limit is $\phi$.  Since direct limits are right exact, it suffices to prove the Lemma in the case that $\mathcal{M}$ and $\mathcal{N}$ are coherent.  Let $i:\operatorname{SSupp }\mathcal{N} \rightarrow Y$ be inclusion.  By naturality of the unit of the adjoint pair $(i^{*},i_{*})$, the diagram
$$
\xymatrix{
f_{*}\mathcal{M} \ar[d] \ar[r]^{f_{*}\phi} & f_{*}\mathcal{N} \ar[d] \\
f_{*}i_{*}i^{*}\mathcal{M} \ar[r]_{f_{*}i_{*}i^{*}\phi} & f_{*}i_{*}i^{*}\mathcal{N}
}
$$
whose verticals are the unit of $(i^{*},i_{*})$, commutes.  Since $\phi$ is an epi, $\operatorname{SSupp }\mathcal{N} \subset \operatorname{SSupp }\mathcal{M}$ so that the left vertical is an epimorphism.  Since $\mathcal{N}$ is rla with respect to $f$, $fi$ is affine by Corollary \ref{cor.sexend}, so the bottom horizontal is an epimorphism.  Finally the right hand vertical is an isomorphism.  Thus $f_{*}\phi$ is an epimorphism as well.
\end{proof}

\begin{definition} \label{def.adip}
Let $X$, $Y$, and $Z$ be schemes, and suppose $\mathcal{E}$ is a ${\mathcal{O}}_{X \times Y}$-module and $\mathcal{F}$ is an ${\mathcal{O}}_{Y \times Z}$-module.  $\mathcal{E}$ and $\mathcal{F}$ have the {\bf affine direct image property}\index{affine direct image property} if the restriction of the projection map $\pr_{13}$ to

$$
\operatorname{SSupp}({\pr_{12}}^{*}\mathcal{E} \otimes {\pr_{23}}^{*}\mathcal{F})
$$
is affine.
\end{definition}

The following result is an immediate consequence of Corollary \ref{cor.sexend}.
\begin{lemma}
If $\mathcal{E}$ is a coherent ${\mathcal{O}}_{X \times Y}$-module and $\mathcal{F}$ is a coherent ${\mathcal{O}}_{Y \times Z}$-module, then $\mathcal{E}$ and $\mathcal{F}$ have the affine direct image property \index{affine direct image property}if the restriction of the projection map $\pr_{13}$ to

$$
\operatorname{Supp}({\pr_{12}}^{*}\mathcal{E} \otimes {\pr_{23}}^{*}\mathcal{F})
$$
is affine.
\end{lemma}

Of course, we have the following
\begin{lemma} \label{lem.rlaadip}
If $\mathcal{E}$ is an ${\mathcal{O}}_{X \times Y}$-bimodule and $\mathcal{F}$ is an ${\mathcal{O}}_{Y \times Z}$-bimodule, then $\pr_{12}^{*}\mathcal{E} \otimes \pr_{23}^{*}\mathcal{F}$ is rla \index{relatively locally affine}with respect to $\pr_{13}$.  In particular, if $\mathcal{E}$ and $\mathcal{F}$ are coherent, then $\mathcal{E}$ and $\mathcal{F}$ have the affine direct image property\index{affine direct image property}.
\end{lemma}

\begin{proof}
Since $X \times Y$ and $Y \times Z$ are noetherian, $\mathcal{E}$ and $\mathcal{F}$ are direct limits of coherent submodules.  Furthermore, since $\pr_{12}^{*}$, $\pr_{23}^{*}$ and $-\otimes-$ commute with direct limits, $\pr_{12}^{*}\mathcal{E} \otimes \pr_{23}{*}\mathcal{F}$ is a direct limit of submodules of the form $\pr_{12}^{*}\mathcal{E}'\otimes \pr_{23}^{*}\mathcal{F}'$, where $\mathcal{E}'$ and $\mathcal{F}'$ are coherent submodules of $\mathcal{E}$ and $\mathcal{F}$, respectively.  Thus, if $\mathcal{M}$ is a coherent submodule of $\pr_{12}^{*}\mathcal{E} \otimes \pr_{23}^{*}\mathcal{F}$, then the fact that $X \times Y \times Z$ is noetherian implies $\mathcal{M}$ is a submodule of $\pr_{12}^{*}\mathcal{E}' \otimes \pr_{23}^{*}\mathcal{F}'$ for coherent submodules $\mathcal{E}'$ and $\mathcal{F}'$ of $\mathcal{E}$ and $\mathcal{F}$ respectively.  In particular, it suffices to prove the Lemma in the case that $\mathcal{E}$ and $\mathcal{F}$ are coherent.

We next show that $\pr_{13}$ restricted to $\operatorname{Supp }\pr_{12}^{*}\mathcal{E} \otimes \pr_{23}^{*}\mathcal{F}$ is affine.  Let $U \subset X$ and $V \subset Z$ be affine open sets.  Since $\mathcal{E}$ and $\mathcal{F}$ are coherent bimodules, $U \times Y \cap \operatorname{Supp }\mathcal{E}$ and $Y \times V \cap \operatorname{Supp }\mathcal{F}$ are affine open subsets of $\operatorname{Supp }\mathcal{E}$ and $\operatorname{Supp }\mathcal{F}$ respectively.

Since $\mathcal{E}$ and $\mathcal{F}$ are coherent,
$$
\operatorname{Supp }\pr_{12}^{*}\mathcal{E} \otimes \pr_{23}^{*}\mathcal{F} = \operatorname{Supp }\pr_{12}^{*}\mathcal{E} \cap \operatorname{Supp }\pr_{23}^{*}\mathcal{F}.
$$
In addition, $\operatorname{Supp }\pr_{12}^{*}\mathcal{E}=\pr_{12}^{-1}(\operatorname{Supp }\mathcal{E})$ and $\operatorname{Supp }\pr_{23}^{*}\mathcal{F}=\pr_{23}^{-1}(\operatorname{Supp }\mathcal{F})$.  Thus, $\pr_{13}^{-1}(U \times V)$ restricted to $\operatorname{Supp }\pr_{12}^{*}\mathcal{E} \otimes \pr_{23}^{*}\mathcal{F}$ is
\begin{equation} \label{eqn.sepeee}
(U \times Y \cap \operatorname{Supp }\mathcal{E}) \times V \cap U \times (Y \times V \cap \operatorname{Supp }\mathcal{F}).
\end{equation}
The latter set is the intersection of affine open sets.  Since $\operatorname{Supp }\pr_{12}^{*}\mathcal{E} \otimes \pr_{23}^{*}\mathcal{F}$ is separated, (\ref{eqn.sepeee}) is affine.

By the previous lemma, $\mathcal{E}$ and $\mathcal{F}$ have the affine direct image property.  It follows that $\pr_{13}$ restricted to $\operatorname{Supp }\mathcal{M}$ is affine as desired.
\end{proof}

\begin{corollary} \cite[Proposition 2.6, p.443]{translation}
The tensor product of bimodules is right exact.\index{tensor product!of bimodules!is right exact}
\end{corollary}

\begin{proof}
By Lemma \ref{lem.rlaadip}, Lemma \ref{lem.exact} applies, so that the result follows from the fact that ordinary tensor products and inverse image functors are right exact.
\end{proof}

\begin{lemma} \label{cor.tensor}
Let $X$ be a scheme, and suppose
$$
\xymatrix{
0 \ar[r] & \mathcal{K} \ar[r]^{\kappa} & \mathcal{M} \ar[r]^{\phi} & \mathcal{M}' \ar[r] & 0
}
$$
and
$$
\xymatrix{
0 \ar[r] & \mathcal{L} \ar[r]^{\gamma} & \mathcal{N} \ar[r]^{\psi} & \mathcal{N}' \ar[r] & 0
}
$$
are short exact sequences of ${\mathcal{O}}_{X}$-modules.  Then
$$
\operatorname{ker }\phi \otimes \psi = (\kappa \otimes 1)(\mathcal{K} \otimes \mathcal{N})+(1 \otimes \gamma)(\mathcal{M} \otimes \mathcal{L}).
$$
\end{lemma}

\begin{proof}
By the naturality of the tensor product, we have a commutative diagram
$$
\xymatrix{
\mathcal{K}\otimes \mathcal{N} \ar[r]^{\kappa \otimes 1} & \mathcal{M} \otimes \mathcal{N} \ar[r]^{\phi\otimes 1} \ar[d]_{1 \otimes \psi} & \mathcal{M}'\otimes \mathcal{N} \ar[r] \ar[d]^{1 \otimes \psi} & 0 \\
& \mathcal{M} \otimes \mathcal{N}' \ar[r]_{\phi \otimes 1} & \mathcal{M}' \otimes \mathcal{N}'
}
$$
whose top row is exact.  Thus $\kappa \otimes 1 (\mathcal{K} \otimes \mathcal{N}) \subset \operatorname{ker }\phi \otimes \psi$.  Similarly, $1 \otimes \gamma (\mathcal{M} \otimes \mathcal{L}) \subset \operatorname{ker }\phi \otimes \psi$.  Thus
$$
\kappa \otimes 1 (\mathcal{K} \otimes \mathcal{N}) + 1 \otimes \gamma (\mathcal{M} \otimes \mathcal{L}) \subset \operatorname{ker }\phi \otimes \psi
$$
defines an inclusion between these sheaves.  Thus, to show equality, we need only show equality locally.  But local equality follows from \cite[formula 7, p.267]{comalg2}.
\end{proof}

\begin{corollary} \label{cor.bimodtensor}\index{tensor product!of bimodules}
Let $X$, $Y$ and $Z$ be schemes and suppose
$$
\xymatrix{
0 \ar[r] & \mathcal{K} \ar[r]^{\kappa} & \mathcal{M} \ar[r]^{\phi} & \mathcal{M}' \ar[r] & 0
}
$$
is a short exact sequence of ${\mathcal{O}}_{X \times Y}$-modules,
$$
\xymatrix{
0 \ar[r] & \mathcal{L} \ar[r]^{\gamma} & \mathcal{N} \ar[r]^{\psi} & \mathcal{N}' \ar[r] & 0
}
$$
is a short exact sequence of ${\mathcal{O}}_{Y \times Z}$-modules, and $\mathcal{M}$ and $\mathcal{N}$ are bimodules.  Then
$$
\operatorname{ker }\phi \otimes_{{\mathcal{O}}_{Y}} \psi = (\kappa \otimes_{{\mathcal{O}}_{Y}} 1)(\mathcal{K} \otimes_{{\mathcal{O}}_{Y}} \mathcal{N})+(1 \otimes_{{\mathcal{O}}_{Y}} \gamma)(\mathcal{M} \otimes_{{\mathcal{O}}_{Y}} \mathcal{L}).
$$
\end{corollary}

\begin{proof}
We have short exact sequences
$$
\xymatrix{
0 \ar[r] & \pr_{12}^{*}\kappa(\pr_{12}^{*}\mathcal{K}) \ar[r] & \pr_{12}^{*}\mathcal{M} \ar[r]^{\pr_{12}^{*}\phi} & \pr_{12}^{*}\mathcal{M}' \ar[r] & 0
}
$$
and
$$
\xymatrix{
0 \ar[r] & \pr_{23}^{*}\gamma(\pr_{23}^{*}\mathcal{L}) \ar[r] & \pr_{23}^{*}\mathcal{N} \ar[r]^{\pr_{23}^{*}\psi} & \pr_{23}^{*}\mathcal{N}' \ar[r] & 0.
}
$$
Thus, by Lemma \ref{cor.tensor},
$$
\operatorname{ker }(\pr_{12}^{*}\phi \otimes \pr_{23}^{*}\psi) = 1 \otimes \pr_{23}^{*}\gamma (\pr_{12}^{*}\mathcal{M} \otimes \pr_{23}^{*}\mathcal{L}) + \pr_{12}^{*}\kappa \otimes 1 (\pr_{12}^{*}\mathcal{K} \otimes \pr_{23}\mathcal{N}).
$$
Since $\pr_{13*}$ is left exact,
$$
\pr_{13*}\operatorname{ker }(\pr_{12}^{*}\phi \otimes \pr_{23}^{*}\psi) = \operatorname{ker }\pr_{13*}(\pr_{12}^{*}\phi \otimes \pr_{23}^{*}\psi).
$$
Since all modules involved are bimodules, $\pr_{12}^{*}\mathcal{M} \otimes \pr_{23}^{*}\mathcal{L}$, $\pr_{12}^{*}\mathcal{M} \otimes \pr_{23}^{*}\mathcal{N}$, and $\pr_{12}^{*}\mathcal{K} \otimes \pr_{23}^{*}\mathcal{N}$ are rla with respect to $\pr_{13}$.  Thus, by Lemma \ref{lem.exact} and a variant of \cite[Proposition 2.2 (2), p. 439]{translation},
$$
\pr_{13*}(1 \otimes \pr_{23}^{*}\gamma (\pr_{12}^{*}\mathcal{M} \otimes \pr_{23}^{*}\mathcal{L})) + \pr_{13*}(\pr_{12}^{*}\kappa \otimes 1 (\pr_{12}^{*}\mathcal{K} \otimes \pr_{23}\mathcal{N}))=
$$
$$
\pr_{13*}(1 \otimes \pr_{23}^{*}\gamma) (\pr_{13*}(\pr_{12}^{*}\mathcal{M} \otimes \pr_{23}^{*}\mathcal{L})) +
$$
$$
\pr_{13*}(\pr_{12}^{*}\kappa \otimes 1) (\pr_{13*}(\pr_{12}^{*}\mathcal{K} \otimes \pr_{23}\mathcal{N}))
$$
as desired.
\end{proof}

\section{Bimodule algebras}\index{bimodule algebra|(}\index{bimodule algebra!graded|(}
We review the definition of bimodule algebra and modules over a bimodule algebra from \cite{translation}.

\begin{definition} \label{def.bimodalg} \cite[Definition 3.1, p.449]{translation}
\begin{enumerate}
\item{An {\bf ${\mathcal{O}}_{X}$-bimodule algebra $\mathcal{B}$} \index{bimodule algebra|textbf}is an algebra object in the category of ${\mathcal{O}}_{X}$-bimodules.}
\item{A {\bf (right) $\mathcal{B}$-module $\mathcal{M}$}\index{B module@$\mathcal{B}$-module|textbf} is an ${\mathcal{O}}_{X}$-module, together with a right ${\mathcal{O}}_{X}$-linear multiplication map $\mu_{\mathcal{M}}:\mathcal{M} \otimes_{{\mathcal{O}}_{X}} \mathcal{B} \rightarrow \mathcal{M}$ satisfying the unit and associative axioms.}
\end{enumerate}
\end{definition}
We elaborate on item (1), above.  $\mathcal{B}$ is an algebra object in the category of ${\mathcal{O}}_{X}$-bimodules if $d:X \rightarrow X \times X$ is diagonal and there exist ${\mathcal{O}}_{X}$-linear maps $\upsilon:d_{*}{\mathcal{O}}_{X} \rightarrow \mathcal{B}$ and $\mu:\mathcal{B} \otimes \mathcal{B} \rightarrow \mathcal{B}$ satisfying the following axioms:
$
\mbox{({\it associativity}) the diagram}
$
\begin{equation} \label{eqn.associativity}
\xymatrix{
\mathcal{B} \otimes \mathcal{B} \otimes \mathcal{B} \ar[r]^{\mu \otimes \mathcal{B}} \ar[d]_{\mathcal{B} \otimes \mu} & \mathcal{B} \otimes \mathcal{B} \ar[d]^{\mu} \\
\mathcal{B} \otimes \mathcal{B} \ar[r]_{\mu} & \mathcal{B}
}
\end{equation}
commutes,
$
\mbox{({\it left unit}) the diagram}
$
\begin{equation} \label{eqn.lunit}
\xymatrix{
d_{*}{\mathcal{O}}_{X} \otimes \mathcal{B} \ar[rr]^{\upsilon \otimes \mathcal{B}} \ar[dr]_{{\mathcal{O}}_{\mu}} & & \mathcal{B} \otimes \mathcal{B}  \ar[dl]^{\mu} \\
& \mathcal{B} &
}
\end{equation}
commutes, and finally,

$
\mbox{({\it right unit}) the diagram}
$
\begin{equation} \label{eqn.runit}
\xymatrix{
\mathcal{B} \otimes d_{*}{\mathcal{O}}_{X}  \ar[rr]^{\mathcal{B} \otimes \upsilon} \ar[dr]_{{}_{\mu}{\mathcal{O}}} & & \mathcal{B} \otimes \mathcal{B} \ar[dl]^{\mu} \\
& \mathcal{B} &
}
\end{equation}
commutes.  The associativity and unit axiom for a right $\mathcal{B}$-module are analogous to associativity and the right unit axiom above.  This construction can be graded:
\begin{definition} \cite[p.453]{translation} \label{def.graded}
\begin{enumerate}
\item{A {\bf graded ${\mathcal{O}}_{X}$-bimodule algebra $\mathcal{B}$}\index{bimodule algebra!graded|textbf} is an ${\mathcal{O}}_{X}$-bimodule algebra, equipped with a decomposition $\underset{n \in \mathbb{Z}}{\oplus}{\mathcal{B}}_{n}$ as ${\mathcal{O}}_{X}$-bimodules such that $\mu=\underset{n \in \mathbb{Z}}{\oplus}\mu_{n,m}$ with $\mu_{m,n}$ mapping ${\mathcal{B}}_{m} \otimes_{{\mathcal{O}}_{X}} {\mathcal{B}}_{n} \rightarrow {\mathcal{B}}_{m+n}$ and with $\upsilon (d_{*}{\mathcal{O}}_{X}) \subset {\mathcal{B}}_{0}$.}

\item{A {\bf graded (right) $\mathcal{B}$-module $\mathcal{M}$}\index{B module@$\mathcal{B}$-module!graded|textbf} is equipped with a decomposition $\mathcal{M}=\oplus {\mathcal{M}}_{n}$ and $\mu_{\mathcal{M}}=\underset{m,n}{\oplus}(\mu_{\mathcal{M}})_{m,n}$ with $(\mu_{\mathcal{M}})_{m,n}:{\mathcal{M}}_{m} \otimes {\mathcal{B}}_{n} \rightarrow {\mathcal{M}}_{m+n}$}
\end{enumerate}
\end{definition}

\begin{definition}
Let $\mathcal{B}$ be a graded ${\mathcal{O}}_{X}$-bimodule.  Then $\mathcal{I} \subset \mathcal{B}$ is a { \bf graded two-sided ideal of $\mathcal{B}$}\index{graded two-sided ideal|textbf} if
\begin{itemize}
\item{}
$$
\mathcal{I} = \bigoplus_{i \geq 0}\mathcal{I} \cap {\mathcal{B}}_{i}
$$
and
\item{}
$$
\mu_{ij}({\mathcal{I}}_{i} \otimes {\mathcal{B}}_{j}) \subset {\mathcal{I}}_{i+j}
$$
and
$$
\mu_{ij}({\mathcal{B}}_{i} \otimes {\mathcal{I}}_{j}) \subset {\mathcal{I}}_{i+j}
$$
\end{itemize}
where ${\mathcal{I}}_{i} = \mathcal{I} \cap {\mathcal{B}}_{i}$.
\end{definition}

\begin{lemma}
Suppose $\mathcal{B}$ is a graded ${\mathcal{O}}_{X}$-bimodule algebra and $\mathcal{I}$ is a graded ideal of $\mathcal{B}$.  Then $\mathcal{B}/\mathcal{I}$ inherits a graded bimodule algebra structure from $\mathcal{B}$.
\end{lemma}

\begin{proof}
We define the multiplication for $\mathcal{B}/\mathcal{I}$.  The proof that this multiplication defines a bimodule algebra structure is straightforward but tedious, so we omit it.  Since $\mathcal{I}$ is a graded ideal, we note that ${\mathcal{I}}_{i} \otimes {\mathcal{B}}_{j} + {\mathcal{B}}_{i} \otimes {\mathcal{I}}_{j}$ is in the kernel of the composition
$$
\xymatrix{
{\mathcal{B}}_{i} \otimes {\mathcal{B}}_{j} \ar[r]^{\mu_{ij}} & {\mathcal{B}}_{i+j} \ar[r]^{\pi} & {\mathcal{B}}_{i+j}/{\mathcal{I}}_{i+j}.
}
$$
By Corollary \ref{cor.bimodtensor}, we have an isomorphism
$$
\frac{{\mathcal{B}}_{i} \otimes {\mathcal{B}}_{j}}{{\mathcal{I}}_{i} \otimes {\mathcal{B}}_{j} + {\mathcal{B}}_{i} \otimes {\mathcal{I}}_{j}} \rightarrow \frac{{\mathcal{B}}_{i}}{{\mathcal{I}}_{i}} \otimes \frac{{\mathcal{B}}_{j}}{{\mathcal{I}}_{j}}.
$$
Give $\mathcal{B}/\mathcal{I}$ the grading $(\mathcal{B}/\mathcal{I})_{i} = {\mathcal{B}}_{i}/{\mathcal{I}}_{i}$ and define the multiplication on $\mathcal{B}/\mathcal{I}$ as the composition
$$
\frac{{\mathcal{B}}_{i}}{{\mathcal{I}}_{i}} \otimes \frac{{\mathcal{B}}_{j}}{{\mathcal{I}}_{j}} \rightarrow \frac{{\mathcal{B}}_{i} \otimes {\mathcal{B}}_{j}}{{\mathcal{I}}_{i} \otimes {\mathcal{B}}_{j} + {\mathcal{B}}_{i} \otimes {\mathcal{I}}_{j}} \rightarrow \frac{{\mathcal{B}}_{i+j}}{{\mathcal{I}}_{i+j}}.
$$
\end{proof}
We omit the routine proof of the following

\begin{lemma} \label{lem.inheritmod}
Suppose ${\mathcal{B}}$ is a graded ${\mathcal{O}}_{X}$-bimodule, $\mathcal{I} \subset \mathcal{B}$ is a graded ideal and $\mathcal{M}$ is a graded $\mathcal{B}$-module such that, for all $i$, the composition
$$
\xymatrix{
{\mathcal{M}}_{i} \otimes {\mathcal{I}}_{j} \ar[r] & {\mathcal{M}}_{i} \otimes {\mathcal{B}}_{j} \ar[r]^{\mu_{ij}} & {\mathcal{M}}_{i+j}
}
$$
is zero.  Then $\mathcal{M}$ inherits a $\mathcal{B}/\mathcal{I}$-module structure.\index{B module@$\mathcal{B}$-module!graded}
\end{lemma}

\begin{definition}
Let $Y$ be a scheme.  An ${\mathbb{N}}$-graded ${\mathcal{O}}_{Y}$-bimodule algebra {\bf $\mathcal{B}$ is generated in degree one} \index{bimodule algebra!graded!generated in degree one|textbf}if, for all $j \geq 1$,
$$
\mu_{1,j}:{\mathcal{B}}_{1} \otimes_{{\mathcal{O}}_{Y}} {\mathcal{B}}_{j} \rightarrow {\mathcal{B}}_{j+1}
$$
and
$$
\mu_{j,1}:{\mathcal{B}}_{j} \otimes_{{\mathcal{O}}_{Y}} {\mathcal{B}}_{1} \rightarrow {\mathcal{B}}_{j+1}
$$
are epimorphisms.

An ${\mathbb{N}}$-graded right $\mathcal{B}$-module {\bf $\mathcal{M}$ is generated in degree zero}\index{B module@$\mathcal{B}$-module!graded!generated in degree zero|textbf} if, for all $i \geq 0$,
$$
\mu_{i,1}:{\mathcal{M}}_{i} \otimes_{{\mathcal{O}}_{Y}} {\mathcal{B}}_{0} \rightarrow {\mathcal{M}}_{i+1}
$$
is an epimorphism.
\end{definition}

\begin{example}\label{example.dave} (Quotients of Tensor Algebras).  \cite[Defintion 2.11, p.22]{qrs}  Let $\mathcal{F}$ be a coherent ${\mathcal{O}}_{X}$-bimodule.  The graded bimodule algebra $T(\mathcal{F})$, the {\bf tensor algebra of $\mathcal{F}$},\index{tensor algebra|textbf} is defined as follows: $(T(\mathcal{F}))_{i} = {\mathcal{F}}^{\otimes i}$ for $i \geq 1$ and $(T(\mathcal{F}))_{0} = d_{*}{\mathcal{O}}_{X}$,
$$
\mu_{ij}:{\mathcal{F}}^{\otimes i} \otimes {\mathcal{F}}^{\otimes j} \rightarrow {\mathcal{F}}^{\otimes i+j}
$$
is just the usual isomorphism of tensor products, and the unit map is the identity map.

If $\mathcal{I}$ is a graded ideal of $T(\mathcal{F})$,  then $T(\mathcal{F})/\mathcal{I}$ inherits a bimodule algebra structure from $T(\mathcal{F})$.  It is easily seen that $T(\mathcal{F})/\mathcal{I}$ is generated in degree one.
\end{example}

\begin{lemma} \label{lem.gen1}
Let $Y$ be a scheme and let $\mathcal{B}$ be an ${\mathbb{N}}$-graded ${\mathcal{O}}_{Y}$-bimodule algebra generated in degree one.  Let $\mathcal{M}$ be an ${\mathbb{N}}$-graded (right) $\mathcal{B}$-module generated in degree zero, with multiplication maps $\mu_{ij}^{\mathcal{M}}$.  If, for each $i \geq 0$, there exists an ${\mathcal{O}}_{Y}$-module morphism $\psi_{i}:{\mathcal{M}}_{i} \rightarrow {\mathcal{M}}_{i}$ such that the diagram
$$
\xymatrix{
{\mathcal{M}}_{i} \otimes {\mathcal{B}}_{1} \ar[r]^{\mu_{i,1}^{\mathcal{M}}} \ar[d]_{\psi_{i}\otimes {\mathcal{B}}_{1}} & {\mathcal{M}}_{i+1} \ar[d]^{\psi_{i+1}} \\
{\mathcal{M}}_{i} \otimes {\mathcal{B}}_{1} \ar[r]_{\mu_{i,1}^{\mathcal{M}}} & {\mathcal{M}}_{i+1}
}
$$
commutes, then the map $\psi:\mathcal{M} \rightarrow \mathcal{M}$ whose $i$th component is $\psi_{i}$, is a $\mathcal{B}$-module map.
\end{lemma}

\begin{proof}
To prove the lemma, we need to show that for all $j \geq 1$, the diagram
\begin{equation} \label{eqn.gen1}
\xymatrix{
{\mathcal{M}}_{i} \otimes {\mathcal{B}}_{j} \ar[r]^{\mu_{i,j}^{\mathcal{M}}} \ar[d]_{\psi_{i}\otimes {\mathcal{B}}_{j}} & {\mathcal{M}}_{i+j} \ar[d]^{\psi_{i+j}} \\
{\mathcal{M}}_{i} \otimes {\mathcal{B}}_{j} \ar[r]_{\mu_{i,j}^{\mathcal{M}}} & {\mathcal{M}}_{i+j}
}
\end{equation}
commutes.  We proceed by induction on $j$.  By hypothesis, (\ref{eqn.gen1}) commutes when $j=1$.  Suppose (\ref{eqn.gen1}) commutes when $j$ is replaced by $j-1$.  We must show (\ref{eqn.gen1}) commutes.  To prove this fact, we claim that all faces of the cube
$$
\xymatrix{
& {\mathcal{M}}_{i}\otimes {\mathcal{B}}_{1} \otimes {\mathcal{B}}_{j-1} \ar[rr] \ar[dd] & & {\mathcal{M}}_{i+1} \otimes {\mathcal{B}}_{j-1} \ar[dd] \\
{\mathcal{M}}_{i}\otimes {\mathcal{B}}_{1} \otimes {\mathcal{B}}_{j-1} \ar[dd] \ar[ur] \ar[rr] & & {\mathcal{M}}_{i+1} \otimes {\mathcal{B}}_{j-1} \ar[dd] \ar[ur] & \\
& {\mathcal{M}}_{i} \otimes {\mathcal{B}}_{j} \ar[rr] & & {\mathcal{M}}_{i+j} \\
{\mathcal{M}}_{i} \otimes {\mathcal{B}}_{j} \ar[rr] \ar[ur] & & {\mathcal{M}}_{i+j} \ar[ur] &
}
$$
whose maps into the page are tensor products of $\psi$ and whose other maps are tensors of multiplications, commute.  To prove this, we need only show that all faces except possibly the bottom commute, since by the hypothesis that $\mathcal{B}$ and ${\mathcal{M}}$ are generated in degree one, the left vertical maps are epis.  The top and right face commute by the induction hypothesis, the left face commutes by naturality properties of the tensor product, and the front and back face commute by the associativity of $\mu^{\mathcal{M}}$.  The result follows.
\end{proof}

\begin{lemma} \label{lem.smalldef}
Let $n$ be a positive integer.  Let $Y$ be a scheme and suppose $\mathcal{B}$ is an ${\mathbb{N}}$-graded ${\mathcal{O}}_{Y}$-bimodule algebra such that
\begin{itemize}
\item{the unit map $\upsilon:d_{*}{\mathcal{O}}_{Y} \rightarrow {\mathcal{B}}_{0}$ is an isomorphism and}
\item{the multiplication map $\mu_{1,j}$ is an isomorphism for $j+1 \leq n$.}
\end{itemize}
Suppose, further, that for $0 \leq i \leq n$, ${\mathcal{M}}_{i}$ is an ${\mathcal{O}}_{Y}$-module and, for $0 \leq i < n$ there exists a map $\nu_{i}:{\mathcal{M}}_{i} \otimes_{{\mathcal{O}}_{Y}} {\mathcal{B}}_{1} \rightarrow {\mathcal{M}}_{i+1}$.  Then there exists a unique multiplication map $\mu^{\mathcal{M}}$ such that $(\mu^{\mathcal{M}})_{i}=\nu_{i}$, making $\mathcal{M} = \oplus_{i=0}^{n}{\mathcal{M}}_{i}$ a truncated ${\mathcal{B}}$-module of length $n+1$.
\end{lemma}

\begin{proof}
We define $(\mu^{\mathcal{M}})_{i,j}$ inductively.  For $0 \leq i \leq n$, define $(\mu^{\mathcal{M}})_{i,0}$ as the composition
$$
\xymatrix{
{\mathcal{M}}_{i}\otimes_{{\mathcal{O}}_{X}} {\mathcal{B}}_{0} \ar[rr]^{{\mathcal{M}}_{i} \otimes_{{\mathcal{O}}_{X}} \upsilon^{-1}} & & {\mathcal{M}}_{i}\otimes_{{\mathcal{O}}_{X}} d_{*}{\mathcal{O}}_{X} \ar[r] & {\mathcal{M}}_{i},
}
$$
where the rightmost map is right multiplication.  Suppose, for $0 \leq i \leq n$ and $0 \leq j' < j$ with $i+j \leq n$, $(\mu^{\mathcal{M}})_{i,j'}$ is defined.  For $0 \leq i \leq n$ such that $i+j \leq n$, define $(\mu^{\mathcal{M}})_{i,j}$ as the composition
$$
\xymatrix{
{\mathcal{M}}_{i}\otimes_{{\mathcal{O}}_{X}} {\mathcal{B}}_{1} \otimes_{{\mathcal{O}}_{X}} {\mathcal{B}}_{j-1} \ar[rr]^{\nu_{i} \otimes_{{\mathcal{O}}_{X}} {\mathcal{B}}_{j-1}} & & {\mathcal{M}}_{i+1} \otimes_{{\mathcal{O}}_{X}} {\mathcal{B}}_{j-1} \ar[d]^{\mu_{i+1,j-1}^{\mathcal{M}}} \\
{\mathcal{M}}_{i}\otimes_{{\mathcal{O}}_{X}} {\mathcal{B}}_{j} \ar[u]^{{\mathcal{M}}_{i} \otimes_{{\mathcal{O}}_{X}} \mu_{1,j-1}^{-1}} & & {\mathcal{M}}_{i+j}
}
$$
For all combinations of $i$ and $j$ such that $\mu_{i,j}$ has not been defined above, let $(\mu^{\mathcal{M}})_{i,j}=0$.  From this definition it is obvious that $\mu^{\mathcal{M}}$ satisfies the unit axiom.  We show that $\mu^{\mathcal{M}}$ is associative by showing that, if $i,j,k$ are nonnegative integers, then the diagram
\begin{equation} \label{eqn.stacks}
\xymatrix{
{\mathcal{M}}_{i}\otimes {\mathcal{B}}_{j} \otimes {\mathcal{B}}_{k} \ar[d]_{\mu_{i,j}^{\mathcal{M}}} \ar[r]^{{\mathcal{M}}_{i} \otimes \mu_{j,k}} & {\mathcal{M}}_{i} \otimes {\mathcal{B}}_{j+k} \ar[d]^{\mu_{i,j+k}^{\mathcal{M}}} \\
{\mathcal{M}}_{i+j} \otimes {\mathcal{B}}_{k} \ar[r]_{\mu_{i+j,k}^{\mathcal{M}}} & {\mathcal{M}}_{i+j+k}
}
\end{equation}
commutes.  If $i+j+k > n$, then both routes of the diagram are $0$.  Otherwise, the diagram can be decomposed into a collection of pairs of diagrams
$$
\xymatrix{
{\mathcal{M}}_{i'}\otimes {\mathcal{B}}_{j'} \otimes {\mathcal{B}}_{k'} \ar[rr]^{{\mathcal{M}}_{i'}\otimes \mu_{j',k'}} \ar[d]_{{\mathcal{M}}_{i'} \otimes \mu_{1,j'-1}^{-1} \otimes {\mathcal{B}}_{k'}} & & {\mathcal{M}}_{i'}\otimes {\mathcal{B}}_{j'+ k'} \ar[d]^{{\mathcal{M}}_{i'}\otimes \mu_{1,j'+k'-1}^{-1}} \\
{\mathcal{M}}_{i'} \otimes {\mathcal{B}}_{1} \otimes {\mathcal{B}}_{j'-1} \otimes {\mathcal{B}}_{k'} \ar[d]_{\nu_{i'}\otimes {\mathcal{B}}_{j'-1} \otimes {\mathcal{B}}_{k'}} \ar[rr]^{{\mathcal{M}}_{i'}\otimes {\mathcal{B}}_{1} \otimes \mu_{j'-1,k'}} & & {\mathcal{M}}_{i'} \otimes {\mathcal{B}}_{1} \otimes {\mathcal{B}}_{j'+k'-1} \ar[d]^{\nu_{i'} \otimes {\mathcal{B}}_{j'-1+k'}} \\
{\mathcal{M}}_{i'+1} \otimes {\mathcal{B}}_{j'-1} \otimes {\mathcal{B}}_{k'} \ar[rr]_{{\mathcal{M}}_{i'+1} \otimes \mu_{j'-1,k'}} & & {\mathcal{M}}_{i'+1} \otimes {\mathcal{B}}_{j'+k'-1}
}
$$
stacked vertically.  The top square commutes by associativity of $\mu$, while the bottom square commutes by functoriality of the tensor product.  Thus, (\ref{eqn.stacks}) commutes.  The fact that $\mu^{\mathcal{M}}$ is unique follows from its construction.
\end{proof}\index{bimodule algebra|)}\index{bimodule algebra!graded|)}

\section{Lifting structures}
Let $X$ and $Y$ be schemes.  We show the associativity of the bimodule tensor product and the left and right scalar multiplication maps are indexed, relying heavily on results from the previous chapter.  We then use this fact to show that if $f:V \rightarrow U$ is a morphism in $\sf{S}$, $\tilde{f}=f \times \id_{X}:V \times X \rightarrow U \times X$, $\mathcal{B}$ is an ${\mathcal{O}}_{U \times X}$-bimodule algebra and $\mathcal{M}$ is a $\mathcal{B}$-module, then ${\tilde{f}}^{2*}\mathcal{B}$ inherits a bimodule algebra structure from $\mathcal{B}$ and ${\tilde{f}}^{*}\mathcal{M}$ inherits a ${\tilde{f}}^{2*}\mathcal{B}$-module structure from $\mathcal{M}$.

\begin{lemma} \label{example.bimods}
For any pair of schemes $X_{1}$ and $X_{2}$ over a base scheme $W$, let ${\sf{Bimod}}_{W}(X_{1}-X_{2})$ denote the category of ${\mathcal{O}}_{W}$-central ${\mathcal{O}}_{X_{1}}-{\mathcal{O}}_{X_{2}}$-bimodules and let ${\sf{bimod}}_{W}(X_{1}-X_{2})$ denote the category of coherent ${\mathcal{O}}_{W}$-central ${\mathcal{O}}_{X_{1}}-{\mathcal{O}}_{X_{2}}$-bimodules.  If $X$ and $Y$ are schemes, then the assignments $U \mapsto {\sf{Bimod}}_{U}(U \times X-U \times Y)$ and $U \mapsto {\sf{Coh}}U \times X$, where $U \in {\sf{S}}$, define $\sf{S}$-indexed categories $\mathfrak{Bimod(X-Y)}$ and ${\mathfrak{Coh X}}$.
\end{lemma}

\begin{proof}
We note that ${\sf{Bimod}}_{U}(U \times X-U \times Y)$ is a full subcategory ${\sf{Qcoh}}((U \times X)\times_{U}(U \times Y))$ and ${\sf{Coh}}U \times X$ is a full subcategory of ${\sf{Qcoh}}U \times X$.  Since the assignments $U \mapsto {\sf{Qcoh}}((U \times X)\times_{U}(U \times Y))$ and $U \mapsto {\sf{Qcoh}}U \times X$ define indexed categories as in Example \ref{example.spaces}, the assertion follows.
\end{proof}

\begin{lemma} \label{lem.important}
Let $I$ and $J$ be finite sets and let $\{A_{i}\}_{i \in I}$ and $\{B_{j}\}_{j \in J}$ be collections of schemes.  For all morphisms $f:V \rightarrow U$ in $\sf{S}$ let $f_{A} = \times_{i}(f \times \id_{A_{i}})$ and let $f_{B} = \times_{j}(f \times \id_{B_{j}})$.  Suppose, for all $U \in \sf{S}$, there exist morphisms
$$
g^{U}:\times_{i}(U \times A_{i})_{U} \rightarrow \times_{j}(U \times B_{j})_{U}
$$
such that, if $f:V \rightarrow U$ is a morphism in $\sf{S}$, then the diagram of schemes
\begin{equation} \label{eqn.cartwheelfirst}
\xymatrix{
\times_{i}(V \times A_{i})_{V} \ar[rr]^{f_{A}} \ar[d]_{g^{V}} & & \times_{i}(U \times A_{i})_{U} \ar[d]^{g^{U}} \\
\times_{j}(V \times B_{j})_{V} \ar[rr]_{f_{B}} & & \times_{j}(U \times B_{j})_{U}
}
\end{equation}
is a pullback.  Let ${\sf{A}}^{U}= {\sf{Qcoh }}\times_{i}(U \times A_{i})_{U}$ and ${\sf{B}}^{U} = {\sf{Qcoh }}\times_{j}(U \times B_{j})_{U}$.

\begin{enumerate}
\item{}
The assignment $f \mapsto f_{A}^{*},f_{B}^{*}$ makes $\mathfrak{A}$ and $\mathfrak{B}$ $\sf{S}$-indexed, naturally adjointed categories,

\item{}
If
$$
\Phi:g_{*}^{U}f_{A*} \Longrightarrow f_{B*}g_{*}^{V}
$$
is equality and
$$
\Lambda^{g_{*}}:f_{B}^{*}g_{*}^{U} \Longrightarrow g_{*}^{V}f_{A}^{*}
$$
is the dual of the square
\begin{equation} \label{eqn.superfirst}
\xymatrix{
{\sf{A}}^{V} \ar@<1ex>[r]^{f_{A*}} \ar[d]_{g_{*}^{V}} \ar@{=>}[dr]|{\Phi}\hole & {\sf{A}}^{U} \ar@<1ex>[l]^{f_{A}^{*}} \ar[d]^{g_{*}^{U}}\\
{\sf{B}}^{V} \ar@<1ex>[r]^{f_{B*}}  & {\sf{B}}^{U} \ar@<1ex>[l]^{f_{B}^{*}}.
}
\end{equation}
then the canonical natural isomorphism
$$
(\Lambda^{g^{*}})^{-1}:g^{V*}f_{B}^{*} \Longrightarrow f_{A}^{*}g^{U*}
$$
is dual to the square
\begin{equation} \label{eqn.supersecond}
\xymatrix{
{\sf{A}}^{U} \ar@<3ex>[r]^{g_{*}^{U}} \ar[d]_{f_{A}^{*}} \ar@{=>}[dr]|{\Lambda^{g_{*}}}\hole & {\sf{B}}^{U} \ar@<-1ex>[l]^{{g}^{U*}} \ar[d]^{f_{B}^{*}}\\
{\sf{A}}^{V} \ar@<-1ex>[r]^{g_{*}^{V}}  & {\sf{B}}^{V} \ar@<3ex>[l]^{g^{V*}}.
}
\end{equation}
and
\item{}
$$
\Lambda^{g^{*}}:f_{A}^{*}g^{U*} \Longrightarrow g^{V*}f_{B}^{*}
$$
and
$$
\Lambda^{g_{*}}:f_{B}^{*}g^{U}_{*} \Longrightarrow g^{V}_{*}f_{A}^{*}
$$
gives $g_{*}$ and $g^{*}$ the structure of weakly indexed functors.
\end{enumerate}
\end{lemma}

\begin{proof}
To prove 1, we note that $\mathfrak{A}$ and $\mathfrak{B}$ are indexed as in Example \ref{example.spaces} and $\mathfrak{A}$ and $\mathfrak{B}$ are adjointed as in Example \ref{example.schemes}.

We now prove 2.  Since (\ref{eqn.supersecond}) is the rotate of (\ref{eqn.superfirst}) (Proposition \ref{prop.rotation}), we know that if $\Psi$ is the dual of (\ref{eqn.supersecond}), the diagram
\begin{equation} \label{eqn.superthird}
\xymatrix{
g^{V*}f_{B}^{*}g_{*}^{U}f_{A*}f_{A}^{*}g^{U*} \ar@{=>}[rrr]^{=} & & & g^{V*}f_{B}^{*}f_{B*}g^{V}_{*}f_{A}^{*}g^{U*} \ar@{=>}[d] \\
g^{V*}f_{B}^{*} \ar@{=>}[u] \ar@{=>}[rrr]_{\Psi} & & & f_{A}^{*}g^{U*}
}
\end{equation}
whose verticals are various units and counits, commutes by Proposition \ref{prop.rotation}.  We show that the top route in (\ref{eqn.superthird}) equals $(\Lambda^{g^{*}})^{-1}$.  The diagram
$$
\xymatrix{
g^{V*}f_{B}^{*} \ar@{=>}[r] \ar@{=>}[d]_{=} & g^{V*}f_{B}^{*}g_{*}^{U}f_{A*}g^{V*}f_{B}^{*} \ar@{=>}[r] \ar@{=>}[d]_{\id*(\Lambda^{g^{*}})^{-1}}  & g^{V*}f_{B}^{*} \ar@{=>}[d]^{(\Lambda^{g^{*}})^{-1}} \\
g^{V*}f_{B}^{*} \ar@{=>}[r] & g^{V*}f_{B}^{*}g_{*}^{U}f_{A*}f_{A}^{*}g^{U*} \ar@{=>}[r] & f_{A}^{*}g^{U*}
}
$$
whose unlabeled arrows are various units and counits, commutes.  The right square commutes by naturality of $(\Lambda^{g^{*}})^{-1}$, while the left square commutes by the universal property of the unit of an adjoint pair.  Since the top route is $(\Lambda^{g^{*}})^{-1}$ while the bottom route is the top of (\ref{eqn.superthird}), $\Psi=(\Lambda^{g^{*}})^{-1}$.

To prove 3, let $\Psi_{f}$ denote the dual 2-cell of the square associated to (\ref{eqn.cartwheelfirst}).  For each $f:V \rightarrow U$, let
$$
\Phi_{f}:g^{U}_{*}f_{A*} \Longrightarrow f_{B*}g^{V}_{*}
$$
be equality and let $H^{U} = g_{*}^{U}$.  Suppose $g:W \rightarrow V$ is a morphism in $\sf{S}$.  Then it is easy to check that
$$
\Phi_{fg}=(f_{*}*\Phi_{g}) \circ (\Phi_{f}*g_{*}).
$$
Thus, the hypothesis of Proposition \ref{prop.indexedfunct} are satisfied when we define $H^{U}$ and $\Phi_{f}$ as above.  We conclude that $\Psi_{f}=\Lambda^{g_{*}}:f_{B}^{*}g^{U}_{*} \Longrightarrow g_{*}^{V}f_{A}^{*}$ gives $\{g^{U}_{*}\}$ an indexed structure.  The fact that $\Lambda^{g^{*}}$ gives $g^{*}$ an indexed structure can be checked locally.
\end{proof}
The following three results are immediate consequences of the Lemma.
\begin{corollary} \label{cor.im1}
Suppose $I$ is a finite set with $J \subset I$ and $\{X_{i}\}_{i \in I}$ is a collection of schemes.  For all $U \in \sf{S}$, suppose $g^{U}$ is the projection
$$
{\pr_{J}}^{U}:\times_{i \in I}(U \times X_{i})_{U} \rightarrow \times_{j \in J}(U \times X_{j})_{U},
$$
${\sf{A}}^{U}={\sf{Qcoh}}\times_{i}(U \times X_{i})$, and ${\sf{B}}^{U}={\sf{Qcoh}}\times_{j}(U \times X_{j})$.  Then the assignments $U \mapsto {\pr_{J}}^{U}_{*}$ and $U \mapsto {\pr_{J}}^{U*}$ define (weakly) indexed functors
$$
\mathfrak{\pr_{J*}}:\mathfrak{A} \rightarrow \mathfrak{B}
$$
and
$$
\mathfrak{\pr_{J}^{*}}:\mathfrak{B} \rightarrow \mathfrak{C}
$$
as in Lemma \ref{lem.important}.
\end{corollary}

\begin{corollary} \label{cor.im2}
Suppose $X$ is a scheme.  For all $U \in \sf{S}$, suppose $g^{U}$ is the diagonal
$$
d^{U}:U \times X \rightarrow (U \times X)_{U}^{2},
$$
${\sf{A}}^{U}={\sf{Qcoh}}U \times X$ and ${\sf{B}}^{U}={\sf{Qcoh}}(U \times X)_{U}^{2}$.  Then the assignments $U \mapsto d^{U}_{*}$ and $U \mapsto d^{U*}$ define indexed functors
$$
\mathfrak{d_{*}}:\mathfrak{A} \rightarrow \mathfrak{B}
$$
and
$$
\mathfrak{d^{*}}:\mathfrak{B} \rightarrow \mathfrak{A}
$$
as in Lemma \ref{lem.important}.
\end{corollary}

\begin{corollary} \label{cor.im3}
Suppose $X$ and $Y$ are schemes.  For all $U \in \sf{S}$, suppose $g^{U}$ is the morphism
$$
e^{U}=(d \times \id_{Y})^{U}:(U \times X)\times_{U}(U \times Y) \rightarrow (U \times X)_{U}^{2} \times_{U}(U \times Y),
$$
${\sf{A}}^{U}={\sf{Qcoh}}((U \times X)\times(U \times Y))$ and ${\sf{B}}^{U}={\sf{Qcoh}}((U \times X)_{U}^{2}\times_{U}(U \times Y))$.  Then the assignments $U \mapsto e^{U}_{*}$ and $U \mapsto e^{U*}$ define indexed functors
$$
\mathfrak{e_{*}}:\mathfrak{A} \rightarrow \mathfrak{B}
$$
and
$$
\mathfrak{e^{*}}:\mathfrak{B} \rightarrow \mathfrak{A}
$$
as in Lemma \ref{lem.important}.
\end{corollary}

\begin{lemma} \label{lem.isommonoid}
Suppose $X$, $Y$ and $Z$ are schemes.  For all $U \in \sf{S}$, let
$$
\pr_{i,j}^{U}:(U \times X) \times_{U} (U \times Y) \times_{U} (U \times Z) \rightarrow U \times X, U \times Y, U \times Z
$$
and
$$
\pr_{i}^{U}:(U \times X) \times_{U} (U \times Y) \rightarrow U \times X, U \times Y
$$
be projection maps, let $(- \otimes_{{\mathcal{O}}_{Y}}-)^{U}$ denote the functor
$$
\pr_{13*}^{U}(\pr_{12}^{*U}- \otimes \pr_{23}^{U*}-):{\sf{Bimod}}_{U}(U \times X-U \times Y) \times {\sf{Bimod}}_{U}(U \times Y-U \times Z) \rightarrow
$$
$$
{\sf{Bimod}}_{U}(U \times X - U \times Z)
$$
and let $(- \otimes_{{\mathcal{O}}_{X}}-)^{U}$ denote the functor
$$
\pr_{2*}^{U}(\pr_{1}^{*U}- \otimes -):{\sf{Coh}}_{U}(U \times X) \times {\sf{bimod}}_{U}(U \times X-U \times Y) \rightarrow {\sf{Coh}}(U \times Y).
$$
Then the assignments $U \mapsto (- \otimes_{{\mathcal{O}}_{Y}}-)^{U}$ and $U \mapsto (- \otimes_{{\mathcal{O}}_{X}}-)^{U}$ define indexed functors
$$
- \otimes_{{\mathcal{O}}_{\mathfrak{Y}}}-:\mathfrak{Bimod(X-Y)} \times \mathfrak{Bimod(Y-Z)} \rightarrow \mathfrak{Bimod(X-Z)}
$$
and
$$
- \otimes_{{\mathcal{O}}_{\mathfrak{X}}}-:\mathfrak{Coh X} \times \mathfrak{bimod(X-Y)} \rightarrow \mathfrak{Coh Y}.
$$
\end{lemma}\index{tensor product!of bimodules!is indexed}\index{tensor product!of a module and a bimodule!is indexed}

\begin{proof}
We prove the first assertion.  The second assertion is proven in a similar manner.  We need only note that, by Corollary \ref{cor.im1} and Example \ref{example.indexcheck}, $- \otimes_{{\mathcal{O}}_{\mathfrak{Y}}}-$ is the composition of indexed functors so that $- \otimes_{{\mathcal{O}}_{\mathfrak{Y}}}-$ inherits an indexed structure.  Specifically, suppose $f:V \rightarrow U$ is a morphism in $\sf{S}$.  Let
$$
f_{XY}=(f \times \id_{X}) \times (f\times \id_{Y}):(V \times X) \times_{V} (V \times Y) \rightarrow (U \times X) \times_{U}(U \times Y)
$$
and define $f_{YZ}$, $f_{XZ}$ and $f_{XYZ}$ similarly.  There is an isomorphism
\begin{equation} \label{eqn.deltanew}
\Delta:f_{XZ}^{*}(- \otimes_{{\mathcal{O}}_{Y}}-)^{U} \Longrightarrow (f_{XY}^{*}- \otimes_{{\mathcal{O}}_{Y}}f_{YZ}^{*}-)
\end{equation}
as follows:
\begin{align} \label{eqn.delta}
f_{XZ}^{*}(- \otimes_{{\mathcal{O}}_{Y}}-)^{U} & = f_{XZ}^{*}\pr_{13*}^{U}(\pr_{12}^{*U}- \otimes \pr_{23}^{U*}-) \\
& \cong \pr_{13*}^{V}f_{XYZ}^{*}(\pr_{12}^{*U}- \otimes \pr_{23}^{U*}-) \nonumber \\
& \cong \pr_{13*}^{V}(f_{XYZ}^{*}\pr_{12}^{*U}- \otimes f_{XYZ}^{*}\pr_{23}^{U*}-) \nonumber \\
& \cong \pr_{13*}^{V}(\pr_{12}^{V*}f_{XY}^{*}- \otimes \pr_{23}^{V*}f_{YZ}^{*}-) \nonumber \\
& = (f_{XY}^{*}- \otimes_{{\mathcal{O}}_{Y}}f_{YZ}^{*}-)^{V} \nonumber.
\end{align}
\end{proof}

\begin{corollary} \label{cor.projy}
Let $I$ and $J$ be finite sets and let $\{A_{i}\}_{i \in I}$ and $\{B_{j}\}_{j \in J}$ be collections of schemes.  Suppose, for all $U \in \sf{S}$, there exist morphisms
$$
g^{U}:\times_{i}(U \times A_{i})_{U} \rightarrow \times_{j}(U \times B_{j})_{U}
$$
such that, if $f:V \rightarrow U$ is a morphism in $\sf{S}$, then the diagram of schemes
$$
\xymatrix{
\times_{i}(V \times A_{i})_{V} \ar[rr]^{ \times_{i}(f \times \id_{A_{i}})} \ar[d]_{g^{V}} & & \times_{i}(U \times A_{i})_{U} \ar[d]^{g^{U}} \\
\times_{j}(V \times B_{j})_{V} \ar[rr]_{\times_{j}(f \times \id_{B_{j}})} & & \times_{j}(U \times B_{j})_{U}
}
$$
is a pullback.  Then the projection formula (\ref{prop.canon} (1)) \index{projection formula!is weakly indexed}
\begin{equation} \label{eqn.projy}
\Psi^{U}:g^{U}_{*}-\otimes- \Longrightarrow g_{*}^{U}(- \otimes g^{U*}-)
\end{equation}
which is a natural transformation between functors from
$$
{\sf{Qcoh}}(\times_{i}(U \times A_{i})_{U})  \times {\sf{Qcoh}}(\times_{j}(U \times B_{j})_{U})
$$
to
$$
{\sf{Qcoh}}\times_{j}(U \times B_{j})_{U}
$$
defines an indexed natural transformation between weakly indexed functors.
\end{corollary}

\begin{proof}
We note that the assignments $U \mapsto g_{*}^{U}$ and $U \mapsto g^{U*}$ define weakly indexed natural transformations by Lemma \ref{lem.important}.  Now, $\Psi^{U}$ is given by the composition:
$$
\xymatrix{
g^{U}_{*}- \otimes - \ar@{=>}[r] & g^{U}_{*}g^{U*}(g^{U}_{*}- \otimes -) \ar@{=>}[r] & g^{U}_{*}(g^{U*}g^{U}_{*}- \otimes g^{U*}-) \ar@{=>}[r] &
}
$$
$$
g^{U}_{*}(- \otimes g^{U*}-)
$$
where the first and last maps are products of id's, units, and counits, and the second map is induced by the commutativity of $g^{U*}$ and the tensor product.  Since the commutativity of $g^{U*}$ and the tensor product is indexed (this can be checked locally), it suffices to show that, if $(\eta^{U}, \epsilon^{U})$ are the unit and counit of $(g^{U*},g_{*}^{U})$, the collections $\{\eta^{U}\}$ and $\{\epsilon^{U}\}$ define indexed natural transformations.  The square
$$
\xymatrix{
{\sf{Qcoh}}\times_{i}(U \times A_{i})_{U} \ar@<1ex>[r]^{g^{U}_{*}} \ar[d]_{\operatorname{id}} \ar@{=>}[dr]|{\operatorname{id}}\hole & {\sf{Qcoh}}\times_{j}(U \times B_{j})_{U} \ar@<1ex>[l]^{{g^{U}}^{*}} \ar[d]^{\operatorname{id}}\\
{\sf{Qcoh}}\times_{i}(U \times A_{i})_{U} \ar@<1ex>[r]^{g^{U}_{*}}  & {\sf{Qcoh}}\times_{j}(U \times B_{j})_{U} \ar@<1ex>[l]^{g^{U*}}.
}
$$
defines a square of indexed categories by Lemma \ref{lem.applicat}.  By the remark following Definition \ref{def.squareof}, $\{\eta^{U}\}$ and $\{\epsilon^{U}\}$ define indexed natural transformations.  Thus, $\Psi^{U}$ is a composition of indexed functors so it is indexed.
\end{proof}
In \cite{lipman}, Lipman describes the importance of
compatibility\index{compatibilities} results in eliminating the
noetherian hypothesis from the statement of Grothendieck duality.
He then gives examples of diagrams \cite[2.2.1, 2.2.2]{lipman}
which arise in the theory of Grothendieck duality, alluding to the
complications involved in showing they commute.  The fact that
\cite[2.2.1]{lipman} commutes (in the category of quasi-coherent
modules over schemes) follows immediately from the fact that the
projection formula is a weakly indexed natural transformation
(which follows from a mild generalization of Corollary
\ref{cor.projy}).  It seems likely that \cite[2.2.2]{lipman}
commutes for a similar reason.  In fact, since the composition of
indexed natural transformations is indexed, one can build up,
using only diagrams of the form (\ref{eqn.lippy}) on page
\pageref{eqn.lippy}, a large number of intimidating diagrams which
must commute.  This suggests that it would be fruitful to pursue a
coherence theorem of the form suggested by Lipman in the context
of squares of indexed categories\index{square of indexed
categories}.

\begin{corollary} \label{cor.hithere}
Keep the notation as in Lemma \ref{lem.important}.  Let $K$ and $L$ be finite sets and let $\{C_{k}\}_{k \in K}$ and $\{D_{l}\}_{l \in L}$ be collections of schemes.  For all $f:V \rightarrow U$ in $\sf{S}$, let $f_{C}=\times_{k}(f \times \id_{C_{k}})$ and let $f_{D}=\times_{l}(f \times \id_{D_{l}})$.  Suppose, for all $U \in \sf{S}$, there exist morphisms
$$
h^{U}:\times_{k}(U \times C_{k})_{U} \rightarrow \times_{l}(U \times D_{l})_{U}
$$
such that, if $f:V \rightarrow U$ is a morphism in $\sf{S}$, then the diagram of schemes

$$
\xymatrix{
\times_{k}(V \times C_{k})_{V} \ar[rr]^{f_{C}} \ar[d]_{h^{V}} & & \times_{k}(U \times C_{k})_{U} \ar[d]^{h^{U}} \\
\times_{l}(V \times D_{l})_{V} \ar[rr]_{f_{D}} & & \times_{l}(U \times D_{l})_{U}
}
$$
is a pullback.  Let ${\sf{C}}^{U}= {\sf{Qcoh }}\times_{k}(U \times C_{k})_{U}$ and ${\sf{D}}^{U} = {\sf{Qcoh }}\times_{l}(U \times D_{l})_{U}$.

Suppose
$$
\Phi:g_{*}^{U}f_{A*} \Longrightarrow f_{B*}g_{*}^{V}
$$
is equality and
$$
\Lambda^{g_{*}}:f_{B}^{*}g_{*}^{U} \Longrightarrow g_{*}^{V}f_{A}^{*}
$$
is the dual of the square
$$
\xymatrix{
{\sf{A}}^{V} \ar@<1ex>[r]^{f_{A*}} \ar[d]_{g_{*}^{V}} \ar@{=>}[dr]|{\Phi}\hole & {\sf{A}}^{U} \ar@<1ex>[l]^{f_{A}^{*}} \ar[d]^{g_{*}^{U}}\\
{\sf{B}}^{V} \ar@<1ex>[r]^{f_{B*}}  & {\sf{B}}^{U} \ar@<1ex>[l]^{f_{B}^{*}}.
}
$$
Suppose $\Phi':h_{*}^{U}f_{C*} \Longrightarrow f_{D*}h_{*}^{V}$ is equality and
$$
\Lambda^{h_{*}}:f_{D}^{*}h_{*}^{U} \Longrightarrow h_{*}^{V}f_{C}^{*}
$$
is dual to $\Phi'$ as above.  If $\Lambda^{g^{*}}:f_{A}^{*}g^{U*} \Longrightarrow  g^{V*}f_{B}^{*}$ and $\Lambda^{h^{*}}:f_{C}^{*}h^{U*} \Longrightarrow h^{V*}f_{D}^{*}$ are the usual isomorphisms then
$$
\Lambda^{h^{*}}:f_{C}^{*}h^{U*} \Longrightarrow h^{V*}f_{D}^{*}
$$
and
$$
\Lambda^{h_{*}}:f_{D}^{*}h^{U}_{*} \Longrightarrow h^{V}_{*}f_{C}^{*}
$$
gives $h_{*}$ and $h^{*}$ the structure of weakly indexed functors.  Furthermore, if the sets $\{H^{U} \}$ and $\{I^{U}\}$ define indexed functors and
\begin{equation} \label{eqn.appl}
\xymatrix{
{\sf{A}}^{U} \ar@<3ex>[r]^{{g_{*}}^{U}} \ar[d]_{H^{U}} \ar@{=>}[dr]|{\Phi^{U}}\hole & {\sf{B}}^{U} \ar@<-1ex>[l]^{{g}^{U*}} \ar[d]^{I^{U}}\\
{\sf{C}}^{U} \ar@<-1ex>[r]^{{h_{*}}^{U}}  & {\sf{D}}^{U} \ar@<3ex>[l]^{{h^{*}}^{U}}
}
\end{equation}
is a square, then $\{ \Phi^{U} \}$ is indexed if and only if $\{\Psi^{U}\}$ is indexed.
\end{corollary}

\begin{proof}
The first two assertions follow from Lemma \ref{lem.important}.  By Lemma \ref{lem.important} and Lemma \ref{lem.applicat}, the data (\ref{eqn.appl}) defines a square of indexed categories.  The result now follows from Corollary \ref{cor.reduce}.
\end{proof}

\begin{corollary} \label{cor.indexcan}
Let $K$ be finite a set, $I, L \subset K$, $J \subset I, L$ and let $\{X_{k}\}_{k \in K}$ be a collection of schemes.  Suppose, for all $U \in \sf{S}$
$$
{\pr^{I}_{J}}^{U}:\times_{i \in I}(U \times X_{i})_{U} \rightarrow \times_{j \in J}(U \times X_{j})_{U}
$$
is projection, ${\pr^{L}_{J}}^{U}$, ${\pr^{K}_{I}}^{U}$ and ${\pr^{K}_{L}}^{U}$ are defined similarly, and
\begin{equation} \label{eqn.localya}
\xymatrix{
\times_{k}(U \times X_{k})_{U} \ar[rr]^{{\pr^{K}_{I}}^{U}} \ar[d]_{{\pr^{K}_{L}}^{U}} & & \times_{i}(U \times X_{i})_{U} \ar[d]^{{\pr^{I}_{J}}^{U}} \\
\times_{l}(U \times X_{l})_{U} \ar[rr]_{{\pr^{L}_{J}}^{U}} & & \times_{j}(U \times X_{j})_{U}
}
\end{equation}
is a pullback.  Then the dual to (\ref{eqn.localya}),
\begin{equation} \label{eqn.indexcan}
\Psi^{U}:{\pr^{L}_{J}}^{U*}{\pr^{I}_{J}}^{U}_{*} \Longrightarrow {\pr^{L}_{K}}^{U}_{*}{\pr^{K}_{I}}^{U*}
\end{equation}
which is a natural transformation between functors from
$$
{\sf{Qcoh}}\times_{i}(U \times X_{i})_{U}
$$
to
$$
{\sf{Qcoh}}\times_{l}(U \times X_{l})_{U},
$$
defines an indexed natural transformation between weakly indexed functors.
\end{corollary}

\begin{proof}
For $U \in \sf{S}$, let ${\sf{A}}^{U}={\sf{rla}}_{{\pr^{I}_{J}}^{U}}$, ${\sf{B}}^{U} = {\sf{Qcoh }}\times_{j}(U \times X_{j})_{U}$, ${\sf{C}}^{U}={\sf{rla}}_{{\pr^{K}_{L}}^{U}}$ and ${\sf{D}}^{U} = {\sf{Qcoh }}\times_{l}(U \times X_{l})_{U}$.  If $g^{U}={\pr^{I}_{J}}^{U}$, $h^{U}={\pr^{K}_{L}}^{U}$, $H^{U}={\pr^{K}_{I}}^{U*}$ and $I^{U}={\pr^{L}_{J}}^{U*}$ then, by Example \ref{example.indexcheck}, Corollary \ref{cor.im1} and Lemma \ref{lem.important}, the hypothesis of Corollary \ref{cor.hithere} are satisfied.  Thus, to show the data $\{\Psi^{U} \}$ defines an indexed natural transformation, it suffices to show that, if
$$
\Phi^{U}:{\pr^{I}_{J}}^{U}_{*}{\pr^{K}_{I}}^{U}_{*}={\pr^{L}_{J}}^{U}_{*}{\pr^{K}_{L}}^{U}_{*},
$$
then the data $\{\Phi^{U}\}$ defines a natural transformation.  This is trivial, and the assertion follows.
\end{proof}

\begin{proposition} \label{prop.assindex}
Suppose $X$, $Y$ and $Z$ are schemes.  Then associativity of the tensor product (Proposition \ref{prop.tensor}) defines indexed natural transformations\index{associativity of tensor products is indexed}
$$
(-\otimes_{{\mathcal{O}}_{\mathfrak{Y}}} \otimes -)\otimes_{{\mathcal{O}}_{\mathfrak{Z}}}- \Longrightarrow -\otimes_{{\mathcal{O}}_{\mathfrak{Y}}} \otimes (-\otimes_{{\mathcal{O}}_{\mathfrak{Z}}}-)
$$
between indexed functors from
$$
\mathfrak{Bimod(X-Y)} \times \mathfrak{Bimod(Y-Z)} \times \mathfrak{Bimod(Z-T)}
$$
to
$$
\mathfrak{Bimod(X-T)}
$$
and
$$
(-\otimes_{{\mathcal{O}}_{\mathfrak{X}}} \otimes -)\otimes_{{\mathcal{O}}_{\mathfrak{Y}}}- \Longrightarrow -\otimes_{{\mathcal{O}}_{\mathfrak{X}}} \otimes (-\otimes_{{\mathcal{O}}_{\mathfrak{Y}}}-)
$$
between indexed functors from
$$
\mathfrak{Coh X} \times \mathfrak{bimod(X-Y)} \times \mathfrak{bimod(Y-Z)}
$$
to
$$
\mathfrak{Coh Z}.
$$
\end{proposition}

\begin{proof}
Each functor above is a composition of indexed functors, hence inherits an indexed structure.  Since the associativity natural transformation is a composition of natural transformations of the form (\ref{eqn.projy}) (on page \pageref{eqn.projy}) and (\ref{eqn.indexcan}), the result follows from Corollary \ref{cor.projy} and Corollary \ref{cor.indexcan}
\end{proof}

\begin{proposition} \label{prop.unitindex}
Keep the notation as in Corollary \ref{cor.im2} and Lemma \ref{lem.isommonoid}.  For all $U \in \sf{S}$, the assignment
$$
U \mapsto \pr_{13*}^{U}(\pr_{12}^{U*}-\otimes \pr_{23}^{U*}d_{*}^{U}{\mathcal{O}}_{U \times Y})
$$
defines an indexed functor
$$
- \otimes_{{\mathcal{O}}_{\mathfrak{Y}}}\mathfrak{d_{*}}{\mathcal{O}}_{\mathfrak{Y}}:\mathfrak{Bimod(X-Y)} \rightarrow \mathfrak{Bimod(X-Y)}.
$$
Similarly, there are indexed functors
$$
\mathfrak{d_{*}}{\mathcal{O}}_{\mathfrak{X}} \otimes_{{\mathcal{O}}_{\mathfrak{X}}}-:\mathfrak{Bimod(X-Y)} \rightarrow \mathfrak{Bimod(X-Y)}
$$
and
$$
- \otimes_{{\mathcal{O}}_{\mathfrak{X}}}\mathfrak{d_{*}}{\mathcal{O}}_{\mathfrak{X}}:\mathfrak{Qcoh X} \rightarrow \mathfrak{Qcoh X}.
$$
Furthermore, the left and right linear multiplication maps (Proposition \ref{prop.scalar}) define indexed natural isomorphisms\index{linear multiplication map!is indexed}
$$
- \otimes_{{\mathcal{O}}_{\mathfrak{Y}}}\mathfrak{d_{*}}{\mathcal{O}}_{\mathfrak{Y}} \Longrightarrow \operatorname{id}_{\mathfrak{Bimod(X-Y)}},
$$
$$
\mathfrak{d_{*}}{\mathcal{O}}_{\mathfrak{X}} \otimes_{{\mathcal{O}}_{\mathfrak{X}}} \Longrightarrow \operatorname{id}_{\mathfrak{Bimod(X-Y)}},
$$
and
$$
- \otimes_{{\mathcal{O}}_{\mathfrak{X}}}\mathfrak{d_{*}}{\mathcal{O}}_{\mathfrak{X}} \Longrightarrow \operatorname{id}_{\mathfrak{Qcoh X}}.
$$
\end{proposition}

\begin{proof}
The fact that the first three assignments define indexed functors follows from Corollary \ref{cor.im2}, Lemma \ref{lem.isommonoid} and Example \ref{example.indexcheck}.  The fact that the given families of natural isomorphisms are indexed follows from an argument similar to that given in Proposition \ref{prop.assindex} so we omit it.
\end{proof}

\begin{theorem} \label{theorem.lift}
Let $f:V \rightarrow U$ be a morphism in $\sf{S}$, let $\tilde{f}=f \times \id_{X}:V \times X \rightarrow U \times X$ and suppose the left composite of the isomorphism
$$
\lambda:\tilde{f}^{2*}d^{U}_{*}{\mathcal{O}}_{U \times X} \rightarrow d^{V}_{*}\tilde{f}^{*}{\mathcal{O}}_{U \times X} \rightarrow d^{V}_{*}{\mathcal{O}}_{V \times X}
$$
comes from the indexed structure of $\mathfrak{d}$.  Suppose the data $(\mathcal{B},\mu,\upsilon)$ consisting of an ${{\mathcal{O}}_{U\times X}}$-bimodule $\mathcal{B}$, a morphism $\mu:\mathcal{B} \otimes \mathcal{B} \rightarrow \mathcal{B}$ and a morphism $d_{*}^{U}{\mathcal{O}}_{U\times X} \rightarrow \mathcal{B}$ defines a bimodule algebra.  Then the data $(\tilde{f}^{2*}\mathcal{B}, \tilde{f}^{2*}\mu \delta^{-1}, \tilde{f}^{2*}\upsilon \lambda^{-1})$, where $\delta$ is induced by the natural isomorphism (\ref{eqn.deltanew}) defines an ${\mathcal{O}}_{V \times X}$-bimodule algebra.\index{bimodule algebra}\index{bimodule algebra!graded}

Furthermore, if the data $(\mathcal{M}, \nu)$ consisting of an ${{\mathcal{O}}_{U \times X}}$-module $\mathcal{M}$ and a morphism $\nu:\mathcal{M} \otimes \mathcal{B} \rightarrow \mathcal{M}$ defines a right $\mathcal{B}$-module, then the data $(\tilde{f}^{*}\mathcal{M}, {\tilde{f}}^{*}\nu \gamma^{-1})$, where $\gamma$ is induced by the indexed structure of the appropriate tensor product, defines an $\tilde{f}^{2*}\mathcal{B}$-module.\index{B module@$\mathcal{B}$-module}\index{B module@$\mathcal{B}$-module!graded}

\end{theorem}

\begin{proof}
To show that a bimodule algebra lifts, one needs to check the associativity and left and right unit axioms.  To show that a module lifts, one needs to check the associativity and right unit axiom.  A straightforward verification shows that these axioms hold by Proposition \ref{prop.assindex} and Proposition \ref{prop.unitindex}.
\end{proof}
The above construction can obviously be graded.

\section{The definition of $\Gamma_{n}$}
Assume $X$ is an $S$-scheme and $\mathcal{B}$ is an ${\mathcal{O}}_{X^{2}}$-module.  Suppose $U \in \sf{S}$ has structure map $f:U \rightarrow S$, let $\tilde{f}=f \times \id_{X}:U \times X \rightarrow S \times X$ and let ${\mathcal{B}}^{U}=\tilde{f}^{2*}\mathcal{B}$.  Before reading further in this chapter, the reader is advised to review the last two sentences of Section 1.6.
\begin{definition} \label{def.family}
Let $\mathcal{B}$ be a graded ${{\mathcal{O}}_{X}}$-bimodule algebra and let $U$ be an affine scheme.  A {\bf family of $\mathcal{B}$-point modules \index{family!of point modules}parameterized by $\mathbf{U}$} or a {\bf $\mathbf{U}$-family over $\mathcal{B}$} is a graded ${\mathcal{B}}^{U}$-module ${\mathcal{M}}_{0} \oplus {\mathcal{M}}_{1} \oplus \cdots $ generated in degree zero such that, for each $i\geq 0$ there exists a map
$$
q_{i}: U \rightarrow X
$$
and an invertible ${\mathcal{O}}_{U}$-module ${\mathcal{L}}_{i}$ with ${\mathcal{L}}_{0} \cong {\mathcal{O}}_{U}$ and
$$
{\mathcal{M}}_{i} \cong (\id_{U} \times q_{i})_{*} {\mathcal{L}}_{i}.
$$
A {\bf family of truncated $\mathcal{B}$-point modules of length $\mathbf{n+1}$ parameterized by $\mathbf{U}$}, \index{family!of truncated point modules}or a {\bf truncated $\mathbf{U}$-family of length $\mathbf{n+1}$} is  a graded ${\mathcal{B}}^{U}$-module ${\mathcal{M}}_{0} \oplus {\mathcal{M}}_{1} \oplus \cdots $generated in degree zero, such that for each $i\geq 0$ there exists a map
$$
q_{i}: U \rightarrow X
$$
and an invertible ${{\mathcal{O}}_{U}}$-module ${\mathcal{L}}_{i}$ with ${\mathcal{L}}_{0} \cong {\mathcal{O}}_{U}$,
$$
{\mathcal{M}}_{i} \cong (\id_{U} \times q_{i})_{*} {\mathcal{L}}_{i}
$$
for $i \leq n$ and ${\mathcal{M}}_{i}=0$ for $i>n$.
\end{definition}
{\it Remark.}  If $f:V \rightarrow U$ is a morphism in $\sf{S}$, $\tilde{f}=f \times \id_{X}:V \times X \rightarrow U \times X$ and ${\mathcal{B}}$ is an ${{\mathcal{O}}_{X}}$-bimodule algebra, then the bimodule algebras $\tilde{f}^{2*}{\mathcal{B}}^{U}$ and ${\mathcal{B}}^{V}$ are canonically isomorphic.  Thus, there is a canonical bijection between isomorphism classes of $\tilde{f}^{2*}{\mathcal{B}}^{U}$-modules and isomorphism classes of ${\mathcal{B}}^{V}$-modules.  We will make use of this fact without further mention.

In the following proposition, we note that the $q_{i}$'s appearing in Definition \ref{def.family} are unique.  That is, if two $U$-families, $\mathcal{M}$ and $\mathcal{M}'$ are isomorphic, then the map $q_{i}$ in the definition of ${\mathcal{M}}_{i}$ {\it equals} the map $q_{i}'$ in the definition of ${\mathcal{M}'}_{i}$.  To prove this result, we need to know the following fact about separated schemes.

\begin{lemma} \label{lem.gsep} \cite[ex. 4.8, p.107]{alggeo}

Let $X$ be a separated scheme, and suppose $U$ is affine.  Then any map
$$
q:U \rightarrow X
$$
is affine.

\end{lemma}

\begin{lemma} \label{lem.sep}

Suppose $R$ is a commutative ring, let $A$ and $B$ be commutative $R$-algebras and suppose $g,h:B \rightarrow A$ are algebra morphisms.  Suppose $M=Ae$ and $N=Af$ are free $A$-modules of rank 1.  Let $M$ and $N$ be $A\otimes B$-modules by letting $a\otimes b \cdot m = ag(b)m$ and letting $a \otimes b \cdot n = ah(b)n$ respectively.  If
$$
\psi:M \rightarrow N
$$
is an $A\otimes B$-module isomorphism, then $g=h$.
\end{lemma}
\begin{proof}
Since $\psi$ is an $A\otimes B$-module isomorphism, it is an $A$-module isomorphism.  Thus, $\psi(e)=u f$ for some unit $u \in A$.  Let $b \in B$.  Then
$$
\psi((1\otimes b) \cdot e)=\psi(g(b)e)=u g(b)f.
$$
But
$$
\psi((1\otimes b) \cdot e)= (1 \otimes b)\cdot \psi(e)=(1 \otimes b) \cdot u f=h(b)u f.
$$
Comparing these results, we find that $h(b)=g(b)$ as desired.
\end{proof}
\begin{proposition} \label{prop.uniqueness}
Let $\mathcal{M}$ and $\mathcal{N}$ be invertible ${\mathcal{O}}_{U}$-modules and suppose

$$
g,h:U \rightarrow X.
$$
If $\psi:(\id_{U} \times g)_{*}\mathcal{M} \rightarrow (\id_{U} \times h)_{*}{\mathcal{N}}$ is an isomorphism, then $g=h$.
\end{proposition}
\begin{proof}
Since $(\id_{U} \times g)_{*}\mathcal{M} \cong (\id_{U} \times h)_{*}{\mathcal{N}}$, their supports are equal.  In particular, $g=h$ point-wise.

We now show that the sheaf components of the maps $g$ and $h$  are equal.  By taking the appropriate open cover of Supp $(\id_{U} \times g)_{*}\mathcal{M} \subset U \times X$, we will reduce our result to the previous lemma.  Let $U' \subset U$ be an open affine set on which both $\mathcal{M}$ and ${\mathcal{N}}$ are free.  Let Spec $B \subset X$ be open affine.  Since $X$ is separated, by Lemma $\ref{lem.sep}$, both $g$ and $h$ are affine.  Thus, $(\id_{U} \times g)^{-1}(U' \times \operatorname{Spec }B) = (\id_{U} \times h)^{-1}(U' \times \operatorname{Spec }B)$ is an affine open set, Spec $A$, of $U$ on which both $\mathcal{M}$ and ${\mathcal{N}}$ are free.  Then $\psi$ gives an $A\otimes B$-module isomorphism

$$
\psi(\operatorname{Spec }A \otimes B): \mathcal{M}(\operatorname{Spec }A)_{A \otimes B} \rightarrow {\mathcal{N}}(\operatorname{Spec }A)_{A \otimes B}.
$$
and the result follows from the previous lemma.
\end{proof}
Suppose $f:V \rightarrow U$ and $X$ is a scheme, and define $\tilde{f}=f \times \id_{X}:V \times X \rightarrow U \times X$.
\begin{lemma} \label{lem.anothersquare}
Suppose $f:V \rightarrow U$ is a morphism of affine schemes and $q:U \rightarrow X$ is a morphism.  Then the dual 2-cells
$$
{\tilde{f}}^{*}(\id_{U}\times q)_{*} \Longrightarrow (\id_{V}\times qf)_{*}f^{*}
$$
and
$$
{\tilde{f}}^{2*}(\id_{U} \times q)^{2}_{*} \Longrightarrow (\id_{V}\times qf)^{2}_{*}f^{*}
$$
induced by the diagrams
$$
\xymatrix{
V \ar[r]^{f} \ar[d]_{\id_{V} \times qf} & U \ar[d]^{\id_{U} \times q} \\
V \times X \ar[r]_{f \times \id_{X}} & U \times X
}
$$
and
$$
\xymatrix{
V \ar[r]^{f} \ar[d]_{(\id_{V} \times qf)^{2}} & U \ar[d]^{(\id_{U} \times q)^{2}} \\
(V \times X)_{V}^{2} \ar[r]_{{\tilde{f}}^{2}} & (U \times X)_{U}^{2}
}
$$
respectively, are isomorphisms.
\end{lemma}

\begin{proof}
The first assertion follows from Proposition \ref{prop.canon} (2) since the diagram
$$
\xymatrix{
V \ar[r]^{f} \ar[d]_{\id_{V} \times qf} & U \ar[d]^{\id_{U} \times q} \\
V \times X \ar[r]_{f \times \id_{X}} & U \times X
}
$$
is a pullback diagram and the vertical arrows are affine by Lemma \ref{lem.gsep}.  The fact that the diagram is a pullback can be verified by a routine argument involving the universal property of the pullback.  The second assertion follows in a similar manner.
\end{proof}
Let $U$, $V$ and $W$ be objects in $\sf{S}$, and let $p:U \rightarrow S$ and $s:W \rightarrow S$ be the structure maps.  Let $\mathcal{B}$ be an ${\mathcal{O}}_{X}$-bimodule algebra, and assume all families are over $\mathcal{B}$.
\begin{lemma} \label{lem.function}
If $\mathcal{M}$ is a $U$-family and $f:V \rightarrow U$ is a morphism, then ${\tilde{f}}^{*}\mathcal{M}$ has a natural $V$-family structure.  Furthermore, if $\mathcal{M} \cong \mathcal{N}$ as ${\mathcal{B}}^{U}$-modules, then $\tilde{f}^{*}\mathcal{M} \cong \tilde{f}^{*}\mathcal{N}$ as ${\mathcal{B}}^{V}$-modules.
\end{lemma}

\begin{proof}
Suppose $\mathcal{M}$ is a $U$-family.  The natural $V$-family structure on ${\tilde{f}}^{*}\mathcal{M}$ is given by Theorem \ref{theorem.lift}.

Next, suppose $\phi:\mathcal{M} \rightarrow \mathcal{N}$ is an isomorphism of ${\mathcal{B}}^{U}$-modules.  By Lemma \ref{lem.isommonoid} and a simple computation, ${\tilde{f}}^{*}\alpha:{\tilde{f}}^{*}\mathcal{M} \rightarrow {\tilde{f}}^{*}\mathcal{N}$ is an isomorphism of ${\mathcal{B}}^{V}$-modules.
\end{proof}
We are now ready to define the functor we will eventually represent.
\begin{proposition} \label{prop.gamman}
The assignment $\Gamma_{n}:{\sf{S}} \rightarrow {\sf{Sets}}$ sending $U$ to
$$
\{ \mbox{isomorphism classes of truncated $U$-families of length n+1} \}
$$
and sending $f:V \rightarrow U$ to the map $\Gamma_{n}(f)$ defined by $\Gamma_{n}(f)[\mathcal{M}] = [{\tilde{f}}^{*}\mathcal{M}]$, is a functor.
\end{proposition}
\begin{proof}
By Lemma \ref{lem.function}, $\Gamma_{n}(-)$ is well defined.  We prove that the assignments in the statement of the proposition are compatible with composition.  Let
$$
\xymatrix{
W \ar[r]^{g} & V \ar[r]^{f} & U
}
$$
be morphisms in $\sf{S}$, let $\tilde{f}=f \times \id_{X}$ and let $\tilde{g}=g \times \id_{X}$.  Then we have a commutative diagram of $S$-schemes
$$
\xymatrix{
(W \times X)_{W}^{2} \ar[r]^{{\tilde{g}}^{2}} \ar[d]_{\pr_{2}^{W}} & (V \times X)_{V}^{2} \ar[r]^{{\tilde{f}}^{2}} \ar[d]^{\pr_{2}^{V}} & (U \times X)_{U}^{2} \ar[d]^{\pr_{2}^{U}} \\
W \times X \ar[r]_{\tilde{g}} & V \times X \ar[r]_{\tilde{f}} & U \times X
}
$$
By naturality of $(\tilde{f}\tilde{g})^{*} \cong {\tilde{g}}^{*}{\tilde{f}}^{*}$, the diagram
$$
\xymatrix{
{(\tilde{f}\tilde{g})}^{*}{\pr_{2}^{U}}_{*}({\pr_{1}}^{U*} \mathcal{M} \otimes {\tilde{p}}^{2*}\mathcal{B}) \ar[r] \ar[d] & {(\tilde{f}\tilde{g})}^{*}\mathcal{M} \ar[d] \\
{\tilde{g}}^{*}{\tilde{f}}^{*}{\pr_{2}^{U}}_{*}({\pr_{1}}^{U*} \mathcal{M} \otimes {\tilde{p}}^{2*}\mathcal{B}) \ar[r] & {\tilde{g}}^{*}{\tilde{f}}^{*}\mathcal{M}
}
$$
commutes.  Thus
$$
\Gamma_{n}(fg)[\mathcal{M}] = [{(\tilde{f}\tilde{g})}^{*}\mathcal{M}] = [{\tilde{g}}^{*}{\tilde{f}}^{*}\mathcal{M}] = \Gamma_{n}(g)\Gamma_{n}(f)[\mathcal{M}].
$$
\end{proof}
\begin{definition} \label{def.gamma}\index{functor of flat families of truncated $\mathcal{B}$-point modules}
The functor {\bf $\Gamma_{n}(-)$} is the {\bf functor of flat families of truncated $\mathcal{B}$-point modules of length $\mathbf{n+1}$}.
\end{definition}

\chapter{Compatibility with Descent} In this chapter, we define the
notion of a functor being compatible with descent and of a functor
being locally determined by a subfunctor.  When $P$ is a subscheme
of a projectivization, we find a subfunctor
$\operatorname{Hom}_{S}^{\Fr}(-,P)$ of
$\operatorname{Hom}_{S}(-,P)$ which locally determines
$\operatorname{Hom}_{S}(-,P)$.  We then find a subfunctor
${\Gamma}_{n}^{\Fr}(-)$ of ${\Gamma}_{n}(-)$ which locally
determines ${\Gamma}_{n}(-)$.  This allows us to reduce the
question of the representability of $\Gamma_{n}(-)$ to the
equivalence of $\Gamma_{n}^{\Fr}(-)$ and
$\operatorname{Hom}_{S}^{\Fr}(-,P)$.

\section{Local determination of a functor by a subfunctor}

\begin{definition} \label{def.compwithdescent}
A contravariant functor $F:{\sf{S}} \rightarrow {\sf{Sets}}$ is {\bf compatible with descent}\index{compatibility with descent|textbf} if, for every faithfully flat morphism of affine $S$-schemes $f:V \rightarrow U$, and for $f_{i}:V \times_{U} V \rightarrow V$ the projection maps, the sequence

$$
\xymatrix{
0 \ar[r] & F(U) \ar[r]^{F(f)} & F(V) \ar@<1ex>[r]^{F(f_{1})} \ar@<-1ex>[r]_{F(f_{2})} & F(V \times_{U} V)
}
$$
is exact.  In other words, $F$ is compatible with descent if $F$ is a sheaf in the faithfully flat topology on $\sf{S}$.
\end{definition}
The motivation behind the following technical definition is provided by Example \ref{example.grothendieck2}.  Readers familiar with Grothendieck's algebraic characterization of maps into projectivizations (Proposition \ref{prop.groth}) are encouraged to read Definition \ref{def.freemorph}, Corollary \ref{cor.defofsub} and Example \ref{example.grothendieck2} before reading Definition \ref{def.technical}.
\begin{definition} \label{def.technical}
Let $F:{\sf{S}} \rightarrow {\sf Sets}$ be a contravariant functor and suppose $F' \subset F$ is a subfunctor.  $F$ is {\bf locally determined by $F'$}\index{local determination of a functor by a subfunctor|textbf}\index{local determination of a functor by a subfunctor|(} if
\begin{itemize}
\item{$F$ is compatible with descent, and}
\item{for all $U \in {\sf S}$ and all finite subsets $\{x_{i}\} \subset F(U)$, there exists a $V \in {\sf S}$ and a faithfully flat map $f:V \rightarrow U$, such that $F(f)(x_{i}) \in F'(V)$ for all $i$.}
\end{itemize}
\end{definition}
{\it Remark.}  It follows readily from the above definition that $F$ is locally determined by itself if and only if $F$ is compatible with descent.
\begin{proposition} \label{prop.subfunctor}
Suppose $F$ and $G$ are locally determined by $F'$ and $G'$ respectively.  Then every natural transformation $\Omega:F' \Longrightarrow G'$ extends to a natural transformation $\overline{\Omega}:F \Longrightarrow G$ such that if $\Omega$ is monic, then $\overline{\Omega}$ is monic and if $\Omega$ is epic, then $\overline{\Omega}$ is epic.
\end{proposition}

\begin{proof}
Suppose $\Omega: F' \Longrightarrow  G'$, $U \in {\sf S}$ and $x \in F(U)$.  By the hypothesis that $F$ is locally determined by $F'$, there exists an affine $S$-scheme $V$ and a faithfully flat morphism $f:V \rightarrow U$ such that $F(f)(x) \in F'(V)$.  By the hypothesis that $F$ and $G$ are compatible with descent, the diagram
\begin{equation} \label{eqn.nat}
\xymatrix{
0 \ar[r] & F(U) \ar[r]^{F(f)} & F(V) \ar@<1ex>[r]^{F(f_{1})} \ar@<-1ex>[r]_{F(f_{2})} \ar@{-->}[d]_{\Omega_{V}} & F(V \times_{U} V) \ar@{-->}[d]^{\Omega_{V \times_{U}V}} \\
0 \ar[r] & G(U) \ar[r]^{G(f)} & G(V) \ar@<1ex>[r]^{G(f_{1})} \ar@<-1ex>[r]_{G(f_{2})} & G(U \times_{V} U)
}
\end{equation}
has exact rows and commutes, where the vertical columns are defined only for elements of the subfunctor $F' \subset F$.  Thus, $F(f)(x) \in F'(V)$ maps to an element, $y \in G(V)$ such that $G(f_{1})(y)=G(f_{2})(y)$.  Thus, there is a unique element $z \in G(U)$ such that $G(f)(z)=y$.  In this way, we get an assignment
$$
\overline{\Omega}_{U}:F(U) \rightarrow G(U).
$$
We show that this assignment is well defined, i.e. if we choose another $V' \in {\sf S}$ and locally finite, faithfully flat map $f':V' \rightarrow U$ such that $F(f')(x) \in F'(V')$, then we may use the same argument above with $V'$ in place of $V$ and still find that $x$ is assigned to $z$.  Let $V' \in {\sf S}$ and let $f':V' \rightarrow U$ be a faithfully flat map such that $F(f')(x) \in F'(V')$.  Let $V''=V \times_{U}V'$.  Then the diagram
$$
\xymatrix{
& F(V) \ar[rr] \ar@{-->}[dd] & & F(V'') \ar@{-->}[dd] \\
F(U) \ar[ur] \ar[rr] & & F(V') \ar@{-->}[dd] \ar[ur] & \\
& G(V) \ar[rr] & & G(V'') \\
G(U) \ar[rr] \ar[ur] & & G(V') \ar[ur] &
}
$$
whose verticals are the maps $\Omega$, commutes:  the top and bottom commute since $F$ and $G$ are functors, and the sides commute when restricted to $F'$ by the naturality of $\Omega$.  Thus, if $x \in F(U)$ gets sent to $w$ and $w'$ in $F(V)$ and $F(V')$ respectively, and if $w$ and $w'$ get sent to $y$ and $y'$ in $G(V)$ and $G(V')$ respectively, then $y$ and $y'$ get sent to the same element, $y'' \in G(V'')$.  But $y$, and hence $y''$, is the image of $z \in G(U)$.  But $y'$ also maps to $y''$.  Since the projection $V'' \rightarrow V'$ is faithfully flat, the induced map $G(V') \rightarrow G(V'')$ is an injection so that $y'$ must also be the image of $z$.  Thus, $\overline{\Omega}_{U}$ is well defined.

We show that $\overline{\Omega}$ is natural.  Let $g:U' \rightarrow U$ be a morphism in ${\sf S}$.  Let $x \in F(U)$.  Let $V \in {\sf S}$ such that there exists a faithfully flat map $f:V \rightarrow U$ such that $F(f)(x) \in F'(V)$.  Let $V'$ be the pullback of
$$
\xymatrix{
& V \ar[d]^{f} \\
U' \ar[r]_{g} & U
}
$$
and let $f':V' \rightarrow U'$ be the induced morphism, which is faithfully flat since $f$ is.  Since $F(f)(x) \in F'(V)$, $F(f')F(g)(x) \in F'(V')$, the left and right side of (\ref{eqn.cubeit}) commute with respect to $x$ by construction of $\overline{\Omega}$.  Since $F$ and $G$ are functors, the top and bottom of $(\ref{eqn.cubeit})$ commute.  Since $\Omega$ is natural, the diagram
$$
\xymatrix{
F'(V) \ar[d]_{\Omega_{V}} \ar[r] & F'(V') \ar[d]^{\Omega_{V'}} \\
G'(V) \ar[r] &  G'(V')
}
$$
commutes.  Thus, the cube
\begin{equation} \label{eqn.cubeit}
\xymatrix{
& F(V) \ar[dd] \ar[rr] & & F(V') \ar[dd] \\
F(U) \ar[dd] \ar[ur] \ar[rr] & & F(U') \ar[ur] \ar[dd] & \\
& G(V)  \ar[rr] & & G(V') \\
G(U) \ar[rr] \ar[ur] & & G(U') \ar[ur] &
}
\end{equation}
whose verticals are $\overline{\Omega}$, commutes with respect to $x$.  Since the maps going ``into'' the page are injections, the front face of the diagram commutes.  We conclude that $\overline{\Omega}$ is natural.

We next show that $\overline{\Omega}_{U}$ is monic if $\Omega_{U}$ is monic.  If $x$ and $x' \in F(U)$ are distinct, then by hypothesis, there exists a $V \in {\sf S}$ and a faithfully flat map $f:V \rightarrow U$ such that $F(f)(x)$ and $F(f)(x')$ are distinct and both in $F'(V)$.  A standard diagram chase in diagram (\ref{eqn.nat}) shows that $\overline{\Omega}_{U}(x) \neq \overline{\Omega}_{U}(x')$.

Finally, we show that $\overline{\Omega}_{U}$ is epic if $\Omega_{U}$ is epic.  Let $y \in G(U)$.  Choose a $V$, and $f:V \rightarrow U$, as above such that $G(f)(y) \in G'(V)$.  Then there exists a $w \in F'(V)$ which gets mapped to $G(f)(y)$.  Furthermore, since $G(f_{1}f)(y)=G(f_{2}f)(y)$, $F(f_{1}f)(w)=F(f_{2}f)(w)$.  Thus, by exactness of Diagram $\ref{eqn.nat}$, there exists a $x \in F(U)$ such that $F(f)(x)=w$.  As above, $\overline{\Omega}_{U}(x)=y$, proving that $\overline{\Omega}_{U}$ is epic.
\end{proof}\index{local determination of a functor by a subfunctor|)}

In order to apply Proposition $\ref{prop.subfunctor}$ to a pair of functors $F$ and $G$, we first need to show that $F$ and $G$ are compatible with descent.  This happens to be the case with $\operatorname{Hom}_{S}(-,Y)$ and $\Gamma_{n}(-)$:
\begin{theorem} \label{theorem.comp1} \cite[Proposition 28.3, p.72]{etale}
If $P$ is an $S$-scheme, then the functor
$$
\operatorname{Hom}_{S}(-,P)
$$
is compatible with descent.
\end{theorem}
At the end of this chapter, we will prove the following

\begin{theorem} \label{theorem.comp2}
If $\mathcal{B}$ is generated in degree one, then $\Gamma_{n}(-)$ is compatible with descent.
\end{theorem}

\section{An algebraic description of maps into projectivizations}
In order to describe the functor $\operatorname{Hom}_{S}^{\Fr}(-,P)$ which locally determines $\operatorname{Hom}_{S}(-,P)$, we must first characterize maps into projective bundles.  In particular, we review a correspondence, due to Grothendieck, between morphisms into projective bundles and algebraic data.

For any scheme $X$ and any ${\mathcal{O}}_{X}$-module $\mathcal{M}$, let ${\mathbb{S}}_{X}(\mathcal{M})$ denote the {\bf symmetric algebra} of $\mathcal{M}$ \index{symmetric algebra}(\cite[ex. 5.16, p.127]{alggeo}).  If $X=\mbox{Spec }A$, then we sometimes also write ${\mathbb{S}}_{A}(\mathcal{F}(\mbox{Spec }A))$ instead of ${\mathbb{S}}_{X}(\mathcal{F})$.

\begin{proposition} \label{prop.groth} \cite[Proposition 4.2.3, p. 73]{ega2}
Let $q: Y \rightarrow Z$ be a morphism of schemes.  Then, given an ${\mathcal{O}}_{Z}$-module $\mathcal{E}$ there is a bijective correspondence between the set of $Z$-morphisms $r:Y \rightarrow {\mathbb{P}}_{Z}(\mathcal{E})$, and the set of equivalence classes of pairs $(\mathcal{L}, \phi)$ composed of an invertible ${\mathcal{O}}_{Y}$-module $\mathcal{L}$ and a epimorphism $\phi:q^{*}(\mathcal{E}) \rightarrow \mathcal{L}$, where two pairs $(\mathcal{L}, \phi)$ and $({\mathcal{L}}', \phi ')$ are equivalent if there exists an ${\mathcal{O}}_{Y}$-module isomorphism $\tau: \mathcal{L} \rightarrow {\mathcal{L}}'$ such that

$$
\xymatrix
{
q^{*}(\mathcal{E}) \ar[r]^{\phi} \ar[dr]_{\phi '} & \mathcal{L}  \ar[d]^{\tau} \\
& \mathcal{L} '
}
$$
commutes.
\end{proposition}

\begin{proof}
While we do not give a complete proof of this result, we will describe the bijection between $\operatorname{Hom }_{Z}(-,{\mathbb{P}}_{Z}(\mathcal{E}))$ and the set of classes of pairs $(\mathcal{L}, \phi)$.  Let $q:Y \rightarrow Z$ and let $r:Y \rightarrow {\mathbb{P}}_{Z}(\mathcal{E})$ be a $Z$-morphism.  If $p:{\mathbb{P}}_{Z}(\mathcal{E}) \rightarrow Z$ is the structure morphism and if ${\mathcal{O}}$ is the structure sheaf of ${\mathbb{P}}_{Z}(\mathcal{E})$, then let $\alpha:p^{*}\mathcal{E} \rightarrow {\mathcal{O}}(1)$ be the canonical epimorphism \cite[Proposition 4.1.6, p. 72]{ega2}.  If $\phi$ is the composition
$$
\xymatrix{
q^{*}\mathcal{E} \ar[r]^{\cong} & r^{*}p^{*}\mathcal{E} \ar[r]^{r^{*}\alpha} & r^{*}{\mathcal{O}}(1).
}
$$
then $r$ corresponds to $(r^{*}{\mathcal{O}}(1),\phi)$.

Conversely, let $(\mathcal{L}, \phi)$  be a pair consisting of an invertible ${\mathcal{O}}_{Y}$-module $\mathcal{L}$ and a epimorphism $\phi:q^{*}(\mathcal{E}) \rightarrow \mathcal{L}$.  Applying the functor ${\mathbb{P}}_{Y}(-)$ to the epimorphism $\phi$ gives us a closed immersion, so that the composition
\begin{equation} \label{eqn.projective}
\xymatrix{
Y \ar[r]^{\cong} & {\mathbb{P}}_{Y}(\mathcal{L}) \ar[r]^{{\mathbb{P}}_{Y}(\phi)} & {\mathbb{P}}_{Y}(q^{*}\mathcal{E}).
}
\end{equation}
is a closed immersion.  Composing the map (\ref{eqn.projective}) with the map
$$
{\mathbb{P}}_{Y}(q^{*}\mathcal{E}) \rightarrow {\mathbb{P}}_{Z}(\mathcal{E})
$$
(\cite[4.1.3.1, p. 71]{ega2}) gives us the desired $Z$-morphism $r:Y \rightarrow {\mathbb{P}}_{Z}(\mathcal{E})$.
\end{proof}

When the morphism $q$ in Proposition \ref{prop.groth} is a closed immersion, we have the following result.
\begin{lemma} \label{lem.bij2}
Suppose $f:Y \rightarrow Z$ is a closed immersion of schemes.  If $\mathcal{F}$ is an ${\mathcal{O}}_{Y}$-module and $\mathcal{G}$ is an ${\mathcal{O}}_{Z}$-module, then the adjoint bijection
$$
{\operatorname{Hom}}_{Y}(f^{*}\mathcal{G},\mathcal{F}) \cong {\operatorname{Hom}}_{Z}(\mathcal{G},f_{*}\mathcal{F})
$$
is a bijection on epimorphisms.
\end{lemma}

\begin{proof}
Let $\psi:f^{*}\mathcal{G} \rightarrow \mathcal{F}$ be an epimorphism.  Since $f$ is a closed immersion, $f$ is affine so that $f_{*}$ is exact.  Thus, $f_{*}\psi$ is an epimorphism.  In addition, the map
$$
\mathcal{G} \rightarrow f_{*}f^{*}\mathcal{G}
$$
is epi, a fact we check locally.  Let $z \in Z$ such that there is a $y \in Y$ with $f(y)=z$.  Since $f$ is a closed immersion, there exists an ideal $I \subset {\mathcal{O}}_{Z,z}$ such that ${\mathcal{O}}_{Y,y} \cong {\mathcal{O}}_{Z,z}/I$.  Thus,
$$
(f_{*}f^{*}\mathcal{G})_{z} \cong (f^{*}\mathcal{G})_{y} \cong {\mathcal{G}}_{z}\otimes_{{\mathcal{O}}_{Z,z}} {\mathcal{O}}_{Y,y} \cong  {\mathcal{G}}_{z}\otimes_{{\mathcal{O}}_{Z,z}}{\mathcal{O}}_{Z,z}/I.
$$
Conversely, suppose $\phi:\mathcal{G} \rightarrow f_{*}\mathcal{F}$ is an epi. Since $f$ is affine, $f^{*}f_{*}\mathcal{F} \rightarrow \mathcal{F}$ is epi.  Thus, the composition
$$
\xymatrix{
f^{*}\mathcal{G} \ar[r]^{f^{*}\phi} & f^{*}f_{*}\mathcal{F} \ar[r] &  \mathcal{F}
}
$$
is an epi, as desired.
\end{proof}

\begin{lemma}
Suppose $f:Y \rightarrow Z$ is a closed immersion of schemes.  Let $\mathcal{F}$ and $\mathcal{F}'$ be ${\mathcal{O}}_{Y}$-modules and let $\mathcal{G}$ be an ${\mathcal{O}}_{Z}$-module.  Call two epis $\psi:f^{*}\mathcal{G} \rightarrow \mathcal{F}$ and $\psi':f^{*}\mathcal{G} \rightarrow \mathcal{F}'$ equivalent if there exists an isomorphism $\tau:\mathcal{F} \rightarrow \mathcal{F}'$ such that the diagram
$$
\xymatrix{
f^{*}\mathcal{G} \ar[r]^{\psi} \ar[dr]_{\psi '} & \mathcal{F} \ar[d]^{\tau} \\
& \mathcal{F}'
}
$$
commutes.  Call two epis $\phi:\mathcal{G} \rightarrow f_{*}\mathcal{F}$ and $\phi':\mathcal{G} \rightarrow f_{*}\mathcal{F}'$ equivalent if there exists an isomorphism $\nu:f_{*}\mathcal{F} \rightarrow f_{*}\mathcal{F}'$ such that the diagram
$$
\xymatrix{
\mathcal{G} \ar[r]^{\phi} \ar[dr]_{\phi '} & f_{*}\mathcal{F} \ar[d]^{\nu} \\
& f_{*}\mathcal{F}'
}
$$
commutes.  Then the adjoint bijection
$$
{\operatorname{Hom}}_{Y}(f^{*}\mathcal{G},\mathcal{F}) \cong {\operatorname{Hom}}_{Z}(\mathcal{G},f_{*}\mathcal{F})
$$
induces a bijection between equivalence classes of epimorphisms.
\end{lemma}

\begin{proof}
Suppose two epis $\psi:f^{*}\mathcal{G} \rightarrow \mathcal{F}$ and $\psi':f^{*}\mathcal{G} \rightarrow \mathcal{F}'$ are equivalent.  Then there exists a $\tau:\mathcal{F} \rightarrow \mathcal{F}'$ such that the diagram
$$
\xymatrix{
\mathcal{G} \ar[r] & f_{*}f^{*}\mathcal{G} \ar[r]^{f_{*}\psi} \ar[dr]_{f_{*}\psi '} & f_{*}\mathcal{F} \ar[d]^{f_{*}\tau} \\
& & f_{*}\mathcal{F}'
}
$$
commutes.  Thus, the adjoint bijection maps classes to classes by Lemma $\ref{lem.bij2}$.  Conversely, suppose two epis $\phi:\mathcal{G} \rightarrow f_{*}\mathcal{F}$ and $\phi':\mathcal{G} \rightarrow f_{*}\mathcal{F}'$ are equivalent.  Then there exists an isomorphism $\nu:f_{*}\mathcal{F} \rightarrow f_{*}\mathcal{F}'$ such that the diagram
$$
\xymatrix{
\mathcal{G} \ar[r]^{\phi} \ar[dr]_{\phi '} & f_{*}\mathcal{F} \ar[d]^{\nu} \\
 & f_{*}\mathcal{F}'
}
$$
commutes.  Let $\epsilon:f^{*}f_{*} \Longrightarrow \id$ be the counit of the pair $(f^{*},f_{*})$.  Since $f$ is a closed immersion, $\epsilon$ is an isomorphism.  If $\phi$ and $\phi'$ have right adjuncts to $\psi$ and $\psi'$, then $\psi$ and $\psi'$ are equivalent via the isomorphism
$$
\xymatrix{
\mathcal{F} \ar[r]^{\epsilon^{-1}} & f^{*}f_{*}\mathcal{F} \ar[r]^{f^{*}\nu} & f^{*}f_{*}\mathcal{F}' \ar[r]^{\epsilon} & \mathcal{F}.
}
$$
\end{proof}
The next fact now follows easily from Proposition \ref{prop.groth}.
\begin{corollary} \label{cor.groth}
Let $f: Y \rightarrow Z$ be a closed immersion.  Given an ${\mathcal{O}}_{Z}$-module $\mathcal{G}$ there is a bijective correspondence between the set of $Z$-morphisms $r:Y \rightarrow {\mathbb{P}}_{Z}(\mathcal{G})$, and the set of equivalence classes of pairs $(f_{*}\mathcal{L}, \phi)$ composed of the direct image of an invertible ${\mathcal{O}}_{Y}$-module $\mathcal{L}$ and an epimorphism $\phi:\mathcal{G} \rightarrow f_{*}\mathcal{L}$, where two pairs $(f_{*}\mathcal{L}, \phi)$ and $(f_{*}{\mathcal{L}}', \phi ')$ are equivalent if there exists an ${\mathcal{O}}_{Y}$-module isomorphism $\nu: f_{*}\mathcal{L} \rightarrow f_{*}{\mathcal{L}}'$ such that
$$
\xymatrix
{
\mathcal{G} \ar[r]^{\phi} \ar[dr]_{\phi '} & f_{*}\mathcal{L}  \ar[d]^{\nu} \\
& f_{*}\mathcal{L} '
}
$$
commutes.
\end{corollary}

\begin{corollary} \label{cor.kernel}
Let $f:Y \rightarrow Z$ be a closed immersion, let $\mathcal{G}$ be an ${\mathcal{O}}_{Z}$-module, and let $\mathcal{L}$ and $\mathcal{L}'$ be invertible ${\mathcal{O}}_{Y}$-modules.  If two epimorphisms $\phi:\mathcal{G} \rightarrow f_{*}\mathcal{L}$ and $\phi':\mathcal{G} \rightarrow f_{*}\mathcal{L}'$ have the same kernel, then they correspond, under the bijection of Corollary $\ref{cor.groth}$, to the same morphism $r:Y \rightarrow {\mathbb{P}}_{Z}(\mathcal{G})$.
\end{corollary}

\begin{proof}
This fact follows from the previous corollary in light of the existence of the commutative diagram with exact rows
$$
\xymatrix{
0 \ar[r] & \mbox{ker} \phi \ar[r] \ar[d]_{=} & \mathcal{G} \ar[d]^{=} \ar[r]^{\phi} & f_{*}\mathcal{L} \ar[r] & 0 \\
0 \ar[r] & \mbox{ker} \phi \ar[r]  & \mathcal{G}  \ar[r]_{\phi'} & f_{*}\mathcal{L}' \ar[r] & 0
}
$$
\end{proof}

\begin{lemma} \label{lem.proj}
If $f:Y \rightarrow Z$ is a morphism of schemes and $\mathcal{C}$ is an ${\mathcal{O}}_{Z}$-module, then there exists a map
\begin{equation} \label{eqn.proj}
\xymatrix{
{\mathbb{P}}_{Y}(f^{*}\mathcal{C}) \ar[r] & {\mathbb{P}}_{Z}(\mathcal{C})
}
\end{equation}
which induces an isomorphism
$$
{\mathbb{P}}_{Y}(f^{*}\mathcal{C}) \cong {\mathbb{P}}_{Z}(\mathcal{C}) \times_{Z}Y.
$$
Furthermore, if $\mathcal{A}$ and $\mathcal{B}$ are ${\mathcal{O}}_{Z}$-modules and $\phi:\mathcal{A} \rightarrow \mathcal{B}$ is an epi, then the diagram
\begin{equation} \label{eqn.com}
\xymatrix{
{\mathbb{P}}_{Y}(f^{*}\mathcal{B}) \ar[d] \ar[rrr]^{{\mathbb{P}}_{Y}(f^{*}\phi)} & & & {\mathbb{P}}_{Y}(f^{*}\mathcal{A}) \ar[d] \\
{\mathbb{P}}_{Z}(\mathcal{B}) \times_{Z} Y \ar[d] \ar[rrr]^{{\mathbb{P}}_{Z}(\phi) \times \id_{Y}} & & & {\mathbb{P}}_{Z}(\mathcal{A}) \times_{Z}Y \ar[d] \\
{\mathbb{P}}_{Z}(\mathcal{B}) \ar[rrr]^{{\mathbb{P}}_{Z}(\phi)} & & & {\mathbb{P}}_{Z}(\mathcal{A})
}
\end{equation}
commutes.
\end{lemma}

\begin{proof}
The construction of (\ref{eqn.proj}) and the proof of the first assertion is given in \cite[Proposition 3.5.3, p.62]{ega2}.  We outline the construction of (\ref{eqn.proj}).  Let Spec $A \subset Z$ be an affine open set and suppose Spec $A' \subset f^{-1}(\operatorname{Spec }A)$.  Then there is a map of graded rings

$$
g_{\mathcal{C}}:{\mathbb{S}}_{A}(\mathcal{C}(U)) \rightarrow {\mathbb{S}}_{A}(\mathcal{C}(U)) \otimes_{A}A' \rightarrow  {\mathbb{S}}_{A'}(\mathcal{C}(U) \otimes_{A}A')
$$
sending a homogeneous element $b$ to the image of $b\otimes 1$ in the right map.  By \cite[ex. 2.14b, p.80]{alggeo}, this map of graded rings induces a map of schemes
\begin{equation} \label{eqn.mapp}
{\mathbb{P}}_{\operatorname{Spec }A'}(f^{*}\mathcal{C}|_{\operatorname{Spec }A'}) \rightarrow {\mathbb{P}}_{\operatorname{Spec }A}(\mathcal{C}|_{\operatorname{Spec }A}).
\end{equation}
The maps ($\ref{eqn.mapp}$) glue to give a map ($\ref{eqn.proj}$).

Returning to a local setting, the diagram of graded rings
$$
\xymatrix{
{\mathbb{S}}_{A'}(\mathcal{B}(U) \otimes_{A}A') & & & {\mathbb{S}}_{A'}(\mathcal{A}(U) \otimes_{A}A') \ar[lll]_{{\mathbb{S}}_{A'}(\phi(U) \otimes_{A} A')} \\
{\mathbb{S}}_{A}(\mathcal{B}(U))\otimes_{A}A' \ar[u] & & & {\mathbb{S}}_{A}(\mathcal{A}(U))\otimes_{A} A' \ar[u] \ar[lll]^{{\mathbb{S}}_{A}(\phi(U))\otimes_{A} A'} \\
{\mathbb{S}}_{A}(\mathcal{B}(U)) \ar[u] & & & {\mathbb{S}}_{A}(\mathcal{A}(U)) \ar[u] \ar[lll]^{{\mathbb{S}}_{A}(\phi(U))}
}
$$
commutes.  Thus, (\ref{eqn.com}) commutes since it is the gluing of Proj(-) applied to the above diagram.
\end{proof}

\begin{lemma} \label{lem.compat}
Suppose
$$
\xymatrix{
A \ar[r]^{a} & B \ar[r]^{b} & C
}
$$
is a diagram of schemes, $\mathcal{C}$ is an ${\mathcal{O}}_{C}$-module and the vertical and top horizontal arrows in the diagram
$$
\xymatrix{
{\mathbb{P}}_{B}(b^{*}\mathcal{C}) \ar[r] & {\mathbb{P}}_{C}(\mathcal{C})  \\
{\mathbb{P}}_{A}(a^{*}b^{*}\mathcal{C}) \ar[r] \ar[u] & {\mathbb{P}}_{A}((ba)^{*}\mathcal{C}) \ar[u]
}
$$
are the maps (\ref{eqn.proj}).  Then the diagram commutes.
\end{lemma}

\begin{proof}
These maps are all induced locally from maps of graded rings.  These graded ring maps make the diagram commute, as in the previous lemma.
\end{proof}

\begin{lemma} \label{lem.newneed}
Suppose
$$
\xymatrix{
A \ar[r]^{a} & B \ar[r]^{b} & C
}
$$
is a diagram of schemes, $\mathcal{C}$ is an ${\mathcal{O}}_{C}$-module, and $r:A \rightarrow {\mathbb{P}}_{C}(\mathcal{C})$ corresponds to the epimorphism $(ba)^{*}\mathcal{C} \rightarrow \mathcal{L}$ via the correspondence constructed in Proposition \ref{prop.groth}. Then there exists a morphism $s:A \rightarrow {\mathbb{P}}_{B}(b^{*}\mathcal{C})$ such that the diagram
$$
\xymatrix{
A \ar[r]^{s} \ar[dr]_{r} & {\mathbb{P}}_{B}(b^{*}\mathcal{C}) \ar[d] \\
& {\mathbb{P}}_{C}(\mathcal{C})
}
$$
whose right vertical is (\ref{eqn.proj}), commutes.
\end{lemma}

\begin{proof}
If $r$ corresponds to the epimorphism $\phi:(ba)^{*}\mathcal{C} \rightarrow \mathcal{L}$, then the map $s$ corresponds to the epimorphism
$$
\xymatrix{
a^{*}b^{*}\mathcal{C} \ar[r]^{\cong} & (ba)^{*}\mathcal{C} \ar[r]^{\phi} & \mathcal{L}.
}
$$
Thus, there is a commutative diagram
$$
\xymatrix{
(ba)^{*}\mathcal{C} \ar[d] \ar[r]^{\phi} & \mathcal{L} \\
a^{*}b^{*}\mathcal{C} \ar[ur] &
}
$$
Since ${\mathbb{P}}_{A}$ is functorial, the diagram
$$
\xymatrix{
A \ar[r] & {\mathbb{P}}_{A}(\mathcal{L}) \ar[dr] \ar[r] & {\mathbb{P}}_{A}((ba)^{*}\mathcal{C}) \ar[r] & {\mathbb{P}}_{C}(\mathcal{C}) \\
& & {\mathbb{P}}_{A}(a^{*}b^{*}\mathcal{C}) \ar[u] \ar[r] & {\mathbb{P}}_{B}(b^{*}\mathcal{C}) \ar[u].
}
$$
commutes by Lemma \ref{lem.compat}.  Since the top route is $r$ while the bottom route is $s$ composed with the right vertical (\ref{eqn.proj}), the assertion follows.
\end{proof}

\begin{lemma} \label{lem.algproj}
Suppose $f:Y \rightarrow Z$ is a morphism of schemes, $\mathcal{C}$ is an ${\mathcal{O}}_{Z}$-module, $W$ be a scheme and there is a map $r:W \rightarrow {\mathbb{P}}_{Y}(f^{*}\mathcal{C})$.  Suppose that $r$ projects to a map $q:W \rightarrow Y$.  If $r$ corresponds to an epi of sheaves
$$
\phi:q^{*}f^{*}\mathcal{C} \rightarrow \mathcal{L}
$$
then $r$ composed with the map ${\mathbb{P}}_{Y}(f^{*}\mathcal{C}) \rightarrow {\mathbb{P}}_{Z}(\mathcal{C})$ (\ref{eqn.proj}) corresponds to the epi
\begin{equation} \label{eqn.epii}
\xymatrix{
(fq)^{*}\mathcal{C} \ar[r] & q^{*}f^{*}\mathcal{C} \ar[r]^{\phi} & \mathcal{L}
}
\end{equation}
\end{lemma}

\begin{proof}
By definition of the correspondence in Proposition \ref{prop.groth}, $r$ corresponds to the top row of the diagram
$$
\xymatrix{
W \ar[rr]^{\cong} & & {\mathbb{P}}_{W}(\mathcal{L}) \ar[rr]^{{\mathbb{P}}_{W}(\phi)} & & {\mathbb{P}}_{W}(q^{*}f^{*}\mathcal{C}) \ar[r] \ar[d] & {\mathbb{P}}_{Y}(f^{*}\mathcal{C}) \ar[d] \\
& & & & {\mathbb{P}}_{W}((fq)^{*}\mathcal{C}) \ar[r] & {\mathbb{P}}_{Z}(\mathcal{C})
}
$$
which commutes by Lemma \ref{lem.compat}.  Thus, $r$ composed with the right vertical of the diagram is just the bottom route of the diagram.  By the definition of the correspondence in Proposition \ref{prop.groth}, it is easy to see that the epimorphism (\ref{eqn.epii}) corresponds to the bottom route.
\end{proof}

\section{Free morphisms and free families}\index{free morphism|(}
Let $P$ be a subscheme of a projectivization.  In this section, we define the functors ${\operatorname{Hom}}_{S}^{\Fr}(-,P)$ and ${\Gamma}_{n}^{\Fr}(-)$ and we prove that they locally determine ${\operatorname{Hom}}_{S}(-,P)$ and ${\Gamma}_{n}(-)$, respectively.
\begin{definition} \label{def.freemorph}
Let $q:Y \rightarrow Z$ be a morphism of schemes, and let $\mathcal{E}$ be an ${\mathcal{O}}_{Z}$-module.  A $Z$-morphism $r:Y \rightarrow {\mathbb{P}}_{Z}(\mathcal{E})$ is {\bf free} if, under the bijection of Proposition $\ref{prop.groth}$, $r$ corresponds to an equivalence class containing an element of the form $({\mathcal{O}}_{Y}, \phi)$.  If $i:P \rightarrow {\mathbb{P}}_{Z}(\mathcal{E})$ is a subscheme, {\bf a morphism $f:Y \rightarrow P$ is free} if $if$ is free.
\end{definition}\index{free morphism|textbf}
The following example illustrates the significance of Definition \ref{def.freemorph}.
\begin{example} \label{example.grothendieck}
Suppose $k$ is a field, $Z=\operatorname{Spec }k$, $Y=\operatorname{Spec }R$ is an affine $Z$-scheme and $q:Y \rightarrow Z$ is the structure map.  Suppose, further, that $\mathcal{E}$ is an $n$-dimensional $k$-vector space and $r:Y \rightarrow {\mathbb{P}}_{Z}(\mathcal{E})$ is a free $Z$-morphism which corresponds, under the bijection in Proposition \ref{prop.groth}, to the class of an epimorphism $\phi:q^{*}\mathcal{E} \rightarrow {\mathcal{O}}_{Y}$.  Since $q^{*}\mathcal{E} \cong {\mathcal{O}}_{Y}^{\oplus n}$, $\phi$ is determined by a non-zero $n$-tuple, $(a_{1}, \ldots, a_{n})$, of elements of $R$.  Furthermore, $({\mathcal{O}}_{Y},{\phi}')$ is a pair in the class of $({\mathcal{O}}_{Y},\phi)$ such that $\phi'$ is determined by $(a_{1}',\ldots,a_{n}') \in R^{\oplus n}$, if and only if there exists a unit $a \in R$ such that $(a_{1}',\ldots,a_{n}')=a(a_{1}, \ldots, a_{n})$.  In particular, free $Z$-morphisms $r:Y \rightarrow {\mathbb{P}}_{Z}(\mathcal{E})$ are parameterized by points of ${\mathbb{P}}_{Z}(k^{n})$ with coordinates in $R$.
\end{example}

\begin{lemma} \label{lem.nat1}
Suppose $f:V \rightarrow U$ is a morphism of affine schemes.  If $q:U \rightarrow Z$ is given and $r:U \rightarrow {\mathbb{P}}_{Z}(\mathcal{E})$ is free, then $rf$ is free.
\end{lemma}

\begin{proof}
If $r$ corresponds to the class containing $\phi:q^{*}\mathcal{E} \rightarrow {\mathcal{O}}_{U}$, then $rf$ corresponds to the class containing the composition
$$
\xymatrix{
(qf)^{*}\mathcal{E} \ar[r]^{\cong} & f^{*}q^{*}\mathcal{E} \ar[r]^{f^{*}\phi} & f^{*}{\mathcal{O}}_{U} \ar[r]^{\cong} & {\mathcal{O}}_{V}
}
$$
\cite[4.2.8, p.75]{ega2}.
\end{proof}

We give a definition which appears in the statement of the next result.

\begin{definition} \label{def.freefamily}
A truncated $U$-family $\mathcal{M}$ of length $n+1$ is called {\bf free} \index{family!free, of truncated point modules|textbf}\index{family!free, of truncated point modules|(}if, for $0 \leq i \leq n$, ${\mathcal{L}}_{i}$ in the definition of ${\mathcal{M}}_{i}$ (\ref{def.family}) is isomorphic to ${\mathcal{O}}_{U}$.
\end{definition}

\begin{lemma} \label{lem.nat2}
Let $f: V \rightarrow U$ be a morphism in $\sf{S}$ and suppose $\mathcal{M}$ is a free truncated $U$-family of length $n+1$.  Then $\Gamma_{n}(f)[\mathcal{M}]$ is an isomorphism class of free truncated families parameterized by $V$.
\end{lemma}

\begin{proof}
Let $\tilde{f}=f\times \id_{X}:V \times X \rightarrow U \times X$.  Suppose ${\mathcal{M}}_{i} \cong (\id_{U} \times q_{i})_{*}{\mathcal{O}}_{U}$ for some map $q_{i}:U \rightarrow X$.  Since $\Gamma_{n}(f)[\mathcal{M}]=[{\tilde{f}}^{*}\mathcal{M}]$ (Proposition \ref{prop.gamman}), in order to prove the result, we need only note that
$$
({\tilde{f}}^{*}\mathcal{M})_{i} \cong (\id_{V} \times q_{i}f)_{*}f^{*}{\mathcal{O}}_{U} \cong (\id_{V} \times q_{i}f)_{*}{\mathcal{O}}_{V}.
$$
\end{proof}

\begin{corollary} \label{cor.defofsub}
Let $\mathcal{E}$ be an ${\mathcal{O}}_{Z}$-module and let $i:P \rightarrow {\mathbb{P}}_{Z}(\mathcal{E})$ be a subscheme.  The assignments
$$
\operatorname{Hom}_{S}^{\Fr}(-,P):\sf{S} \rightarrow {\sf Sets}
$$
defined by
$$
U \mapsto  \mbox{ subset of $\operatorname{Hom}_{S}(U,P)$ of free morphisms}
$$
and
$$
\Gamma_{n}^{\Fr}(-):\sf{S} \rightarrow {\sf Sets}
$$
defined by
$$
U \mapsto \mbox{ subset of $\Gamma_{n}(U)$ of free truncated $U$-families of length $n+1$},
$$
naturally induce subfunctors of $\operatorname{Hom}_{S}(-,P)$ and $\Gamma_{n}(-)$, respectively\index{Gammanfr@$\Gamma_{n}^{\Fr}$|textbf}.
\end{corollary}

\begin{proof}
By Lemma $\ref{lem.nat1}$ and Lemma $\ref{lem.nat2}$, the assignments inherit the property of being functors from $\Gamma_{n}(-)$ and $\operatorname{Hom}_{S}(-,Y)$.
\end{proof}

\begin{lemma} \label{lem.invertible}
\cite[ex. 4.12a, p.136]{comalg}.  If $U = \operatorname{Spec} R$ is an affine scheme and ${\mathcal{L}}_{1}, \ldots, {\mathcal{L}}_{n}$ are invertible ${\mathcal{O}}_{U}$-modules then there exists a finite set of $R$-module generators of $R$, $\{s_{j} \}$, such that each ${\mathcal{L}}_{i}$ is free on each $\operatorname{Spec}R_{s_{j}}$.
\end{lemma}
The following example illustrates the idea behind the proof that $\operatorname{Hom}_{S}(-,P)$ is locally determined by $\operatorname{Hom}_{S}^{\Fr}(-,P)$.
\begin{example} \label{example.grothendieck2}
Suppose $Y$ is a noetherian affine $Z$-scheme and $q:Y \rightarrow Z$ is the structure map.  Suppose, further, that $\mathcal{E}$ is an ${\mathcal{O}}_{Z}$-module  and $r:Y \rightarrow {\mathbb{P}}_{Z}(\mathcal{E})$ is a $Z$-morphism which corresponds, under the bijection in Proposition \ref{prop.groth}, to the class of an epimorphism $\phi:q^{*}\mathcal{E} \rightarrow \mathcal{L}$.  Let $\{U_{1}, \ldots, U_{m}\}$ be an affine open cover of $Y$.  Then the map
$$
f:\underset{i=1}{\overset{m}{\Pi}}U_{i} \rightarrow Y
$$
is faithfully flat.  Since $\operatorname{Hom}_{S}(-,{\mathbb{P}}_{Z}(\mathcal{E}))$ is compatible with descent, $r$ is determined by morphisms $r_{i}:U_{i} \rightarrow {\mathbb{P}}_{Z}(\mathcal{E})$.  This is a restatement of the fact that any morphism $r:Y \rightarrow {\mathbb{P}}_{Z}(\mathcal{E})$ can be constructed by gluing compatible morphisms on the cover $\{U_{1}, \ldots, U_{m}\}$.

If, however, $\mathcal{L}$ is free on the affine open cover $\{U_{1}, \ldots, U_{m}\}$, then the maps $r_{i}$ are free morphisms.  By Lemma \ref{lem.invertible}, one can always find a finite affine open cover of $Y$ on which $\mathcal{L}$ is free.  Thus, any $Z$-morphism $r:Y \rightarrow {\mathbb{P}}_{Z}(\mathcal{E})$ can be constructed by gluing {\it free} morphisms.    Since, as illustrated by Example \ref{example.grothendieck}, free morphisms are easier to work with than general morphisms, this observation will be useful to us.

More generally, if $x_{1}, \ldots, x_{n}$ are $Z$-morphisms from $Y$ to ${\mathbb{P}}_{Z}(\mathcal{E})$, the fact that there exists an affine open cover $\{U_{1}, \ldots, U_{m}\}$ of $Y$ such that $x_{j}$ is determined by free morphisms $x_{ij}:U_{i} \rightarrow {\mathbb{P}}_{Z}(\mathcal{E})$ for each $1 \leq j \leq n$ is exactly the fact that $\operatorname{Hom}_{S}(-,{\mathbb{P}}_{Z}(\mathcal{E}))$ is locally determined by $\operatorname{Hom}_{S}^{\Fr}(-,{\mathbb{P}}_{Z}(\mathcal{E}))$.
\end{example}

\begin{lemma}
If $\mathcal{E}$ is an ${\mathcal{O}}_{Z}$-module and $i:P \rightarrow {\mathbb{P}}_{Z}(\mathcal{E})$ is a subscheme then $\operatorname{Hom}_{S}(-,P)$ is locally determined by $\operatorname{Hom}_{S}^{\Fr}(-,P)$ and $\Gamma_{n}(-)$ is locally determined by $\Gamma_{n}^{\Fr}(-)$.
\end{lemma}
\begin{proof}
The proof of the first assertion was outlined in the previous example, so we only prove the second assertion.  By theorem $\ref{theorem.comp2}$, $\Gamma_{n}$ is compatible with descent.  Let $U$=Spec $R$ be an affine noetherian scheme, let $J$ be a finite index set, and let $[{\mathcal{M}}^{j}] \in \Gamma_{n}(U)$ for $j \in J$.  Let ${\mathcal{L}}_{ij}$ be the invertible ${\mathcal{O}}_{U}$-module appearing in the $i$-th component of ${\mathcal{M}}^{j}$.  By Lemma $\ref{lem.invertible}$, there exists a finite generating set $\{r_{k}\} \subset R$ such that ${\mathcal{L}}_{ij}$ is free on Spec $R_{r_{k}}$ for each $i$, $j$, and $k$.  If $Q = \Pi_{k}R_{r_{k}}$, $\operatorname{Spec}Q$=$V$ and $f:V \rightarrow U$, then $f$ is faithfully flat, and $\Gamma_{n}(f)([{\mathcal{M}}^{j}])$ is the class of a free $V$-family.
\end{proof}

\begin{corollary} \label{cor.freefunct}
Let $i:P \rightarrow {\mathbb{P}}_{Z}(\mathcal{E})$ be a subscheme.  In order to prove that $\Gamma_{n}(-)$ is represented by $P$, it suffices to show that $\Gamma_{n}^{\Fr}(-) \cong \operatorname{Hom}_{S}^{\Fr}(-,P)$.
\end{corollary}\index{free morphism|)}\index{family!free, of truncated point modules|)}

\section{The proof that $\Gamma_{n}$ is compatible with descent}
We now assume, for the remainder of this section that $\mathcal{B}$ is generated in degree one.  We use this assumption explicitly in only one step of our proof, Corollary \ref{cor.extend}.  Unless otherwise stated, we assume all families are truncated, of length $n$.  Under these hypothesis, we prove that the functor $\Gamma_{n}(-)$ is compatible with descent (Theorem \ref{theorem.comp2}).  Our proof employs some basic results from descent theory, which we now review.

\subsection{Descent theory}
The following definitions and Theorem \ref{theorem.descent} are taken from \cite[Section 6.1]{neron}.

\begin{definition} \label{def.descent}
Suppose
$$
h: Y' \rightarrow Y,
$$
$Y'' =  Y' \times_{Y} Y'$, and $Y'''=Y' \times_{Y} Y' \times_{Y} Y'$.  Let $h_{i}$ and $h_{ij}$ be the projection maps from $Y''$ and $Y'''$ respectively.  Let $\mathcal{F}$ be an ${\mathcal{O}}_{Y'}$-module.  A {\bf covering datum of $\mathcal{F}$}\index{covering datum|textbf} is an isomorphism $\phi:h_{1}^{*}\mathcal{F} \rightarrow h_{2}^{*}\mathcal{F}$.  A covering datum $\phi$ of $\mathcal{F}$ is called a {\bf descent datum} if the diagram

\begin{equation} \label{eqn.datum}
\xymatrix{
h_{12}^{*}h_{1}^{*}\mathcal{F} \ar[r]^{h_{12}^{*}\phi} \ar[d] & h_{12}^{*}h_{2}^{*}\mathcal{F} \ar[r] & h_{23}^{*}h_{1}^{*}\mathcal{F} \ar[r]^{h_{23}^{*}\phi} & h_{23}^{*}h^{*}\mathcal{F} \ar[d] \\
h_{13}^{*}h_{1}^{*}\mathcal{F} \ar[rrr]_{h_{13}^{*}\phi} & & & h_{13}^{*}h^{*}\mathcal{F}
}
\end{equation}
commutes, where all unlabeled arrows are the natural isomorphisms.  In this case, $\phi$ is said to satisfy the {\bf cocycle condition}\index{cocycle condition|textbf}.  The {\bf category of ${\mathcal{O}}_{Y'}$-modules with descent data}\index{category of oyprime modules with descent data@category of ${\mathcal{O}}_{Y'}$-modules with descent data|textbf} is the category whose objects are pairs $(\mathcal{F},\phi)$ where $\mathcal{F}$ is an ${\mathcal{O}}_{Y'}$-module and $\phi$ is a descent datum on $\mathcal{F}$, and a morphism between $(\mathcal{F},\phi)$ and $(\mathcal{G},\psi)$ is a morphism $\gamma: \mathcal{F} \rightarrow \mathcal{G}$ such that the diagram

\begin{equation} \label{eqn.morph}
\xymatrix{
h_{1}^{*}\mathcal{F} \ar[r]^{\phi} \ar[d]_{h_{1}^{*}\gamma} & h_{2}^{*}\mathcal{F} \ar[d]^{h_{2}^{*}\gamma} \\
h_{1}^{*}\mathcal{G} \ar[r]_{\psi} & h_{2}^{*}\mathcal{G}
}
\end{equation}
commutes.  If this holds, we say $h_{1}^{*}\gamma=h_{2}^{*}\gamma$.
\end{definition}
We will need the following inheritance fact.
\begin{lemma} \label{lem.inher}
Suppose $(\mathcal{F},\phi)$ is an ${\mathcal{O}}_{Y'}$-module with descent datum\index{descent datum}, $\mathcal{G}$ is an ${\mathcal{O}}_{Y'}$-module, and
$$
\psi:{h_{1}}^{*}\mathcal{G} \rightarrow {h_{2}}^{*}\mathcal{G}
$$
is an isomorphism.  If $\mu:\mathcal{F} \rightarrow \mathcal{G}$ is an epi making the diagram
$$
\xymatrix{
{h_{1}}^{*}\mathcal{F} \ar[r]^{\phi} \ar[d]_{{h_{1}}^{*}\mu} & {h_{2}}^{*}\mathcal{F} \ar[d]^{{h_{2}}^{*}\mu} \\
{h_{1}}^{*}\mathcal{G} \ar[r]_{\psi} & {h_{2}}^{*}\mathcal{G}
}
$$
commute, then $\psi$ is a descent datum\index{descent datum}.
\end{lemma}
We now introduce some useful terminology.
\begin{definition} \label{def.candesdat}
Suppose $\mathcal{H}$ is an ${\mathcal{O}}_{Y}$-module.  Then there the natural map
$$
\phi:{h_{1}}^{*}h^{*}\mathcal{H} \rightarrow {h_{2}}^{*}h^{*}\mathcal{H}
$$
induced by the canonical map of functors
$$
{h_{1}}^{*}h^{*} \Longrightarrow (hh_{1})^{*}=(hh_{2})^{*} \Longrightarrow {h_{2}}^{*}h^{*}.
$$
is called the {\bf canonical descent datum on $h^{*}\mathcal{H}$}.\index{descent datum!canonical|textbf}
\end{definition}
We have the following main theorem of descent theory:

\begin{theorem} \label{theorem.descent}
(Grothendieck).  Let $h:Y' \rightarrow Y$ be faithfully flat and quasi-compact.  Then the functor $\mathcal{H} \mapsto h^{*}\mathcal{H}$, which goes from ${\mathcal{O}}_{Y}$-modules to ${\mathcal{O}}_{Y'}$-modules with descent data\index{category of oyprime modules with descent data@category of ${\mathcal{O}}_{Y'}$-modules with descent data}, is an equivalence of categories.
\end{theorem}

\subsection{Modified descent data}
We now prove a theorem which will allow us to work in a slightly different setting then that above.
\begin{definition} \label{def.mdescent}
Retain the notation in Definition \ref{def.descent}.  Let $Z''$ and $Z'''$ be schemes such that there exists isomorphisms
$$
z'':Y'' \rightarrow Z''
$$
and
$$
z''':Y''' \rightarrow Z'''
$$
and maps
$$
e_{i}:Z'' \rightarrow Y'
$$
and
$$
e_{ij}:Z''' \rightarrow Z''
$$
making the diagram
$$
\xymatrix{
Y''' \ar[r]^{h_{ij}} \ar[d]_{z'''} & Y'' \ar[d]^{z''} \ar[r]^{h_{k}} & Y' \ar[d]^{=}  \\
Z''' \ar[r]_{e_{ij}} & Z'' \ar[r]_{e_{k}} & Y'
}
$$
commute.  The {\bf category of ${\mathcal{O}}_{Y'}$-modules with descent data modified by $z''$ and $z'''$}\index{category of oyprime modules with descent data@category of ${\mathcal{O}}_{Y'}$-modules with descent data!modified|textbf} or simply the {\bf category of ${\mathcal{O}}_{Y'}$-modules with modified descent data} when the context is clear, is the category whose objects are pairs $(\mathcal{F}',\phi)$, where $\mathcal{F}'$ is an ${\mathcal{O}}_{Y'}$-module and $\phi:{e_{1}}^{*}\mathcal{F}' \rightarrow {e_{2}}^{*}\mathcal{F}'$ is an isomorphism making the descent datum diagram (\ref{eqn.datum}), modified by replacing all ``$h$'''s by ``$e$'''s, commute, and whose morphisms, $\gamma:(\mathcal{F}',\phi) \rightarrow (\mathcal{G}',\psi)$ make diagram $(\ref{eqn.morph})$, modified by replacing all ``$h$'''s by ``$e$'''s, commute.
\end{definition}
As before, $e_{1}^{*}h^{*}$ is naturally isomorphic to $e_{2}^{*}h^{*}$ providing canonical descent data for any module of the form $h^{*}\mathcal{H}$ where $\mathcal{H}$ is an ${\mathcal{O}}_{Y}$-module.

Since the category of ${\mathcal{O}}_{Y'}$-modules with descent data is so similar to the category of ${\mathcal{O}}_{Y'}$-modules with modified descent data, we expect the categories to be equivalent.

\begin{theorem} \label{theorem.mod}
The functor from ${\mathcal{O}}_{Y'}$-modules with descent data\index{category of oyprime modules with descent data@category of ${\mathcal{O}}_{Y'}$-modules with descent data} to ${\mathcal{O}}_{Y'}$-modules with modified descent data sending an object $(\mathcal{F}',\phi)$ to $(\mathcal{F}',\psi)$, where $\psi$ is the composition

$$
\xymatrix{
{e_{1}}^{*}\mathcal{F}' \ar[r] & {z''}^{*}{h_{1}}^{*}\mathcal{F}' \ar[r]^{{z''}^{*}\phi} & {z''}^{*}{h_{2}}^{*}\mathcal{F}' \ar[r] & {e_{2}}^{*}\mathcal{F}'
}
$$
and sending a morphism to itself, is an equivalence of categories.
\end{theorem}

\begin{proof}
The only thing which is not obvious is that assignment we describe is a functor {\it into} the correct category.  This follows from a routine, but tedious, computation.
\end{proof}

\begin{corollary} \label{cor.descent}
Retain the notation of Definition \ref{def.mdescent}.  If $h:Y' \rightarrow Y$ is faithfully flat and quasi-compact, the functor $\mathcal{H} \mapsto h^{*}\mathcal{H}$, which goes from ${\mathcal{O}}_{Y}$-modules to ${\mathcal{O}}_{Y'}$-modules with descent data modified by $z''$ and $z'''$, is an equivalence of categories.
\end{corollary}
For the remainder of this chapter, we use Corollary
\ref{cor.descent} to prove Theorem \ref{theorem.comp2}.  Before
reading further in this chapter, the reader is advised to review
the last two sentences of Section 1.6.
\subsection{Standard families.}  Let $U$ be an affine, noetherian scheme, and suppose all families are truncated of length $n+1$ over an ${\mathcal{O}}_{X}$-bimodule algebra $\mathcal{B}$.
\begin{definition}
A $U$-family or truncated $U$-family $\mathcal{M}$ is in {\bf standard form} or is {\bf standard}\index{family!standard|textbf} if, for each $i$ such that ${\mathcal{M}}_{i} \neq 0$, ${\mathcal{M}}_{i}=(\id_{U} \times q_{i})_{*} {\mathcal{L}}_{i}$ for some map $q_{i}:U \rightarrow X$ and some invertible sheaf ${\mathcal{L}}_{i}$ on $U$.
\end{definition}
\begin{lemma} \label{lem.standard}
Any truncated $U$-family is canonically isomorphic to a truncated family in standard form.
\end{lemma}

\begin{proof}
Suppose $\mathcal{M}$ be a truncated $U$-family and
\begin{equation} \label{eqn.standard}
{\mathcal{M}}_{i} \cong (\id_{U}\times q_{i})_{*}{\mathcal{L}}_{i}.
\end{equation}
We construct a truncated $U$-family, $\mathcal{P}$, in standard form, such that $\mathcal{M} \cong \mathcal{P}$.  Suppose $\mathcal{M}$ has multiplication $\mu_{ij}$ and let ${\mathcal{P}}_{i}=(\id_{U} \times q_{i})_{*}{\mathcal{L}}_{i}$.  Let $\mu_{ij}': {\mathcal{P}}_{i} \otimes_{{\mathcal{O}}_{U \times X}} {\mathcal{B}}^{U}_{j} \rightarrow {\mathcal{P}}_{i+1}$ be defined as the composition:

$$
\xymatrix{
{\mathcal{P}}_{i} \otimes_{{\mathcal{O}}_{U \times X}} {\mathcal{B}}^{U}_{j} \ar[r]^{\cong} & {\mathcal{M}}_{i} \otimes_{{\mathcal{O}}_{U \times X}} {\mathcal{B}}^{U}_{j} \ar[r]^{\hskip .3in {\mu}_{ij}} & {\mathcal{M}}_{i+j} \ar[r]^{\cong} & {\mathcal{P}}_{i+j}.
}
$$
It is easy to check that $\mathcal{M} \cong \mathcal{P}$ under the maps (\ref{eqn.standard}).
\end{proof}
The following three technical lemmas are used to prove Proposition \ref{prop.preker}.  This proposition is used, in turn, to prove Corollary \ref{cor.extend}.
\begin{lemma} \cite[Lemma 1.31, p.58]{alggeo1} \label{lem.sup1}
Let $Y$ be a scheme and suppose $\mathcal{R}$ is a coherent ${\mathcal{O}}_{Y}$-module.  Suppose $Z \subset Z' \subset Y$ are open sets such that
$$
Z \cap \operatorname{Supp }\mathcal{R} = Z' \cap \operatorname{Supp }\mathcal{R}.
$$
Then the restriction map $r:\mathcal{R}(Z') \rightarrow \mathcal{R}(Z)$ is an isomorphism of ${\mathcal{O}}_{Y}(Z')$-modules.
\end{lemma}

\begin{lemma} \label{lem.sup2}
Suppose $\mathcal{L}$ is an invertible ${\mathcal{O}}_{U}$-module, $q:U \rightarrow X$ is a morphism and $\pr_{i}:(U \times X)_{U}^{2} \rightarrow U \times X$ is the standard projection.  Then
\begin{enumerate}
\item{}
the set
$$
\{q^{-1}(V) \times W | V, W \subset X \mbox{ are affine open} \}
$$
forms an open cover of $U \times X$, and
\item{}
if $V, W \subset X$ are open affine and $V' \subset q^{-1}(V)$ is open then
$$
(U \times X) \times_{U} (V' \times W) \cap \operatorname{Supp } \pr_{1}^{*}(\id_{U} \times q)_{*}\mathcal{L} =
$$
$$
(U \times V) \times_{U} (V' \times W) \cap \operatorname{Supp } \pr_{1}^{*}(\id_{U} \times q)_{*}\mathcal{L}.
$$
\end{enumerate}
\end{lemma}

\begin{proof}
The proof of 1 is trivial.  To prove 2, we first note that the bottom expression is a subset of the top.  Conversely, we note that
$$
\mbox{Supp }\pr_{1}^{*}(\id_{U}\times q)_{*}\mathcal{L} = \{(a,b,a,d)\in(U \times X) \times_{U} (U \times X) | b=q(a) \}.
$$
Thus, an element of the top, $(a,b,a,d)$, must have the property that $a \in V'$, $b=q(a)$, and $d \in W$.  Since $V' \subset q^{-1}(V)$, an element of the top must have second coordinate in $V$.
\end{proof}

\begin{lemma} \label{lem.reduce}
Let $R$ be a commutative ring, let $A$, $B$ and $C$ be $R$-algebras, and suppose there is a map of $R$-algebras
$$
h:B \rightarrow A.
$$
Let $M$ be a free $A$-module of rank one, let $N$ be an $A$-central $A\otimes B-A \otimes C$-bimodule, and let $P$ be any right $A \otimes C$-module.  Give $M$ the $A\otimes B$-module structure it inherits from $h$.  Suppose
$$
\mu:M\otimes_{A \otimes B}N \rightarrow P
$$
is any morphism of right $A \otimes C$-modules.  If
$$
\psi:M \rightarrow M
$$
is an $A \otimes B$-module map, then ker $\mu \subset$ ker $\mu (\psi \otimes N)$.
\end{lemma}
\begin{proof}
Since $\psi$ is an $A\otimes B$-module morphism, $\psi$ is an $A$-module morphism.  Since $M$ is a free $A$-module of rank one, $\psi$ is just multiplication by an element of $A$, $a$.  Suppose $c \in \operatorname{ker }\mu$ has the form $c=\Sigma_{i}m_{i}\otimes n_{i} \in M \otimes N$.  Then
$$
\mu(\psi \otimes N)(c)=\mu(\Sigma_{i}(a\otimes 1) \cdot m_{i} \otimes n_{i})=\mu(\Sigma_{i} m_{i} \otimes (a \otimes 1) \cdot n_{i}) =
$$
$$
\mu(\Sigma_{i} m_{i} \otimes n_{i} \cdot (a \otimes 1))=\mu(\Sigma_{i} m_{i} \otimes n_{i})\cdot (a \otimes 1) = 0.
$$
\end{proof}

\begin{proposition} \label{prop.preker}
Suppose $\mathcal{L}$ is an invertible ${\mathcal{O}}_{U}$-module, $\mathcal{B}$ is an ${\mathcal{O}}_{(U \times X)_{U}^{2}}$-module, $\mathcal{P}$ is an ${\mathcal{O}}_{U\times X}$-module and $q:U \rightarrow X$ is a morphism.  Let $\pr_{i}:(U \times X)_{U}^{2} \rightarrow U \times X$ be the standard projection maps.  If
$$
\mu: \pr_{2*}(\pr_{1}^{*}(\id_{U} \times q)_{*}{\mathcal{L}} \otimes_{{\mathcal{O}}_{X}} \mathcal{B}) \rightarrow \mathcal{P}
$$
and
$$
\psi:(\id_{U} \times q)_{*}{\mathcal{L}} \rightarrow (\id_{U} \times q)_{*}{\mathcal{L}}
$$
are morphisms, then ker $\mu \subset$ ker $\mu (\psi \otimes \mathcal{B})$.
\end{proposition}

\begin{proof}
We reduce the statement of the Proposition to that of Lemma \ref{lem.reduce}.  Suppose $V, W \subset X$ are open affine.  Let $V' \subset q^{-1}(V)$ be an affine open set such that $\mathcal{L}$ restricted to $V'$ is free.  We have
\begin{equation} \label{eqn.bigex1}
\pr_{2*}(\pr_{1}^{*}(\id_{U} \times q)_{*}\mathcal{L} \otimes \mathcal{B}) (V' \times W) =
\end{equation}
$$
\pr_{1}^{*}(\id_{U} \times q)_{*}\mathcal{L} \otimes \mathcal{B} ((U \times X) \times_{U} (V' \times W)).
$$
By Lemma \ref{lem.sup2} part 2, and Lemma \ref{lem.sup1}, (\ref{eqn.bigex1}) is isomorphic, as ${\mathcal{O}}_{U \times X}(V' \times W)$-modules, to
\begin{equation} \label{eqn.bigex2}
\pr_{1}^{*}(\id_{U} \times q)_{*}\mathcal{L} \otimes \mathcal{B} ((U \times V) \times_{U} (V' \times W)).
\end{equation}
Since $(U \times V) \times_{U} (V' \times W)$ is affine, (\ref{eqn.bigex2}) equals
$$
\pr_{1}^{*}(\id_{U} \times q)_{*}\mathcal{L} ((U \times V) \times_{U} (V' \times W)) \otimes \mathcal{B} ((U \times V) \times_{U} (V' \times W))
$$
which, in turn, is easily seen to be naturally isomorphic to
$$
(\id_{U} \times q)_{*}\mathcal{L}(V' \times V) \otimes_{{\mathcal{O}}_{U \times X}(V' \times V)} \mathcal{B}((V' \times V) \times_{U} (V' \times W)).
$$
The result now follows from Lemma \ref{lem.reduce}.
\end{proof}

\begin{corollary} \label{cor.extend}
Let $\mathcal{P}$ be a $U$-family in standard form\index{family!standard}.  Any ${\mathcal{O}}_{U \times X}$-module automorphism, $\psi_{0}$, of ${\mathcal{P}}_{0}$ can be extended to a ${\mathcal{B}}^{U}$-module automorphism of $\mathcal{P}$.
\end{corollary}

\begin{proof}
Let $\mu_{i,1}:{\mathcal{P}}_{i} \otimes {\mathcal{B}}_{1} \rightarrow {\mathcal{P}}_{i+1}$ be the multiplication map, and let $\psi_{0}'$ be the restriction of $\psi_{0} \otimes {\mathcal{B}}_{1}^{U}$ to ker $\mu_{0,1}$.  By Proposition \ref{prop.preker}, the following diagram, whose rows are exact
$$
\xymatrix{
0 \ar[r] & \mbox{ker }\mu_{0,1}  \ar[r] \ar[d]_{\psi_{0}'} & {\mathcal{P}}_{0} \otimes {\mathcal{B}}_{1}^{U} \ar[r]^{\mu_{0}} \ar[d]^{\psi_{0} \otimes {\mathcal{B}}_{1}^{U}} & {\mathcal{P}}_{1} \ar[r] & 0 \\
0 \ar[r] & \mbox{ker }\mu_{0,1}  \ar[r] & {\mathcal{P}}_{0} \otimes {\mathcal{B}}_{1}^{U} \ar[r]^{\mu_{0}} & {\mathcal{P}}_{1} \ar[r] & 0
}
$$
commutes.  Thus, there exists a unique ${\mathcal{O}}_{U \times X}$-module automorphism, $\psi_{1}$ of ${\mathcal{P}}_{1}$, completing the diagram.  Continuing in this way, we have an ${\mathcal{O}}_{U \times X}$-module automorphism $\psi$ whose $i$th component is $\psi_{i}$, with the property that the diagram
$$
\xymatrix{
{{\mathcal{P}}_{i} \otimes {\mathcal{B}}_{1}^{U}} \ar[r]^{\mu_{i}} \ar[d]_{\psi_{i} \otimes {\mathcal{B}}_{1}^{U}} & {\mathcal{P}}_{i+1} \ar[d]^{\psi_{i+1}} \\
{\mathcal{P}}_{i} \otimes {\mathcal{B}}_{1}^{U} \ar[r]_{\mu_{i}} & {\mathcal{P}}_{i+1}
}
$$
commutes.  By Lemma \ref{lem.gen1}, $\psi$ is a ${\mathcal{B}}^{U}$-module automorphism, as desired.
\end{proof}

\begin{corollary} \label{cor.zero}
If $\psi: \mathcal{M} \rightarrow \mathcal{N}$ is an isomorphism of $U$-families and if
$$
\gamma:{\mathcal{M}}_{0} \rightarrow {\mathcal{N}}_{0}
$$
is an isomorphism of ${\mathcal{O}}_{U\times X}$-modules, then there exists an isomorphism $\theta: \mathcal{M} \rightarrow \mathcal{N}$ of $U$-families such that $\theta_{0}=\gamma$.
\end{corollary}

\begin{proof}
Suppose $\mathcal{N}$ is canonically isomorphic to the $U$-family in standard form $\mathcal{P}$.  Now, $\psi_{0} \gamma^{-1}:{\mathcal{N}}_{0} \rightarrow {\mathcal{N}}_{0}$, so gives an automorphism ${\mathcal{P}}_{0} \rightarrow {\mathcal{P}}_{0}$.  By Corollary \ref{cor.extend}, this automorphism extends to an automorphism of $\mathcal{P}$, which gives an automorphism $\delta$ of $\mathcal{N}$, such that $\delta_{0}=\psi_{0} \gamma^{-1}$.  Now let $\theta=\delta^{-1} \psi$.  Then $\theta_{0}=\gamma \psi_{0}^{-1} \psi_{0}=\gamma$, as desired.
\end{proof}

\subsection{$\Gamma_{n}$ is compatible with descent}
We now set up notation which we will use for the rest of this
chapter. Let $p:U \rightarrow S$ and $q:V \rightarrow S$ be
affine, noetherian schemes, and let $f:V \rightarrow U$ be a
faithfully flat morphism.  Let $f_{i}: V_{U}^{2} \rightarrow V$
and $f_{ij}: V_{U}^{3} \rightarrow V_{U}^{2}$ be projection maps.
For ease of reading, we alter our notation convention for the rest
of the chapter by letting $e=f \times \id_{X}$, $e_{i}=f_{i}
\times \id_{X}$, and $e_{ij}=f_{ij} \times \id_{X}$.  We have the
following commuting diagram:
$$
\xymatrix{
V \times_{U} V \times_{U} V \times X \ar[d] \ar@<2ex>[r]^{e_{12}} \ar[r]^{e_{23}} \ar@<-2ex>[r]^{e_{13}} & V \times_{U} V \times X \ar[d] \ar@<1ex>[r]^{e_{1}} \ar@<-1ex>[r]^{e_{2}} & V \times X \ar[d] \ar[r]^{e} & U \times X \ar[d] \\
V \times_{U} V \times_{U} V \ar@<2ex>[r]^{f_{12}} \ar[r]^{f_{23}} \ar@<-2ex>[r]^{f_{13}} & V \times_{U} V \ar@<1ex>[r]^{f_{1}} \ar@<-1ex>[r]^{f_{2}} & V \ar[r]^{f} & U
}
$$
Since $f$ is faithfully flat and locally finite, so is $e$, and if
$$
z'':Y''=(V \times X)_{U\times X}^{2} \rightarrow V_{U}^{2} \times X=Z'',
$$
$$
z''':Y'''=(V \times X)_{U\times X}^{3} \rightarrow V_{U}^{3} \times X=Z'''
$$
and
$$
h_{ij}:(V \times X)_{U\times X}^{3} \rightarrow (V \times X)_{U \times X}^{2}
$$
are projections, then the diagram
$$
\xymatrix{
Y''' \ar[r]^{h_{ij}} \ar[d]_{z'''} & Y'' \ar[d]^{z''} \ar[r]^{h_{k}} & Y' \ar[d]^{=} \\
Z''' \ar[r]_{e_{ij}} & Z'' \ar[r]_{e_{k}} & Y'
}
$$
commutes.  Thus, Corollary \ref{cor.descent} may be employed to prove the following

\begin{proposition} \label{prop.descentc}
With the notation above, the functor $\mathcal{H} \mapsto e^{*}\mathcal{H}$ which goes from ${\mathcal{O}}_{U \times X}$-modules to ${\mathcal{O}}_{V \times X}$-modules with modified descent data,\index{category of oyprime modules with descent data@category of ${\mathcal{O}}_{Y'}$-modules with descent data!modified} is an equivalence of categories.
\end{proposition}
We now reintroduce some notation.  Let $\tilde{p}=p \times \id_{X}:U \times X \rightarrow S \times X$, and for any $U \in \sf{S}$, let $\pr_{i}^{U}:(U \times X)_{U}^{2} \rightarrow U \times X$ be the standard projection maps.  Let $\mathcal{F}$ be an ${\mathcal{O}}_{X}$-bimodule.  To each $U \in \sf{S}$, let
$$
F^{U} = \pr_{2*}^{U}(\pr_{1}^{U*}(-) \otimes {\tilde{p}}^{2*}\mathcal{F}),
$$
and suppose, for $f:V \rightarrow U$ in $\sf{S}$,
\begin{equation} \label{eqn.theta}
\Lambda_{f}:{\tilde{f}}^{*}F^{U} \Longrightarrow F^{V}{\tilde{f}}^{*}
\end{equation}
is the isomorphism defined in Example \ref{example.functors}.  This assignment defines an indexed functor $\mathfrak{F}:\mathfrak{X} \rightarrow \mathfrak{X}$, where $\mathfrak{X}$ is the $\sf{S}$-indexed category defined in Example \ref{example.spaces}, such that ${\sf{X}}^{U}={\sf{Mod }}(U \times X)$.

\subsection{Injectivity of $\Gamma_{n}(U) \rightarrow \Gamma_{n}(V)$}
\begin{lemma}
Let $\mathcal{M}$ and $\mathcal{N}$ be truncated $U$-families.  Suppose $\theta: {e}^{*}\mathcal{M} \rightarrow {e}^{*}\mathcal{N}$ is a ${\mathcal{B}}^{V}$-module isomorphism.  Suppose ${e_{1}}^{*}\theta_{0}={e_{2}}^{*}\theta_{0}$.  Then ${e_{1}}^{*}\theta={e_{2}}^{*}\theta$.
\end{lemma}

\begin{proof}
Since ${e_{1}}^{*}\theta_{0}={e_{2}}^{*}\theta_{0}$, there exists, by Proposition $\ref{prop.descentc}$, a $\psi_{0}:{\mathcal{M}}_{0} \rightarrow {\mathcal{N}}_{0}$ such that ${e}^{*}\psi_{0}=\theta_{0}$.

Since there is a morphism of functors ${e}^{*}F^{U} \Longrightarrow F^{V}{e}^{*}$ as in (\ref{eqn.theta}), the diagram
\begin{equation} \label{eqn.topguy}
\xymatrix{
{e}^{*}F^{U}{\mathcal{M}}_{0} \ar[r]^{{e}^{*}F^{U}\psi_{0}} \ar[d] & {e}^{*}F^{U}{\mathcal{N}}_{0} \ar[d] \\
F^{V}{e}^{*}{\mathcal{M}}_{0} \ar[r]_{F^{V}{e}^{*}\psi_{0}}  & F^{V}{e}^{*}{\mathcal{N}}_{0}
}
\end{equation}
commutes, and the columns are isomorphisms (since they are just the maps (\ref{eqn.theta})).  Since $\theta$ is a $\mathcal{B}$-module map, the diagram

\begin{equation} \label{eqn.bottomguy}
\xymatrix{
{{e}_{i}}^{*}F^{V}{e}^{*}{\mathcal{M}}_{0} \ar[r]^{{{e}_{i}}^{*}F^{V}\theta_{0}} \ar[d] & {{e}_{i}}^{*}F^{V}{e}^{*}{\mathcal{N}}_{0} \ar[d] \\
{{e}_{i}}^{*}{e}^{*}{\mathcal{M}}_{1} \ar[r]_{{{e}_{i}}^{*}\theta_{1}} & {{e}_{i}}^{*}{e}^{*}{\mathcal{M}}_{1}
}
\end{equation}
commutes.  Applying $e_{i}^{*}$ to (\ref{eqn.topguy}) and stacking (\ref{eqn.topguy}) above (\ref{eqn.bottomguy}) gives a commutative diagram:
$$
\xymatrix{
{{e}_{i}}^{*}{e}^{*}F^{U}{\mathcal{M}}_{0} \ar[rrr]^{{{e}_{i}}^{*}{e}^{*}F^{U}\psi_{0}} \ar[d] & & & {{e}_{i}}^{*}{e}^{*}F^{U}{\mathcal{N}}_{0} \ar[d] \\
{{e}_{i}}^{*}{e}^{*}{\mathcal{M}}_{1} \ar[rrr]_{{{e}_{i}}^{*}\theta_{1}} & & & {{e}_{i}}^{*}{e}^{*}{\mathcal{M}}_{1}.
}
$$
Since multiplication $F^{V}e^{*}{\mathcal{M}}_{0} \rightarrow e^{*}{\mathcal{M}}_{1}$ is defined as the composition
$$
F^{V}e^{*}{\mathcal{M}}_{0} \rightarrow e^{*}F^{U}{\mathcal{M}}_{0} \rightarrow e^{*}{\mathcal{M}}_{1}
$$
the verticals of this diagram are ${{e}_{i}}^{*}{e}^{*}$ applied to the multiplication, $\mu$, of the $\mathcal{B}$-modules $\mathcal{M}$ and $\mathcal{N}$.  This shows that the left and right face of
$$
\xymatrix{
& & {e_{1}}^{*}{e}^{*}F^{U}{\mathcal{M}}_{0} \ar[r] \ar[dll] \ar[dd] & {e_{1}}^{*}{e}^{*}F^{U}{\mathcal{N}}_{0} \ar[dd] \ar[dll] \\
{e_{2}}^{*}e^{*}F^{U}{\mathcal{M}}_{0} \ar[r] \ar[dd] & {e_{2}}^{*}{e}^{*}F^{U}{\mathcal{N}}_{0} \ar[dd] & & \\
& & {e_{1}}^{*}{e}^{*}{\mathcal{M}}_{1} \ar[dll] \ar[r] & {e_{1}}^{*}e^{*}{\mathcal{N}}_{1} \ar[dll] \\
{e_{2}}^{*}{e}^{*}{\mathcal{M}}_{1} \ar[r] & {e_{2}}^{*}{e}^{*}{\mathcal{N}}_{1} & &
}
$$
commute.  Since ${{e}_{1}}^{*}{e}^{*}\mu={{e}_{2}}^{*}{e}^{*}\mu$ and
$$
{e_{1}}^{*}{e}^{*}F^{U}\psi_{0}={e_{2}}^{*}{e}^{*}F^{U}\psi_{0},
$$
the top and sides of the cube commute as well.  Since the vertical maps are epis, the bottom is commutative so ${{e}_{1}}^{*}\theta_{1}={{e}_{2}}^{*}\theta_{1}$.  Using this fact and repeating the above argument, we find that ${e_{1}}^{*}\theta={e_{2}}^{*}\theta$.
\end{proof}
Let us prove part of Theorem $\ref{theorem.comp1}$.
\begin{proposition}
With the notation as in the beginning of this section, $\Gamma_{n}(f):\Gamma_{n}(U) \rightarrow \Gamma_{n}(V)$ is an injection.
\end{proposition}
\begin{proof}
Let $\mathcal{M}$ and $\mathcal{N}$ be $U$-families.  By Lemma $\ref{lem.standard}$, we may assume $\mathcal{M}$ and $\mathcal{N}$ are in standard form.  Suppose $[{e}^{*}\mathcal{M}]=[{e}^{*}\mathcal{N}]$, i.e. suppose ${e}^{*}\mathcal{M} \cong {e}^{*}\mathcal{N}$ as $V$-families.  By Corollary $\ref{cor.zero}$, there exists a $\theta: {e}^{*}\mathcal{M} \cong {e}^{*}\mathcal{N}$ such that $\theta_{0}$ is the identity.  Thus ${e_{1}}^{*}\theta_{0}={e_{2}}^{*}\theta_{0}$.  By the previous lemma, ${e_{1}}^{*}\theta={e_{2}}^{*}\theta$.  Thus, by Proposition $\ref{prop.descentc}$, there exists an isomorphism $\delta:\mathcal{M} \rightarrow \mathcal{N}$ such that ${e}^{*}\delta=\theta$.  Finally, we must show that $\delta$ is a $\mathcal{B}$-module morphism.  If $\mu'$ denotes the multiplication maps of the modules $e^{*}\mathcal{M}$ and $e^{*}\mathcal{N}$, then we know
$$
\xymatrix{
F^{V}{e}^{*}{\mathcal{M}}_{i} \ar[d]_{\mu_{i}'} \ar[r]^{F^{V}\theta_{i}} & F^{V}{e}^{*}{\mathcal{N}}_{i} \ar[d]^{\mu_{i}'} \\
{e}^{*}{\mathcal{M}}_{i+1} \ar[r]_{\theta_{i+1}} & {e}^{*}{\mathcal{N}}_{i+1}
}
$$
commutes since $\theta$ is a $\mathcal{B}$-module morphism.  By naturality of the indexing map (\ref{eqn.theta}), the diagram
$$
\xymatrix{
{e}^{*}F^{U}{\mathcal{M}}_{i} \ar[d] \ar[r]^{{e}^{*}F^{U}\delta_{i}} & {e}^{*}F^{U}{\mathcal{N}}_{i} \ar[d] \\
F^{V}{e}^{*}{\mathcal{M}}_{i} \ar[r]_{F^{V}\theta_{i}} & F^{V}{e}^{*}{\mathcal{N}}_{i}
}
$$
commutes, and when this diagram is placed above the previous diagram, the columns are ${e}^{*}$ applied to the multiplication map of the $\mathcal{B}$-modules $\mathcal{M}$ and $\mathcal{N}$ by definition of $\mu'$.  Since the composition of the previous two diagrams commute, and it equals the application of the functor ${e}^{*}$ to a diagram of ${\mathcal{O}}_{U \times X}$-modules, the diagram commutes after removing $e^{*}$ by Proposition $\ref{prop.descentc}$.  Thus $\delta$ is a $\mathcal{B}$-module morphism.
\end{proof}

\subsection{Completing the proof that $\Gamma_{n}(-)$ is compatible with descent.}

\begin{lemma} \label{lem.rigid}
Let $\mathcal{M}$ be a $V$-family such that ${e_{1}}^{*}\mathcal{M} \cong {e_{2}}^{*}\mathcal{M}$ as ${\mathcal{B}}^{V \times_{U} V}$-modules.  Then, for some $q:U \rightarrow X$, ${\mathcal{M}}_{0} \cong {e}^{*}(\id_{U} \times q)_{*} {\mathcal{O}}_{U}$.
\end{lemma}
\begin{proof}
Without loss of generality, we may assume $\mathcal{M}$ is in standard form.  Thus, for some $q':V \rightarrow X$, ${\mathcal{M}}_{0}=(\id_{V} \times q')_{*}{\mathcal{O}}_{V}$.  We may now compute

\begin{equation} \label{eqn.lasty}
\xymatrix{
{e_{1}}^{*}{\mathcal{M}}_{0}={e_{1}}^{*}(\id_{V} \times q')_{*}{\mathcal{O}}_{V} \ar[r] & (\id_{V \times_{U} V} \times q'f_{1})_{*}{f_{1}}^{*}{\mathcal{O}}_{V}
}
\end{equation}
where the right morphism, the dual 2-cell induced by the diagram of schemes
$$
\xymatrix{
V \times_{U}V \ar[r]^{f_{1}} \ar[d]_{\id_{V \times_{U}V} \times q'f_{1}} & V \ar[d]^{\id_{V} \times q'} \\
V \times_{U}V \times X \ar[r]_{e_{1}} & V \times X,
}
$$
is an isomorphism by Lemma $\ref{lem.anothersquare}$.  Furthermore,

$$
{e_{1}}^{*}(\id_{V} \times q')_{*}{\mathcal{O}}_{V} \cong {e_{2}}^{*}(\id_{V} \times q')_{*}{\mathcal{O}}_{V} \cong (\id_{V \times_{U} V} \times q'f_{2})_{*}{f_{2}}^{*}{\mathcal{O}}_{V}
$$
where the first isomorphism exists by hypothesis and the second is an isomorphism as in (\ref{eqn.lasty}).  Thus, by Proposition $\ref{prop.uniqueness}$, $q'f_{1}=q'f_{2}$ so by Theorem $\ref{theorem.comp1}$, there exists a $q:U \rightarrow X$ such that $qf=q'$.  Since $f^{*}{\mathcal{O}}_{U} \cong {\mathcal{O}}_{V}$,
$$
{\mathcal{M}}_{0} = (\id_{V} \times qf)_{*}{\mathcal{O}}_{V} \cong {e}^{*}(\id_{U} \times q)_{*}{\mathcal{O}}_{U},
$$
where the last isomorphism is the dual 2-cell induced by the diagram
$$
\xymatrix{
V \ar[r]^{f} \ar[d]_{\id_{V} \times qf} & U \ar[d]^{\id_{U} \times q} \\
V \times X \ar[r]_{e} & U \times X.
}
$$
\end{proof}
\begin{lemma} \label{lem.mzero}
Let $\mathcal{M}$ be a $V$-family, and suppose ${e_{1}}^{*}\mathcal{M} \cong {e_{2}}^{*}\mathcal{M}$ as ${\mathcal{B}}^{V\times_{U}V}$-modules.  Then there exists a $V$-family $\mathcal{M}' \cong \mathcal{M}$ and an isomorphism $\theta:{e_{1}}^{*}\mathcal{M}' \rightarrow {e_{2}}^{*}\mathcal{M}'$ such that $\theta_{0}$ is a descent datum.\index{descent datum}
\end{lemma}

\begin{proof}
By Lemma $\ref{lem.rigid}$, there is a $V$-family $\mathcal{M}' \cong \mathcal{M}$ with ${\mathcal{M}'}_{0} = {e}^{*}(\id_{U} \times q)_{*}{\mathcal{O}}_{U}$.  Thus $e_{1}^{*}{\mathcal{M}}_{0} \cong e_{1}^{*}{\mathcal{M}} \cong e_{2}^{*}\mathcal{M} \cong e_{2}^{*}{\mathcal{M}}_{0}$.  Since the canonical isomorphism

$$
{e_{1}}^{*}{e}^{*}(\id_{U} \times q)_{*}{\mathcal{O}}_{U} \cong {e_{2}}^{*}{e}^{*}(\id_{U} \times q)_{*}{\mathcal{O}}_{U}
$$
is a descent datum, the result follows from Corollary $\ref{cor.zero}$.
\end{proof}

\begin{lemma} \label{lem.hardnat}
Let $\mathcal{L}$ be an invertible ${\mathcal{O}}_{U}$-module and let $\mathcal{N}=(\id_{U} \times q)_{*}{\mathcal{L}}$.  Then the diagram
$$
\xymatrix{
{e_{1}}^{*}{e}^{*}F^{U}\mathcal{N} \ar[r] \ar[d] & F^{V \times_{U} V}{e_{1}}^{*}{e}^{*}\mathcal{N} \ar[d]  \\
{e_{2}}^{*}{e}^{*}F^{U}\mathcal{N} \ar[r] & F^{V \times_{U} V}{e_{2}}^{*}{e}^{*}\mathcal{N}
}
$$
whose horizontal maps are just compositions of maps (\ref{eqn.theta}) (on page \pageref{eqn.theta}), commutes.
\end{lemma}
\begin{proof}
We must show that the diagram
$$
\xymatrix{
{e_{1}}^{*}{e}^{*}F^{U}\mathcal{N} \ar[r] \ar[d] & F^{V \times_{U} V}{e_{1}}^{*}{e}^{*}\mathcal{N} \ar[d]  \\
(e{e_{1}})^{*}F^{U}\mathcal{N} \ar[r] \ar[d] & F^{V \times_{U} V}({ee_{1}})^{*}\mathcal{N} \ar[d]  \\
({ee_{2}})^{*}F^{U}\mathcal{N} \ar[r] \ar[d] & F^{V \times_{U} V}({ee_{2}})^{*}\mathcal{N} \ar[d] \\
{e_{2}}^{*}{e}^{*}F^{U}\mathcal{N} \ar[r] & F^{V \times_{U} V}{e_{2}}^{*}{e}^{*}\mathcal{N}
}
$$
commutes, where the vertical maps are the natural maps.  Since $F$ is indexed by Lemma \ref{lem.isommonoid}, the top and bottom squares commute.  Since the middle verticals are equalities, the middle square also commutes and the assertion follows.
\end{proof}

\begin{proposition}
Let $\mathcal{M}$ be a $V$-family, and suppose there is a ${\mathcal{B}}^{V \times_{U}V}$-module isomorphism $\theta:{e_{1}}^{*}\mathcal{M} \rightarrow {e_{2}}^{*}\mathcal{M}$.  Then there exists a $U$-family $\mathcal{N}$ such that ${e}^{*}\mathcal{N} \cong \mathcal{M}$ as ${\mathcal{B}}^{V}$-modules.
\end{proposition}

\begin{proof}
By Lemma $\ref{lem.mzero}$, we may assume that there exists a $q:U \rightarrow X$ such that ${\mathcal{M}}_{0}={e}^{*}(\id_{U}\times q)_{*}{\mathcal{O}}_{U}$ and that $\theta_{0}$ is a canonical descent datum.  By Lemma $\ref{lem.hardnat}$ and the fact that $\theta$ is a ${\mathcal{B}}^{V}$-module morphism, the diagram
$$
\xymatrix{
{e_{1}}^{*}{e}^{*}F^{U}(\id_{U} \times q)_{*}{\mathcal{O}}_{U} \ar[r] \ar[d] & {e_{2}}^{*}{e}^{*}F^{U}(\id_{U} \times q)_{*}{\mathcal{O}}_{U} \ar[d] \\
F^{V \times_{U} V}{e_{1}}^{*}{\mathcal{M}}_{0} \ar[r]^{F^{V \times_{U} V}\theta_{0}} \ar[d]_{\mu_{0}} & F^{V \times_{U} V}{e_{2}}^{*}{\mathcal{M}}_{0} \ar[d]^{\mu_{0}} \\
{e_{1}}^{*}{\mathcal{M}}_{1} \ar[r]_{\theta_{1}} & {e_{2}}^{*}{\mathcal{M}}_{1}
}
$$
commutes, where the top vertical maps are indexing maps (\ref{eqn.theta}) on page \pageref{eqn.theta}.  Let ${\mathcal{N}}_{0}=(\id_{U} \times q)_{*}{\mathcal{O}}_{U}$.  Expanding the bottom square, we get a diagram
$$
\xymatrix{
{e_{1}}^{*}{e}^{*}F^{U}{\mathcal{N}}_{0} \ar[rr] \ar[d] & & {e_{2}}^{*}{e}^{*}F^{V}{\mathcal{N}}_{0} \ar[d] \\
F^{V \times_{U} V}{e_{1}}^{*}{e}^{*}{\mathcal{N}}_{0} \ar[rr]^{F^{V \times_{U} V}\theta_{0}} \ar[d] & & F^{V \times_{U} V}{e_{2}}^{*}{e}^{*}{\mathcal{N}}_{0} \ar[d] \\
{e_{1}}^{*}F^{V}{e}^{*}{\mathcal{N}}_{0} \ar[d]_{{e_{1}}^{*}\mu} & & {e_{2}}^{*}F^{V}{e}^{*}{\mathcal{N}}_{0} \ar[d]^{{e_{2}}^{*}\mu} \\
{e_{1}}^{*}{\mathcal{M}}_{1} \ar[rr]_{\theta_{1}} &  & {e_{2}}^{*}{\mathcal{M}}_{1}
}
$$
where, again, the top two pairs of vertical maps are indexing maps (\ref{eqn.theta}).  But by definition of these vertical maps, this diagram equals
$$
\xymatrix{
{e_{1}}^{*}{e}^{*}F^{U}{\mathcal{N}}_{0} \ar[r] \ar[d] & {e_{2}}^{*}{e}^{*}F^{U}{\mathcal{N}}_{0} \ar[d] \\
{e_{1}}^{*}F^{V}{e}^{*}{\mathcal{N}}_{0} \ar[d]_{{e_{1}}^{*}\mu} & {e_{2}}^{*}F^{V}{e}^{*}{\mathcal{N}}_{0} \ar[d]^{{e_{2}}^{*}\mu} \\
{e_{1}}^{*}{\mathcal{M}}_{1} \ar[r]_{\theta_{1}} & {e_{2}}^{*}{\mathcal{M}}_{1}
}
$$
whose vertical maps are ${e_{i}}^{*}$ applied to a morphism
$$
{e}^{*}F^{U}{\mathcal{N}}_{0} \rightarrow {\mathcal{M}}_{1}.
$$
Since the top morphism of this diagram is a descent datum, $\theta_{1}$ is a descent datum by Lemma $\ref{lem.inher}$ and there is an isomorphism $\psi:{\mathcal{M}}_{1} \rightarrow {e}^{*}{\mathcal{N}}_{1}$ for some ${\mathcal{O}}_{U \times X}$-module ${\mathcal{N}}_{1}$ by Proposition $\ref{prop.descentc}$.  Therefore, we have a commutative diagram
$$
\xymatrix{
{e_{1}}^{*}{e}^{*}F^{U}{\mathcal{N}}_{0} \ar[r] \ar[d] & {e_{2}}^{*}{e}^{*}F^{U}{\mathcal{N}}_{0} \ar[d] \\
{e_{1}}^{*}F^{V}{e}^{*}{\mathcal{N}}_{0} \ar[d]_{{e_{1}}^{*}\psi\mu} & {e_{2}}^{*}F^{V}{e}^{*}{\mathcal{N}}_{0} \ar[d]^{{e_{2}}^{*}\psi\mu} \\
{e_{1}}^{*}{e}^{*}{\mathcal{N}}_{1} \ar[r] & {e_{2}}^{*}{e}^{*}{\mathcal{N}}_{1}
}
$$
where the bottom horizontal map equals $({e_{2}}^{*}\psi)\theta_{1}({e_{1}}^{*}\psi)^{-1}$ and where the vertical maps are ${e_{i}}^{*}$ applied to a morphism
$$
\phi:{e}^{*}F^{U}{\mathcal{N}}_{0} \rightarrow {e}^{*}{\mathcal{N}}_{1}.
$$
Since $\theta_{1}$ is a descent datum, so is $({e_{2}}^{*}\psi)\theta_{1}({e_{1}}^{*}\psi)^{-1}$, so that by Proposition $\ref{prop.descentc}$, there exists a $\mu':F^{U}{\mathcal{N}}_{0} \rightarrow {\mathcal{N}}_{1}$ such that ${e}^{*}\mu'=\phi$.  Thus we have a commutative diagram
$$
\xymatrix{
{e_{1}}^{*}{e}^{*}F^{U}{\mathcal{N}}_{0} \ar[r] \ar[d]_{{e_{1}}^{*}{e}^{*}\mu'} & {e_{2}}^{*}{e}^{*}F^{U}{\mathcal{N}}_{0} \ar[d]^{{e_{2}}^{*}{e}^{*}\mu'} \\
{e_{1}}^{*}{e}^{*}{\mathcal{N}}_{1} \ar[r] & {e_{2}}^{*}{e}^{*}{\mathcal{N}}_{1}.
}
$$
Since $e_{i}$ and $e$ are faithfully flat, ${e_{i}}^{*}$ and ${e}^{*}$ are exact.  Thus, ker ${e_{i}}^{*}{e}^{*}\mu' = {e_{i}}^{*}{e}^{*}\operatorname{ker } \mu'$ and $\mu'$ is epi since $e^{*}\mu'$ is epi.  Since $\mu'$ is an epi, there exists a {\it unique} morphism
$$
{e_{1}}^{*}{e}^{*}{\mathcal{N}}_{1} \rightarrow {e_{2}}^{*}{e}^{*}{\mathcal{N}}_{1}
$$
making the diagram commute.  Since the top morphism is the canonical one, we conclude that $\theta_{1}$ is also the canonical morphism (Definition \ref{def.candesdat}).  The proof follows by induction.
\end{proof}

Thus, we have succeeded in proving Theorem $\ref{theorem.comp1}$.

\chapter{The Representation of $\Gamma_{n}$ for Low $n$} In this
chapter, we assume, unless otherwise stated, that $X$ is a
separated, noetherian scheme, $\mathcal{E}$ is a coherent
${\mathcal{O}}_{X}$-bimodule, and there is a graded ideal
$\mathcal{I} \subset T(\mathcal{E})$ whose first nonzero component
occurs in degree $m > 1$.  Let $\mathcal{B} =
T(\mathcal{E})/\mathcal{I}$.  We represent $\Gamma_{n}$ for $n <
m$ (Proposition \ref{prop.zeron} and Theorem \ref{theorem.rep1}).
Before reading the remainder of this chapter, the reader is
advised to review the last two sentences of Section 1.6.
\section{The representation of $\Gamma_{0}$}
\begin{proposition} \label{prop.zeron}
$\Gamma_{0}(-) \cong \operatorname{Hom}_{S}(-,X)$.
\end{proposition}

\begin{proof}
Let $U$ be an affine scheme.  Then
$$
\Gamma_{0}(U) = \{\mbox{isomorphism classes of truncated $U$-families of length 1} \}.
$$
By Proposition $\ref{prop.uniqueness}$, these classes are indexed by elements of $\operatorname{Hom}_{S}(U,X)$:  if $\mathcal{M}$ is a truncated $U$-family of length one, $\mathcal{M} \cong (\id_{U} \times q)_{*}{\mathcal{O}}_{U}$ for a unique map $q:U \rightarrow X$.  This gives a bijection
$$
\sigma_{U}:\Gamma_{0}(U) \rightarrow \operatorname{Hom}_{S}(U,X)
$$
defined by $\sigma_{U}([\mathcal{M}])=q$.  We show $\sigma$ is natural.  Let $q:U \rightarrow X$, so that $[(\id_{U} \times q)_{*}{\mathcal{O}}_{U}] \in \Gamma_{0}(U)$ and let $f:V \rightarrow U$ be a morphism of affine $S$-schemes.  Then
$$
\Gamma_{0}(f)[(\id_{U} \times q)_{*}{\mathcal{O}}_{U}] = [(\id_{V} \times qf)_{*}{\mathcal{O}}_{V}]
$$
so that
$$
\xymatrix{
\operatorname{Hom}_{S}(U,X) \ar[d]_{- \circ f} \ar[r]^{\sigma_{U}} & \Gamma_{0}(U) \ar[d]^{\Gamma_{0}(f)}\\
\operatorname{Hom}_{S}(V,X) \ar[r]_{\sigma_{V}} & \Gamma_{0}(V)
}
$$
commutes.  Thus, $\sigma$ defines an isomorphism
$$
\Sigma:\Gamma_{0}(-) \Longrightarrow \operatorname{Hom}_{S}(-,X).
$$
\end{proof}

\section{The representation of $\Gamma_{n}$ for $0 < n < m$}
Before stating the main result of this chapter (Theorem
\ref{theorem.rep1}), we introduce some notation.  Suppose $Y$ and
$Z$ are separated, noetherian schemes, $\mathcal{E}$ is an
${\mathcal{O}}_{X \times Y}$-module and $\mathcal{F}$ is an
${\mathcal{O}}_{Y \times Z}$-module.  Suppose also that $\pr_{i}:X
\times Y \rightarrow X,Y$, $\pr_{i}:Y \times Z \rightarrow Y, Z$
(same notation!) and $\pr_{ij}:X \times Y \times Z \rightarrow X
\times Y, X \times Z, Y \times Z$ are projection maps, and
$p_{\mathcal{E}}$ and $p_{\mathcal{F}}$ are the structure maps of
${\mathbb{P}}_{X \times Y}(\mathcal{E})$ and ${\mathbb{P}}_{Y
\times Z}(\mathcal{F})$ respectively.  We denote by
$$
{\mathbb{P}}_{X \times Y}(\mathcal{E}) \otimes_{Y} {\mathbb{P}}_{Y \times Z}(\mathcal{F})
$$
the fiber product in the following diagram:
$$
\xymatrix
{
{\mathbb{P}}_{X \times Y}(\mathcal{E}) \otimes_{Y} {\mathbb{P}}_{Y \times Z}(\mathcal{F}) \ar@{-->}[rr]^{q_{\mathcal{F}}} \ar@{-->}[d]_{q_{\mathcal{E}}} & & {\mathbb{P}}_{Y \times Z}(\mathcal{F}) \ar[d]^{\pr_{1} \circ p_{\mathcal{F}}} \\
{\mathbb{P}}_{X \times Y}(\mathcal{E}) \ar[rr]_{\pr_{2} \circ p_{\mathcal{E}}} & & Y.
}
$$
Since $X \times Y \times Z$ is the fiber product of the following diagram:
$$
\xymatrix
{
& Y \times Z \ar[d]^{\pr_{1}} \\
X \times Y \ar[r]_{\pr_{2}} & Y
}
$$
we have the commutative diagram
$$
\xymatrix
{
& {\mathbb{P}}_{X \times Y}(\mathcal{E}) \otimes_{Y} {\mathbb{P}}_{Y \times Z}(\mathcal{F})   \ar[dr]^{q_{\mathcal{F}}} \ar[dl]_{q_{\mathcal{E}}}  & \\
{\mathbb{P}}_{X \times Y}(\mathcal{E}) \ar[dd]^{p_{\mathcal{E}}} & & {\mathbb{P}}_{Y \times Z}(\mathcal{F}) \ar[dd]_{p_{\mathcal{F}}} \\
& X \times Y \times Z \ar[dl]^{\pr_{12}} \ar[dr]_{\pr_{23}} & \\
X \times Y \ar[dr]_{\pr_{2}} & & Y \times Z \ar[dl]^{\pr_{1}} \\
& Y &.
}
$$
By the universal property of the fiber product, there is a map $t:{\mathbb{P}}_{X \times Y}(\mathcal{E}) \otimes_{Y} {\mathbb{P}}_{Y \times Z}(\mathcal{F}) \rightarrow X \times Y \times Z$, making
$$
\xymatrix
{
& {\mathbb{P}}_{X \times Y}(\mathcal{E}) \otimes_{Y} {\mathbb{P}}_{Y \times Z}(\mathcal{F}) \ar[dr]^{q_{\mathcal{F}}} \ar[dl]_{q_{\mathcal{E}}} \ar@{-->}[dd]^{t} & \\
{\mathbb{P}}_{X \times Y}(\mathcal{E}) \ar[dd]^{p_{\mathcal{E}}} & & {\mathbb{P}}_{Y \times Z}(\mathcal{F}) \ar[dd]_{p_{\mathcal{F}}} \\
& X \times Y \times Z \ar[dl]^{\pr_{12}} \ar[dr]_{\pr_{23}} & \\
X \times Y \ar[dr]_{\pr_{2}} & & Y \times Z \ar[dl]^{\pr_{1}} \\
& Y &.
}
$$
commute.  Thus, ${\mathbb{P}}_{X \times Y}(\mathcal{E}) \otimes_{Y} {\mathbb{P}}_{Y \times Z}(\mathcal{F})$ is an $X \times Z$-scheme.

{\it Remark.}  Using the tensor product notation to denote a fiber product of spaces becomes useful when $X=Y=Z$.
\begin{lemma}
Retain the notation above.  If $T$ is a scheme and ${\mathcal{G}}$ is an ${\mathcal{O}}_{Z \times T}$-module, then there is an isomorphism
$$
({\mathbb{P}}_{X \times Y}(\mathcal{E}) \otimes_{Y} {\mathbb{P}}_{Y \times Z}(\mathcal{F})) \otimes_{Z} {\mathbb{P}}_{Z \times T}(\mathcal{G}) \rightarrow {\mathbb{P}}_{X \times Y}(\mathcal{E}) \otimes_{Y} ({\mathbb{P}}_{Y \times Z}(\mathcal{F}) \otimes_{Z} {\mathbb{P}}_{Z \times T}(\mathcal{G})).
$$
\end{lemma}
\begin{proof}
The proof is a tedious exercise in the use of the universal property of the fibre product.
\end{proof}

{\it Remark}:  The same result holds if we replace the various projective bundles by arbitrary $X \times Y$, $Y\times Z$, and $Z \times W$ schemes.
By the lemma, the space $({\mathbb{P}}_{X^{2}}(\mathcal{F}))^{\otimes n}$ is well defined.

\begin{theorem} \label{theorem.rep1}
For $0< n$, there is a monomorphism
$$
\Phi:\Gamma_{n}(-) \Longrightarrow \operatorname{Hom}_{S}(-,({\mathbb{P}}_{X^{2}}(\mathcal{E}))^{\otimes n})
$$
which is an isomorphism for $n < m$.
\end{theorem}
A more detailed description of $\Phi$ is given in Proposition
\ref{prop.trans}.  For the remainder of this chapter, let $U$ be
an affine, noetherian scheme and suppose
$$
\pr_{i}: (U \times X)_{U}^{2} \rightarrow U \times X
$$
are the standard projection maps.  Let $q:U \rightarrow X$, and let $\mathcal{M}$ be an ${\mathcal{O}}_{U \times X}$-module isomorphic to $(\id_{U} \times q)_{*}{\mathcal{O}}_{U}$.  Let $Z = \operatorname{SSupp }\pr_{1}^{*}\mathcal{M}$, and let $i:Z \rightarrow (U \times X)_{U}^{2}$ be the inclusion map.

\begin{lemma}
Let $R$ be a commutative ring and let $A$, $B$, and $C$ be commutative $R$-algebras.  Let $q:B \rightarrow A$, and let $A$ have a $B$-module structure via the rule
$$
b \cdot a = q(b)a.
$$
In what follows, all unadorned tensor products are over $R$.  Then the $A \otimes C$-module map
$$
f:A \otimes C \rightarrow (A \otimes B) \otimes_{A} (A \otimes C) /\operatorname{ann }(((A \otimes B) \otimes_{A}(A \otimes C))\otimes_{A \otimes B}A)
$$
sending $a\otimes c$ to the image of $a \otimes 1 \otimes 1 \otimes c$ in
$$
(A \otimes B) \otimes_{A} (A \otimes C) /\operatorname{ann }(((A \otimes B) \otimes_{A}(A \otimes C))\otimes_{A \otimes B}A),
$$
is an isomorphism.
\end{lemma}

\begin{proof}
The map $f$ is an epimorphism since, if $a \in A$, $b \in B$, and $c \in C$, then $a\otimes b \otimes 1 \otimes c - q(b) a \otimes 1 \otimes 1 \otimes c \in\operatorname{ann }(((A \otimes B) \otimes_{A}(A \otimes C))\otimes_{A \otimes B}A)$.

We now show the map is a monomorphism.  First, if $A \otimes C$ is given an $(A \otimes B) \otimes_{A} (A \otimes C)$-module structure defined by the rule
$$
(a \otimes b \otimes 1 \otimes c) \cdot (c' \otimes a')=cc' \otimes aq(b)a'
$$
for $a, a' \in A$, $b \in B$ and $c \in C$, then it is easy to see that there is a natural $(A \otimes B) \otimes_{A} (A \otimes C)$-module isomorphism
\begin{equation} \label{eqn.nate}
((A \otimes B) \otimes_{A} (A \otimes C)) \otimes_{A \otimes B} A \rightarrow A \otimes C.
\end{equation}
Let $I$ be a finite set.  For all $i \in I$, let $a_{i} \in A$ and $c_{i} \in C$.  Suppose $\Sigma_{i}a_{i}\otimes 1 \otimes 1 \otimes c_{i} \in \operatorname{ann } ((A \otimes B)\otimes_{A} (A \otimes C) \otimes_{A \otimes B} A)$.  Then, under the isomorphism (\ref{eqn.nate}), $\Sigma_{i}a_{i}\otimes 1 \otimes 1 \otimes c_{i} \in \operatorname{ann } A\otimes C$.  Thus, $\Sigma_{i}a_{i} \otimes c_{i} \in A \otimes C$ equals $0$, proving that $f$ is a monomorphism.
\end{proof}

\begin{proposition} \label{prop.isomp}
The map $\pr_{2}i:Z \rightarrow U \times X$ is an isomorphism.
\end{proposition}

\begin{proof}
First, we show that as a map of topological spaces,
$$
\pr_{2}i:Z \rightarrow U \times X
$$
is an isomorphism.  This follows easily from the fact that, since
$$
Z = \{(a,b,a,c) \in (U \times X)_{U}^{2} | b=q(a) \},
$$
the continuous map
$$
(\id_{U} \times q) \times (\id_{U} \times \id_{X}):U \times X \rightarrow Z \subset (U \times X)_{U}^{2}
$$
is inverse to $\pr_{2}i$.

We next show that $\pr_{2}i$ is an isomorphism of schemes.  Let $\mathcal{I} = \mbox{ann }\pr_{1}^{*}\mathcal{M}$.  We must show that if
$$
\pr_{2}^{\#}:{\mathcal{O}}_{U \times X} \rightarrow {\pr_{2}}_{*}{\mathcal{O}}_{(U \times X)_{U}^{2}}
$$
is the sheaf morphism from the map $\pr_{2}$, and if
$$
\xymatrix{
i^{\#}:{\mathcal{O}}_{(U \times X)_{U}^{2}} \ar[r]^{\delta} & {\mathcal{O}}_{(U \times X)_{U}^{2}}/\mathcal{I} \ar[r] & i_{*}i^{-1}({\mathcal{O}}_{(U \times X)_{U}^{2}}/\mathcal{I})
}
$$
is the sheaf morphism from the map $i$, where the first arrow is the usual quotient map, then the sheaf morphism
$$
{\pr_{2}}_{*}i^{\#} \circ \pr_{2}^{\#}:{\mathcal{O}}_{U \times X} \rightarrow {\pr_{2}}_{*}i_{*}i^{-1}({\mathcal{O}}_{(U \times X)_{U}^{2}}/\mathcal{I})
$$
is an isomorphism of ${\mathcal{O}}_{U \times X}$-modules.  Since the rightmost composite of $i^{\#}$ is an isomorphism (which can be checked pointwise), in order to show ${\pr_{2}}_{*}i^{\#} \circ \pr_{2}^{\#}$ is an isomorphism, it suffices to show
$$
{\pr_{2}}_{*}\delta \circ \pr_{2}^{\#}: {\mathcal{O}}_{U \times X} \rightarrow \pr_{2*}({\mathcal{O}}_{(U\times X)_{U}^{2}}/\mathcal{I})
$$
is an isomorphism.  Let $\operatorname{Spec }B$ and $\operatorname{Spec }C$ be affine open subsets of $X$ and let $\operatorname{Spec }A =q^{-1}(\operatorname{Spec }B)$.  The composition of
$$
\pr_{2*} \delta \circ \pr_{2}^{\#}(\operatorname{Spec }A \otimes C): {\mathcal{O}}_{U \times X}(\operatorname{Spec }A \otimes C) \rightarrow \pr_{2*}({\mathcal{O}}_{(U\times X)_{U}^{2}}/\mathcal{I})(\operatorname{Spec }A \otimes C) =
$$
$$
({\mathcal{O}}_{(U\times X)_{U}^{2}}/\mathcal{I})((\operatorname{Spec }A \times X) \times_{\operatorname{Spec }A} (\operatorname{Spec }A \otimes C))
$$
with the $A \otimes C$-module isomorphism
$$
({\mathcal{O}}_{(U\times X)_{U}^{2}}/\mathcal{I})((\operatorname{Spec }A \times X) \times_{\operatorname{Spec }A} (\operatorname{Spec }A \otimes C)) \rightarrow
$$
$$
({\mathcal{O}}_{(U\times X)_{U}^{2}}/\mathcal{I})(\operatorname{Spec }((A \otimes B) \otimes_{A} (A \otimes C)))
$$
given by restriction (isomorphic by Lemma \ref{lem.sup1} and Lemma \ref{lem.sup2}) is just the map of $A \otimes C$-modules
$$
f:A \otimes C \rightarrow (A \otimes B) \otimes_{A} (A \otimes C) /\operatorname{ann }(((A \otimes B) \otimes_{A}(A \otimes C))\otimes_{A \otimes B}A)
$$
sending $a\otimes c$ to the image of $a \otimes 1 \otimes 1 \otimes c$ in
$$
(A \otimes B) \otimes_{A} (A \otimes C) /\operatorname{ann }(((A \otimes B) \otimes_{A}(A \otimes C))\otimes_{A \otimes B}A).
$$
The fact that this map is an isomorphism follows from the previous lemma.
\end{proof}

\begin{proposition} \label{prop.adjointness}
Let $\mathcal{N}$ be an ${\mathcal{O}}_{U\times X}$-module, and let $\mathcal{F}$ be an ${\mathcal{O}}_{(U \times X)_{U}^{2}}$-module.  Then $(\pr_{2}i)_{*}$ has a right adjoint $(\pr_{2}i)^{!}$ and there is a natural bijection

\begin{equation} \label{eqn.adjoints}
\psi: \operatorname{Hom}_{U\times X}(\pr_{2*}(\pr_{1}^{*}\mathcal{M} \otimes \mathcal{F}),\mathcal{N}) \rightarrow
\end{equation}
$$
\operatorname{Hom}_{(U\times X)_{U}^{2}}(\mathcal{F},\mathcal{H}\it{om}_{{\mathcal{O}}_{(U\times X)_{U}^{2}}}(\pr_{1}^{*}\mathcal{M},i_{*}(\pr_{2}i)^{!}\mathcal{N})).
$$
\end{proposition}

\begin{proof}
The first assertion follows from the previous lemma.  To prove the second assertion, we note that the functors
$$
(\pr_{2}i)_{*}\circ i^{*} \circ (\pr_{1}^{*}\mathcal{M} \otimes -):{\sf{Qcoh }} (U\times X)_{U}^{2} \rightarrow {\sf{Qcoh }} U \times X
$$
and
$$
\mathcal{H}\it{om}_{{\mathcal{O}}_{(U \times X)_{U}^{2}}}(\pr_{1}^{*}\mathcal{M},-) \circ i_{*} \circ (\pr_{2}i)^{!}:{\sf{Qcoh }}U\times X \rightarrow {\sf{Qcoh }}(U \times X)_{U}^{2}
$$
are the left adjoint and right adjoint, respectively, of an adjoint pair.  The assertion then follows from the fact that
$$
i_{*}i^{*}(\pr_{1}^{*}\mathcal{M} \otimes \mathcal{F}) \cong \pr_{1}^{*}\mathcal{M}\otimes \mathcal{F}.
$$
\end{proof}
We prove Theorem \ref{theorem.rep1} by first proving the following three preliminary results.

\begin{lemma} \label{lem.unit}
Suppose $r:U \rightarrow X$, $\mathcal{N} \cong (\id_{U} \times r)_{*}{\mathcal{O}}_{U}$,
$$
F:{\sf{Qcoh }} (U \times X)_{U}^{2} \rightarrow {\sf{Qcoh }} U \times X
$$
is the functor ${\pr_{2}}_{*}\circ (\pr_{1}^{*}{\mathcal{M}} \otimes -)$, and
$$
G:{\sf{Qcoh }} U \times X \rightarrow {\sf{Qcoh }}(U \times X)_{U}^{2}
$$
is the functor $\mathcal{H}\it{om}_{{\mathcal{O}}_{(U \times X)_{U}^{2}}}(\pr_{1}^{*}\mathcal{M},-) \circ i_{*} \circ ({\pr_{2}i})^{!}$.  Then $(F,G)$ form an adjoint pair such that if
$$
\eta:\id_{{\sf{Qcoh }} (U \times X)_{U}^{2}} \Longrightarrow GF
$$
is the unit map and
$$
\epsilon:FG \Longrightarrow \id_{{\sf{Qcoh }} U \times X}
$$
is the counit map, then $\eta_{G\mathcal{N}}$, and ${\epsilon}_{\mathcal{N}}$ are isomorphisms.
\end{lemma}

\begin{proposition} \label{prop.simpleiso}
Let $r:U \rightarrow X$, let $\mathcal{N} \cong (\id_{U} \times r)_{*}{\mathcal{O}}_{U}$.  Then there is an isomorphism of ${\mathcal{O}}_{(U \times X)_{U}^{2}}$-modules
$$
\mathcal{H}\it{om}_{{\mathcal{O}}_{(U \times X)_{U}^{2}}}(\pr_{1}^{*}\mathcal{M},i_{*}(\pr_{2}i)^{!}\mathcal{N}) \rightarrow ((\id_{U} \times q)\times (\id_{U} \times r))_{*}{\mathcal{O}}_{U}.
$$
\end{proposition}

\begin{proposition} \label{prop.bij}
Let $\mathcal{F}$ be an ${\mathcal{O}}_{(U \times X)_{U}^{2}}$-module and let $\mathcal{N}$ be an ${\mathcal{O}}_{U \times X}$-module.  Then the bijection $(\ref{eqn.adjoints})$,
$$
\psi:\operatorname{Hom}_{U\times X}({\pr_{2}}_{*}(\pr_{1}^{*}\mathcal{M} \otimes \mathcal{F}),\mathcal{N}) \rightarrow
$$
$$
\operatorname{Hom}_{(U\times X)_{U}^{2}}(\mathcal{F},\mathcal{H}\it{om}_{{\mathcal{O}}_{(U\times X)_{U}^{2}}}(\pr_{1}^{*}\mathcal{M},i_{*}(\pr_{2}i)^{!}\mathcal{N})).
$$
is a bijection on epimorphisms.
\end{proposition}

\subsection{The proof of Lemma \ref{lem.unit}}
The next result follows easily from the definitions, so we omit the proof.
\begin{lemma} \label{lem.imm}
Suppose $i:W \rightarrow Y$ is a closed immersion of schemes and $\mathcal{A}$ and $\mathcal{B}$ are ${\mathcal{O}}_{W}$-modules.  Then there is a natural isomorphism of ${\mathcal{O}}_{Y}$-modules
$$
i_{*}\mathcal{H}\it{om}_{W}(\mathcal{A},\mathcal{B}) \rightarrow \mathcal{H}\it{om}_{Y}(i_{*}\mathcal{A},i_{*}\mathcal{B}).
$$
\end{lemma}

\begin{lemma} \label{lem.hom}
Suppose $\mathcal{P}$ is an ${\mathcal{O}}_{(U \times X)_{U}^{2}}$-module such that the natural morphism
$$
\mathcal{P} \rightarrow i_{*}i^{*}\mathcal{P}
$$
is an isomorphism.  Then there is an ${\mathcal{O}}_{Z}$-module isomorphism
$$
i^{*}\mathcal{H}\it{om}_{{\mathcal{O}}_{(U \times X)_{U}^{2}}}(\pr_{1}^{*}\mathcal{M},\mathcal{P}) \rightarrow \mathcal{H}\it{om}_{{\mathcal{O}}_{Z}} (i^{*}\pr_{1}^{*}\mathcal{M},i^{*}\mathcal{P}).
$$
\end{lemma}

\begin{proof}
By Lemma \ref{lem.imm},
\begin{equation} \label{eqn.hom}
i_{*}\mathcal{H}\it{om}_{{\mathcal{O}}_{Z}} (i^{*}\pr_{1}^{*}\mathcal{M},i^{*}\mathcal{P}) \cong \mathcal{H}\it{om}_{{\mathcal{O}}_{(U \times X)_{U}^{2}}} (i_{*}i^{*}\pr_{1}^{*}\mathcal{M},i_{*}i^{*}\mathcal{P}),
\end{equation}
while, since
$$
i_{*}i^{*}\pr_{1}^{*}\mathcal{M} \cong \pr_{1}^{*}\mathcal{M},
$$
and
$$
i_{*}i^{*}\mathcal{P} \cong \mathcal{P},
$$
the latter sheaf is isomorphic to
$$
\mathcal{H}\it{om}_{{\mathcal{O}}_{(U \times X)_{U}^{2}}} (\pr_{1}^{*}\mathcal{M},\mathcal{P}).
$$
Thus, applying $i^{*}$ to (\ref{eqn.hom}), we have a sheaf isomorphism
$$
i^{*}i_{*}\mathcal{H}\it{om}_{{\mathcal{O}}_{Z}} (i^{*}\pr_{1}^{*}\mathcal{M},i^{*}\mathcal{P}) \rightarrow i^{*}\mathcal{H}\it{om}_{{\mathcal{O}}_{(U \times X)_{U}^{2}}} (\pr_{1}^{*}\mathcal{M},\mathcal{P}).
$$
The assertion follows from the fact that
$$
i^{*}i_{*}\mathcal{H}\it{om}_{{\mathcal{O}}_{Z}} (i^{*}\pr_{1}^{*}\mathcal{M},i^{*}\mathcal{P}) \cong \mathcal{H}\it{om}_{{\mathcal{O}}_{Z}} (i^{*}\pr_{1}^{*}\mathcal{M},i^{*}\mathcal{P}).
$$
\end{proof}

\begin{lemma} \label{lem.free}
The ${\mathcal{O}}_{Z}$-module $i^{*}\pr_{1}^{*}\mathcal{M}$ is free, of rank one.
\end{lemma}

\begin{proof}
We need to show that
$$
i^{*}\pr_{1}^{*}\mathcal{M} \cong i^{-1}({\mathcal{O}}_{(U \times X)_{U}^{2}}/\operatorname{ann }\pr_{1}^{*}\mathcal{M})
$$
as $i^{-1}({\mathcal{O}}_{(U \times X)_{U}^{2}}/\operatorname{ann }\pr_{1}^{*}\mathcal{M})$-modules.  Since $i^{*}i_{*} \cong \id_{Z}$ we need only show that
$$
i_{*}i^{*}\pr_{1}^{*}\mathcal{M} \cong i_{*}i^{-1}({\mathcal{O}}_{(U \times X)_{U}^{2}}/\operatorname{ann }\pr_{1}^{*}\mathcal{M})
$$
as $i_{*}i^{-1}({\mathcal{O}}_{(U \times X)_{U}^{2}}/\operatorname{ann }\pr_{1}^{*}\mathcal{M})={\mathcal{O}}_{(U \times X)_{U}^{2}}/\operatorname{ann }\pr_{1}^{*}\mathcal{M}$-modules.  But $i_{*}i^{*}\pr_{1}^{*}\mathcal{M} \cong \pr_{1}^{*}\mathcal{M}$.  Thus, we need only show that there is an ${\mathcal{O}}_{(U \times X)_{U}^{2}}/\operatorname{ann }\pr_{1}^{*}\mathcal{M}$-module isomorphism
$$
\pr_{1}^{*}\mathcal{M} \cong {\mathcal{O}}_{(U \times X)_{U}^{2}}/\operatorname{ann }\pr_{1}^{*}\mathcal{M}.
$$
We construct a map locally and show that it glues together.  We first observe that the set
$$
\mathcal{C} = \{ \operatorname{Spec }(A \otimes B) \otimes_{A} (A \otimes C) | q^{-1}(\operatorname{Spec }B)=\operatorname{Spec }A, \mbox{ and } \operatorname{Spec }B, \operatorname{Spec }C \subset X \}
$$
is an open cover of $\operatorname{Supp }\pr_{1}^{*}\mathcal{M}$ by Lemma \ref{lem.sup2}.  Thus, to define a map of sheaves $\pr_{1}^{*}\mathcal{M} \rightarrow {\mathcal{O}}_{(U \times X)_{U}^{2}}/\operatorname{ann }\pr_{1}^{*}\mathcal{M}$, it suffices to define, for each $V \in \mathcal{C}$, a map $f_{V}:\pr_{1}^{*}\mathcal{M}(V) \rightarrow {\mathcal{O}}_{(U \times X)_{U}^{2}}/\operatorname{ann }\pr_{1}^{*}\mathcal{M}(V)$ and then to show that if $V'$ is also in $\mathcal{C}$, then $f_{V}(V\cap V')=f_{V'}(V \cap V')$.

Suppose $V = \operatorname{Spec }(A \otimes B) \otimes_{A} (A \otimes C) \in \mathcal{C}$.  Since
$$
{\mathcal{O}}_{(U \times X)_{U}^{2}}(\operatorname{Spec }(A \times B)\otimes_{A}(A \otimes C))=(A \otimes B) \otimes_{A} (A \otimes C)
$$
and
$$
\pr_{1}^{*}\mathcal{M}(\operatorname{Spec }(A\otimes B) \otimes_{A} (A \otimes C))=(A \otimes B) \otimes_{A} (A \otimes C) \otimes_{A \otimes B}A,
$$
there is an epimorphism
\begin{equation} \label{eqn.epimorph}
{\mathcal{O}}_{(U \times X)_{U}^{2}}(\operatorname{Spec }(A\otimes B) \otimes_{A} (A\otimes C)) \rightarrow \pr_{1}^{*}\mathcal{M}(\operatorname{Spec }(A\otimes B) \otimes_{A} (A \otimes C))
\end{equation}
sending
$$
a \otimes b \otimes  1 \otimes c \in (A \otimes B) \otimes_{A} (A \otimes C)
$$
to
$$
(a \otimes b \otimes 1 \otimes c) \otimes 1 \in (A \otimes B \otimes C)\otimes_{A \otimes B}A.
$$
The kernel of this epimorphism is exactly ann $(A \otimes B) \otimes_{A} (A \otimes C)\otimes_{A \otimes B}A$.  Thus, the map (\ref{eqn.epimorph}) gives an isomorphism
$$
f_{V}:({\mathcal{O}}_{(U \times X)_{U}^{2}}/\operatorname{ann }\pr_{1}^{*}\mathcal{M})(\operatorname{Spec }(A \otimes B) \times_{A} (A \otimes C)) \rightarrow
$$
$$
\pr_{1}^{*}\mathcal{M}(\operatorname{Spec }(A\otimes B) \otimes_{A} (A \otimes C)).
$$
We show that these isomorphisms glue to give an isomorphism of sheaves.  Suppose $V' = \operatorname{Spec }(A' \otimes B') \otimes_{A'} (A' \otimes C') \in \mathcal{C}$. Then as before, we have an isomorphism
$$
f_{V'}:({\mathcal{O}}_{(U \times X)_{U}^{2}}/\operatorname{ann }\pr_{1}^{*}\mathcal{M})(\operatorname{Spec }(A'\otimes B') \otimes_{A'}(A' \otimes C')) \rightarrow
$$
$$
\pr_{1}^{*}\mathcal{M}(\operatorname{Spec }(A'\otimes B') \otimes_{A'} (A' \otimes C')).
$$
We must show
$$
f_{V}(\operatorname{Spec }(A\otimes B) \otimes_{A} (A\otimes C) \cap \operatorname{Spec }(A' \otimes B') \otimes_{A'} (A'\otimes C') )=
$$
$$
{f'}_{V'}(\operatorname{Spec }(A\otimes B) \otimes_{A} (A \otimes C) \cap \operatorname{Spec }(A' \otimes B') \otimes_{A'} (A' \otimes C')).
$$
Since $U$ and $X$ are separated, the intersection of two affine open subsets of $U$ is affine, and the intersection of two affine open subsets of $X$ is affine.  Suppose Spec $A \cap $ Spec $A'$ = Spec $A''$, Spec $B \cap $ Spec $B'$ = Spec $B''$, and Spec $C \cap $ Spec $C'$ = Spec $C''$.  Then we need to show that
$$
f_{V}(\operatorname{Spec }(A'' \otimes B'') \otimes_{A''}(A'' \otimes C''))=
$$
$$
f_{V'}(\operatorname{Spec }(A'' \otimes B'') \otimes_{A''}(A'' \otimes C'')).
$$
This follows immediately from the definition of $f_{V}$ and $f_{V'}$.
\end{proof}

\begin{lemma} \label{lem.i}
Let $F:\sf{A} \rightarrow \sf{B}$ be a functor and suppose $(F,G)$ is an adjoint pair.  If $\mathcal{N}$ is an object of $\sf{B}$ such that $\epsilon_{\mathcal{N}}:FG\mathcal{N} \rightarrow \mathcal{N}$ is an isomorphism, then $\eta_{G\mathcal{N}}$ is an isomorphism.
\end{lemma}

\begin{proof}
Suppose $\epsilon_{\mathcal{N}}$ is an isomorphism.  Since the composition
$$
\xymatrix{
G\mathcal{N} \ar@{=>}[r]^{\eta_{G\mathcal{N}}} & GFG\mathcal{N} \ar@{=>}[r]^{G*\epsilon_{\mathcal{N}}} & G\mathcal{N}
}
$$
is the identity, and since $G*\epsilon_{\mathcal{N}}$ is an isomorphism,
$$
\eta_{G\mathcal{N}}=(G*\epsilon_{\mathcal{N}})^{-1}
$$
so that $\eta_{G\mathcal{N}}$ is an isomorphism.
\end{proof}

\begin{lemma}
Let $r:U \rightarrow X$, let $\mathcal{N} \cong (\id_{U} \times r)_{*}{\mathcal{O}}_{U}$, let
$$
F:{\sf{Qcoh }} (U \times X)_{U}^{2} \rightarrow {\sf{Qcoh }} U \times X
$$
be the functor ${\pr_{2}}_{*}\circ (\pr_{1}^{*}{\mathcal{M}} \otimes -)$, and let
$$
G:{\sf{Qcoh }} U \times X \rightarrow {\sf{Qcoh }}(U \times X)_{U}^{2}
$$
be the functor $\mathcal{H}\it{om}_{{\mathcal{O}}_{(U \times X)_{U}^{2}}}(\pr_{1}^{*}\mathcal{M},-) \circ i_{*} \circ ({\pr_{2}i})^{!}$.  Then $(F,G)$ form an adjoint pair such that if
$$
\eta:\id_{{\sf{Qcoh }} (U \times X)_{U}^{2}} \Longrightarrow GF
$$
is the unit map and
$$
\epsilon:FG \Longrightarrow \id_{{\sf{Qcoh }} U \times X}
$$
is the counit map, then $\eta_{G\mathcal{N}}$ and ${\epsilon}_{\mathcal{N}}$ are isomorphisms.
\end{lemma}

\begin{proof}
By Lemma \ref{lem.i}, in order to prove that $\eta_{G\mathcal{N}}$ is an isomorphism, it suffices to show that $\epsilon_{\mathcal{N}}$ is an isomorphism.  To complete the proof of the lemma, we show that $\epsilon_{\mathcal{N}}$ is an isomorphism.  Let
$$
F':{\sf{Qcoh }} (U \times X)_{U}^{2} \rightarrow {\sf{Qcoh }}(U \times X)_{U}^{2}
$$
be the functor $\pr_{1}^{*}{\mathcal{M}} \otimes -$, and let
$$
G':{\sf{Qcoh }} (U \times X)_{U}^{2} \rightarrow {\sf{Qcoh }}(U \times X)_{U}^{2}
$$
be the functor $\mathcal{H}\it{om}_{{\mathcal{O}}_{(U \times X)_{U}^{2}}}(\pr_{1}^{*}\mathcal{M},-)$.  If $\epsilon':F'G' \rightarrow \id$ is the counit of the adjoint pair $(F',G')$, then $\epsilon:FG \rightarrow \id$ is given by the composition
$$
\xymatrix{
\pr_{2*}F'G'i_{*}(\pr_{2}i)^{!} \ar@{=>}[rrr]^{\pr_{2*}*\epsilon'*i_{*}(\pr_{2}i)^{!}} & & & (\pr_{2}i)_{*}(\pr_{2}i)^{!} \ar@{=>}[r] & \id
}
$$
where the rightmost map is the counit map of the adjoint pair $((\pr_{2}i)_{*},(\pr_{2}i)^{!})$.  Since $\pr_{2}i$ is an isomorphism (Proposition \ref{prop.isomp}), to show that $\epsilon_{\mathcal{N}}$ is an isomorphism, it suffices to show that ${\epsilon'}_{i_{*}(\pr_{2}i)^{!}\mathcal{N}}$ is an isomorphism.  This map is defined on an affine open set $V \subset (U \times X)_{U}^{2}$ by sending $a \otimes \beta$ in
$$
F'G'i_{*}(\pr_{2}i)^{!}\mathcal{N}(V)=\pr_{1}^{*}\mathcal{M}(V) \otimes \operatorname{Hom}_{{\mathcal{O}}_{(U \times X)_{U}^{2}}(V)}(\pr_{1}^{*}\mathcal{M}(V),i_{*}(\pr_{2}i)^{!}\mathcal{N}(V))
$$
to $\beta(a) \in i_{*}({\pr_{2}i})^{!}\mathcal{N}$.  To prove that ${\epsilon'}_{i_{*}(\pr_{2}i)^{!}\mathcal{N}}$ is an isomorphism, we claim it suffices to prove that $i^{*}{\epsilon'}_{i_{*}(\pr_{2}i)^{!}\mathcal{N}}$ is an isomorphism.  For, suppose this is the case.  Then, if ${\eta}'':\id \Longrightarrow i_{*}i^{*}$ is the unit, the top map in the commutative diagram
\begin{equation} \label{eqn.adjj}
\xymatrix{
i_{*}i^{*}F'G'i_{*}(\pr_{2}i)^{!}\mathcal{N} \ar[rr]^{i_{*}i^{*}\epsilon'} & & i_{*}i^{*}i_{*}(\pr_{2}i)^{!}\mathcal{N} \\
F'G'i_{*}(\pr_{2}i)^{!}\mathcal{N} \ar[rr]_{\epsilon'} \ar[u]^{\eta''} & & i_{*}(\pr_{2}i)^{!}\mathcal{N} \ar[u]_{\eta''}
}
\end{equation}
is an isomorphism.  Since ${\eta''}_{\pr_{1}^{*}\mathcal{M}}$ is an isomorphism, ${\eta''}_{\pr_{1}^{*}\mathcal{M} \otimes \mathcal{H}}$ is an isomorphism for all $\mathcal{H}$ in ${\sf{Qcoh }}(U\times X)_{U}^{2}$.  Thus, ${\eta''}_{F'G'i_{*}(\pr_{2}i)^{!}\mathcal{N}}$ is an isomorphism.  Since ${\eta''}_{i_{*}(\pr_{2}i)^{!}\mathcal{N}}$ is also an isomorphism by Lemma \ref{lem.i}, the top and sides of (\ref{eqn.adjj}) are isomorphisms, so the bottom of (\ref{eqn.adjj}) is also an isomorphism.

To complete the proof, we show that
$$
i^{*}{\epsilon'}_{i_{*}(\pr_{2}i)^{!}\mathcal{N}}:i^{*}(\pr_{1}^{*}\mathcal{M} \otimes \mathcal{H}\it{om}_{{\mathcal{O}}_{(U \times X)_{U}^{2}}}(\pr_{1}^{*}\mathcal{M},i_{*}(\pr_{2}i)^{!}\mathcal{N})) \rightarrow i^{*}i_{*}(\pr_{2}i)^{!}\mathcal{N}
$$
is an isomorphism.  By Lemma \ref{lem.hom}, there is an isomorphism
$$
\psi:i^{*}\mathcal{H}\it{om}_{{\mathcal{O}}_{(U \times X)_{U}^{2}}}(\pr_{1}^{*}\mathcal{M},i_{*}(\pr_{2}i)^{!}\mathcal{N})) \rightarrow \mathcal{H}\it{om}_{{\mathcal{O}}_{Z}}(i^{*}\pr_{1}^{*}\mathcal{M},i^{*}i_{*}(\pr_{2}i)^{!}\mathcal{N})).
$$
Thus, there is a morphism
$$
\xymatrix{
i^{*}\pr_{1}^{*}\mathcal{M} \otimes \mathcal{H}\it{om}_{{\mathcal{O}}_{Z}}(i^{*}\pr_{1}^{*}\mathcal{M},i^{*}i_{*}(\pr_{2}i)^{!}\mathcal{N}) \ar[d]^{i^{*}\pr_{1}^{*}\mathcal{M} \otimes \psi^{-1}} & \\
i^{*}\pr_{1}^{*}\mathcal{M} \otimes i^{*}\mathcal{H}\it{om}_{{\mathcal{O}}_{(U \times X)_{U}^{2}}}(\pr_{1}^{*}\mathcal{M},i_{*}(\pr_{2}i)^{!}\mathcal{N}) \ar[d] &  \\
i^{*}(\pr_{1}^{*}\mathcal{M} \otimes \mathcal{H}\it{om}_{{\mathcal{O}}_{(U \times X)_{U}^{2}}}(\pr_{1}^{*}\mathcal{M},i_{*}(\pr_{2}i)^{!}\mathcal{N})) \ar[r]^{i^{*}{\epsilon'}} & i^{*}i_{*}(\pr_{2}i)^{!}\mathcal{N}.
}
$$
Since $i^{*}\pr_{1}^{*}\mathcal{M}$ is a free rank one ${\mathcal{O}}_{Z}$-module by Lemma \ref{lem.free}, there is a natural isomorphism
$$
\xymatrix{
i^{*}\pr_{1}^{*}\mathcal{M} \otimes \mathcal{H}\it{om}_{{\mathcal{O}}_{Z}}(i^{*}\pr_{1}^{*}\mathcal{M},i^{*}i_{*}({\pr_{2}i})^{!}\mathcal{N}) \ar[d] & & \\
\mathcal{H}\it{om}_{{\mathcal{O}}_{Z}}(i^{*}\pr_{1}^{*}\mathcal{M},i^{*}i_{*}({\pr_{2}i})^{!}\mathcal{N}) \ar[r] & i^{*}i_{*}(\pr_{2}i)^{!}\mathcal{N}.
}
$$
An easy local computation shows that the diagram
$$
\xymatrix{
i^{*}\pr_{1}^{*}\mathcal{M} \otimes \mathcal{H}\it{om}_{{\mathcal{O}}_{Z}}(i^{*}\pr_{1}^{*}\mathcal{M},i^{*}i_{*}({\pr_{2}i})^{!}\mathcal{N}) \ar[dr]^{(i^{*}\epsilon') \circ \psi^{-1}}  \ar[d] &  & \\
\mathcal{H}\it{om}_{{\mathcal{O}}_{Z}}(i^{*}\pr_{1}^{*}\mathcal{M},i^{*}i_{*}({\pr_{2}i})^{!}\mathcal{N}) \ar[r] & i^{*}i_{*}(\pr_{2}i)^{!}\mathcal{N}.
}
$$
commutes.  Thus, $i^{*}\epsilon'$, and hence $\epsilon'$, is an isomorphism.
\end{proof}

\subsection{The proof of Proposition \ref{prop.simpleiso}}

\begin{lemma} \label{cor.homi}
If $r:U \rightarrow X$, and if $\mathcal{N} \cong (\id_{U} \times r)_{*}{\mathcal{O}}_{U}$ then there is an isomorphism of ${\mathcal{O}}_{Z}$-modules
$$
i^{*}\mathcal{H}\it{om}_{{\mathcal{O}}_{(U \times X)_{U}^{2}}}(\pr_{1}^{*}\mathcal{M},i_{*}(\pr_{2}i)^{!}\mathcal{N}) \rightarrow (\pr_{2}i)^{!}\mathcal{N}.
$$
\end{lemma}

\begin{proof}
By Lemma $\ref{lem.i}$, the natural morphism
$$
i_{*}(\pr_{2}i)^{!}\mathcal{N} \rightarrow i_{*}i^{*}i_{*}(\pr_{2}i)^{!}\mathcal{N}
$$
is an isomorphism.  Thus, by Lemma $\ref{lem.hom}$, there is an isomorphism
$$
i^{*}\mathcal{H}\it{om}_{{\mathcal{O}}_{(U \times X)_{U}^{2}}}(\pr_{1}^{*}\mathcal{M},i_{*}(\pr_{2}i)^{!}\mathcal{N}) \rightarrow \mathcal{H}\it{om}_{{\mathcal{O}}_{Z}} (i^{*}\pr_{1}^{*}\mathcal{M},i^{*}i_{*}(\pr_{2}i)^{!}\mathcal{N}).
$$
By Lemma $\ref{lem.free}$,
$$
\mathcal{H}\it{om}_{{\mathcal{O}}_{Z}} (i^{*}\pr_{1}^{*}\mathcal{M},i^{*}i_{*}(\pr_{2}i)^{!}\mathcal{N}) \cong i^{*}i_{*}(\pr_{2}i)^{!}\mathcal{N},
$$
and this last sheaf is naturally isomorphic to $(\pr_{2}i)^{!}\mathcal{N}$.
\end{proof}

\begin{lemma} \label{lem.*!}
Let $f:Y \rightarrow Z$ be an isomorphism of schemes.  Then the functors
$$
f^{*}: \sf{Qcoh }Z \rightarrow \sf{Qcoh }Y
$$
and
$$
f^{!}: \sf{Qcoh }Y \rightarrow \sf{Qcoh }Z
$$
are naturally equivalent.
\end{lemma}

\begin{proof}
Since $f$ is an isomorphism, $f$ is affine and there are adjoint pairs $(f^{*},f_{*})$ and $(f_{*},f^{!})$ with counits $\epsilon$ and $\epsilon'$, respectively, which are isomorphisms.  Thus, we have a natural isomorphism
$$
\xymatrix{
f^{!} \ar@{=>}[r]^{{\epsilon}^{-1}*f^{!}} & f^{*}f_{*}f^{!} \ar@{=>}[r]^{f^{*}*{\epsilon'}} & f^{*}.
}
$$
\end{proof}

\begin{proposition}
Let $r:U \rightarrow X$, let $\mathcal{N} \cong (\id_{U} \times r)_{*}{\mathcal{O}}_{U}$.  Then there is an isomorphism of ${\mathcal{O}}_{(U \times X)_{U}^{2}}$-modules
$$
\mathcal{H}\it{om}_{{\mathcal{O}}_{(U \times X)_{U}^{2}}}(\pr_{1}^{*}\mathcal{M},i_{*}(\pr_{2}i)^{!}\mathcal{N}) \rightarrow ((\id_{U} \times q)\times (\id_{U} \times r))_{*}{\mathcal{O}}_{U}.
$$
\end{proposition}

\begin{proof}
By Lemma \ref{cor.homi}, there is an isomorphism
$$
i^{*}\mathcal{H}\it{om}_{{\mathcal{O}}_{(U \times X)_{U}^{2}}}(\pr_{1}^{*}\mathcal{M},i_{*}(\pr_{2}i)^{!}\mathcal{N}) \rightarrow (\pr_{2}i)^{!}\mathcal{N}.
$$
By Lemma \ref{lem.unit},
$$
i_{*}i^{*}\mathcal{H}\it{om}_{{\mathcal{O}}_{(U \times X)_{U}^{2}}}(\pr_{1}^{*}\mathcal{M},i_{*}(\pr_{2}i)^{!}\mathcal{N}) \cong \mathcal{H}\it{om}_{{\mathcal{O}}_{(U \times X)_{U}^{2}}}(\pr_{1}^{*}\mathcal{M},i_{*}(\pr_{2}i)^{!}\mathcal{N}),
$$
so that
$$
\mathcal{H}\it{om}_{{\mathcal{O}}_{(U \times X)_{U}^{2}}}(\pr_{1}^{*}\mathcal{M},i_{*}(\pr_{2}i)^{!}\mathcal{N}) \cong i_{*}(\pr_{2}i)^{!}\mathcal{N}.
$$
By Lemma $\ref{lem.*!}$,
\begin{equation} \label{eqn.fruit}
i_{*}(\pr_{2}i)^{!}\mathcal{N} \cong i_{*}i^{*}{\pr_{2}}^{*}\mathcal{N} \cong i_{*}i^{*}{\pr_{2}}^{*}(\id_{U} \times r)_{*}{\mathcal{O}}_{U}.
\end{equation}
We simplify (\ref{eqn.fruit}) by defining a map $a$ which makes the diagram
\begin{equation} \label{eqn.alpha}
\xymatrix{
& & & \operatorname{SSupp }(\pr_{1}^{*}(\id_{U} \times q)_{*}{\mathcal{O}}_{U}) \ar[d]^{i} \\
U \ar[rrr]_{(\id_{U} \times q) \times (\id_{U} \times r)} \ar@{-->}[urrr]^{a}  & & & (U \times X)_{U}^{2}.
}
\end{equation}
commute.  For $p \in U$, define $a(p)=(p, q(p), p, r(p)) \in Z$.  Let $\mathcal{I} = \mbox{ann }\pr_{1}^{*}\mathcal{M}$.  We now define the sheaf component of $a$,
$$
a^{\#}:i^{-1}({\mathcal{O}}_{(U \times X)_{U}^{2}}/\mathcal{I}) \rightarrow a_{*}{\mathcal{O}}_{U}.
$$
To define $a^{\#}$, we note that the sheaf map
$$
(\id_{U} \times q) \times (\id_{U} \times r)^{\#}:{\mathcal{O}}_{(U \times X)_{U}^{2}} \rightarrow ((\id_{U}\times q) \times (\id_{U} \times r))_{*}{\mathcal{O}}_{U}
$$
has $\mathcal{I}$ in its kernel.  This follows readily by examining the morphism locally.  Thus, there is a sheaf map
$$
({\mathcal{O}}_{(U \times X)_{U}^{2}}/\mathcal{I}) \rightarrow ((\id_{U}\times q) \times (\id_{U} \times r))_{*} {\mathcal{O}}_{U}.
$$
Applying $i^{-1}$ to this map, we have a map
$$
a^{\#}:i^{-1}({\mathcal{O}}_{(U \times X)_{U}^{2}}/\mathcal{I}) \rightarrow i^{-1}((\id_{U}\times q) \times (\id_{U} \times r))_{*} {\mathcal{O}}_{U}.
$$
The assertion follows by noting that
$a_{*}=i^{-1}((\id_{U} \times q) \times (\id_{U} \times r))_{*}.$  With this definition of $a$, it is easy to see that (\ref{eqn.alpha}) commutes.

We next show that the diagram of schemes
$$
\xymatrix{
U \ar[rrr]^{{\id}_{U}} \ar[d]_{a} & & & U \ar[d]^{\id_{U} \times r}\\
\mbox{SSupp }(\pr_{1}^{*}(\id_{U} \times q)_{*}{\mathcal{O}}_{U}) \ar[rrr]_{\pr_{2}i} & & & U \times X
}
$$
is a pullback diagram.  In particular, we show that
$$
a \times \id_{U}:U \rightarrow Z \times_{U \times X} U
$$
is an isomorphism.  Let $\pi_{2}:Z \times_{U \times X} U \rightarrow U$ be projection.  It is clear that the composition
$$
\xymatrix{
U \ar[rr]^{a \times \id_{U}} & & Z \times_{U \times X} U \ar[rr]^{\pi_{2}} & & U
}
$$
is the identity map.  Thus, in order to show $a \times \id_{U}$ is an isomorphism, it suffices to show that the composition
$$
\xymatrix{
Z \times_{U \times X} U \ar[rr]^{\pi_{2}} & & U \ar[rr]^{a \times \id_{U}} & & Z\times_{U \times X} U
}
$$
is the identity.  This is true at the level of spaces.  We must show that the map of sheaves
$$
{\mathcal{O}}_{Z \times_{U \times X}U} \rightarrow (a \times \id_{U})_{*}{\mathcal{O}}_{U} \rightarrow (a \times \id_{U})_{*}{\pi_{2}}_{*}{\mathcal{O}}_{Z \times_{U \times X}U}
$$
is the identity map, where the left map comes from the scheme map $a \times \id_{U}$ and the right map comes from the scheme map $\pi_{2}$.  A local computation bears out this fact.  Thus, by Lemma $\ref{lem.anothersquare}$, $i^{*}{\pr_{2}}^{*}(\id_{U} \times r)_{*} \cong a_{*}$.  Applying $i_{*}$ to this isomorphism yields the isomorphism of functors
$$
i_{*}i^{*}{\pr_{2}}^{*}(\id_{U} \times r)_{*} \cong i_{*}a_{*}=((\id_{U} \times q) \times (\id_{U} \times r))_{*}.
$$
\end{proof}

\subsection{The proof of Proposition \ref{prop.bij}}

\begin{proposition}
Let $\mathcal{F}$ be an ${\mathcal{O}}_{(U \times X)_{U}^{2}}$-module and let $\mathcal{N}$ be an ${\mathcal{O}}_{U \times X}$-module.  Then the bijection $(\ref{eqn.adjoints})$,
$$
\psi:\operatorname{Hom}_{U\times X}({\pr_{2}}_{*}(\pr_{1}^{*}\mathcal{M} \otimes \mathcal{F}),\mathcal{N}) \rightarrow
$$
$$
\operatorname{Hom}_{(U\times X)_{U}^{2}}(\mathcal{F},\mathcal{H}\it{om}_{{\mathcal{O}}_{(U\times X)_{U}^{2}}}(\pr_{1}^{*}\mathcal{M},i_{*}(\pr_{2}i)^{!}\mathcal{N})).
$$
is a bijection on epimorphisms.
\end{proposition}

\begin{proof}
Let
$$
F: {\sf{Qcoh }} (U \times X)_{U}^{2} \rightarrow {\sf{Qcoh }} U \times X
$$
be the functor
$$
(\pr_{2}i)_{*}\circ i^{*} \circ (\pr_{1}^{*}\mathcal{M} \otimes -),
$$
and let
$$
G: {\sf{Qcoh }} U \times X \rightarrow {\sf{Qcoh }} (U \times X)_{U}^{2}
$$
be the functor
$$
\mathcal{H}\it{om}_{{\mathcal{O}}_{(U \times X)_{U}^{2}}}(\pr_{1}^{*}\mathcal{M},-) \circ i_{*} \circ (\pr_{2}i)^{!}.
$$
Let $\eta:\id \rightarrow GF$ be the unit of the adjoint pair $(F,G)$.  By definition, if $f$ is an element of the left hand side of (\ref{eqn.adjoints}), then $\psi(f)= Gf \eta_{\mathcal{F}}$.  Thus, if $f$ is an epimorphism, then to show $\psi(f)$ is an epimorphism, it suffices to show that $\eta_{\mathcal{F}}$ and $Gf$ are epimorphisms.

We first show that the unit map
$$
\eta_{\mathcal{F}}:\mathcal{F} \rightarrow GF\mathcal{F}
$$
is an epimorphism of ${\mathcal{O}}_{(U \times X)_{U}^{2}}$-modules.  Let
$$
F': {\sf{Qcoh }} (U \times X)_{U}^{2} \rightarrow {\sf{Qcoh }} (U \times X)_{U}^{2}
$$
be the functor
$$
\pr_{1}^{*}\mathcal{M} \otimes - ,
$$
and let
$$
G': {\sf{Qcoh }} (U \times X)_{U}^{2} \rightarrow {\sf{Qcoh }} (U \times X)_{U}^{2}
$$
be the functor
$$
\mathcal{H}\it{om}_{{\mathcal{O}}_{(U \times X)_{U}^{2}}}(\pr_{1}^{*}\mathcal{M},-).
$$
Then $(F',G')$ form an adjoint pair, with unit, $\eta'$.  We show that $\eta'_{\mathcal{F}}$ is an epimorphism, which will later help us show $\eta_{\mathcal{F}}$ is also an epimorphism.  We first note that, to prove
$$
{\eta'}_{\mathcal{F}}:\mathcal{F} \rightarrow \mathcal{H}\it{om}_{{\mathcal{O}}_{(U \times X)_{U}^{2}}}(\pr_{1}^{*}\mathcal{M},\pr_{1}^{*}\mathcal{M}\otimes \mathcal{F})
$$
is an epimorphism, it suffices to show $i^{*}\eta'$ is an epi.  For, supposing that $i^{*}\eta'$ is an epi, $i_{*}i^{*} \eta'$ must also be an epi since $i_{*}$ is exact.  Since $i$ is a closed immersion, the unit map ${\eta''}_{\mathcal{F}}:\mathcal{F} \rightarrow i_{*}i^{*}\mathcal{F}$ is an epi.  Since the unit map ${\eta''}_{G'F'\mathcal{F}}$ is an isomorphism, the commutativity of the diagram
$$
\xymatrix{
\mathcal{F} \ar[rr]^{{\eta'}_{\mathcal{F}}} \ar[d]_{{\eta''}_{\mathcal{F}}} & &  G'F'\mathcal{F} \ar[d]^{{\eta''}_{G'F'\mathcal{F}}} \\
i_{*}i^{*}\mathcal{F} \ar[rr]_{i_{*}i^{*}{\eta'}_{\mathcal{F}}} & & i_{*}i^{*}G'F'\mathcal{F}
}
$$
implies ${\eta'}_{\mathcal{F}}$ is an epi.  Therefore, in order to show ${\eta'}_{\mathcal{F}}$ is epi, it suffices to show $i^{*}{\eta'}_{\mathcal{F}}$ is epi.  We prove this.  By Lemmas $\ref{lem.hom}$ and $\ref{lem.free}$, there is a natural isomorphism
$$
\delta:i^{*}\mathcal{H}\it{om}_{{\mathcal{O}}_{(U \times X)_{U}^{2}}}(\pr_{1}^{*}\mathcal{M},\pr_{1}^{*}\mathcal{M} \otimes \mathcal{F}) \rightarrow i^{*}\mathcal{F}
$$
such that $\delta i^{*} \eta' = \operatorname{id}_{\mathcal{F}}$.  Thus, $i^{*}\eta'$ is an epi so that $\eta'$ is an epi.

Next, by Proposition $\ref{prop.isomp}$, the composition of unit maps
$$
\pr_{1}^{*}\mathcal{M} \otimes \mathcal{F} \rightarrow i_{*}i^{*}(\pr_{1}^{*}\mathcal{M} \otimes \mathcal{F}) \rightarrow i_{*}(\pr_{2}i)^{!}(\pr_{2}i)_{*}i^{*}(\pr_{1}^{*}\mathcal{M} \otimes \mathcal{F})
$$
is an isomorphism.  Thus, we have a natural isomorphism
$$
\mathcal{H}\it{om}_{{\mathcal{O}}_{(U \times X)_{U}^{2}}}(\pr_{1}^{*}\mathcal{M},\pr_{1}^{*}\mathcal{M} \otimes \mathcal{F}) \rightarrow GF\mathcal{F}
$$
making the diagram
$$
\xymatrix{
& \mathcal{H}\it{om}_{{\mathcal{O}}_{(U \times X)_{U}^{2}}}(\pr_{1}^{*}\mathcal{M},\pr_{1}^{*}\mathcal{M} \otimes \mathcal{F}) \ar[d] \\
\mathcal{F} \ar[ur]^{\eta'_{\mathcal{F}}} \ar[r]_{\eta_{\mathcal{F}}} & GF\mathcal{F}
}
$$
commute.  In particular, the unit map $\eta_{\mathcal{F}}:\mathcal{F} \rightarrow GF \mathcal{F}$ is an epimorphism.

Finally, suppose
$$
f:(\pr_{2}i)_{*}i^{*}(\pr_{1}^{*}\mathcal{M} \otimes \mathcal{F}) \rightarrow \mathcal{N}
$$
is an epimorphism.  We must show that the induced map of sheaves,
$$
Gf:G(\pr_{2}i)_{*}i^{*}(\pr_{1}^{*}\mathcal{M} \otimes \mathcal{F}) \rightarrow G\mathcal{N}
$$
is also an epimorphism.  Since, by Proposition $\ref{prop.isomp}$, $i_{*}$ and $(\pr_{2}i)^{!}$ are exact, we note that $i_{*}(\pr_{2}i)^{!}f$ is an epimorphism.  To show $Gf$ is an epimorphism, we need only show that $i^{*}Gf$ is an epimorphism.  To this end, we note that $\eta:\id \Longrightarrow i_{*}i^{*}$ applied to both the domain and codomain of $i_{*}(\pr_{2}i)^{!}f$ is an isomorphism.  Thus, by Lemmas $\ref{lem.hom}$ and $\ref{lem.free}$, $i^{*}G'$ applied to $i_{*}(\pr_{2}i)^{!}f$ is an epimorphism.  Since $i^{*}G=i^{*}G'i_{*}(\pr_{2}i)^{!}$, $i^{*}Gf$ is an epimorphism so that $Gf$ is an epimorphism as desired.
\end{proof}

\subsection{The proof of Theorem \ref{theorem.rep1}}
If $\mathcal{G}$ is an ${\mathcal{O}}_{X^{2}}$-module, and $p:(U \times X)_{U}^{2} \rightarrow X^{2}$ is projection, then let ${\mathcal{G}}^{U} = p^{*}\mathcal{G}$.  If $q_{i},q_{j}:U \rightarrow X$ are morphisms and $d:U \rightarrow U \times U$ is the diagonal morphism, then let
$$
q_{ij} = (\id_{U} \times q_{i} \times \id_{U} \times q_{j})\circ(d \times d) \circ d:U \rightarrow (U \times X)_{U}^{2}.
$$
Suppose $\mathcal{M}$ is a free truncated $U$-family of length $n+1$ with multiplication maps $\mu_{i,j}$.  For $0 \leq k \leq n$, let
$$
i_{k}:\mbox{SSupp }\pr_{1}^{*}{\mathcal{M}}_{k} \rightarrow (U \times X)_{U}^{2}
$$
be inclusion and let
$$
F^{\mathcal{M}}_{k}= (\pr_{2}i_{k})_{*} i_{k}^{*} (\pr_{1}^{*}{\mathcal{M}}_{k} \otimes_{{\mathcal{O}}_{(U \times X)_{U}^{2}}} - )
$$
and
$$
G^{\mathcal{M}}_{k}= {\mathcal{H}}{\it om}_{{\mathcal{O}}_{(U \times X)_{U}^{2}}}(\pr_{1}^{*}{\mathcal{M}}_{k},-) i_{k*} (\pr_{2}i_{k})^{!}.
$$
Finally, for $i,j \geq 0$, let the right adjunct of $\mu_{i,j} \in \operatorname{Hom }_{U \times X}(F^{\mathcal{M}}_{i}{\mathcal{E}}^{U \otimes j},{\mathcal{M}}_{i+j})$ be denoted $\operatorname{Ad }\mu_{i,j} \in \operatorname{Hom }_{(U \times X)_{U}^{2}}({\mathcal{E}}^{U \otimes j},G^{\mathcal{M}}_{i}{\mathcal{M}}_{i+j})$.

\begin{proposition} \label{prop.trans}
For $1 \leq n$, there is a natural transformation
$$
\Phi:\Gamma_{n}(-) \Longrightarrow \operatorname{Hom}_{S}(-,{\mathbb{P}}_{X^{2}}(\mathcal{E})^{\otimes n}).
$$
Furthermore, suppose $U$ is a noetherian affine scheme and
$$
j:\operatorname{Hom}_{S}(U,{\mathbb{P}}_{X^{2}}(\mathcal{E})^{\otimes n}) \rightarrow \operatorname{Hom}_{S}(U,{\mathbb{P}}_{(U \times X)_{U}^{2}}({\mathcal{E}}^{U})^{\otimes n})
$$
is the map sending $f$ to the composition
$$
\xymatrix{
U \ar[r]^{\hskip -.4in f \times \id_{U}} & {\mathbb{P}}_{X^{2}}(\mathcal{E})^{\otimes n} \times_{S} U \ar[r]^{\cong} & {\mathbb{P}}_{(U \times X)_{U}^{2}}({\mathcal{E}}^{U})^{\otimes n}.
}
$$
Then $j \circ \Phi_{U}([\mathcal{M}])$ corresponds to the $n$ epimorphisms
$$
\xymatrix{
{\mathcal{E}}^{U} \ar[r]^{\hskip -.1in \operatorname{Ad }\mu_{i,1}} & G_{i}{\mathcal{M}}_{i+1} \ar[r]^{\cong} & q_{i,i+1*}{\mathcal{O}}_{U}
}
$$
whose rightmost map is given by Proposition \ref{prop.simpleiso}.
\end{proposition}

\begin{proof}
We construct a map of sets
\begin{equation} \label{eqn.prim}
\phi_{U}:\{\mbox{truncated free $U$-families of length $n+1$} \} \rightarrow \operatorname{Hom}_{S}^{\Fr}(U,{\mathbb{P}}_{X^{2}}(\mathcal{E})^{\otimes n}).
\end{equation}\index{family!free, of truncated point modules|(}
We then show the $\phi_{U}$ does not distinguish between isomorphic families, so $\phi_{U}$ induces a set map
$$
\phi_{U}^{\Fr}:\Gamma_{n}^{\Fr}(U) \rightarrow \operatorname{Hom}_{S}^{\Fr}(U,{\mathbb{P}}_{X^{2}}(\mathcal{E})^{\otimes n}).
$$
We next show that $\phi_{U}^{\Fr}$ is natural, i.e. defines a natural transformation
$$
\Phi^{\Fr}:\Gamma_{n}^{\Fr}(-) \Longrightarrow \operatorname{Hom}_{S}^{\Fr}(-,{\mathbb{P}}_{X^{2}}(\mathcal{E})^{\otimes n}).
$$
By Corollary \ref{cor.freefunct}, $\Phi^{\Fr}$ extends to a natural transformation
$$
\Phi:\Gamma_{n}(-) \Longrightarrow \operatorname{Hom}_{S}(-,{\mathbb{P}}_{X^{2}}(\mathcal{E})^{\otimes n})
$$
as desired.

{\it Step 1:  the construction of a set map
$$
\phi_{U}:\{\mbox{truncated free $U$-families of length $n+1$} \} \rightarrow \operatorname{Hom}_{S}^{\Fr}(U,{\mathbb{P}}_{X^{2}}(\mathcal{E})^{\otimes n}).
$$
}  We define $\phi_{U}$ (\ref{eqn.prim}).  For $0 \leq j \leq n$, let $q_{j}:U \rightarrow X$ be a morphism such that
$$
{\mathcal{M}}_{j} \cong (\id_{U} \times q_{j})_{*}{\mathcal{O}}_{U}.
$$
and let
$$
\mu_{j,1}:F^{\mathcal{M}}_{j}{\mathcal{E}}^{U} \rightarrow {\mathcal{M}}_{j+1}
$$
be the ${\mathcal{B}}^{U}$-module multiplication map of $\mathcal{M}$.  Since $\mu_{j,1}$ is an epi by definition of $\mathcal{M}$, by Proposition \ref{prop.bij} $\mu_{j,1}$ corresponds to an epi
$$
\gamma_{j}:  {\mathcal{E}}^{U} \rightarrow G_{j}^{\mathcal{M}}{\mathcal{M}}_{j+1}.
$$
By Proposition $\ref{prop.simpleiso}$, this epimorphism can be composed with an isomorphism
$$
G_{j}^{\mathcal{M}}{\mathcal{M}}_{j+1} \rightarrow q_{j,j+1*}{\mathcal{O}}_{U}
$$
to give an epimorphism
$$
{\gamma'}_{j}:{\mathcal{E}}^{U} \rightarrow q_{j,j+1*}{\mathcal{O}}_{U}.
$$
By the adjointness of the pair $(q_{j,j+1}^{*},q_{j,j+1*})$, and by the affinity of $q_{j,j+1}$, this epimorphism corresponds to an epimorphism
\begin{equation} \label{eqn.sillness}
q_{j,j+1}^{*}p^{*}\mathcal{E} = q_{j,j+1}^{*}{\mathcal{E}}^{U} \rightarrow {\mathcal{O}}_{U}
\end{equation}
By Proposition \ref{prop.groth}, this epimorphism corresponds to a map
$$
r_{j}:U \rightarrow {\mathbb{P}}_{(U\times X)_{U}^{2}}({\mathcal{E}}^{U}) \cong {\mathbb{P}}_{X^{2}}(\mathcal{E}) \times_{X^{2}}(U \times X)_{U}^{2} \rightarrow {\mathbb{P}}_{X^{2}}(\mathcal{E})
$$
whose projection to the base is $q_{j} \times q_{j+1}$.  Thus, the maps $(r_{0}, \ldots, r_{n-1})$ give us a map
$$
r:U \rightarrow {\mathbb{P}}_{X^{2}}(\mathcal{E})^{\otimes n}
$$
whose projection to the base is $(q_{0},q_{n})$.  We define $\phi_{U}(\mathcal{M})=r$.  Since the left-most composite of $r_{j}$ is a $U$-morphism, the second assertion of the Proposition holds.

{\it Step 2: $\phi_{U}$ induces a set map
$$
\phi_{U}^{\Fr}:\Gamma_{n}^{\Fr}(U) \rightarrow \operatorname{Hom}_{S}^{\Fr}(U,{\mathbb{P}}_{X^{2}}(\mathcal{E})^{\otimes n}).
$$
}
We need to show $\phi_{U}$ does not distinguish between isomorphic families.  Suppose $\mathcal{N}$ is another truncated $U$-family of length $n+1$ and $\theta: \mathcal{M} \rightarrow \mathcal{N}$ is an isomorphism of ${\mathcal{B}}^{U}$-modules.  Suppose, as above, that $\mathcal{N}$ determines an epimorphism
$$
\delta_{j}:{\mathcal{E}}^{U} \rightarrow \mathcal{H}\it{om}_{{\mathcal{O}}_{(U\times X)_{U}^{2}}}(\pr_{1}^{*}{\mathcal{N}}_{j},i_{j*}(\pr_{2}i)^{!}{\mathcal{N}}_{j+1})).
$$
and
$$
{\delta'}_{j}:{\mathcal{E}}^{U} \rightarrow  q_{j,j+1*}{\mathcal{O}}_{U}.
$$
We claim that, to show $\phi_{U}(\mathcal{M})=\phi_{U}(\mathcal{N})$, it suffices to show $\operatorname{ker }\gamma_{j}=\operatorname{ker }\delta_{j}$.  If this held, then there would exist an isomorphism
$$
\tau_{j}:q_{j,j+1*}{\mathcal{O}}_{U} \rightarrow q_{j,j+1*}{\mathcal{O}}_{U}
$$
such that
\begin{equation} \label{eqn.tau}
\tau_{j} {\gamma'}_{j} = {\delta'}_{j}.
\end{equation}
But since $q_{j,j+1}$ is a closed immersion, the counit map
$$
q_{j,j+1}^{*}q_{j,j+1*}{\mathcal{O}}_{U} \rightarrow {\mathcal{O}}_{U}
$$
is an isomorphism.  Thus, applying the functor $q_{j,j+1*}$ to equation (\ref{eqn.tau}) allows us to deduce that $\phi_{U}(\mathcal{M}) = \phi_{U}(\mathcal{N})$.  We conclude that, in order to show $\phi_{U}(\mathcal{M})=\phi_{U}(\mathcal{N})$, it suffices to show that $\operatorname{ker }\gamma_{j}=\operatorname{ker }\delta_{j}$.  Let $c_{j}:\operatorname{ker }\delta_{j} \rightarrow {\mathcal{E}}^{U}$.  We show $\gamma_{j} c_{j}=0$.  Since, for any $U$-family $\mathcal{P}$, $(F_{j}^{\mathcal{P}},G_{j}^{\mathcal{P}})$ is an adjoint pair, the diagram
\begin{equation} \label{eqn.biggle}
\xymatrix{
\operatorname{Hom}_{U \times X}(F_{j}^{\mathcal{N}}{\mathcal{E}}^{U},{\mathcal{M}}_{j+1}) \ar[d]_{- \circ F_{j}^{\mathcal{N}}c_{j}}  \ar[r] & \operatorname{Hom}_{U \times X}({\mathcal{E}}^{U},G_{j}^{\mathcal{N}}{\mathcal{M}}_{j+1}) \ar[d]^{- \circ c_{j}} \\
\operatorname{Hom}_{U \times X}(F_{j}^{\mathcal{N}}\operatorname{ker }\delta_{j},{\mathcal{M}}_{j+1}) \ar[r] & \operatorname{Hom}_{U \times X}(\operatorname{ker }\delta_{j},G_{j}^{\mathcal{N}}{\mathcal{M}}_{j+1})
}
\end{equation}
commutes.  If
$$
\nu_{j}:F_{j}^{\mathcal{N}}{\mathcal{E}}^{U}={\mathcal{N}}_{j} \otimes_{{\mathcal{O}}_{U \times X}} {\mathcal{E}}^{U} \rightarrow {\mathcal{N}}_{j+1}
$$
is ${\mathcal{B}}^{U}$-module multiplication, then
$$
\theta_{j+1}^{-1}\nu_{j} \in \operatorname{Hom}_{U \times X}(F_{j}^{\mathcal{N}}{\mathcal{E}}^{U},{\mathcal{M}}_{j+1})
$$
goes, via the top of (\ref{eqn.biggle}), to the map
$$
\xymatrix{
{\mathcal{E}}^{U} \ar[r] & G_{j}^{\mathcal{N}}F_{j}^{\mathcal{N}}{\mathcal{E}}^{U} \ar[r]^{G_{j}^{\mathcal{N}}\nu} & G_{j}^{\mathcal{N}}{\mathcal{N}}_{j+1} \ar[r]^{G_{j}^{\mathcal{N}}{\theta_{j+1}^{-1}}} & G_{j}^{\mathcal{N}}{\mathcal{M}}_{j+1}
}
$$
where the left-most map is the unit of the pair $(F_{j}^{\mathcal{N}},G_{j}^{\mathcal{N}})$.  The left-most two maps compose to give $\delta_{j}$.  Thus, via the top route of (\ref{eqn.biggle}), $\theta_{j+1}^{-1}\nu_{j}$ goes to zero.  By the commutativity of (\ref{eqn.biggle}),
\begin{equation} \label{eqn.nicce}
\theta_{j+1}^{-1}\nu_{j} \circ F_{j}^{\mathcal{N}}c_{j} \in \operatorname{Hom}_{U \times X}(F_{j}^{\mathcal{N}}\operatorname{ker }\delta_{j},{\mathcal{M}}_{j+1})
\end{equation}
equals zero.  In addition, since $\theta$ is a $\mathcal{B}$-module morphism, the diagram
$$
\xymatrix{
{\mathcal{M}}_{j} \otimes_{{\mathcal{O}}_{U \times X}}{\mathcal{E}}^{U} \ar[r]^{\mu_{j}} \ar[d]_{\theta_{j} \otimes_{{\mathcal{O}}_{U \times X}} {\mathcal{E}}^{U}} & {\mathcal{M}}_{j+1} \ar[d]^{\theta_{j+1}} \\
{\mathcal{N}}_{j} \otimes_{{\mathcal{O}}_{U \times X}}{\mathcal{E}}^{U} \ar[r]_{\nu_{j}} & {\mathcal{N}}_{j+1}
}
$$
commutes so that
\begin{equation} \label{eqn.ya}
{\theta_{j+1}}^{-1}\nu_{j}(\theta_{j} \otimes_{{\mathcal{O}}_{U \times X}} {\mathcal{E}}^{U}) = \mu_{j}.
\end{equation}
We next note that $\pr_{2*}(\pr_{1}^{*}\theta_{j} \otimes -):F_{j}^{\mathcal{M}} \Longrightarrow F_{j}^{\mathcal{N}}$ is a natural equivalence, so that the diagram
\begin{equation} \label{eqn.equiv}
\xymatrix{
\operatorname{Hom}_{U \times X}(F_{j}^{\mathcal{N}}{\mathcal{E}}^{U},{\mathcal{M}}_{j+1}) \ar[d]_{- \circ F_{j}^{\mathcal{N}}c_{j}} \ar[rr]^{ - \circ {\pr_{2}}_{*}(\pr_{1}^{*}\theta_{j} \otimes {\mathcal{E}}^{U})} & & \operatorname{Hom}_{U \times X}(F_{j}^{\mathcal{M}}{\mathcal{E}}^{U},{\mathcal{M}}_{j+1}) \ar[d]^{- \circ F_{j}^{\mathcal{M}}c_{j}} \\
\operatorname{Hom}_{U \times X}(F_{j}^{\mathcal{N}}\operatorname{ker }\delta_{j},{\mathcal{M}}_{j+1})  \ar[rr] & & \operatorname{Hom}_{U \times X}(F_{j}^{\mathcal{M}}\operatorname{ker }\delta_{j},{\mathcal{M}}_{j+1})
}
\end{equation}
commutes by functoriality of $\operatorname{Hom}_{U \times X}(-,{\mathcal{M}}_{j+1})$.  Thus, since $\theta_{j+1}^{-1}\nu_{j} \circ F_{j}^{\mathcal{N}}c_{j} = 0$ by (\ref{eqn.nicce}), and since ${\theta_{j+1}}^{-1}\nu_{j}(\theta_{j} \otimes_{{\mathcal{O}}_{U \times X}} {\mathcal{E}}^{U}) = \mu_{j}$ by (\ref{eqn.ya}), the commutativity of (\ref{eqn.equiv}) implies
\begin{equation} \label{eqn.gotit}
\mu_{j}F_{j}^{\mathcal{M}}c_{j}=0.
\end{equation}
Finally, by the adjointness of $(F_{j}^{\mathcal{M}},G_{j}^{\mathcal{M}})$, the diagram
$$
\xymatrix{
\operatorname{Hom}_{U \times X}(F_{j}^{\mathcal{M}}{\mathcal{E}}^{U},{\mathcal{M}}_{j+1}) \ar[d]_{- \circ F_{j}^{\mathcal{M}}c_{j}} \ar[r] & \operatorname{Hom}_{U \times X}({\mathcal{E}}^{U},G_{j}^{\mathcal{M}}{\mathcal{M}}_{j+1}) \ar[d]^{- \circ c_{j}} \\
\operatorname{Hom}_{U \times X}(F_{j}^{\mathcal{M}}\operatorname{ker }\delta_{j},{\mathcal{M}}_{j+1}) \ar[r] & \operatorname{Hom}_{U \times X}(\operatorname{ker }\delta_{j},G_{j}^{\mathcal{M}}{\mathcal{M}}_{j+1})
}
$$
commutes.  Since $\mu_{j}$, an element in the upper left, goes to $\gamma_{j}$ on the upper right, (\ref{eqn.gotit}) implies $\gamma_{j}c_{j} = 0$ as desired.  We conclude $\operatorname{ker} \gamma_{j} \subset \operatorname{ker} \delta_{j}$.  A similar argument shows $\operatorname{ker} \delta_{j} \subset \operatorname{ker} \gamma_{j}$ so $\operatorname{ker} \gamma_{j} = \operatorname{ker} \delta_{j}$.  Thus, $\phi_{U}$ induces a map
$$
\phi^{\Fr}_{U}:\Gamma_{n}^{\Fr}(-) \rightarrow \operatorname{Hom}_{k}^{\Fr}(-,{\mathbb{P}}_{X^{2}}(\mathcal{F})^{\otimes n}).
$$

{\it Step 3:  compatibilities related to the proof that $\Phi^{\Fr}$ is natural}.  Let $V$ be an affine scheme, suppose $f:V \rightarrow U$ is a morphism, and let $\tilde{f}=f\times \id_{X}:V \times X \rightarrow U \times X$.  Let $\mathcal{M}$ be a ${\mathcal{B}}^{U}$-module, and give $\tilde{f}^{*}\mathcal{M}$ its ${\mathcal{B}}^{V}$-module structure (Theorem \ref{theorem.lift}).  Let $F_{j}^{U} = F_{j}^{\mathcal{M}}$ and $G_{j}^{U} = G_{j}^{\mathcal{M}}$.  Let $i_{j}^{V}: \operatorname{SSupp }\pr_{1}^{V*}\tilde{f}^{*}{\mathcal{M}}_{j} \rightarrow (V \times X)_{V}^{2}$ be inclusion.  Finally, let
$$
F_{j}^{V}:{\sf{Qcoh }} (V \times X)_{V}^{2} \rightarrow {\sf{Qcoh }} (V\times X)_{V}^{2}
$$
be the functor $\pr_{2*}^{V}(\pr_{1}^{V*}{\tilde{f}^{*}\mathcal{M}}_{j} \otimes -)$, and let
$$
G_{j}^{V}:{\sf{Qcoh }} V \times X \rightarrow {\sf{Qcoh }} (V \times X)_{V}^{2}
$$
be the functor $\mathcal{H}\it{om}_{{\mathcal{O}}_{(V \times X)_{V}^{2}}}(\pr_{1}^{*}\tilde{f}^{*}\mathcal{M}_{j},-) \circ {i}_{j*}^{V} \circ ({\pr_{2}^{V}i_{j}^{V}})^{!}$.  Since $(F_{j}^{U},G_{j}^{U})$ and $(F_{j}^{V},G_{j}^{V})$ are adjoint pairs, if $\Psi^{1}$ is the composition of isomorphisms
$$
F_{j}^{V}\tilde{f}^{2*}= \pr_{2*}^{V}(\pr_{1}^{V*}{\tilde{f}^{*}\mathcal{M}}_{j} \otimes \tilde{f}^{2*}-) \Longrightarrow \pr_{2*}^{V}\tilde{f}^{2*}(\pr_{1}^{U*}{\mathcal{M}}_{j} \otimes -) \Longrightarrow
$$
$$
\tilde{f}^{*}\pr_{2*}(\pr_{1}^{*}{\mathcal{M}}_{j} \otimes -) = \tilde{f}^{*}F_{j}^{U},
$$
where the second nontrivial map is the isomorphism from Proposition \ref{prop.canon}, (2), there exists, by Lemma \ref{prop.duality}, a natural transformation $\Phi^{1}:\tilde{f}^{2*}G_{j}^{U} \Longrightarrow G_{j}^{V}\tilde{f}^{*}$ such that $\Psi^{1}$ is the dual 2-cell to the square
$$
\xymatrix{
{\sf{Qcoh }} U \times X \ar[d]_{\tilde{f}^{*}} \ar@<1ex>[r]^{G_{j}^{U}} \ar@{=>}[dr]|{\Phi^{1}}\hole & {\sf{Qcoh }} (U \times X)_{U}^{2} \ar@<1ex>[l]^{F_{j}^{U}} \ar[d]^{{\tilde{f}}^{2*}} \\
{\sf{Qcoh }} V \times X  \ar@<1ex>[r]^{G_{j}^{V}} & {\sf{Qcoh }} (V \times X)^{2}_{V} \ar@<1ex>[l]^{F^{V}_{j}}  .
}
$$
By Corollary \ref{cor.adcom2}, the diagram
\begin{equation} \label{eqn.digge1}
\xymatrix{
\mbox{Hom}_{U \times X}(F_{j}^{U}{\mathcal{E}}^{U},{\mathcal{M}}_{j+1}) \ar[r] \ar[d]_{\tilde{f}^{*}(-)} & \mbox{Hom}_{(U \times X)_{U}^{2}}({\mathcal{E}}^{U},G_{j}^{U}{\mathcal{M}}_{j+1}) \ar[d]^{\tilde{f}^{2*}(-)} \\
\mbox{Hom}_{V \times X}(\tilde{f}^{*}F_{j}^{U}{\mathcal{E}}^{U},\tilde{f}^{*}{\mathcal{M}}_{j+1}) \ar[d]_{- \circ \Psi^{1}_{{\mathcal{E}}^{U}}} & \mbox{Hom}_{(V \times X)_{V}^{2}}(\tilde{f}^{2*}{\mathcal{E}}^{U},\tilde{f}^{2*}G_{j}^{U}{\mathcal{M}}_{j+1}) \ar[d]^{\Phi^{1}_{{\mathcal{M}}_{j+1}} \circ -} \\
\mbox{Hom}_{V \times X}(F_{j}^{V}\tilde{f}^{2*}{\mathcal{E}}^{U},\tilde{f}^{*}{\mathcal{M}}_{j+1}) \ar[r] & \mbox{Hom}_{(V \times X)_{V}^{2}}(\tilde{f}^{2*}{\mathcal{E}}^{U},G_{j}^{V}\tilde{f}^{*}{\mathcal{M}}_{j+1})
}
\end{equation}
with horizontal maps the adjoint isomorphisms, commutes.  Let $qf_{j,j+1}$ denote the closed immersion $(\id_{V} \times q_{j}f) \times (\id_{V} \times q_{j+1}f):V \rightarrow (V \times X)_{V}^{2}$.  Since $(q_{j,j+1}^{*},q_{j,j+1*})$ and $(qf_{j,j+1}^{*},qf_{j,j+1*})$ are adjoint pairs, if $\Psi^{2}$ is the canonical isomorphism $(qf_{j,j+1})^{*}\tilde{f}^{2*} \Longrightarrow f^{*}q_{j,j+1}^{*}$ induced by the commutative diagram of schemes
$$
\xymatrix{
V \ar[rr]^{qf_{j,j+1}} \ar[d]_{f} & & (V \times X)_{V}^{2} \ar[d]^{\tilde{f}^{2}} \\
U \ar[rr]_{q_{j,j+1}} & & (U \times X)_{U}^{2}
}
$$
then there exists, by Lemma \ref{prop.duality}, a natural transformation $\Phi^{2}:\tilde{f}^{2*}q_{j,j+1*} \Longrightarrow qf_{j,j+1*}f^{*}$ such that $\Psi^{2}$ is dual to the square
$$
\xymatrix{
{\sf{Qcoh }} U \ar[d]_{\tilde{f}^{*}} \ar@<1ex>[rr]^{q_{j,j+1*}}  \ar@{=>}[drr]|{\Phi^{2}}\hole  & & {\sf{Qcoh }} (U \times X)_{U}^{2} \ar@<1ex>[ll]^{q_{j,j+1}^{*}} \ar[d]^{\tilde{f}^{2*}} \\
{\sf{Qcoh }} V  \ar@<1ex>[rr]^{qf_{j,j+1*}} & & {\sf{Qcoh }} (V \times X)_{V}^{2}  \ar@<1ex>[ll]^{qf_{j,j+1}^{*}}.
}
$$
Thus, again by Corollary \ref{cor.adcom2}, the diagram
\begin{equation} \label{eqn.digge2}
\xymatrix{
\operatorname{Hom}_{(U \times X)_{U}^{2}}({\mathcal{E}}^{U},q_{j,j+1*}{\mathcal{O}}_{U}) \ar[d]_{\tilde{f}^{2*}(-)} \ar[r] & \operatorname{Hom}_{U}(q_{j,j+1}^{*}{\mathcal{E}}^{U},{\mathcal{O}}_{U}) \ar[d]^{f^{*}(-)} \\
\operatorname{Hom}_{(V \times X)_{V}^{2}}(\tilde{f}^{2*}{\mathcal{E}}^{U},\tilde{f}^{2*}q_{j,j+1*}{\mathcal{O}}_{U}) \ar[d]_{\Phi^{2}_{{\mathcal{O}}_{U}} \circ -} & \operatorname{Hom}_{V}(f^{*}q_{j,j+1}^{*}{\mathcal{E}}^{U},f^{*}{\mathcal{O}}_{U}) \ar[d]^{- \circ \Psi^{2}_{{\mathcal{E}}^{U}}} \\
\operatorname{Hom}_{(V \times X)_{V}^{2}}(\tilde{f}^{2*}{\mathcal{E}}^{U},q_{j,j+1}f_{*}f^{*}{\mathcal{O}}_{U}) \ar[r] & \operatorname{Hom}_{V}(qf_{j,j+1}^{*}\tilde{f}^{2*}{\mathcal{E}}^{U},f^{*}{\mathcal{O}}_{U})
}
\end{equation}
with horizontal maps the adjoint isomorphisms, commutes.

Finally, if
$$
\zeta:G_{j}^{U}{\mathcal{M}}_{j+1} \rightarrow q_{j,j+1*}{\mathcal{O}}_{U}
$$
is the isomorphism in Proposition $\ref{prop.simpleiso}$, then the diagram
\begin{equation} \label{eqn.middle}
\xymatrix{
\operatorname{Hom}_{(U \times X)_{U}^{2}}({\mathcal{E}}^{U},G_{j}^{U}{\mathcal{M}}_{j+1}) \ar[r]^{\zeta \circ -} \ar[d]_{\tilde{f}^{2*}(-)} & \operatorname{Hom}_{(U \times X)_{U}^{2}}({\mathcal{E}}^{U},q_{j,j+1*}{\mathcal{O}}_{U}) \ar[d]^{\tilde{f}^{2*}(-)} \\
\operatorname{Hom}_{(V \times X)_{V}^{2}}(\tilde{f}^{2*}{\mathcal{E}}^{U},{\tilde{f}}^{2*}G_{j}^{U}{\mathcal{M}}_{j+1}) \ar[r]^{\tilde{f}^{2*}\zeta \circ -} \ar[d]_{\Phi^{1}_{{\mathcal{M}}_{j+1}} \circ -} & \operatorname{Hom}_{(V \times X)_{V}^{2}}(\tilde{f}^{2*}{\mathcal{E}}^{U},\tilde{f}^{2*}q_{j,j+1*}{\mathcal{O}}_{U}) \ar[d]^{\Phi^{2}_{{\mathcal{O}}_{U}} \circ -} \\
\operatorname{Hom}_{(V \times X)_{V}^{2}}(\tilde{f}^{2*}{\mathcal{E}}^{U},G_{j}^{V}\tilde{f}^{*}{\mathcal{M}}_{j+1}) & \operatorname{Hom}_{(V \times X)_{V}^{2}}(\tilde{f}^{2*}{\mathcal{E}}^{U},qf_{j,j+1*}f^{*}{\mathcal{O}}_{U})
}
\end{equation}
commutes by the functoriality of $\tilde{f}^{2*}$.

{\it Step 4:  $\Phi^{1}_{{\mathcal{M}}_{j+1}}$ and $\Phi^{2}$ are isomorphisms}.  We first show that $\Phi^{2}$ is an isomorphism and use this result to show that $\Phi^{1}_{{\mathcal{M}}_{j+1}}$ is an isomorphism.  By Proposition \ref{prop.duality}, $\Phi^{2}$ equals the composition
$$
{\tilde{f}}^{2*}q_{j,j+1*} \Longrightarrow qf_{j,j+1*}qf_{j,j+1}^{*}{\tilde{f}}^{2*}q_{j,j+1*} \Longrightarrow
$$
$$
qf_{j,j+1*}f^{*}q_{j,j+1}^{*}q_{j,j+1*} \Longrightarrow qf_{j,j+1*}f^{*}
$$
where the first map is the unit of the pair $(qf_{j,j+1}^{*}qf_{j,j+1*})$, the second map is the isomorphism $qf_{j,j+1*}*\Psi^{2}*q_{j,j+1*}$, and the last map is the counit of the pair $(q_{j,j+1}^{*},q_{j,j+1*})$.  Since $q_{j,j+1}$ is a closed immersion, this last map is an isomorphism.  Thus, to show that $\Phi^{2}$ is an isomorphism, we need only show that the unit
\begin{equation} \label{eqn.unite}
{\tilde{f}}^{2*}q_{j,j+1*} \Longrightarrow qf_{j,j+1*}qf_{j,j+1}^{*}{\tilde{f}}^{2*}q_{j,j+1*}
\end{equation}
is an isomorphism.  Since, by Lemma \ref{lem.anothersquare}, there is an isomorphism ${\tilde{f}}^{2*}q_{j,j+1*} \Longrightarrow qf_{j,j+1*}f^{*}$, and since the unit
$$
qf_{j,j+1*}f^{*} \Longrightarrow qf_{j,j+1*}qf_{j,j+1}^{*}qf_{j,j+1*}f^{*}
$$
is an isomorphism ($qf_{j,j+1}$ is a closed immersion), by naturality of the unit, (\ref{eqn.unite}) is an isomorphism.  Thus, $\Phi^{2}$ is an isomorphism.

We now show that $\Phi^{1}_{{\mathcal{M}}_{j+1}}$ is an isomorphism.  Since $\Psi^{1}$ is an isomorphism it suffices, by Proposition \ref{prop.duality}, to show that the unit
\begin{equation} \label{eqn.unitya}
{\tilde{f}}^{2*}G_{j}^{U}{\mathcal{M}}_{j+1} \rightarrow G_{j}^{V}F_{j}^{V}{\tilde{f}}^{2*}G_{j}^{U}{\mathcal{M}}_{j+1}
\end{equation}
and the counit
\begin{equation} \label{eqn.counitya}
F_{j}^{U}G_{j}^{U}{\mathcal{M}}_{j+1} \rightarrow {\mathcal{M}}_{j+1}
\end{equation}
are isomorphisms.  The counit (\ref{eqn.counitya}) is an isomorphism by Lemma \ref{lem.unit}.  Thus, we only need to show that (\ref{eqn.unitya}) is an isomorphism.  By Lemma \ref{lem.unit}, the unit
$$
G_{j}^{V}{\tilde{f}}^{*}{\mathcal{M}}_{j+1} \rightarrow G_{j}^{V}F_{j}^{V}G_{j}^{V}{\tilde{f}}^{*}{\mathcal{M}}_{j+1}
$$
is an isomorphism.  In addition, by Proposition \ref{prop.simpleiso}, there are isomorphisms
$$
\xymatrix{
G_{j}^{V}{\tilde{f}}^{*}{\mathcal{M}}_{j+1} \ar[r]^{\cong} & qf_{j,j+1*}f^{*}{\mathcal{O}}_{U} \ar[r]^{(\Phi^{2}_{{\mathcal{O}}_{U}})^{-1}} & {\tilde{f}}^{2}q_{j,j+1*}{\mathcal{O}}_{U} \ar[r]^{\cong} & {\tilde{f}}^{2*}G_{j}^{U}{\mathcal{M}}_{j+1}.
}
$$
Thus, by naturality of the unit, (\ref{eqn.unitya}) is an isomorphism.

{\it Step 5:  completion of the proof that $\Phi^{\Fr}$ is natural}.  The top row of the left to right combination of diagrams (\ref{eqn.digge1}) (on page \pageref{eqn.digge1}), (\ref{eqn.middle}), and (\ref{eqn.digge2}) sends $\mu_{j,1}$, the $j$th component of multiplication of the ${\mathcal{B}}^{U}$-module $\mathcal{M}$, to an epimorphism $q_{j,j+1}^{*}{\mathcal{E}}^{U} \rightarrow {\mathcal{O}}_{U}$ corresponding to $r_{j}:U \rightarrow {\mathbb{P}}_{X^{2}}(\mathcal{E})$.

We show the right hand vertical of the combined diagram sends the epimorphism corresponding to $r_{j}$ to the epimorphism corresponding to $r_{j}f$.  Suppose $\phi:q_{j,j+1}^{*}{\mathcal{E}}^{U} \rightarrow {\mathcal{O}}_{U}$ is an epimorphism corresponding to $r_{j}$.  Then, in particular, $\phi$ corresponds, by Proposition \ref{prop.groth}, to a map $g:U \rightarrow {\mathbb{P}}_{(U \times X)_{U}^{2}}({\mathcal{E}}^{U})$ which, when composed with the projection ${\mathbb{P}}_{(U \times X)_{U}^{2}}({\mathcal{E}}^{U}) \rightarrow {\mathbb{P}}_{X^{2}}(\mathcal{E})$ ((\ref{eqn.proj}), defined on page \pageref{eqn.proj}), equals $r_{j}$.  The right hand vertical of (\ref{eqn.digge2}) sends $\phi$ to the composition
\begin{equation} \label{eqn.original}
\xymatrix{
qf_{j,j+1}^{*}\tilde{f}^{2*}{\mathcal{E}}^{U} \ar[r] & (\tilde{f}^{2}qf_{j,j+1})^{*}{\mathcal{E}}^{U}=(q_{j,j+1}f)^{*}{\mathcal{E}}^{U} \ar[r] & f^{*}q_{j,j+1}^{*}{\mathcal{E}}^{U} \ar[r]^{f^{*}\phi} & f^{*}{\mathcal{O}}_{U}.
}
\end{equation}
The composition
$$
\xymatrix{
(q_{j,j+1}f)^{*}{\mathcal{E}}^{U} \ar[r] & f^{*}q_{j,j+1}^{*}{\mathcal{E}}^{U} \ar[r]^{f^{*}\phi} & f^{*}{\mathcal{O}}_{U}
}
$$
appearing on the right hand side of (\ref{eqn.original}) corresponds to the morphism
$$
\xymatrix{
V \ar[r]^{f} & U \ar[rrr]^{g} & & & {\mathbb{P}}_{(U \times X)_{U}^{2}}({\mathcal{E}}^{U})
}
$$
(\cite[4.2.8, p.75]{ega2}).  Thus, by Lemma \ref{lem.newneed}, (\ref{eqn.original}) corresponds to the top morphism in the commutative diagram
$$
\xymatrix{
V \ar[rrrr] \ar[d]_{=} & & & & {\mathbb{P}}_{(V \times X)_{V}^{2}}(\tilde{f}^{2*}{\mathcal{E}}^{U}) \ar[d] \\
V \ar[r]_{f} & U \ar[rrr]_{g} & & & {\mathbb{P}}_{(U \times X)_{U}^{2}}({\mathcal{E}}^{U})
}
$$
whose right vertical is the map (\ref{eqn.proj}), which is defined  on page \pageref{eqn.proj}.  Thus, the top route of the combined diagram sends $\mu$ to $\Phi_{U}^{\Fr}([\mathcal{M}]) \circ f$ as desired.

Since $\Phi^{1}_{{\mathcal{M}}_{j+1}}$ and $\Phi^{2}$ are isomorphisms, the bottom route of the combined diagrams exists and sends $\mu_{j,1}$ to an epimorphism
$$
\upsilon:qf_{j,j+1}^{*}{\tilde{f}}^{2*}{\mathcal{E}}^{U} \rightarrow f^{*}{\mathcal{O}}_{U}.
$$
Since the combined diagram commutes $\upsilon$ also corresponds to $r_{j}f$.  But the left hand vertical of (\ref{eqn.digge1}), on page \pageref{eqn.digge1}, sends $\mu_{j,1}$ to $\mu_{j,1}'$, the $j$th component of multiplication of the ${\mathcal{B}}^{V}$-module structure on ${\tilde{f}}^{*}\mathcal{M}$ inherited from $\mathcal{M}$ (Lemma \ref{lem.function}).  If $\zeta':G_{j}^{V}{\tilde{f}}^{*}{\mathcal{M}} \rightarrow qf_{j,j+1*}f^{*}{\mathcal{O}}_{U}$ is the isomorphism constructed in Proposition \ref{prop.simpleiso}, then we claim the epimorphism $\upsilon$ and the map constructed by applying the functor $qf_{j,j+1}^{*}$ to the map
\begin{equation} \label{eqn.lastyee}
\xymatrix{
{\tilde{f}}^{2*}{\mathcal{E}}^{U} \ar[r] & G_{j}^{V}F_{j}^{V}{\tilde{f}}^{2*}{\mathcal{E}}^{U} \ar[r]^{G_{j}^{V}\mu_{j,1}'} & G_{j}^{V}{\tilde{f}}^{*}{\mathcal{M}}_{j+1} \ar[r]^{\zeta'} & qf_{j,j+1*}f^{*}{\mathcal{O}}_{U}
}
\end{equation}
correspond to the same map $r_{j}f$ under the correspondence defined in Proposition \ref{prop.groth}.  Since the application of $qf_{j,j+1}^{*}$ to (\ref{eqn.lastyee}) corresponds to $\Phi_{V}^{\Fr}(\Gamma_{n}(f)[\mathcal{M}])$ under the correspondence of Proposition \ref{prop.groth}, this would suffice to prove the proposition.  To prove the claim, let
$$
\upsilon':\tilde{f}^{2*}{\mathcal{E}}^{U} \rightarrow qf_{j,j+1*}\tilde{f}^{*}{\mathcal{O}}_{U}
$$
be the left adjunct of $\upsilon$.  Then $\upsilon'$ equals the composition
$$
\xymatrix{
{\tilde{f}}^{2*}{\mathcal{E}}^{U} \ar[r] & G_{j}^{V}F_{j}^{V}{\tilde{f}}^{2*}{\mathcal{E}}^{U} \ar[r]^{G_{j}^{V}\mu_{j,1}'} & G_{j}^{V}{\tilde{f}}^{*}{\mathcal{M}}_{j+1} \ar[r]^{\beta} & qf_{j,j+1*}f^{*}{\mathcal{O}}_{U}
}
$$
for some isomorphism $\beta:G_{j}^{V}{\tilde{f}}^{*}{\mathcal{M}}_{j+1} \rightarrow qf_{j,j+1*}f^{*}{\mathcal{O}}_{U}$.  Thus (\ref{eqn.lastyee}) equals the composition of $\upsilon'$ with an isomorphism.  In particular, the application of $qf_{j,j+1}^{*}$ to (\ref{eqn.lastyee}) corresponds, via the map in Proposition \ref{prop.groth}, to $r_{j}f$.  We conclude that
$$
\Phi_{V}^{\Fr}(\Gamma_{n}(f)[\mathcal{M}])=\Phi_{U}^{\Fr}([\mathcal{M}]) \circ f
$$
as desired.
\end{proof}

\begin{proposition} \label{prop.refreport}
For $1 \leq n$, the transformation
$$
\Phi:\Gamma_{n}(-) \Longrightarrow \operatorname{Hom}_{S}(-,{\mathbb{P}}_{X^{2}}(\mathcal{E})^{\otimes n})
$$
defined in the proof of Proposition \ref{prop.trans} is a monomorphism of functors.  For $1 \leq n < m$, $\Phi$ is an equivalence.
\end{proposition}

\begin{proof}
Fix notation as in the previous proposition.  By Corollary \ref{cor.freefunct}, in order to prove the Proposition, it suffices to show that $\Phi$ restricted to $\Gamma_{n}^{\Fr}(-)$ is a monomorphism or equivalence to $\operatorname{Hom}_{S}^{\Fr}(-,{\mathbb{P}}_{X^{2}}(\mathcal{E})^{\otimes n})$, depending on $n$.

We first show $\Phi$ is injective.  Suppose $\mathcal{M} = \bigoplus_{j=0}^{n} (\id_{U} \times q_{i})_{*}{\mathcal{O}}_{U}$ and ${\mathcal{N}} = \bigoplus_{j=0}^{n} (\id_{U} \times r_{j})_{*}{\mathcal{O}}_{U}$ are truncated $U$-families of length $n+1$ with multiplications $\mu$ and $\nu$, such that $\Phi([\mathcal{M}])=\Phi([\mathcal{N}])$.  We show that $\mathcal{M}$ and $\mathcal{N}$ are isomorphic ${\mathcal{B}}^{U}$-modules.  We first note that the hypothesis implies ${\mathcal{M}}_{j} = {\mathcal{N}}_{j}$.  For, if $p_{\mathcal{E}}:{\mathbb{P}}_{X^{2}}(\mathcal{E}) \rightarrow X^{2}$ is the structure map, $q_{j} \times q_{j+1} = p_{\mathcal{E}} \Phi([\mathcal{M}]) = p_{\mathcal{E}} \Phi([\mathcal{N}]) = r_{j} \times r_{j+1}$.  Thus, for all $0 \leq j \leq n$, there is equality ${\mathcal{M}}_{j} = {\mathcal{N}}_{j}$.

We define a ${\mathcal{B}}^{U}$-module isomorphism $\theta:\mathcal{M} \rightarrow \mathcal{N}$ by induction.  Let $\theta_{0}=\id_{\mathcal{M}}$.  For $0 \leq j \leq n-1$, we proceed to define $\theta_{j+1}$.  Assume, for $0< j' \leq j$, there exists an isomorphism $\theta_{j'}:{\mathcal{M}}_{j'} \rightarrow {\mathcal{N}}_{j'}$ such that the diagram
$$
\xymatrix{
F_{j'-1}^{\mathcal{M}}{\mathcal{E}}^{U} \ar[r]^{\mu_{j'-1}} \ar[d]_{F_{j'-1}\theta_{j'-1}} & {\mathcal{M}}_{j'} \ar[d]^{\theta_{j'}} \\
F_{j'-1}^{\mathcal{N}}{\mathcal{E}}^{U} \ar[r]_{\nu_{j'-1}}  & {\mathcal{M}}_{j'}
}
$$
commutes.  Since $\Phi([\mathcal{M}])=\Phi([\mathcal{N}])$, there must exist an isomorphism $\tau_{j}:{\mathcal{O}}_{U} \rightarrow {\mathcal{O}}_{U}$ such that the diagram
$$
\xymatrix{
&  q_{j,j+1}^{*}G_{j}^{\mathcal{M}}{\mathcal{M}}_{j+1} \ar[r]^{\hskip .3in \cong} & {\mathcal{O}}_{U} \ar[dd]^{\tau} \\
q_{j,j+1}^{*}{\mathcal{E}}^{U} \ar[ur]^{\alpha} \ar[dr]_{\beta} & & \\
&  q_{j,j+1}^{*}G_{j}^{\mathcal{N}}{\mathcal{N}}_{j+1} \ar[r]_{\hskip .3in \cong} & {\mathcal{O}}_{U}
}
$$
commutes (Proposition \ref{prop.groth}), where the rightmost isomorphisms are the composition
$$
q_{j,j+1}^{*}G_{j}^{\mathcal{M}}{\mathcal{M}}_{j+1} \rightarrow q_{j,j+1}^{*}q_{j,j+1*}{\mathcal{O}}_{U} \rightarrow {\mathcal{O}}_{U},
$$
the left map of this composition is the isomorphism of Proposition \ref{prop.simpleiso}, and $\alpha$ and $\beta$ are epimorphisms which correspond, by Proposition \ref{prop.groth}, to $\Phi([\mathcal{M}])$ and $\Phi([\mathcal{N}])$, respectively.  Thus, there exists an isomorphism $\tau':q_{j,j+1}^{*}G_{j}^{\mathcal{M}}{\mathcal{M}}_{j+1} \rightarrow q_{j,j+1}^{*}G_{j}^{\mathcal{N}}{\mathcal{N}}_{j+1}$ such that
$$
\xymatrix{
&  q_{j,j+1}^{*}G_{j}^{\mathcal{M}}{\mathcal{M}}_{j+1} \ar[dd]^{\tau'} \\
q_{j,j+1}^{*}{\mathcal{E}}^{U} \ar[ur]^{\alpha} \ar[dr]_{\beta} &  \\
&  q_{j,j+1}^{*}G_{j}^{\mathcal{N}}{\mathcal{N}}_{j+1}
}
$$
commutes.  Applying the functor $q_{j,j+1*}$ to this diagram gives a commutative diagram
$$
\xymatrix{
& & q_{j,j+1*}q_{j,j+1}^{*}G_{j}^{\mathcal{M}}{\mathcal{M}}_{j+1} \ar[dd]^{q_{j,j+1*}\tau'} \\
{\mathcal{E}}^{U} \ar[r] & q_{j,j+1*}q_{j,j+1}^{*}{\mathcal{E}}^{U} \ar[ur]^{q_{j,j+1*}\alpha} \ar[dr]_{q_{j,j+1*}\beta} &  \\
& &  q_{j,j+1*}q_{j,j+1}^{*}G_{j}^{\mathcal{N}}{\mathcal{N}}_{j+1}
}
$$
where the left map is the unit map, $\eta$, of the pair $(q_{j,j+1}^{*},q_{j,j+1*})$.  Since $\eta_{G_{j}^{\mathcal{M}}{\mathcal{M}}_{j+1}}$ is an isomorphism by Lemma \ref{lem.unit}, we have a commutative diagram
$$
\xymatrix{
{\mathcal{E}}^{U} \ar[r] \ar[d]_{=} & G_{j}^{\mathcal{M}}{\mathcal{M}}_{j+1} \ar[d]^{\cong} \\
{\mathcal{E}}^{U} \ar[r]  & G_{j}^{\mathcal{N}}{\mathcal{N}}_{j+1}.
}
$$
Applying $F_{j}^{\mathcal{M}}$ to this diagram, gives us a commutative diagram
\begin{equation} \label{eqn.comm1}
\xymatrix{
F_{j}^{\mathcal{M}}{\mathcal{E}}^{U} \ar[r] \ar[d]_{=} & F_{j}^{\mathcal{M}}G_{j}^{\mathcal{M}}{\mathcal{M}}_{j+1} \ar[d]^{\cong}  \\
F_{j}^{\mathcal{M}}{\mathcal{E}}^{U} \ar[r] & F_{j}^{\mathcal{M}}G_{j}^{\mathcal{N}}{\mathcal{N}}_{j+1}.
}
\end{equation}
Since $\pr_{2*}(\pr_{1}^{*}\theta_{j}\otimes -)$ is a natural equivalence between $F_{j}^{\mathcal{M}}$ and $F_{j}^{\mathcal{N}}$, the diagram
\begin{equation} \label{eqn.comm2}
\xymatrix{
F_{j}^{\mathcal{M}}{\mathcal{E}}^{U} \ar[r] \ar[d]_{\pr_{2*}(\pr_{1}^{*}\theta_{j} \otimes {\mathcal{E}}^{U})} & F_{j}^{\mathcal{M}}G_{j}^{{\mathcal{N}}_{j+1}} \ar[d]^{\pr_{2*}(\pr_{1}^{*}\theta_{j} \otimes G_{j}^{\mathcal{N}})}{\mathcal{N}}_{j+1}   \\
F_{j}^{\mathcal{N}}{\mathcal{E}}_{U} \ar[r] & F_{j}^{\mathcal{N}}G_{j}^{\mathcal{N}}{\mathcal{N}}_{j+1}.
}
\end{equation}
commutes also.
Stacking (\ref{eqn.comm1}) above (\ref{eqn.comm2}) gives a commutative diagram
$$
\xymatrix{
F_{j}^{\mathcal{M}}{\mathcal{E}}^{U} \ar[r] \ar[d]_{\cong} & F_{j}^{\mathcal{M}}G_{j}^{\mathcal{M}}{\mathcal{M}}_{j+1} \ar[d]^{\cong} \ar[r] & {\mathcal{M}}_{j+1} \\
F_{j}^{\mathcal{N}}{\mathcal{E}}^{U} \ar[r] & F_{j}^{\mathcal{N}}G_{j}^{\mathcal{N}}{\mathcal{N}}_{j+1} \ar[r] & {\mathcal{N}}_{j+1}.
}
$$
where the rightmost maps are the counits of the pair $(F_{j}^{\mathcal{M}},G_{j}^{\mathcal{M}})$.  By Lemma \ref{lem.unit}, these counits are isomorphism, so there exists an isomorphism $\theta_{j+1}:{\mathcal{M}}_{j+1} \rightarrow {\mathcal{N}}_{j+1}$ such that the diagram
$$
\xymatrix{
F_{j}^{\mathcal{M}}{\mathcal{E}}^{U} \ar[r] \ar[d]_{\cong} & {\mathcal{M}}_{j+1} \ar[d]^{\theta_{j+1}} \\
F_{j}^{\mathcal{N}}{\mathcal{E}}^{U} \ar[r] & {\mathcal{N}}_{j+1}.
}
$$
commutes.  Since the top morphism is $\mu_{j}$ while the bottom morphism is $\nu_{j}$ by the definition of $\Phi$, we have succeeded, by Lemma \ref{lem.gen1}, in constructing a ${\mathcal{B}}^{U}$-module isomorphism $\theta:\mathcal{M} \rightarrow \mathcal{N}$, as desired.  Thus $\Phi$ is injective.

Now, suppose $1 \leq n < m$.  We show $\Phi$ is surjective.  Suppose $r:U \rightarrow {\mathbb{P}}_{X^{2}}(\mathcal{E})^{\otimes n}$ is a free $S$-morphism.  Then, on the $j$th tensor, $r$ projects to $(q_{j-1},q_{j})$.  Let ${\mathcal{M}}_{j} = (\id_{U} \times q_{j})_{*}{\mathcal{O}}_{U}$.  We give $\mathcal{M} = \oplus_{j=0}^{n} {\mathcal{M}}_{j}$ a ${\mathcal{B}}^{U}$-module structure such that $\Phi([\mathcal{M}])=r$.  By Proposition \ref{prop.groth}, there exists, for $0 \leq j < n$, epimorphisms
$$
q_{j,j+1}^{*}{\mathcal{E}}^{U} \rightarrow {\mathcal{O}}_{U}.
$$
By the adjointness of the pair $(q_{j,j+1}^{*},q_{j,j+1*})$, this epimorphism corresponds to a map
\begin{equation} \label{eqn.next}
{\mathcal{E}}^{U} \rightarrow q_{j,j+1*}{\mathcal{O}}_{U}
\end{equation}
which is an epimorphism since $q_{j,j+1}$ is a closed immersion.  By Proposition \ref{prop.simpleiso}, there is an isomorphism
$$
{\mathcal{H}}\it{om}_{{\mathcal{O}}_{(U\times X)_{U}^{X}}}(\pr_{1}^{*}{\mathcal{M}}_{q_{j}},i_{*}(\pr_{2}i)^{!}{\mathcal{M}}_{j+1}) \rightarrow q_{j,j+1*}{\mathcal{O}}_{U}
$$
so that (\ref{eqn.next}) corresponds to an epimorphism
$$
{\mathcal{E}}^{U} \rightarrow {\mathcal{H}}\it{om}_{{\mathcal{O}}_{(U\times X)_{U}^{X}}}(\pr_{1}^{*}{\mathcal{M}}_{q_{j}},i_{*}(\pr_{2}i)^{!}{\mathcal{M}}_{j+1}).
$$
By Lemma $\ref{lem.unit}$, this epimorphism corresponds to an epimorphism
$$
\mu_{j}:F_{j}^{\mathcal{M}}{\mathcal{E}}^{U} \rightarrow {\mathcal{M}}_{j+1}.
$$
Since $n < m$, the multiplication maps $\nu_{l,l'}:{\mathcal{B}}_{l} \otimes {\mathcal{B}}_{l'} \rightarrow {\mathcal{B}}_{l+l'}$ are isomorphisms for $l+l' \leq n$.  Thus, by Lemma \ref{lem.smalldef}, there is a unique ${\mathcal{B}}^{U}$-module multiplication, $\mu^{\mathcal{M}}$, for $\mathcal{M}$ such that $(\mu^{\mathcal{M}})_{i,1}=\mu_{i}$.  One may now easily check that $\Phi([\mathcal{M}])=r$.
\end{proof}\index{family!free, of truncated point modules|)}

\chapter{The Bimodule Segre Embedding}\index{bimodule Segre
embedding|(}

Let $X$, $Y$ and $Z$ be separated, noetherian schemes, let $\mathcal{E}$ be a coherent ${\mathcal{O}}_{X \times Y}$-bimodule and let $\mathcal{F}$ be a coherent ${\mathcal{O}}_{Y \times Z}$-bimodule.  Suppose that $\pr_{i}:X \times Y \rightarrow X,Y$, $\pr_{i}:Y \times Z \rightarrow Y, Z$ are projections (same notation!), and $p_{\mathcal{E}}$ and $p_{\mathcal{F}}$ are the structure maps of ${\mathbb{P}}_{X \times Y}(\mathcal{E})$ and ${\mathbb{P}}_{Y \times Z}(\mathcal{F})$ respectively.  We denote by
$$
{\mathbb{P}}_{X \times Y}(\mathcal{E}) \otimes_{Y} {\mathbb{P}}_{Y \times Z}(\mathcal{F})
$$
the fiber product in the following diagram:
\begin{equation} \label{eqn.thoughtno}
\xymatrix
{
{\mathbb{P}}_{X \times Y}(\mathcal{E}) \otimes_{Y} {\mathbb{P}}_{Y \times Z}(\mathcal{F}) \ar@{-->}[rr]^{q_{\mathcal{F}}} \ar@{-->}[d]_{q_{\mathcal{E}}} & & {\mathbb{P}}_{Y \times Z}(\mathcal{F}) \ar[d]^{\pr_{1} \circ p_{\mathcal{F}}} \\
{\mathbb{P}}_{X \times Y}(\mathcal{E}) \ar[rr]_{\pr_{2} \circ p_{\mathcal{E}}} & & X.
}
\end{equation}
In this chapter, we construct a closed immersion
$$
s:{\mathbb{P}}_{X\times Y}(\mathcal{E}) \otimes_{Y} {\mathbb{P}}_{Y \times Z}(\mathcal{F}) \rightarrow {\mathbb{P}}_{X \times Z}({\mathcal{E}} \otimes_{{\mathcal{O}}_{Y}} {\mathcal{F}}),
$$
the {\it bimodule Segre embedding}\index{bimodule Segre embedding}, and prove it is natural in a suitable sense.  If $\mathcal{I} \subset T(\mathcal{E})$ is an ideal which has first nonzero component in degree $m > 0$, and $\mathcal{B} = T(\mathcal{E})/\mathcal{I}$, then we will show that, for $n \geq m$, $\Gamma_{n}$ is represented by the pullback of the diagram
$$
\xymatrix{
& {{\mathbb{P}}_{X^{2}}(\mathcal{E})}^{\otimes n} \ar[d]^{s} \\
{\mathbb{P}}_{X^{2}}({\mathcal{E}}^{n}/{\mathcal{I}}_{n}) \ar[r]  & {\mathbb{P}}_{X^{2}}({\mathcal{E}}^{n}).
}
$$
(Theorem \ref{theorem.bigone}).  After stating our main theorem (Theorem \ref{theorem.bigsegre}) which describes the properties of $s$ we need in order to prove $\Gamma_{n}$ is representable for large $n$, we prove these properties one by one, making consistent use of the algebraic description of $s$ due to Proposition \ref{prop.groth}.

We describe notation we use in this chapter.  Suppose $W$ is a scheme, $d:W \rightarrow W \times W$ is the diagonal morphism, $q_{1}:W \rightarrow X$, $q_{2}:W \rightarrow Y$ and $q_{3}:W \rightarrow Z$ are morphisms, $q_{12} = (q_{1} \times q_{2})\circ d:W \rightarrow X \times Y$, $q_{23} = (q_{2} \times q_{3}) \circ d:W \rightarrow Y \times Z$ and $q_{123} = (q_{1} \times q_{2} \times q_{3})\circ ({\operatorname{id}}_{W}\times d)\circ d:W \rightarrow X \times Y \times Z$.  Finally, suppose $q_{12}^{W} = (\id_{W} \times q_{1} \times \id_{W} \times q_{2})\circ (d \times d)\circ d: W \rightarrow (W \times X) \times_{W} (W \times Y)$.  Let $q_{23}^{W}$ and $q_{123}^{W}$ be defined similarly.

\section{Statement of the main theorem}
Before we can state the main result of this chapter, we must give a number of technical definitions.  Let $U$ be a scheme, Let $X'$, $Y'$ and $Z'$ be $U$-schemes, and let $\pr_{ij}:X \times Y \times Z \rightarrow X,Y,Z$ and $\pr_{ij}':X'\times_{U} Y' \times_{U} Z' \rightarrow X',Y',Z'$ be the standard projection maps.

\begin{definition}
With the above notation, maps $p_{xy}:X' \times_{U} Y' \rightarrow X \times Y$, $p_{yz}:Y' \times_{U} Z' \rightarrow Y \times Z$, $p_{xz}:X' \times_{U} Z' \rightarrow X \times Z$ and $p_{xyz}:X' \times_{U} Y' \times_{U} Z' \rightarrow X \times Y \times Z$ are said to be {\bf overlap compatible}\index{overlap compatibility|textbf} if the diagrams
\begin{equation} \label{eqn.firsty}
\xymatrix{
X' \times_{U} Y' \times_{U} Z' \ar[d]_{\pr_{12}'} \ar[r]^{p_{xyz}} & X \times Y \times Z \ar[d]^{\pr_{12}} &  X' \times_{U} Y' \times_{U} Z' \ar[r]^{p_{xyz}} \ar[d]_{\pr_{23}'} & X \times Y \times Z \ar[d]^{\pr_{23}} \\
X' \times_{U} Y' \ar[r]_{p_{xy}} & X \times Y &  Y' \times_{U} Z' \ar[r]_{p_{yz}} & Y \times Z
}
\end{equation}
commute, and
\begin{equation} \label{eqn.second}
\xymatrix{
X' \times_{U} Y' \times_{U} Z' \ar[d]_{\pr_{13}'} \ar[r]^{p_{xyz}} & X \times Y \times Z \ar[d]^{\pr_{13}} \\
X' \times_{U} Z' \ar[r]_{p_{xz}} & X \times Z
}
\end{equation}
is a pullback.
\end{definition}
The following Lemma is similar to Lemma \ref{lem.isommonoid}.

\begin{lemma} \label{lem.genisommonoid}
Suppose $\mathcal{E}$ is an ${\mathcal{O}}_{X \times Y}$-module and $\mathcal{F}$ is an ${\mathcal{O}}_{Y \times Z}$-module such that $\mathcal{E}$ and $\mathcal{F}$ have the affine direct image property\index{affine direct image property}.  Suppose $p_{xy}$, $p_{yz}$, $p_{xz}$ and $p_{xyz}$ are overlap compatible.  Then $p_{xy}^{*}\mathcal{E}$ and $p_{yz}^{*}\mathcal{F}$ have the affine direct image property \index{affine direct image property}and there is a natural isomorphism
\begin{equation} \label{eqn.itoldyou}
p_{xy}^{*}\mathcal{E} \otimes_{{\mathcal{O}}_{Y'}} p_{yz}^{*}\mathcal{F} \rightarrow p_{xz}^{*}(\mathcal{E} \otimes_{{\mathcal{O}}_{Y}} \mathcal{F}).
\end{equation}
\end{lemma}

\begin{proof}
Since the diagrams (\ref{eqn.firsty}) commute, $\pr^{'*}_{12}p_{xy}^{*} \cong p_{xyz}^{*}p_{12}^{*}$ and $\pr_{23}^{'*}p_{yz}^{*} \cong p_{xyz}^{*}\pr_{23}^{*}$.  Using these isomorphisms, we have
\begin{align*}
p_{xy}^{*}\mathcal{E} \otimes_{{\mathcal{O}}_{Y'}} p_{yz}^{*}\mathcal{F} & = \pr_{13*}^{'}(\pr_{12}^{'*}p_{xy}^{*}\mathcal{E}\otimes \pr_{23}^{'*}p_{yz}^{*}\mathcal{F}) \\
& \cong \pr_{13*}^{'}(p_{xyz}^{*}\pr_{12}^{*}\mathcal{E}\otimes p_{xyz}^{*}\pr_{23}^{*}\mathcal{F}) \\
& \cong \pr_{13*}^{'}p_{xyz}^{*}(\pr_{12}^{*}\mathcal{E}\otimes \pr_{23}^{*}\mathcal{F}) \\
& \cong p_{xz}^{*}\pr_{13*}(\pr_{12}^{*}\mathcal{E} \otimes \pr_{23}^{*}\mathcal{F}) \\
& = p_{xz}^{*}(\mathcal{E} \otimes_{{\mathcal{O}}_{Y}} \mathcal{F}),
\end{align*}
where the last nontrivial isomorphism is an application of Proposition \ref{prop.canon}, which may be invoked since $\mathcal{E}$ and $\mathcal{F}$ have the affine direct image property and since (\ref{eqn.second}) is a pullback.

To prove the first assertion, we must show that
$$
\pr_{12}^{'*}p_{xy}^{*}\mathcal{E}\otimes \pr_{23}^{'*}p_{yz}^{*}\mathcal{F}
$$
is rla with respect to $\pr_{13}'$.  But
$$
\pr_{12}^{'*}p_{xy}^{*}\mathcal{E}\otimes \pr_{23}^{'*}p_{yz}^{*}\mathcal{F} \cong p_{xyz}^{*}(\pr_{12}^{*}\mathcal{E}\otimes \pr_{23}^{*}\mathcal{F})
$$
as above, and the right hand side is rla with respect to $\pr_{13}'$ by Proposition \ref{prop.canon} since (\ref{eqn.second}) is a pullback.
\end{proof}

\begin{example} \label{example.refrep}
Suppose $U$ be an affine noetherian scheme, $X'=U \times X$, $Y'=U \times Y$, $Z'=U \times Z$ and the maps $p_{xyz}$, $p_{xy}$, $p_{xz}$ and $p_{yz}$ are projections.  Then $p_{xyz}$, $p_{xy}$, $p_{xz}$ and $p_{yz}$ are overlap compatible and the isomorphism (\ref{eqn.itoldyou}) gives the tensor product of bimodules its indexed structure (Lemma \ref{lem.isommonoid}).
\end{example}

Let
\begin{equation} \label{eqn.alphaone}
\Omega_{1}:\pr_{12}^{*}q_{12*} = \pr_{12}^{*}\pr_{12*}q_{123*} \Longrightarrow q_{123*}
\end{equation}
and
\begin{equation} \label{eqn.alphatwo}
\Omega_{2}:\pr_{23}^{*}q_{23*} = \pr_{23}^{*}\pr_{23*}q_{123*} \Longrightarrow q_{123*}
\end{equation}
be natural transformations whose last composite is a counit map.  Let
\begin{equation} \label{eqn.betaone}
\Upsilon:q_{123*}- \otimes_{{\mathcal{O}}_{X \times Y \times Z}} q_{123*}- \Longrightarrow q_{123*}q_{123}^{*}(q_{123*}- \otimes_{{\mathcal{O}}_{X \times Y \times Z}} q_{123*}-) \Longrightarrow
\end{equation}
$$
q_{123*}(q_{123}^{*}q_{123*}- \otimes_{{\mathcal{O}}_{W}} q_{123}^{*}q_{123*}-) \Longrightarrow q_{123*}(- \otimes_{{\mathcal{O}}_{W}} -)
$$
be the natural transformation whose first composite is a unit and whose last composite is a counit.

\begin{definition}
Suppose $\mathcal{E}$ and $\mathcal{F}$ have the affine direct image property, and let
$$
f:{\mathbb{P}}_{X\times Y}(\mathcal{E}) \otimes_{Y} {\mathbb{P}}_{Y \times Z}(\mathcal{F}) \rightarrow {\mathbb{P}}_{X \times Z}({\mathcal{E}} \otimes_{{\mathcal{O}}_{Y}} {\mathcal{F}})
$$
be a morphism.
\begin{itemize}
\item{}
$f$ is {\bf functorial}\index{functoriality|textbf} if whenever $\mathcal{E}'$ is an ${\mathcal{O}}_{X \times Y}$-module such that there is an epimorphism $\upsilon:\mathcal{E} \rightarrow \mathcal{E}'$, then $\mathcal{E}'$ and $\mathcal{F}$ have the affine direct image property\index{affine direct image property}, the map
$$
\upsilon \otimes_{{\mathcal{O}}_{Y}} \mathcal{F}:\mathcal{E} \otimes_{{\mathcal{O}}_{Y}} \mathcal{F} \rightarrow \mathcal{E}' \otimes_{{\mathcal{O}}_{Y}} \mathcal{F}
$$
is an epimorphism, and the diagram
\begin{equation}
\xymatrix{
{\mathbb{P}}_{X \times Y}(\mathcal{E}') \otimes_{Y} {\mathbb{P}}_{Y \times Z}(\mathcal{F}) \ar[rrr]^{{\mathbb{P}}_{X \times Y}(\upsilon) \times \id} \ar[d]_{f} & & & {\mathbb{P}}_{X \times Y}(\mathcal{E}) \otimes_{Y} {\mathbb{P}}_{Y \times Z}(\mathcal{F}) \ar[d]^{f} \\
{\mathbb{P}}_{X \times Z}(\mathcal{E}' \otimes_{{\mathcal{O}}_{Y}} \mathcal{F}) \ar[rrr]_{{\mathbb{P}}_{X \times Z}(\upsilon \otimes_{{\mathcal{O}}_{Y}} \mathcal{F})}  & & & {\mathbb{P}}_{X \times Z}(\mathcal{E} \otimes_{{\mathcal{O}}_{Y}} \mathcal{F})
}
\end{equation}
commutes.

\item{}
$f$ is {\bf compatible with base change}\index{compatibility with base change|textbf} if whenever $p_{xy}$, $p_{yz}$, $p_{xz}$ and $p_{xyz}$ are overlap compatible, the diagram
$$
\xymatrix{
{\mathbb{P}}_{X' \times_{U} Y'}(p_{xy}^{*}\mathcal{E}) \otimes_{Y'} {\mathbb{P}}_{Y' \times_{U} Z'}(p_{yz}^{*}\mathcal{F}) \ar[dd] \ar[r]^{\hskip .3in f} & {\mathbb{P}}_{X' \times_{U} Z'}(p_{xy}^{*}\mathcal{E} \otimes_{{\mathcal{O}}_{Y'}} p_{yz}^{*}\mathcal{F}) \ar[d] \\
& {\mathbb{P}}_{X' \times_{U} Z'}(p_{xz}^{*}(\mathcal{E} \otimes_{{\mathcal{O}}_{Y}} \mathcal{F})) \ar[d] \\
{\mathbb{P}}_{X \times Y}(\mathcal{E}) \otimes_{Y} {\mathbb{P}}_{Y \times Z}(\mathcal{F})  \ar[r]_{f} & {\mathbb{P}}_{X \times Z}(\mathcal{E} \otimes_{{\mathcal{O}}_{Y}} \mathcal{F})
}
$$
whose left and bottom right vertical arrows are defined in Lemma \ref{lem.proj}, and whose top right vertical arrow is induced by the isomorphism defined in Lemma \ref{lem.genisommonoid}, commutes.

\item{}
$f$ is {\bf associative}\index{associativity|textbf} if whenever $T$ is a scheme and $\mathcal{G}$ is an ${\mathcal{O}}_{Z \times T}$-bimodule, the diagram
$$
\xymatrix{
({\mathbb{P}}_{X \times Y}(\mathcal{E}) \otimes_{Y} {\mathbb{P}}_{Y \times Z}(\mathcal{F})) \otimes_{Z} {\mathbb{P}}_{Z \times T}(\mathcal{G}) \ar[r]^{\cong} \ar[d]_{f \times \id} & {\mathbb{P}}_{X \times Y}(\mathcal{E}) \otimes_{Y} ({\mathbb{P}}_{Y \times Z}(\mathcal{F}) \otimes_{Z} {\mathbb{P}}_{Z \times T}(\mathcal{G})) \ar[d]^{\id \times f} \\
{\mathbb{P}}_{X \times Z}(\mathcal{E}\otimes_{{\mathcal{O}}_{Y}} \mathcal{F}) \otimes_{Z} {\mathbb{P}}_{Z \times T}(\mathcal{G})  \ar[d]_{f} & {\mathbb{P}}_{X \times Y}(\mathcal{E}) \otimes_{Y} {\mathbb{P}}_{Y \times Z'}(\mathcal{F} \otimes_{{\mathcal{O}}_{Z}} \mathcal{G})  \ar[d]^{f} \\
{\mathbb{P}}_{X \times T}((\mathcal{E}\otimes_{{\mathcal{O}}_{Y}} \mathcal{F})\otimes_{{\mathcal{O}}_{Z}} \mathcal{G}) \ar[r]_{\cong} & {\mathbb{P}}_{X \times T}(\mathcal{E} \otimes_{{\mathcal{O}}_{Y}}(\mathcal{F} \otimes_{{\mathcal{O}}_{Z}} \mathcal{G}))
}
$$
whose bottom row is induced by the associativity isomorphism (\ref{prop.tensor}), commutes.
\end{itemize}
\end{definition}

\begin{theorem} \label{theorem.bigsegre}
Retain the above notation, let $\mathcal{E}$ be an ${{\mathcal{O}}_{X \times Y}}$-module and let $\mathcal{F}$ be an ${{\mathcal{O}}_{Y \times Z}}$-module such that $\mathcal{E}$ and $\mathcal{F}$ have the affine direct image property.  Then there exists a unique map
$$
s:{\mathbb{P}}_{X\times Y}(\mathcal{E}) \otimes_{Y} {\mathbb{P}}_{Y \times Z}(\mathcal{F}) \rightarrow {\mathbb{P}}_{X \times Z}({\mathcal{E}} \otimes_{{\mathcal{O}}_{Y}} {\mathcal{F}})
$$
such that if $r:W \rightarrow {\mathbb{P}}_{X\times Y}(\mathcal{E}) \otimes_{Y} {\mathbb{P}}_{Y \times Z}(\mathcal{F})$ is a morphism whose projection to $X \times Y \times Z$ is $q_{123}$ and $r$ corresponds, via Proposition \ref{prop.groth}, to the pair
$$
\psi_{1}:\mathcal{E} \rightarrow q_{12*}{\mathcal{L}}_{1}
$$
and
$$
\psi_{2}:\mathcal{F} \rightarrow q_{23*}{\mathcal{L}}_{2}
$$
then $s \circ r$ corresponds to the map
\begin{equation} \label{eqn.bigsegre}
\xymatrix{
\mathcal{E} \otimes_{{\mathcal{O}}_{Y}} \mathcal{F} \ar[rr]^{\hskip -.3in \psi_{1} \otimes_{{\mathcal{O}}_{Y}} \psi_{2}} & & q_{12*}{\mathcal{L}}_{1} \otimes_{{\mathcal{O}}_{Y}} q_{23*}{\mathcal{L}}_{2} \ar[rr]^{\hskip .5in \pr_{13*}(\Omega_{1} \otimes \Omega_{2})} & &
}
\end{equation}
$$
\xymatrix{
\pr_{13*}(q_{123*}{\mathcal{L}}_{1} \otimes q_{123*}{\mathcal{L}}_{2}) \ar[rr]^{\hskip .4in \pr_{13*}\Upsilon} & & q_{13*}({\mathcal{L}}_{1} \otimes {\mathcal{L}}_{2}).
}
$$
Furthermore, $s$ is a closed immersion\index{bimodule Segre embedding!is a closed immersion} which is functorial\index{bimodule Segre embedding!is functorial}, compatible with base change\index{bimodule Segre embedding!is compatible with base change} and associative\index{bimodule Segre embedding!is associative}.
\end{theorem}

We call $s$ the {\bf bimodule Segre embedding}\index{bimodule Segre embedding|textbf}.  We construct $s$ as the composition of three maps.  We first construct an isomorphism
$$
s_{1}:{\mathbb{P}}_{X \times Y}(\mathcal{E}) \otimes_{Y} {\mathbb{P}}_{Y \times Z}(\mathcal{F}) \rightarrow {\mathbb{P}}_{X\times Y \times Z}({\pr_{12}}^{*}\mathcal{E}) \times_{X\times Y \times Z}{\mathbb{P}}_{X\times Y \times Z}({\pr_{23}}^{*}\mathcal{F})
$$
(Proposition \ref{prop.firstmap}).  The classical Segre embedding\index{Segre embedding}, whose construction we review (Theorem \ref{prop.segre}), gives us a closed immersion
$$
s_{2}:{\mathbb{P}}_{X\times Y \times Z}({\pr_{12}}^{*}\mathcal{E}) \times_{X\times Y \times Z}{\mathbb{P}}_{X\times Y \times Z}({\pr_{23}}^{*}\mathcal{F}) \rightarrow {\mathbb{P}}_{X \times Y \times Z}({\pr_{12}}^{*}\mathcal{E} \otimes {\pr_{23}}^{*}\mathcal{F})
$$
Finally, since $\mathcal{E}$ and $\mathcal{F}$ have the affine direct image property \index{affine direct image property}by hypothesis, we may construct a closed immersion
$$
s_{3}:{\mathbb{P}}_{X \times Y \times Z}({\pr_{12}}^{*}\mathcal{E} \otimes {\pr_{23}}^{*}\mathcal{F}) \rightarrow  {\mathbb{P}}_{X \times Z}({\pr_{13}}_{*}({\pr_{12}}^{*}\mathcal{E} \otimes {\pr_{23}}^{*}\mathcal{F}))
$$
(Proposition \ref{prop.thirdmap}).

\section{Construction of the bimodule Segre embedding}
We remind the reader of relevant notation.  Let $\mathcal{E}$ be a coherent ${\mathcal{O}}_{X \times Y}$-module and let $\mathcal{F}$ be a coherent ${\mathcal{O}}_{Y \times Z}$-module.  Unless otherwise stated, we assume $\mathcal{E}$ and $\mathcal{F}$ have the affine direct image property.  We denote by
$$
{\mathbb{P}}_{X \times Y}(\mathcal{E}) \otimes_{Y} {\mathbb{P}}_{Y \times Z}(\mathcal{F})
$$
the fiber product in (\ref{eqn.thoughtno}).  By the universal property of the fiber product, there is a map $t:{\mathbb{P}}_{X \times Y}(\mathcal{E}) \otimes_{Y} {\mathbb{P}}_{Y \times Z}(\mathcal{F}) \rightarrow X \times Y \times Z$, making
$$
\xymatrix
{
& {\mathbb{P}}_{X \times Y}(\mathcal{E}) \otimes_{Y} {\mathbb{P}}_{Y \times Z}(\mathcal{F}) \ar[dr]^{q_{\mathcal{F}}} \ar[dl]_{q_{\mathcal{E}}} \ar@{-->}[dd]^{t} & \\
{\mathbb{P}}_{X \times Y}(\mathcal{E}) \ar[dd]^{p_{\mathcal{E}}} & & {\mathbb{P}}_{Y \times Z}(\mathcal{F}) \ar[dd]_{p_{\mathcal{F}}} \\
& X \times Y \times Z \ar[dl]^{\pr_{12}} \ar[dr]_{\pr_{23}} & \\
X \times Y \ar[dr]_{\pr_{2}} & & Y \times Z \ar[dl]^{\pr_{1}} \\
& Y &.
}
$$
commute.  We now proceed to construct and give algebraic descriptions of $s_{1}$, $s_{2}$ and $s_{3}$.
\subsection{$s_{1}:{\mathbb{P}}_{X \times Y}(\mathcal{E}) \otimes_{Y} {\mathbb{P}}_{Y \times Z}(\mathcal{F}) \rightarrow {\mathbb{P}}_{X\times Y \times Z}({\pr_{12}}^{*}\mathcal{E}) \times_{X\times Y \times Z}{\mathbb{P}}_{X\times Y \times Z}({\pr_{23}}^{*}\mathcal{F})$}
The proof of the following result is a tedious application of the universal property of the pullback, so we omit it.
\begin{proposition} \label{prop.firstmap}
The pair $(q_{\mathcal{E}},t)$ induces a map

$$
f_{\mathcal{E}}: {\mathbb{P}}_{X \times Y}(\mathcal{E}) \otimes_{Y}{\mathbb{P}}_{Y \times Z}(\mathcal{F}) \rightarrow  {\mathbb{P}}_{X \times Y}(\mathcal{E}) \times_{X \times Y} X \times Y \times Z
$$
and the pair $(q_{\mathcal{F}},t)$ induces a map
$$
f_{\mathcal{F}}: {\mathbb{P}}_{X \times Y}(\mathcal{E}) \otimes_{Y}{\mathbb{P}}_{Y \times Z}(\mathcal{F}) \rightarrow  {\mathbb{P}}_{Y \times Z}(\mathcal{F}) \times_{Y \times Z} X \times Y \times Z.
$$
The pair $(f_{\mathcal{E}},f_{\mathcal{F}})$ induces an isomorphism
$$
{\mathbb{P}}_{X \times Y}(\mathcal{E}) \otimes_{Y}{\mathbb{P}}_{Y \times Z}(\mathcal{F}) \rightarrow
$$
$$
[{\mathbb{P}}_{X \times Y}(\mathcal{E}) \times_{X \times Y} (X \times Y \times Z)] \times_{X \times Y \times Z}  [{\mathbb{P}}_{Y \times Z}(\mathcal{F}) \times_{Y \times Z} (X \times Y \times Z)]
$$
In particular, there is an isomorphism
\begin{equation} \label{eqn.segfirst}
s_{1}:{\mathbb{P}}_{X \times Y}(\mathcal{E}) \otimes_{Y}{\mathbb{P}}_{Y \times Z}(\mathcal{F}) \rightarrow {\mathbb{P}}_{X \times Y \times Z}({\pr_{12}}^{*}\mathcal{E}) \times_{X \times Y \times Z} {\mathbb{P}}_{X \times Y \times Z}({\pr_{23}}^{*}\mathcal{F}).
\end{equation}
Furthermore, the inverse of $s_{1}$ is induced by the maps
$$
{\mathbb{P}}_{X \times Y \times Z}({\pr_{12}}^{*}\mathcal{E}) \times_{X \times Y \times Z} {\mathbb{P}}_{X \times Y \times Z}({\pr_{23}}^{*}\mathcal{F}) \rightarrow {\mathbb{P}}_{X \times Y \times Z}({\pr_{12}}^{*}\mathcal{E}) \rightarrow {\mathbb{P}}_{X \times Y}(\mathcal{E})
$$
and
$$
{\mathbb{P}}_{X \times Y \times Z}({\pr_{12}}^{*}\mathcal{E}) \times_{X \times Y \times Z} {\mathbb{P}}_{X \times Y \times Z}({\pr_{23}}^{*}\mathcal{F}) \rightarrow {\mathbb{P}}_{X \times Y \times Z}({\pr_{23}}^{*}\mathcal{F}) \rightarrow {\mathbb{P}}_{Y \times Z}(\mathcal{F})
$$
where the second composite of each map is the map constructed in Lemma \ref{lem.proj}.
\end{proposition}
We now describe the algebraic properties of this isomorphism.  We first give an alternate description of the maps $\Omega_{1}$ (\ref{eqn.alphaone}) and $\Omega_{2}$ ((\ref{eqn.alphatwo}) on page \pageref{eqn.alphatwo}).  By Lemma \ref{prop.duality}, the isomorphism
$$
\Pi_{1}:q_{123}^{*}\pr_{12}^{*} \Longrightarrow q_{12}^{*}
$$
is the dual 2-cell to the square
\begin{equation} \label{eqn.square1}
\xymatrix{
{\sf Qcoh} W \ar@<1ex>[r]^{q_{12*}} \ar[d]_{\id} \ar@{=>}[dr]|{\Delta_{1}}\hole & {\sf{Qcoh}} X \times Y \ar[d]^{\pr_{12}^{*}} \ar@<1ex>[l]^{q_{12}^{*}} \\
{\sf Qcoh} W \ar@<1ex>[r]^{q_{123*}}  & {\sf{Qcoh}} X \times Y \times Z  \ar@<1ex>[l]^{q_{123}^{*}}
}
\end{equation}
where $\Delta_{1}$ is the composition
\begin{equation} \label{eqn.alpha1}
\xymatrix{
\pr_{12}^{*}q_{12*} \ar@{=>}[r] & q_{123*}q_{123}^{*}\pr_{12}^{*}q_{12*} \ar@{=>}[rr]^{q_{123*}*\Pi_{1}*q_{12*}} & & q_{123*}q_{12}^{*}q_{12*} \ar@{=>}[r] & q_{123*}.
}
\end{equation}
In a similar fashion, we have morphisms
$$
\Pi_{2}:q_{123}^{*}\pr_{23}^{*} \Longrightarrow q_{23}^{*}
$$
and
\begin{equation} \label{eqn.alpha2}
\Delta_{2}:\pr_{23}^{*}q_{23*} \Longrightarrow q_{123*}.
\end{equation}

\begin{lemma}
Retain the notation above and let $\epsilon$ be the counit of the adjoint pair $(\pr_{12}^{*},\pr_{12*})$.  The map
$$
\Delta_{1}:\pr_{12}^{*}q_{12*}  \Longrightarrow q_{123*}
$$
(\ref{eqn.alpha1}) equals the composition
$$
\xymatrix{
\pr_{12}^{*}q_{12*} \ar[r]^{\hskip -.2in =} &  \pr_{12}^{*}\pr_{12*}q_{123*} \ar@{=>}[r]^{\hskip .3in \epsilon * q_{123*}} & q_{123*}
}
$$
(which is just $\Omega_{1}$, (\ref{eqn.alphaone}) on page \pageref{eqn.alphaone}) and a similar description holds for $\Delta_{2}$ (\ref{eqn.alpha2}), i.e. $\Delta_{2}=\Omega_{2}$ ((\ref{eqn.alphatwo}) on page \pageref{eqn.alphatwo}).
\end{lemma}

\begin{proof}
By definition, $\Delta_{1}$ is the composition
$$
\xymatrix{
\pr_{12}^{*}q_{12*} \ar@{=>}[r] & q_{123*}q_{123}^{*}\pr_{12}^{*}q_{12*} \ar@{=>}[rr]^{\hskip .2in q_{123*}*\Pi_{1}*q_{12*}} & & q_{123*}q_{12}^{*}q_{12*} \ar@{=>}[r] & q_{123*}.
}
$$
The assertion follows from the commutativity of the following diagram, whose unlabeled maps are various units and counits:
$$
\xymatrix{
\pr_{12}^{*}q_{12*} \ar@{=>}[r] \ar@{=>}[d]_{=} & q_{123*}q_{123}^{*}\pr_{12}^{*}q_{12*} \ar@{=>}[rr]^{q_{123*}*\Pi_{1}*q_{12*}} \ar@{=>}[d]^{=} & & q_{123*}q_{12}^{*}q_{12*} \ar@{=>}[dd] \\
\pr_{12}^{*}\pr_{12*}q_{123*} \ar@{=>}[r] \ar@{=>}[d] & q_{123*}q_{123}^{*}\pr_{12}^{*}\pr_{12*}q_{123*}  \ar@{=>}[d] & & \\
q_{123*} \ar@{=>}[r] & q_{123*}q_{123}^{*}q_{123*} \ar@{=>}[rr] & & q_{123*}
}
$$
\end{proof}

\begin{proposition} \label{prop.algfirst}
Suppose
$$
r:W \rightarrow {\mathbb{P}}_{X \times Y}(\mathcal{E}) \otimes_{Y} {\mathbb{P}}_{Y \times Z}(\mathcal{F})
$$
corresponds, via the bijection in Proposition \ref{prop.groth}, to the pair of epis
$$
\phi_{1}:q_{12}^{*}\mathcal{E} \rightarrow {\mathcal{L}}_{1}
$$
and
$$
\phi_{2}:q_{23}^{*}\mathcal{F} \rightarrow {\mathcal{L}}_{2},
$$
with right adjuncts $\psi_{1}$ and $\psi_{2}$.  Then the composition $s_{1} \circ r$ corresponds to the pair of epimorphisms
\begin{equation} \label{eqn.pair1}
\xymatrix{
q_{123}^{*}\pr_{12}^{*}\mathcal{E} \ar[r]^{\hskip .2in \Pi_{1\mathcal{E}}} & q_{12}^{*}\mathcal{E} \ar[r]^{\phi_{1}} & {\mathcal{L}}_{1}
}
\end{equation}

and
\begin{equation} \label{eqn.pair2}
\xymatrix{
q_{123}^{*}{\pr_{23}}^{*}\mathcal{F} \ar[r]^{\hskip .2in \Pi_{2\mathcal{F}}} & q_{23}^{*}\mathcal{F} \ar[r]^{\phi_{2}} & {\mathcal{L}}_{2}
}
\end{equation}
which have right adjuncts
$$
\xymatrix{
\pr_{12}^{*}\mathcal{E} \ar[rr]^{\pr_{12}^{*}\psi_{1}} & & \pr_{12}^{*}q_{12*}{\mathcal{L}}_{1} \ar[r]^{\Omega_{1{\mathcal{L}}_{1}}} & q_{123*}{\mathcal{L}}_{1}
}
$$
and
$$
\xymatrix{
\pr_{23}^{*}\mathcal{F} \ar[rr]^{\pr_{23}^{*}\psi_{2}} & & \pr_{23}^{*}q_{23*}{\mathcal{L}}_{2} \ar[r]^{\Omega_{2{\mathcal{L}}_{2}}} & q_{123*}{\mathcal{L}}_{2}
}
$$
respectively.
\end{proposition}

\begin{proof}
The pair of epis (\ref{eqn.pair1}) and (\ref{eqn.pair2}) correspond to a map
$$
r':W \rightarrow  {\mathbb{P}}_{X \times Y \times Z}({\pr_{12}}^{*}\mathcal{E}) \times_{X \times Y \times Z} {\mathbb{P}}_{X \times Y \times Z}({\pr_{23}}^{*}\mathcal{F})
$$
such that $r'$ composed with $s_{1}^{-1}$ corresponds to the pair $\phi_{1}$ and $\phi_{2}$ (Lemma \ref{lem.algproj}).  Thus, $s_{1}^{-1} \circ r'=r$ so that $r'=s_{1} \circ r$.  The first assertion follows.

Applying Corollary \ref{cor.adcom2} to (\ref{eqn.square1}) gives a commutative diagram
$$
\xymatrix{
\operatorname{Hom}_{W}(q_{12}^{*}\mathcal{E},{\mathcal{L}}_{1}) \ar[r] \ar[dd]_{- \circ \Pi_{1\mathcal{E}}} & \operatorname{Hom}_{X \times Y}(\mathcal{E},q_{12*}{\mathcal{L}}_{1}) \ar[d]^{{\pr_{12}}^{*}} \\
 & \operatorname{Hom}_{X \times Y \times Z}({\pr_{12}}^{*}\mathcal{E},{\pr_{12}}^{*}q_{12*}{\mathcal{L}}_{1}) \ar[d]^{\Omega_{1{\mathcal{L}}_{1}} \circ -} \\
\operatorname{Hom}_{W}(q_{123}^{*}{\pr_{12}}^{*}\mathcal{E},{\mathcal{L}}_{1}) \ar[r] & \operatorname{Hom}_{X \times Y \times Z}({\pr_{12}}^{*}\mathcal{E},q_{123*}{\mathcal{L}}_{1})
}
$$
whence the second assertion.
\end{proof}

\subsection{The classical Segre embedding}\index{Segre embedding|(}
Let $T$ and $W$ be schemes, and suppose $q:W \rightarrow T$ is a morphism.  By Lemma \ref{prop.duality}, the isomorphism
\begin{equation} \label{eqn.newone}
\Theta':q^{*}(- \otimes_{{\mathcal{O}}_{T}}-) \Longrightarrow (- \otimes_{{\mathcal{O}}_{W}}-)q^{*2}
\end{equation}
is dual to the square
\begin{equation} \label{eqn.square2}
\xymatrix{
({\sf Qcoh} W)^{2} \ar@<1ex>[r]^{q_{*}^{2}} \ar@{=>}[dr]|{\Upsilon'}\hole \ar[d]_{- \otimes_{{\mathcal{O}}_{W}}-} & ({\sf{Qcoh}} T)^{2} \ar[d]^{- \otimes_{{\mathcal{O}}_{T}}-} \ar@<1ex>[l]^{q^{*2}} \\
{\sf Qcoh} W \ar@<1ex>[r]^{q_{*}}  & {\sf{Qcoh}} T  \ar@<1ex>[l]^{q^{*}}
}
\end{equation}
where $\Upsilon'$ is the composition
\begin{equation} \label{eqn.newtwo}
\xymatrix{
(- \otimes_{{\mathcal{O}}_{T}}-)q_{*}^{2} \ar@{=>}[r] & q_{*}q^{*}(- \otimes_{{\mathcal{O}}_{T}}-)q_{*}^{2}   \ar@{=>}[r] &
}
\end{equation}
$$
q_{*}(-\otimes_{{\mathcal{O}}_{W}}-)q^{*2}q_{*}^{2} \Longrightarrow  q_{*}(- \otimes_{{\mathcal{O}}_{W}}-),
$$
where the second map is the product $q_{*}*\Theta'*q_{*}^{2}$.

The following result provides motivation for Theorem \ref{theorem.bigsegre}, and is used in its proof.
\begin{theorem} \label{prop.segre}
Let $T$ be a scheme and let $\mathcal{A}$ and $\mathcal{B}$ be ${\mathcal{O}}_{T}$-modules.  Then there exists a unique map

$$
s_{2}:{\mathbb{P}}_{T}(\mathcal{A}) \times_{T} {\mathbb{P}}_{T}(\mathcal{B}) \rightarrow {\mathbb{P}}_{T}({\mathcal{A}} \otimes_{{\mathcal{O}}_{T}} {\mathcal{B}})
$$
such that if $r:W \rightarrow {\mathbb{P}}_{T}(\mathcal{A}) \times_{T} {\mathbb{P}}_{T}(\mathcal{B})$ is a morphism whose projection to $T$ is $q$ and $r$ corresponds, via Proposition \ref{prop.groth}, to the pair
$$
\phi_{1}:q^{*}\mathcal{A} \rightarrow {\mathcal{L}}_{1}
$$
and
$$
\phi_{2}:q^{*}\mathcal{B} \rightarrow {\mathcal{L}}_{2}
$$
with right adjunct $\psi_{1}$ and $\psi_{2}$, respectively, then $s_{2}\circ r$ corresponds to the map
$$
\xymatrix{
q^{*}(\mathcal{A} \otimes \mathcal{B}) \ar[r]^{\Theta'} & q^{*}\mathcal{A}\otimes q^{*}\mathcal{B} \ar[r]^{\phi_{1} \otimes \phi_{2}} & {\mathcal{L}}_{1}\otimes {\mathcal{L}}_{2}.
}
$$
which has right adjunct
$$
\xymatrix{
\mathcal{A} \otimes \mathcal{B} \ar[rr]^{\psi_{1}\otimes \psi_{2}} & & q_{*}{\mathcal{L}}_{1}\otimes q_{*}{\mathcal{L}}_{2} \ar[r]^{\Upsilon'} & q_{*}({\mathcal{L}}_{1} \otimes {\mathcal{L}}_{2}).
}
$$
Furthermore, $s_{2}$ has the following properties:
\begin{enumerate}
\item{$s_{2}$ is a closed immersion,}\index{Segre embedding!is a closed immersion}
\item{$s_{2}$ is functorial:}\index{Segre embedding!is functorial}
If $\mathcal{A}'$ is an  ${\mathcal{O}}_{T}$-module such that there is an epimorphism $\phi:\mathcal{A} \rightarrow \mathcal{A}'$ then the diagram
$$
\xymatrix{
{\mathbb{P}}_{T}(\mathcal{A}') \times_{T} {\mathbb{P}}_{T}(\mathcal{B}) \ar[rrr]^{{\mathbb{P}}_{T}(\phi) \times \id} \ar[d]_{s_{2}} & & & {\mathbb{P}}_{T}(\mathcal{A}) \times_{T} {\mathbb{P}}_{T}(\mathcal{B}) \ar[d]^{s_{2}} \\
{\mathbb{P}}_{T}(\mathcal{A}' \otimes \mathcal{B}) \ar[rrr]  & & & {\mathbb{P}}_{T}(\mathcal{A} \otimes \mathcal{B})
}
$$
commutes, and
\item{$s_{2}$ is compatible with base change:}\index{Segre embedding!is compatible with base change}
If $g:T' \rightarrow T$ is a morphism of schemes, then the diagram
$$
\xymatrix{
{\mathbb{P}}_{T'}(g^{*}\mathcal{A}) \times_{T'} {\mathbb{P}}_{T'}(g^{*}\mathcal{B}) \ar[dd] \ar[r]^{s_{2}} & {\mathbb{P}}_{T'}(g^{*}\mathcal{A} \otimes_{{\mathcal{O}}_{T'}} g^{*}\mathcal{B}) \ar[d] \\
& {\mathbb{P}}_{T'}(g^{*}(\mathcal{A} \otimes_{{\mathcal{O}}_{T}} \mathcal{B})) \ar[d] \\
{\mathbb{P}}_{T}(\mathcal{A}) \times_{T} {\mathbb{P}}_{T}(\mathcal{B})  \ar[r]_{s_{2}} & {\mathbb{P}}_{T}(\mathcal{A} \otimes_{{\mathcal{O}}_{T}} \mathcal{B})
}
$$
whose left and bottom right vertical arrows are defined in Lemma \ref{lem.proj}, commutes.
\end{enumerate}
\end{theorem}

\begin{proof}
We repeat the construction of a map
$$
s_{2}:{\mathbb{P}}_{T}(\mathcal{A}) \times_{T} {\mathbb{P}}_{T}(\mathcal{B}) \rightarrow {\mathbb{P}}_{T}({\mathcal{A}} \otimes_{{\mathcal{O}}_{T}} {\mathcal{B}})
$$
given in \cite[4.3.1, p.76]{ega2}\index{Segre embedding|textbf}.  We show later that $s_{2}$ has the promised characteristic property.

Suppose
$$
u:  {\mathbb{P}}_{T}(\mathcal{A}) \times_{T} {\mathbb{P}}_{T}(\mathcal{B}) \rightarrow T
$$
is the structure map.  We construct $s_{2}$ using Proposition \ref{prop.groth}.  According to Proposition \ref{prop.groth}, in order to find a $T$-morphism from ${\mathbb{P}}_{T}(\mathcal{A}) \times_{T} {\mathbb{P}}_{T}(\mathcal{B})$ to ${\mathbb{P}}_{T}({\mathcal{A}} \otimes_{{\mathcal{O}}_{T}} {\mathcal{B}})$, it suffices to find an invertible sheaf, $\mathcal{L}$ on ${\mathbb{P}}_{T}(\mathcal{A}) \times_{T} {\mathbb{P}}_{T}(\mathcal{B})$, and an epimorphism,
$$
u^{*}(\mathcal{A} {\otimes}_{{\mathcal{O}}_{T}} \mathcal{B}) \rightarrow \mathcal{L}
$$
Let ${\mathcal{L}}_{\mathcal{A}} = {\mathcal{O}}_{{\mathbb{P}}_{T}(\mathcal{A})}$ and ${\mathcal{L}}_{\mathcal{B}} = {\mathcal{O}}_{{\mathbb{P}}_{T}(\mathcal{B})}$ be the structure sheaves on ${\mathbb{P}}_{T}(\mathcal{A})$ and ${\mathbb{P}}_{T}(\mathcal{B})$, respectively.  Let

$$
q_{\mathcal{A}}: {\mathbb{P}}_{T}(\mathcal{A}) \times_{T} {\mathbb{P}}_{T}(\mathcal{B}) \rightarrow {\mathbb{P}}_{T}(\mathcal{A})
$$
and
$$
q_{\mathcal{B}}: {\mathbb{P}}_{T}(\mathcal{A}) \times_{T} {\mathbb{P}}_{T}(\mathcal{B}) \rightarrow {\mathbb{P}}_{T}(\mathcal{B})
$$
be the canonical projection maps.  Define
\begin{equation}
{\mathcal{L}}=q_{\mathcal{A}}^{*}{\mathcal{L}}_{\mathcal{A}} {\otimes} q_{\mathcal{B}}^{*} {\mathcal{L}}_{\mathcal{B}}
\end{equation}
which is an invertible sheaf on ${\mathbb{P}}_{T}(\mathcal{A}) \times_{T} {\mathbb{P}}_{T}(\mathcal{B})$.  Since
$$
u^{*}\mathcal{A} {\otimes} u^{*}\mathcal{B} \cong u^{*}(\mathcal{A} {\otimes}_{{\mathcal{O}}_{T}} \mathcal{B})
$$
it suffices to construct an epimorphism
$$
\rho: u^{*}\mathcal{A} {\otimes} u^{*}\mathcal{B} \rightarrow \mathcal{L}
$$
Now, to construct $\rho$, it suffices to find epis:

\begin{equation} \label{eqn.firstyy}
\rho_{1}: u^{*}\mathcal{A} \rightarrow q_{\mathcal{A}}^{*}{{\mathcal{L}}_{\mathcal{A}}}
\end{equation}
and
\begin{equation} \label{eqn.secondyy}
\rho_{2}: u^{*}\mathcal{B} \rightarrow q_{\mathcal{B}}^{*}{{\mathcal{L}}_{\mathcal{B}}}
\end{equation}
We construct (\ref{eqn.firstyy}).  The construction of (\ref{eqn.secondyy}) is similar, so we omit it.  Let $p_{\mathcal{A}}$ and $p_{\mathcal{B}}$ be the structure maps on ${\mathbb{P}}_{T}(\mathcal{A})$ and ${\mathbb{P}}_{T}(\mathcal{B})$ respectively.  There is a canonical epi
\begin{equation} \label{eqn.newthree}
\delta_{\mathcal{A}}:p_{\mathcal{A}}^{*}\mathcal{A} \rightarrow {\mathcal{L}}_{\mathcal{A}}
\end{equation}
\cite[4.1.6, p.72]{ega2}.  Applying the right exact functor, $q_{\mathcal{A}}^{*}$, to this equation, yields an epi
$$
q_{\mathcal{A}}^{*} p_{\mathcal{A}}^{*}\mathcal{A} \rightarrow q_{\mathcal{A}}^{*}{{\mathcal{L}}_{\mathcal{A}}}
$$

But $p_{\mathcal{A}} q_{\mathcal{A}} = u$.  Thus, there is an epi
$$
u^{*} \mathcal{A} \rightarrow q_{\mathcal{A}}^{*}{{\mathcal{L}}_{\mathcal{A}}}
$$
as desired.

The fact that $s_{2}$ is a closed immersion\index{Segre embedding!is a closed immersion} is proven in \cite[4.3.3, p.77]{ega2}, the fact that $s_{2}$ is functorial\index{Segre embedding!is functorial} is proven in \cite[4.3.4, p.78]{ega2}, and the fact that $s_{2}$ is compatible with base change\index{Segre embedding!is compatible with base change} is proven in \cite[4.3.5, p.78]{ega2}.
\end{proof}

\begin{proposition} \label{prop.algsecond}
Let $T$ be a scheme and let $\mathcal{A}$ and $\mathcal{B}$ be ${\mathcal{O}}_{T}$-modules.  Then the map
$$
s_{2}:{\mathbb{P}}_{T}(\mathcal{A}) \times_{T} {\mathbb{P}}_{T}(\mathcal{B}) \rightarrow {\mathbb{P}}_{T}({\mathcal{A}} \otimes_{{\mathcal{O}}_{T}} {\mathcal{B}})
$$
defined in the proof of Theorem \ref{prop.segre} has the property that if $r:W \rightarrow {\mathbb{P}}_{T}(\mathcal{A}) \times_{T} {\mathbb{P}}_{T}(\mathcal{B})$ is a morphism whose projection to $T$ is $q$ and $r$ corresponds, via Proposition \ref{prop.groth}, to the pair
$$
\phi_{1}:q^{*}\mathcal{A} \rightarrow {\mathcal{L}}_{1}
$$
and
$$
\phi_{2}:q^{*}\mathcal{B} \rightarrow {\mathcal{L}}_{2}
$$
with right adjuncts $\psi_{1}$ and $\psi_{2}$ respectively, then $s_{2}\circ r$ corresponds to the map
$$
\xymatrix{
q^{*}(\mathcal{A} \otimes \mathcal{B}) \ar[r]^{\Theta'} & q^{*}\mathcal{A}\otimes q^{*}\mathcal{B} \ar[r]^{\phi_{1} \otimes \phi_{2}} & {\mathcal{L}}_{1}\otimes {\mathcal{L}}_{2}.
}
$$
with right adjunct
$$
\xymatrix{
\mathcal{A} \otimes \mathcal{B} \ar[r]^{\hskip -.12in \psi_{1} \otimes \psi_{2}} & q_{*}{\mathcal{L}}_{1}\otimes q_{*}{\mathcal{L}}_{2} \ar[r]^{\Upsilon'} & q_{*}({\mathcal{L}}_{1} \otimes {\mathcal{L}}_{2})
}
$$
where $\Theta'$ and $\Upsilon'$ were defined in (\ref{eqn.newone}) on page \pageref{eqn.newone} and in (\ref{eqn.newtwo}) on page \pageref{eqn.newtwo}, respectively.
\end{proposition}

\begin{proof}
Let $u:{\mathbb{P}}_{T}(\mathcal{A}) \times_{T} {\mathbb{P}}_{T}(\mathcal{B}) \rightarrow T$ be projection to the base, and suppose $r:W \rightarrow {\mathbb{P}}_{T}(\mathcal{A}) \times_{T} {\mathbb{P}}_{T}(\mathcal{B})$ projects to the two maps $r_{\mathcal{A}}:W \rightarrow {\mathbb{P}}_{T}(\mathcal{A})$ and $r_{\mathcal{B}}:W \rightarrow {\mathbb{P}}_{T}(\mathcal{B})$.  Thus, we have the following identities between maps:
$$
\xymatrix{
q_{\mathcal{A}}r=r_{\mathcal{A}} & & q_{\mathcal{B}}r=r_{\mathcal{B}} & & u=p_{\mathcal{A}}q_{\mathcal{A}}=p_{\mathcal{B}}q_{\mathcal{B}}
}
$$
As in the proof of Theorem \ref{prop.segre}, $s_{2}$ corresponds to a sheaf morphism
$$
\xymatrix{
u^{*}(\mathcal{A} \otimes \mathcal{B}) \ar[r] & {q_{\mathcal{A}}}^{*} {p_{\mathcal{A}}}^{*}\mathcal{A} \otimes {q_{\mathcal{B}}}^{*} {p_{\mathcal{B}}}^{*}\mathcal{B} \ar[rr]^{q_{\mathcal{A}}^{*}\delta_{\mathcal{A}} \otimes q_{\mathcal{B}}^{*}\delta_{\mathcal{B}}} & &  {q_{\mathcal{A}}}^{*} {\mathcal{L}}_{\mathcal{A}} \otimes  {q_{\mathcal{B}}}^{*} {\mathcal{L}}_{\mathcal{B}}
}
$$
where $\delta_{\mathcal{A}}$ is defined by (\ref{eqn.newthree}) and $\delta_{\mathcal{B}}$ is defined similarly.  Thus, the composition $s_{2} \circ r$ corresponds to the upper route in the following diagram (\cite[4.2.8, p. 75]{ega2}).
\begin{equation} \label{eqn.maps}
\xymatrix{
r^{*}u^{*}(\mathcal{A} \otimes \mathcal{B}) \ar[r] & r^{*}({q_{\mathcal{A}}}^{*} {p_{\mathcal{A}}}^{*}\mathcal{A} \otimes {q_{\mathcal{B}}}^{*} {p_{\mathcal{B}}}^{*}\mathcal{B}) \ar[rr]^{r^{*}(q_{\mathcal{A}}^{*}\delta_{\mathcal{A}} \otimes q_{\mathcal{B}}^{*}\delta_{\mathcal{B}})} \ar[d] & & r^{*}({q_{\mathcal{A}}}^{*} {\mathcal{L}}_{\mathcal{A}} \otimes  {q_{\mathcal{B}}}^{*} {\mathcal{L}}_{\mathcal{B}}) \ar[d] \\
(ur)^{*}(\mathcal{A} \otimes \mathcal{B}) \ar[u] \ar[r] & {r_{\mathcal{A}}}^{*}{p_{\mathcal{A}}}^{*}\mathcal{A} \otimes {r_{\mathcal{B}}}^{*}{p_{\mathcal{B}}}^{*}\mathcal{B} \ar[rr]^{{r_{\mathcal{A}}}^{*}\delta_{\mathcal{A}}\otimes {r_{\mathcal{B}}}^{*}\delta_{\mathcal{B}}} & & {r_{\mathcal{A}}}^{*}{\mathcal{L}}_{\mathcal{A}} \otimes {r_{\mathcal{B}}}^{*}{\mathcal{L}}_{\mathcal{B}}
}
\end{equation}
Since $r^{*}$ commutes naturally with tensor products, the right square commutes, and a local computation bears out the fact that the left square commutes.  Thus, by Proposition \ref{prop.groth}, $s_{2} \circ r$ corresponds to the bottom route of this diagram.

By the construction of the correspondence in Proposition \ref{prop.groth}, the map $q_{\mathcal{A}}:{\mathbb{P}}_{T}(\mathcal{A}) \times_{T} {\mathbb{P}}_{T}(\mathcal{B}) \rightarrow {\mathbb{P}}_{T}(\mathcal{A})$ corresponds to the epimorphism
$$
\xymatrix{
u^{*}\mathcal{A} \ar[r] & q^{*}_{\mathcal{A}}p^{*}_{\mathcal{A}}\mathcal{A} \ar[r]^{q^{*}_{\mathcal{A}}\delta_{\mathcal{A}}} & q_{\mathcal{A}}^{*}{\mathcal{L}}_{\mathcal{A}}.
}
$$
Thus, $r_{\mathcal{A}}=q_{\mathcal{A}}r$ corresponds to the epimorphism
$$
\xymatrix{
\sigma_{\mathcal{A}}:(p_{\mathcal{A}}r_{\mathcal{A}})^{*}\mathcal{A} \ar[r] & {r_{\mathcal{A}}}^{*}{p_{\mathcal{A}}}^{*}\mathcal{A} \ar[r]^{{r_{\mathcal{A}}}^{*}\delta_{\mathcal{A}}} & {r_{\mathcal{A}}}^{*}{\mathcal{L}}_{\mathcal{A}}
}
$$
(\cite[4.2.8, p.75]{ega2}).  But, by hypothesis, $r_{\mathcal{A}}$ also corresponds to an epimorphism
$$
\phi_{1}:(p_{\mathcal{A}}r_{\mathcal{A}})^{*}\mathcal{A} \rightarrow {\mathcal{L}}_{1}
$$
Thus, there exists an isomorphism $\tau_{\mathcal{A}}:{r_{\mathcal{A}}}^{*}{\mathcal{L}}_{\mathcal{A}} \rightarrow {\mathcal{L}}_{1}$ making the diagram
$$
\xymatrix{
q^{*}\mathcal{A} \ar[r]^{\sigma_{\mathcal{A}}} \ar[dr]_{\phi_{1}} & {r_{\mathcal{A}}}^{*}{\mathcal{L}}_{\mathcal{A}} \ar[d]^{\tau_{\mathcal{A}}} \\
& {\mathcal{L}}_{1}
}
$$
commute.  In a similar fashion, there exists a map $\sigma_{\mathcal{B}}$ and an isomorphism $\tau_{\mathcal{B}}$ making the diagram
$$
\xymatrix{
q^{*}\mathcal{B} \ar[r]^{\sigma_{\mathcal{B}}} \ar[dr]_{\phi_{2}} & {r_{\mathcal{B}}}^{*}{\mathcal{L}}_{\mathcal{B}} \ar[d]^{\tau_{\mathcal{B}}} \\
& {\mathcal{L}}_{2}
}
$$
commute.  In particular,
\begin{equation} \label{eqn.segrecounts}
\xymatrix{
q^{*}(\mathcal{A}\otimes \mathcal{B}) \ar[r] & q^{*}\mathcal{A}\otimes q^{*}\mathcal{B} \ar[rr]^{\sigma_{\mathcal{A}} \otimes \sigma_{\mathcal{B}}} \ar[drr]_{\phi_{1} \otimes \phi_{2}} & & {r_{\mathcal{A}}}^{*}{\mathcal{L}}_{\mathcal{A}} \otimes {r_{\mathcal{B}}}^{*}{\mathcal{L}}_{\mathcal{B}} \ar[d]^{\tau_{\mathcal{A} \otimes \tau_{\mathcal{B}}}} \\
& & & {\mathcal{L}}_{1} \otimes {\mathcal{L}}_{2}
}
\end{equation}
commutes so that the two routes of (\ref{eqn.segrecounts}) correspond to the same map
$$
W \rightarrow {\mathbb{P}}_{T}(\mathcal{A} \otimes \mathcal{B}).
$$
But, since the bottom route of (\ref{eqn.maps}) corresponds to $s_{2} \circ r$, and this map is the same as the top of (\ref{eqn.segrecounts}), the bottom of (\ref{eqn.segrecounts}) also corresponds to $s_{2} \circ r$.

To prove the final assertion, we need only note that an application of Corollary \ref{cor.adcom2} to (\ref{eqn.square2}) (on page \pageref{eqn.square2}) gives a commutative diagram
$$
\xymatrix{
\operatorname{Hom}_{({\sf{Mod }}W)^{2}}(q^{*2}(\mathcal{A},\mathcal{B}),({\mathcal{L}}_{1},{\mathcal{L}}_{2})) \ar[r] \ar[d]_{-\otimes_{{\mathcal{O}}_{W}}-} & \operatorname{Hom}_{({\sf{Mod }}T)^{2}}((\mathcal{A},\mathcal{B}),q_{*}^{2}({\mathcal{L}}_{1},{\mathcal{L}}_{2})) \ar[d]^{- \otimes_{{\mathcal{O}}_{T}}-} \\
\operatorname{Hom}_{W}(q^{*}\mathcal{A}\otimes q^{*}\mathcal{B},{\mathcal{L}}_{1}\otimes {\mathcal{L}}_{2}) \ar[d]_{- \circ \Theta'} & \operatorname{Hom}_{T}(\mathcal{A} \otimes \mathcal{B},q_{*}{\mathcal{L}}_{1}\otimes q_{*}{\mathcal{L}}_{2}) \ar[d]^{\Upsilon' \circ -} \\
\operatorname{Hom}_{W}(q^{*}(\mathcal{A} \otimes \mathcal{B}),{\mathcal{L}}_{1}\otimes {\mathcal{L}}_{2}) \ar[r] & \operatorname{Hom}_{T}(\mathcal{A}\otimes \mathcal{B},q_{*}({\mathcal{L}}_{1} \otimes {\mathcal{L}}_{2})).
}
$$
\end{proof}\index{Segre embedding|)}

\subsection{A map $s_{3}:{\mathbb{P}}_{X \times Y \times Z}({\pr_{12}}^{*}\mathcal{E} \otimes {\pr_{23}}^{*}\mathcal{F}) \rightarrow  {\mathbb{P}}_{X \times Z}({\pr_{13}}_{*}({\pr_{12}}^{*}\mathcal{E} \otimes {\pr_{23}}^{*}\mathcal{F}))$}
\begin{proposition} \label{prop.epi}
Let $f:Y \rightarrow Z$ be a morphism of schemes and suppose $\mathcal{A}$ is an ${\mathcal{O}}_{Y}$-module and $i:\operatorname{SSupp }\mathcal{A} \rightarrow Y$ is inclusion.  Suppose also that $fi$ is affine.  Then the counit

$$
\epsilon_{\mathcal{A}}:{f}^{*}{f}_{*}\mathcal{A} \rightarrow \mathcal{A}
$$
is an epi.
\end{proposition}

\begin{proof}
First, we note that the natural map $\mathcal{A} \rightarrow i_{*}i^{*}\mathcal{A}$ is an isomorphism, and, for any ${\mathcal{O}}_{Y}$-module $\mathcal{N}$, the natural map $i_{*}i^{*}\mathcal{N}\rightarrow \mathcal{N}$ is an epimorphism.  Thus, since $i_{*}$ is exact, it suffices to show that the natural map

$$
i^{*}\epsilon_{\mathcal{A}}:i^{*}{f}^{*}{f}_{*}\mathcal{A} \rightarrow i^{*}\mathcal{A}
$$
is an epimorphism.  We construct an epi between these sheaves, and show it is equal to $i^{*}\epsilon_{\mathcal{A}}$.  Let
$$
\Xi_{i,f}:i^{*}{f}^{*} \rightarrow (fi)^{*}
$$
be the isomorphism of functors in Lemma \ref{lem.adjointnew}, and let $\Psi=\Xi_{i,f}*(fi)_{*}$.  Then the map
\begin{equation} \label{eqn.topy}
\xymatrix{
i^{*}{f}^{*}{f}_{*}\mathcal{A} \ar[r] & i^{*}{f}^{*}{f}_{*}i_{*}i^{*}\mathcal{A} \ar[r]^{\Psi_{i^{*}\mathcal{A}}} & (fi)^{*}(fi)_{*}i^{*}\mathcal{A} \ar[r] & i^{*}\mathcal{A}.
}
\end{equation}
is an epi since the final map is an epi \cite[Corollary 1,p.39]{curves}.  We claim that the diagram
\begin{equation} \label{eqn.adjoint}
\xymatrix{
i^{*}{f}^{*}{f}_{*}\mathcal{A} \ar[r] \ar[d]_{i^{*}\epsilon_{\mathcal{A}}} & i^{*}{f}^{*}{f}_{*}i_{*}i^{*}\mathcal{A} \ar[d] \ar[r]^{\Psi} & (fi)^{*}(fi)_{*}i^{*}\mathcal{A} \ar[d] \\
i^{*}\mathcal{A} \ar[r] & i^{*}i_{*}i^{*}\mathcal{A} \ar[r] & i^{*}\mathcal{A}
}
\end{equation}
whose top circuit is (\ref{eqn.topy}), commutes.  The right square commutes by Lemma \ref{lem.adjointnew}, while the left square commutes by naturality of $\epsilon_{\mathcal{A}}$.  Since the bottom row of this diagram is the identity, the assertion follows.
\end{proof}

\begin{lemma} \label{lem.sup}
Let $g:W \rightarrow Z$ be an affine morphism of schemes, and suppose $\mathcal{C}$ is an ${\mathcal{O}}_{W}$-module.  Then $g_{*}\mathcal{C}$ is quasi-coherent and the composition
$$
c:{\mathbb{P}}_{W}(\mathcal{C}) \rightarrow {\mathbb{P}}_{W}(g^{*}g_{*}\mathcal{C}) \rightarrow {\mathbb{P}}_{Z}(g_{*}\mathcal{C})
$$
whose left hand map is induced by the counit of the pair $(g^{*},g_{*})$, and whose right hand map is (\ref{eqn.proj}) (defined on page \pageref{eqn.proj}), is a closed immersion.
\end{lemma}

\begin{proof}
The fact that $g_{*}\mathcal{C}$ is a ${\mathcal{O}}_{W}$-module follows from \cite[II.5.8, p.115]{alggeo}.  Since $g$ is affine, and since the construction of (\ref{eqn.proj}) is local, to prove $c$ is a closed immersion, we may assume that $W$ and $Z$ are affine.  Let $W$=Spec $Q$, $Z$=Spec $R$, and suppose $g:R \rightarrow Q$ is a ring map.  Suppose $C$ is an $Q$-module.  Then, by Lemma \ref{lem.proj}, the map
$$
{\mathbb{P}}_{W}(C_{R} \otimes Q) \rightarrow {\mathbb{P}}_{Z}(C_{R})
$$
is induced by a map of graded rings:
$$
{\mathbb{S}}_{R}(C_{R}) \rightarrow {\mathbb{S}}_{Q}(C_{R} \otimes_{R} Q) \cong {\mathbb{S}}_{R}(C_{R}) \otimes_{R} Q
$$
sending $s_{n}$ to $s_{n} \otimes 1$.  Furthermore, the map
$$
{\mathbb{P}}_{Z}(C) \rightarrow {\mathbb{P}}_{Z}(C_{R}\otimes Q)
$$
is induced by the map of graded rings
$$
{\mathbb{S}}_{Q}(C_{R} \otimes_{R} Q) \rightarrow {\mathbb{S}}_{Q}(C)
$$
given by the canonical epimorphism $C_{R} \otimes_{R} Q \rightarrow C$.  Putting the two graded ring maps together, we have a map
\begin{equation} \label{eq.graded}
{\mathbb{S}}_{R}(C_{R}) \rightarrow {\mathbb{S}}_{Q}(C_{R} \otimes_{R} Q) \rightarrow {\mathbb{S}}_{Q}(C)
\end{equation}
which is a surjection in degrees $\geq 1$.  Thus, the map ${\mathbb{P}}_{W}(C) \rightarrow {\mathbb{P}}_{Z}(C_{R})$ is a closed immersion.
\end{proof}

\begin{lemma} \label{lem.eisenbud} \cite[Lemma 6.4, p.163]{comalg}
Let $M$ and $N$ be $R$-modules, and suppose that $N$ is generated by a family of elements $\{n_{i}\}$.  Every element of $M \otimes_{R}N$ may be written as a finite sum $\Sigma_{i}m_{i}\otimes n_{i}$.  Such an expression is $0$ if and only if there exist elements $m_{j}'$ of $M$ and elements $a_{ij}$ of $R$ such that
$$
\Sigma_{j}a_{ij}m_{j}' = m_{i} \mbox{ for all }i
$$
and
$$
\Sigma_{i}a_{ij}n_{i}=0 \mbox{ in $N$ for all }j.
$$
\end{lemma}

\begin{lemma} \label{lem.isom}
Let $Y$ be a scheme and suppose $\mathcal{A}$ is an ${\mathcal{O}}_{Y}$-module with
$$
i:\operatorname{SSupp }\mathcal{A} \rightarrow Y
$$
inclusion.  Then the morphism (\ref{eqn.proj}) defined on page \pageref{eqn.proj},
$$
{\mathbb{P}}_{\operatorname{SSupp }\mathcal{A}}(i^{*}\mathcal{A}) \rightarrow {\mathbb{P}}_{Y}(\mathcal{A})
$$
is an isomorphism.
\end{lemma}

\begin{proof}
As in the proof of Lemma \ref{lem.sup}, the assertion is local on $Y$.  Thus, suppose $Y=$Spec $R$, and $M$ is an $R$-module with ann $M=I$.  We claim that the map

\begin{equation} \label{eqn.local}
{\mathbb{P}}_{\mbox{Spec }R/I}(M \otimes_{R} R/I) \rightarrow {\mathbb{P}}_{\mbox{Spec }R}(M)
\end{equation}
is an isomorphism.  For, to construct this map, we start with the natural map of sets $M \rightarrow M \otimes R/I$.  This set map is actually a map of $R$-modules.  In fact, it induces a map of $R$-modules

$$
M^{\otimes n} \rightarrow M^{\otimes n} \otimes R/I.
$$
By \ref{lem.eisenbud}, this is an injective map of $R$-modules.  Thus, in degree $n>0$ the induced map on symmetric algebras is an isomorphism, hence the map ($\ref{eqn.local}$) is an isomorphism.
\end{proof}

\begin{proposition} \label{prop.thirdmap}
Let $f:Y \rightarrow Z$ be a morphism of schemes and suppose $\mathcal{A}$ is an ${\mathcal{O}}_{Y}$-module which is rla with respect to $f$\index{relatively locally affine}.  Then the counit
$$
f^{*}f_{*}\mathcal{A} \rightarrow \mathcal{A}
$$
is epi, and the corresponding morphism
$$
d:{\mathbb{P}}_{Y}(\mathcal{A}) \rightarrow {\mathbb{P}}_{Y}(f^{*}f_{*}\mathcal{A}) \rightarrow {\mathbb{P}}_{Z}(f_{*}\mathcal{A})
$$
is a closed immersion.
\end{proposition}

\begin{proof}
The first assertion follows from Proposition \ref{prop.epi}.  Since $i_{*}i^{*}\mathcal{A} \cong \mathcal{A}$, it suffices to prove the next assertion with $i_{*}i^{*}\mathcal{A}$ replacing $\mathcal{A}$.  For, suppose
\begin{equation} \label{eq.im}
{\mathbb{P}}_{Y}(i_{*}i^{*}\mathcal{A}) \rightarrow {\mathbb{P}}_{Y}(f^{*}f_{*}i_{*}i^{*}\mathcal{A}) \rightarrow {\mathbb{P}}_{Z}(f_{*}i_{*}i^{*}\mathcal{A})
\end{equation}
is a closed immersion.  Then by naturality of the counit of the pair $(f^{*},f_{*})$, the left square of
$$
\xymatrix{
{\mathbb{P}}_{Y}(i_{*}i^{*}\mathcal{A}) \ar[r] \ar[d] & {\mathbb{P}}_{Y}(f^{*}f_{*}i_{*}i^{*}\mathcal{A}) \ar[r] \ar[d] & {\mathbb{P}}_{Z}(f_{*}i_{*}i^{*}\mathcal{A}) \ar[d] \\
{\mathbb{P}}_{Y}(\mathcal{A}) \ar[r] & {\mathbb{P}}_{Y}(f^{*}f_{*}\mathcal{A}) \ar[r] & {\mathbb{P}}_{Z}(f_{*}\mathcal{A})
}
$$
commutes.  By Lemma \ref{lem.proj}, the right square commutes.  Thus, in order to prove the proposition, we must show that (\ref{eq.im}) is a closed immersion.  To this end, we note that the diagram

$$
\xymatrix{
{\mathbb{P}}_{Y}(i_{*}i^{*}\mathcal{A}) \ar[r] & {\mathbb{P}}_{Y}(f^{*}f_{*}i_{*}i^{*}\mathcal{A}) \ar[r] & {\mathbb{P}}_{Z}(f_{*}i_{*}i^{*}\mathcal{A}) \\
{\mathbb{P}}_{\operatorname{SSupp }\mathcal{A}}(i^{*}i_{*}i^{*}\mathcal{A}) \ar[r] \ar[u] & {\mathbb{P}}_{\operatorname{SSupp }\mathcal{A}}(i^{*}f^{*}f_{*}i_{*}i^{*}\mathcal{A}) \ar[u] \ar[r] & {\mathbb{P}}_{\operatorname{SSupp }\mathcal{A}}((fi)^{*}(fi)_{*}i^{*}\mathcal{A}) \ar[u] \\
{\mathbb{P}}_{\operatorname{SSupp }\mathcal{A}}(i^{*}\mathcal{A}) \ar[u] \ar[r] &  {\mathbb{P}}_{\operatorname{SSupp }\mathcal{A}}(i^{*}i_{*}i^{*}\mathcal{A}) \ar[u] \ar[r] & {\mathbb{P}}_{\operatorname{SSupp }\mathcal{A}}(i^{*}\mathcal{A}) \ar[u]
}
$$
whose top verticals are (\ref{eqn.proj}), defined on page \pageref{eqn.proj}, and whose bottom verticals are induced by counits of the pairs $(i^{*},i_{*})$, $(f^{*},f_{*})$ and $((fi)^{*},(fi)_{*})$ respectively, commutes.  To prove this fact, we note that the top left square commutes by Lemma \ref{lem.proj}, the top right square commutes by Lemma \ref{lem.compat}, the bottom left square is obviously commutative, and the lower right square commutes by Lemma \ref{lem.adjointnew}.  In addition the left vertical is an isomorphism by Lemma \ref{lem.isom}.  Thus to show the top horizontal map is a closed immersion, it suffices to show that the bottom route of the diagram is a closed immersion.  Since the bottom row is the identity map, we need only show that the right vertical is a closed immersion.  This follows from Lemma \ref{lem.sup} since $\mathcal{A}$ is rla with respect to $f$.
\end{proof}

\begin{lemma} \label{lem.cannon}
Let
$$
\xymatrix{
W \ar[r]^{q} & Y \ar[r]^{f} & Z
}
$$
be a morphism of schemes, let $\epsilon$ be the counit of the pair $(f^{*},f_{*})$, and let $\Xi_{f,q}:q^{*}f^{*} \Longrightarrow (fq)^{*}$ be the isomorphism defined in Lemma \ref{lem.adjointnew}.  Then the dual 2-cell to the diagram of schemes
\begin{equation} \label{eqn.expect}
\xymatrix{
W \ar[d]_{\id} \ar[r]^{q} &  Y \ar[d]^{f} \\
W \ar[r]_{fq} & Z
}
\end{equation}
$$
\Psi:(fq)^{*}f_{*} \Longrightarrow q^{*}
$$
equals the composition
$$
\xymatrix{
(fq)^{*}f_{*} \ar@{=>}[r]^{\Xi_{f,q}^{-1}*f_{*}} & q^{*}f^{*}f_{*} \ar@{=>}[r]^{q^{*}*\epsilon} & q^{*}.
}
$$
\end{lemma}

\begin{proof}
Let $(\eta',\epsilon')$ be the unit and counit of $(q^{*},q_{*})$, and let $\epsilon''$ be the counit of $((fq)^{*},(fq)_{*})$.  Then, by definition, $\Psi$ equals the composition
$$
\xymatrix{
(fq)^{*}f_{*} \ar@{=>}[rr]^{(fq)^{*}f_{*}*\eta'} & & (fq)^{*}f_{*}q_{*}q^{*} \ar@{=>}[r]^{\hskip .3in \epsilon''*q^{*}} & q^{*}.
}
$$
We claim the diagram
\begin{equation} \label{eqn.stary}
\xymatrix{
(fq)^{*}f_{*} \ar@{=>}[r] \ar@{=>}[d]_{\Xi_{f,q}^{-1} * f_{*}} & (fq)^{*}f_{*}q_{*}q^{*} \ar@{=>}[rr] \ar@{=>}[d]^{\Xi_{f,q}^{-1} * (fq)_{*}q^{*}} & & q^{*} \\
q^{*}f^{*}f_{*} \ar@{=>}[r]  & q^{*}f^{*}f_{*}q_{*}q^{*} \ar@{=>}[rr]_{q^{*}*\epsilon''*q_{*}q^{*}} & & q^{*}q_{*}q^{*} \ar@{=>}[u]
}
\end{equation}
whose top row is $\Psi$, commutes.  For, the left square commutes by naturality of $\eta'$, while the right square commutes by Lemma \ref{lem.adjointnew}.  To complete the proof of the assertion, we must show that the bottom of the diagram
$$
\xymatrix{
q^{*}f^{*}f_{*} \ar@{=>}[r] & q^{*}f^{*}f_{*}q_{*}q^{*} \ar@{=>}[r] & q^{*}
}
$$
equals
\begin{equation} \label{eqn.staryy}
\xymatrix{
q^{*}f^{*}f_{*} \ar@{=>}[r]^{q^{*} * \epsilon} & q^{*}.
}
\end{equation}
This follows from the commutativity of
$$
\xymatrix{
q^{*}f^{*}f_{*} \ar@{=>}[d]_{q^{*} * \epsilon} \ar@{=>}[r] & q^{*}f^{*}f_{*}q_{*}q^{*} \ar@{=>}[d]^{q^{*} * \epsilon * q_{*} q^{*}} \ar@{=>}[rr]  & & q^{*}q_{*}q^{*} \ar@{=>}[d] \\
q^{*} \ar@{=>}[r]_{q^{*} * \eta'} &  q^{*}q_{*}q^{*} \ar@{=>}[rr]_{\epsilon' * q^{*}} & & q^{*}
}
$$
since the top route of this diagram is the bottom of (\ref{eqn.stary}) while the bottom route is (\ref{eqn.staryy}).  The right square obviously commutes and the left square commutes by naturality of $\eta'$.
\end{proof}

\begin{lemma} \label{lem.algthird}
Suppose $f:Y \rightarrow Z$ is a morphism of schemes, $\mathcal{A}$ is an ${\mathcal{O}}_{Y}$-module, $i:\operatorname{SSupp }\mathcal{A} \rightarrow Y$ is inclusion and $fi$ is affine.  Suppose further, that $W$ is a scheme and $r:W \rightarrow {\mathbb{P}}_{Y}(\mathcal{A})$ is a morphism with projection $q:W \rightarrow Y$.  Then the counit $\epsilon_{\mathcal{A}}:f^{*}f_{*}\mathcal{A} \rightarrow \mathcal{A}$ is an epimorphism, and if $r$ corresponds to an epimorphism
$$
\phi:q^{*}\mathcal{A} \rightarrow \mathcal{L},
$$
which has right adjunct
$$
\psi:\mathcal{A} \rightarrow q_{*}\mathcal{L}
$$
then the composition of $r$ with the morphism $d$ of Proposition \ref{prop.thirdmap} corresponds to the epimorphism
\begin{equation} \label{eqn.topiii}
\xymatrix{
(fq)^{*}f_{*}\mathcal{A} \ar[r] & q^{*}f^{*}f_{*}\mathcal{A} \ar[r]^{q^{*}\epsilon_{\mathcal{A}}} & q^{*}\mathcal{A} \ar[r]^{\phi} & \mathcal{L}
}
\end{equation}
which has right adjunct
\begin{equation} \label{eqn.topiiir}
\xymatrix{
f_{*}\mathcal{A} \ar[r]^{f_{*}\psi} & f_{*}q_{*}\mathcal{L}.
}
\end{equation}
\end{lemma}

\begin{proof}
The first assertion follows from Proposition \ref{prop.thirdmap}.

To prove the second assertion, we note that the diagram
$$
\xymatrix{
W \ar[r]^{\cong} & {\mathbb{P}}_{W}(\mathcal{L}) \ar[r] & {\mathbb{P}}_{W}(q^{*}\mathcal{A}) \ar[r] \ar[d] & {\mathbb{P}}_{W}(q^{*}f^{*}f_{*}\mathcal{A}) \ar[d] \ar[r] & {\mathbb{P}}_{W}((fq)^{*}f_{*}\mathcal{A}) \ar[d] \\
& & {\mathbb{P}}_{Y}(\mathcal{A}) \ar[r] & {\mathbb{P}}_{Y}(f^{*}f_{*}\mathcal{A}) \ar[r] & {\mathbb{P}}_{Z}(f_{*}\mathcal{A})
}
$$
commutes: the left square commutes by naturality of the counit of the pair $(f^{*},f_{*})$, while the right square commutes by Lemma \ref{lem.compat}.  Now, $d \circ r$ is the bottom circuit of the diagram, so it is also the top circuit.  By applying the correspondence of Proposition \ref{prop.groth} to the epimorphism \ref{eqn.topiii}, it is easy to see that the corresponding map is the top circuit of the diagram.

We next show that the map (\ref{eqn.topiii}) has right adjoint (\ref{eqn.topiiir}).  The dual 2-cell
$$
\Delta:(fq)^{*}f_{*} \Longrightarrow q^{*}
$$
to the square induced by the diagram of schemes
$$
\xymatrix{
W \ar[r]^{q} \ar[d]_{\id} & Y  \ar[d]^{f} \\
W \ar[r]_{fq} & Z
}
$$
equals, by Lemma \ref{lem.cannon}, the first two composites of (\ref{eqn.topiii}).  By Corollary \ref{cor.adcom2}, there is a commutative diagram
$$
\xymatrix{
\operatorname{Hom}_{W}(q^{*}\mathcal{A}, \mathcal{L}) \ar[r] \ar[dd]_{- \circ \Delta} & \operatorname{Hom}_{Y}(\mathcal{A}, q_{*}{\mathcal{L}}) \ar[d]^{f_{*}} \\
 & \operatorname{Hom}_{Z}(f_{*}\mathcal{A}, (fq)_{*}\mathcal{L}) \ar[d]^{\id(-)} \\
\operatorname{Hom}_{W}((fq)^{*}f_{*}\mathcal{A}, \mathcal{L}) \ar[r] & \operatorname{Hom}_{Z}(f_{*}\mathcal{A}, (fq)_{*}\mathcal{L}).
}
$$
whence the result.
\end{proof}

\begin{theorem} \label{theorem.algmap}
Let $\mathcal{E}$ be an ${{\mathcal{O}}_{X \times Y}}$-module and let $\mathcal{F}$ be an ${{\mathcal{O}}_{Y \times Z}}$-module such that $\mathcal{E}$ and $\mathcal{F}$ have the affine direct image property.  Then the bimodule Segre embedding
$$
s:{\mathbb{P}}_{X\times Y}(\mathcal{E}) \otimes_{Y} {\mathbb{P}}_{Y \times Z}(\mathcal{F}) \rightarrow {\mathbb{P}}_{X \times Z}({\mathcal{E}} \otimes_{{\mathcal{O}}_{Y}} {\mathcal{F}})
$$
has the property that if $r:W \rightarrow {\mathbb{P}}_{X\times Y}(\mathcal{E}) \otimes_{Y} {\mathbb{P}}_{Y \times Z}(\mathcal{F})$ is a morphism whose projection to $X \times Y \times Z$ is $q_{123}$ and $r$ corresponds, via Proposition \ref{prop.groth}, to a pair with right adjuncts
$$
\psi_{1}:\mathcal{E} \rightarrow q_{12*}{\mathcal{L}}_{1}
$$
and
$$
\psi_{2}:\mathcal{F} \rightarrow q_{23*}{\mathcal{L}}_{2}
$$
then $s\circ r$ corresponds to the map
$$
\xymatrix{
\mathcal{E} \otimes_{{\mathcal{O}}_{Y}} \mathcal{F} \ar[rr]^{\hskip -.3in \psi_{1} \otimes_{{\mathcal{O}}_{Y}} \psi_{2}} & & q_{12*}{\mathcal{L}}_{1} \otimes_{{\mathcal{O}}_{Y}} q_{23*}{\mathcal{L}}_{2} \ar[rr]^{\hskip .5in \pr_{13*}(\Omega_{1} \otimes \Omega_{2})} & &
}
$$
$$
\xymatrix{
\pr_{13*}(q_{123*}{\mathcal{L}}_{1} \otimes q_{123*}{\mathcal{L}}_{2}) \ar[rr]^{\hskip .4in \pr_{13*}\Upsilon} & & q_{13*}({\mathcal{L}}_{1} \otimes {\mathcal{L}}_{2}).
}
$$
\end{theorem}

\begin{proof}
Suppose $r:W \rightarrow {\mathbb{P}}_{X\times Y}(\mathcal{E}) \otimes_{Y} {\mathbb{P}}_{Y \times Z}(\mathcal{F})$ corresponds to a pair of epimorphism with right adjuncts
$$
\psi_{1}:\mathcal{E} \rightarrow q_{12*}{\mathcal{L}}_{1}
$$
and
$$
\psi_{2}:\mathcal{F} \rightarrow q_{23*}{\mathcal{L}}_{2}.
$$
By Proposition \ref{prop.algfirst}, $s_{1} \circ r$ corresponds to the pair
$$
\xymatrix{
\pr_{12}^{*}\mathcal{E} \ar[rr]^{\pr_{12}^{*}\psi_{1}} & & \pr_{12}^{*}q_{12*}{\mathcal{L}}_{1} \ar[r]^{\Omega_{1{\mathcal{L}}_{1}}} & q_{123*}{\mathcal{L}}_{1}
}
$$
and
$$
\xymatrix{
\pr_{23}^{*}\mathcal{F} \ar[rr]^{\pr_{23}^{*}\psi_{2}} & & \pr_{23}^{*}q_{23*}{\mathcal{L}}_{2} \ar[r]^{\Omega_{2{\mathcal{L}}_{2}}} & q_{123*}{\mathcal{L}}_{2}.
}
$$
By Theorem \ref{prop.segre}, $s_{2} \circ s_{1} \circ r$ corresponds to an epi with right adjunct
$$
\xymatrix{
\pr_{12}^{*}\mathcal{E} \otimes \pr_{23}^{*}\mathcal{F} \ar[rr]^{\hskip -.3in \pr_{12}^{*}\psi_{1} \otimes \pr_{23}^{*}\psi_{2}} & & \pr_{12}^{*}q_{12*}{\mathcal{L}}_{1} \otimes \pr_{23}^{*}q_{23*}{\mathcal{L}}_{2} \ar[r]^{\hskip .6in \Omega_{1} \otimes \Omega_{2}} &
}
$$
$$
\xymatrix{
q_{123*}{\mathcal{L}}_{1} \otimes q_{123*}{\mathcal{L}}_{2} \ar[r]^{\Upsilon} & q_{123*}({\mathcal{L}}_{1} \otimes {\mathcal{L}}_{2}).
}
$$
Finally, by Lemma \ref{lem.algthird}, $s_{3} \circ s_{2} \circ s_{1} \circ r$ corresponds to
$$
\xymatrix{
\mathcal{E} \otimes_{{\mathcal{O}}_{Y}} \mathcal{F} \ar[rr]^{\hskip -.3in \psi_{1} \otimes_{{\mathcal{O}}_{Y}} \psi_{2}} & & q_{12*}{\mathcal{L}}_{1} \otimes_{{\mathcal{O}}_{Y}} q_{23*}{\mathcal{L}}_{2} \ar[rr]^{\hskip .5in \pr_{13*}(\Omega_{1} \otimes \Omega_{2})} & &
}
$$
$$
\xymatrix{
\pr_{13*}(q_{123*}{\mathcal{L}}_{1} \otimes q_{123*}{\mathcal{L}}_{2}) \ar[rr]^{\hskip .4in \pr_{13*}\Upsilon} & & q_{13*}({\mathcal{L}}_{1} \otimes {\mathcal{L}}_{2}).
}
$$
as desired.
\end{proof}

\section{$s$ is functorial}\index{bimodule Segre embedding!is functorial|(}
\begin{theorem} \label{theorem.functorial}
Suppose $\mathcal{E}$ and $\mathcal{E}'$ are ${\mathcal{O}}_{X \times Y}$-modules such that there is an epimorphism $\upsilon:\mathcal{E} \rightarrow \mathcal{E}'$.  Suppose, further, that $\mathcal{F}$ is an ${\mathcal{O}}_{Y \times Z}$-module such that $\mathcal{E}$ and $\mathcal{F}$ have the affine direct image property\index{affine direct image property}.  Then $\mathcal{E}'$ and $\mathcal{F}$ have the direct image property, the map
$$
\upsilon \otimes_{{\mathcal{O}}_{Y}} \mathcal{F}:\mathcal{E} \otimes_{{\mathcal{O}}_{Y}} \mathcal{F} \rightarrow \mathcal{E}' \otimes_{{\mathcal{O}}_{Y}} \mathcal{F}
$$
is an epimorphism, and the diagram
\begin{equation} \label{eqn.functorial}
\xymatrix{
{\mathbb{P}}_{X \times Y}(\mathcal{E}') \otimes_{Y} {\mathbb{P}}_{Y \times Z}(\mathcal{F}) \ar[rrr]^{{\mathbb{P}}_{X \times Y}(\upsilon) \times \id} \ar[d]_{s} & & & {\mathbb{P}}_{X \times Y}(\mathcal{E}) \otimes_{Y} {\mathbb{P}}_{Y \times Z}(\mathcal{F}) \ar[d]^{s} \\
{\mathbb{P}}_{X \times Z}(\mathcal{E}' \otimes_{{\mathcal{O}}_{Y}} \mathcal{F}) \ar[rrr]_{{\mathbb{P}}_{X \times Z}(\upsilon \otimes_{{\mathcal{O}}_{Y}} \mathcal{F})}  & & & {\mathbb{P}}_{X \times Z}(\mathcal{E} \otimes_{{\mathcal{O}}_{Y}} \mathcal{F})
}
\end{equation}
commutes.
\end{theorem}

\begin{proof}
To prove the first assertion, we need only note that since there is an epimorphism $\upsilon:\mathcal{E} \rightarrow \mathcal{E}'$, $\operatorname{Supp }\mathcal{E}' \subset \operatorname{Supp }\mathcal{E}$, so that $\operatorname{Supp }\pr_{12}^{*}\mathcal{E}'\otimes \pr_{23}^{*}\mathcal{F} \subset \operatorname{Supp }\pr_{12}^{*}\mathcal{E} \otimes \pr_{23}^{*}\mathcal{F}$.

We next show that
$$
\upsilon \otimes_{{\mathcal{O}}_{Y}} \mathcal{F}:\mathcal{E} \otimes_{{\mathcal{O}}_{Y}} \mathcal{F} \rightarrow \mathcal{E}' \otimes_{{\mathcal{O}}_{Y}} \mathcal{F}
$$
is an epimorphism.  Since $\pr_{12}^{*}\upsilon \otimes \pr_{23}^{*}\mathcal{F}$ is an epimorphism, we need only note that both $\pr_{12}^{*}\mathcal{E} \otimes \pr_{23}^{*}\mathcal{F}$ and $\pr_{12}^{*}\mathcal{E}' \otimes \pr_{23}^{*}\mathcal{F}$ are rla with respect to $\pr_{13}$, so that the assertion follows from Lemma \ref{lem.exact}.

We now show that (\ref{eqn.functorial}) commutes.  Let $W$ be an $S$-scheme and suppose $r:W \rightarrow {\mathbb{P}}_{X\times Y}(\mathcal{E}') \otimes_{Y} {\mathbb{P}}_{Y \times Z}(\mathcal{F})$ is a morphism whose projection to $X \times Y \times Z$ is $q_{123}$ and $r$ corresponds, via Proposition \ref{prop.groth}, to a pair which has right adjuncts
$$
\psi_{1}:\mathcal{E}' \rightarrow q_{12*}{\mathcal{L}}_{1}
$$
and
$$
\psi_{2}:\mathcal{F} \rightarrow q_{23*}{\mathcal{L}}_{2}.
$$
By Theorem \ref{theorem.algmap}, $s \circ r$ corresponds to the map
$$
\xymatrix{
\mathcal{E}' \otimes_{{\mathcal{O}}_{Y}} \mathcal{F} \ar[rr]^{\hskip -.3in \psi_{1} \otimes_{{\mathcal{O}}_{Y}} \psi_{2}} & & q_{12*}{\mathcal{L}}_{1} \otimes_{{\mathcal{O}}_{Y}} q_{23*}{\mathcal{L}}_{2} \ar[rr]^{\hskip .5in \pr_{13*}(\Omega_{1} \otimes \Omega_{2})} & &
}
$$
$$
\xymatrix{
\pr_{13*}(q_{123*}{\mathcal{L}}_{1} \otimes q_{123*}{\mathcal{L}}_{2}) \ar[rr]^{\hskip .4in \pr_{13*}\Upsilon} & & q_{13*}({\mathcal{L}}_{1} \otimes {\mathcal{L}}_{2}).
}
$$
We claim
\begin{equation} \label{eqn.route}
\xymatrix{
W \ar[d]_{r} & & & \\
{\mathbb{P}}_{X \times Y}(\mathcal{E}') \otimes_{Y} {\mathbb{P}}_{Y \times Z}(\mathcal{F}) \ar[rrr]^{{\mathbb{P}}_{X \times Y}(\upsilon) \times \id} \ar[d]_{s} & & & {\mathbb{P}}_{X \times Y}(\mathcal{E}) \otimes_{Y} {\mathbb{P}}_{Y \times Z}(\mathcal{F}) \ar[d]^{s} \\
{\mathbb{P}}_{X \times Z}(\mathcal{E}' \otimes_{{\mathcal{O}}_{Y}} \mathcal{F}) \ar[rrr]  & & & {\mathbb{P}}_{X \times Z}(\mathcal{E} \otimes_{{\mathcal{O}}_{Y}} \mathcal{F})
}
\end{equation}
commutes.  By Proposition \ref{prop.groth}, the composition $({\mathbb{P}}_{X \times Y}(\upsilon) \times \id) \circ r$ corresponds to a pair of epis with right adjuncts
$$
\xymatrix{
\mathcal{E} \ar[r]^{\upsilon} & \mathcal{E}' \ar[r]^{\hskip -.11in \psi_{1}} &  q_{12*}{\mathcal{L}}_{1}
}
$$
and
$$
\xymatrix{
\mathcal{F} \ar[r]^{\hskip -.11in \psi_{2}} &  q_{23*}{\mathcal{L}}_{2},
}
$$
so that, by Theorem \ref{theorem.algmap}, the right hand route of (\ref{eqn.route}) corresponds to
$$
\xymatrix{
\mathcal{E} \otimes_{{\mathcal{O}}_{Y}} \mathcal{F} \ar[r]^{\upsilon \otimes_{{\mathcal{O}}_{Y}} \mathcal{F}}  &    \mathcal{E}' \otimes_{{\mathcal{O}}_{Y}} \mathcal{F} \ar[rr]^{\hskip -.3in \psi_{1} \otimes_{{\mathcal{O}}_{Y}} \psi_{2}} & & q_{12*}{\mathcal{L}}_{1} \otimes_{{\mathcal{O}}_{Y}} q_{23*}{\mathcal{L}}_{2} \ar[rr]^{\hskip .5in \pr_{13*}(\Omega_{1} \otimes \Omega_{2})} & &
}
$$
$$
\xymatrix{
\pr_{13*}(q_{123*}{\mathcal{L}}_{1} \otimes q_{123*}{\mathcal{L}}_{2}) \ar[rr]^{\hskip .4in \pr_{13*}\Upsilon} & & q_{13*}({\mathcal{L}}_{1} \otimes {\mathcal{L}}_{2}).
}
$$
On the other hand, by Theorem \ref{theorem.algmap}, the left hand route corresponds to the same map.  The assertion follows.
\end{proof}\index{bimodule Segre embedding!is functorial|)}

\section{$s$ is compatible with base change}\index{bimodule Segre embedding!is compatible with base change|(}
\begin{theorem} \label{theorem.nat}
Suppose
$$
\xymatrix{
Y' \ar[d]_{f'} \ar[r]^{p'} & Y \ar[d]^{f} \\
Z' \ar[r]_{p} & Z
}
$$
is a pullback diagram of schemes, with dual 2-cell $\Psi:p^{*}f_{*} \Longrightarrow f^{'}_{*}p^{'*}$.  Let $\mathcal{A}$ be an ${\mathcal{O}}_{Y}$-module which is rla \index{relatively locally affine}with respect to $f$.  Then $\Psi_{\mathcal{A}}$ is an isomorphism and the diagram
\begin{equation} \label{eqn.p}
\xymatrix{
{\mathbb{P}}_{Y'}({p'}^{*}\mathcal{A}) \ar[dd] \ar[r] & {\mathbb{P}}_{Y'}({f'}^{*}{f'}_{*}{p'}^{*}\mathcal{A}) \ar[r] & {\mathbb{P}}_{Z'}({f'}_{*}{p'}^{*}\mathcal{A}) \ar[d]^{{\mathbb{P}}_{Z'}(\Psi_{\mathcal{A}})} \\
& & {\mathbb{P}}_{Z'}({p}^{*}f_{*}\mathcal{A}) \ar[d] \\
{\mathbb{P}}_{Y}(\mathcal{A}) \ar[r] & {\mathbb{P}}_{Y}({f}^{*}f_{*}\mathcal{A}) \ar[r] & {\mathbb{P}}_{Z}({f}_{*}\mathcal{A})
}
\end{equation}
commutes.
\end{theorem}

\begin{proof}
Since $\mathcal{A}$ is rla with respect to $f$, $\Psi_{\mathcal{A}}$ is an isomorphism by Proposition \ref{prop.canon}.  We next show that (\ref{eqn.p}) commutes.  First, by Lemma \ref{lem.compat},
\begin{equation} \label{eqn.anothercompat}
\xymatrix{
{\mathbb{P}}_{Y'}({f'}^{*}p^{*}f_{*}\mathcal{A}) \ar[r] \ar[d] & {\mathbb{P}}_{Z'}(p^{*}f_{*}\mathcal{A}) \ar[d] \\
{\mathbb{P}}_{Y'}((f{p'})^{*}f_{*}\mathcal{A}) \ar[r] & {\mathbb{P}}_{Z}(f_{*}\mathcal{A}) \\
{\mathbb{P}}_{Y'}({p'}^{*}f^{*}f_{*}\mathcal{A}) \ar[u]^{\cong} \ar[r] & {\mathbb{P}}_{Y}(f^{*}f_{*}\mathcal{A}) \ar[u]
}
\end{equation}
commutes.  Thus, the diagram
\begin{equation} \label{eqn.midd}
\xymatrix{
{\mathbb{P}}_{Y'}({f'}^{*}p^{*}f_{*}\mathcal{A}) \ar[r] \ar[d] & {\mathbb{P}}_{Z'}(p^{*}f_{*}\mathcal{A}) \ar[d] \\
{\mathbb{P}}_{Y'}({p'}^{*}f^{*}f_{*}\mathcal{A}) \ar[r] & {\mathbb{P}}_{Z}(f_{*}\mathcal{A})
}
\end{equation}
whose left vertical is the left vertical of (\ref{eqn.anothercompat}) commutes.  We claim the diagram
$$
\xymatrix{
{\mathbb{P}}_{Y'}({p'}^{*}\mathcal{A}) \ar[r] \ar[d]_{=} & {\mathbb{P}}_{Y}(\mathcal{A}) \ar[r] & {\mathbb{P}}_{Y}(f^{*}f_{*}\mathcal{A}) \ar[d] \\
{\mathbb{P}}_{Y'}({p'}^{*}\mathcal{A}) \ar[r] \ar[dd]_{=} & {\mathbb{P}}_{Y'}({p'}^{*}f^{*}f_{*}\mathcal{A}) \ar[r] & {\mathbb{P}}_{Z}(f_{*}\mathcal{A}) \\
& {\mathbb{P}}_{Y'}({f'}^{*}p^{*}f_{*}\mathcal{A}) \ar[u] \ar[r] & {\mathbb{P}}_{Z'}(p^{*}f_{*}\mathcal{A}) \ar[u] \\
{\mathbb{P}}_{Y'}({p'}^{*}\mathcal{A}) \ar[r] & {\mathbb{P}}_{Y'}({f'}^{*}{f'}_{*}{p'}^{*}\mathcal{A}) \ar[r] \ar[u] & {\mathbb{P}}_{Z'}({f'}_{*}{p'}^{*}\mathcal{A}) \ar[u]
}
$$
commutes.  Since the outer circuit is (\ref{eqn.p}), our result will follow.  The lower left rectangle commutes by Proposition \ref{prop.adcom}.  The middle square equals (\ref{eqn.midd}), and the lower right square commutes by Lemma \ref{lem.proj}.  It remains to show the top rectangle commutes.  The top is a composition of the following subdiagrams
$$
\xymatrix{
{\mathbb{P}}_{Y'}({p'}^{*}\mathcal{A}) \ar[r] \ar[d]_{=} & {\mathbb{P}}_{Y}(\mathcal{A}) \ar[r] & {\mathbb{P}}_{Y}(f^{*}f_{*}\mathcal{A}) \ar[d] \\
{\mathbb{P}}_{Y'}({p'}^{*}\mathcal{A}) \ar[r]  & {\mathbb{P}}_{Y'}({p'}^{*}f^{*}f_{*}\mathcal{A}) \ar[r] \ar[ur] & {\mathbb{P}}_{Z}(f_{*}\mathcal{A}).
}
$$
where the diagonal morphism is (\ref{eqn.proj}) (defined on page \pageref{eqn.proj}).  The left subdiagram commutes by Lemma \ref{lem.proj}, while the right subdiagram commutes by Lemma \ref{lem.compat}.  The assertion follows.
\end{proof}

\begin{theorem} \label{theorem.basechanger}
$s$ is compatible with base change.  That is, if $p_{xy}$, $p_{yz}$, $p_{xz}$ and $p_{xyz}$ are overlap compatible, then the diagram
\begin{equation} \label{eqn.basechanger}
\xymatrix{
{\mathbb{P}}_{X' \times_{U} Y'}(p_{xy}^{*}\mathcal{E}) \otimes_{Y'} {\mathbb{P}}_{Y' \times_{U} Z'}(p_{yz}^{*}\mathcal{F}) \ar[dd] \ar[r]^{\hskip .3in s} & {\mathbb{P}}_{X' \times_{U} Z'}(p_{xy}^{*}\mathcal{E} \otimes_{{\mathcal{O}}_{Y'}} p_{yz}^{*}\mathcal{F}) \ar[d] \\
& {\mathbb{P}}_{X' \times_{U} Z'}(p_{xz}^{*}(\mathcal{E} \otimes_{{\mathcal{O}}_{Y}} \mathcal{F})) \ar[d] \\
{\mathbb{P}}_{X \times Y}(\mathcal{E}) \otimes_{Y} {\mathbb{P}}_{Y \times Z}(\mathcal{F})  \ar[r]_{s} & {\mathbb{P}}_{X \times Z}(\mathcal{E} \otimes_{{\mathcal{O}}_{Y}} \mathcal{F})
}
\end{equation}
whose left and bottom right vertical arrows are defined in Lemma \ref{lem.proj}, and whose top right vertical arrow is induced by the isomorphism defined in Lemma \ref{lem.genisommonoid}, commutes.
\end{theorem}

\begin{proof}
We construct two commutative diagrams which have the property that the composite diagram which one gets by placing these diagrams side by side has outer circuit equal to (\ref{eqn.basechanger}).  To simplify the presentation, we omit subscripts on the projectivization functors.  We claim Figure 1, on the following page,
\begin{sidewaysfigure}
$$
\xymatrix{
\mathbb{P}(p_{xy}^{*}\mathcal{E}) \otimes_{Y'} \mathbb{P}(p_{yz}^{*}\mathcal{F}) \ar[r]^{\hskip -.4in s_{1}} \ar[ddd] & \mathbb{P}(\pr_{12}^{'*}p_{xy}^{*}\mathcal{E}) \times_{X' \times_{U} Y' \times_{U} Z'} \mathbb{P}(\pr_{23}^{'*}p_{yz}^{*}\mathcal{F}) \ar[r]^{\hskip .4in s_{2}} \ar[d] & \mathbb{P}(\pr_{12}^{'*}p_{xy}^{*}\mathcal{E} \otimes \pr_{23}^{'*}p_{yz}^{*}\mathcal{F}) \ar[d] \\
& \mathbb{P}(p_{xyz}^{*}\pr_{12}^{*}\mathcal{E}) \times_{X' \times_{U} Y' \times_{U} Z'} \mathbb{P}(p_{xyz}^{*}\pr_{23}^{*}\mathcal{F}) \ar[r]^{\hskip .4in s_{2}} \ar[dd] & \mathbb{P}(p_{xyz}^{*}\pr_{12}^{*}\mathcal{E} \otimes p_{xyz}^{*}\pr_{23}^{*}\mathcal{F}) \ar[d] \\
& & \mathbb{P}(p_{xyz}^{*}(\pr_{12}^{*}\mathcal{E} \otimes \pr_{23}^{*}\mathcal{F})) \ar[d] \\
\mathbb{P}(\mathcal{E}) \otimes_{Y} \mathbb{P}(\mathcal{F}) \ar[r]_{\hskip -.4in s_{1}} & \mathbb{P}(\pr_{12}^{*}\mathcal{E}) \times_{X \times Y \times Z} \mathbb{P}(\pr_{23}^{*}\mathcal{F}) \ar[r]_{\hskip .4in s_{2}} & \mathbb{P}(\pr_{12}^{*}\mathcal{E} \otimes \pr_{23}^{*}\mathcal{F})
}
$$
\caption{}
\end{sidewaysfigure}
whose unlabeled maps are either induced by isomorphisms of the various modules being projectivized, or by maps (\ref{eqn.proj}) (defined on page \pageref{eqn.proj}), commutes.  For, the upper right square commutes by Proposition \ref{prop.segre} (2), while the lower right square commutes by Proposition \ref{prop.segre} (3).  We next show the left hand square commutes.  To show the left hand square commutes, it obviously suffices to show that the square
$$
\xymatrix{
\mathbb{P}(p_{xy}^{*}\mathcal{E}) \otimes_{Y'} \mathbb{P}(p_{yz}^{*}\mathcal{F})  \ar[dd] & \mathbb{P}(\pr_{12}^{'*}p_{xy}^{*}\mathcal{E}) \times_{X' \times_{U} Y' \times_{U} Z'} \mathbb{P}(\pr_{23}^{'*}p_{yz}^{*}\mathcal{F}) \ar[d] \ar[l]_{\hskip -.4in s_{1}^{-1}} \\
& \mathbb{P}(p_{xyz}^{*}\pr_{12}^{*}\mathcal{E}) \times_{X' \times_{U} Y' \times_{U} Z'} \mathbb{P}(p_{xyz}^{*}\pr_{23}^{*}\mathcal{F}) \ar[d] \\
\mathbb{P}(\mathcal{E}) \otimes_{Y} \mathbb{P}(\mathcal{F})  & \mathbb{P}(\pr_{12}^{*}\mathcal{E}) \times_{X \times Y \times Z} \mathbb{P}(\pr_{23}^{*}\mathcal{F}) \ar[l]^{\hskip -.4in s_{1}^{-1}}
}
$$
commutes.  Since each route of this diagram is a map into the fibre product $\mathbb{P}(\mathcal{E}) \otimes_{Y} \mathbb{P}(\mathcal{F})$, each route is determined by a pair of maps, one into ${\mathbb{P}}(\mathcal{E})$ and one into ${\mathbb{P}}(\mathcal{F})$.  Thus, by the description of $s_{1}^{-1}$ given in Proposition \ref{prop.firstmap}, it suffices to show the diagram
$$
\xymatrix{
\mathbb{P}(p_{xy}^{*}\mathcal{E}) \ar[ddd] & \mathbb{P}(\pr_{12}^{'*}p_{xy}^{*}\mathcal{E}) \ar[l] \ar[d]^{\cong} \\
& \mathbb{P}((\pr_{12}p_{xyz})^{*}\mathcal{E})= \mathbb{P}((\pr_{xy}p_{12}')^{*}\mathcal{E}) \ar[d] \ar[ldd] \\
& \mathbb{P}(p_{xyz}^{*}\pr_{12}^{*}\mathcal{E}) \ar[d] \\
\mathbb{P}(\mathcal{E}) & \mathbb{P}(\pr_{12}^{*}\mathcal{E}) \ar[l]
}
$$
whose verticals are maps of the form (\ref{eqn.proj}) and whose diagonal is of the form (\ref{eqn.proj}), commutes (a diagram like that above containing ${\mathbb{P}}(\mathcal{F})$ commutes in a similar manner).  But the upper and lower circuits commute by Lemma \ref{lem.compat}.

Let $\epsilon$ be the counit of the adjoint pair $(\pr_{13}^{*},\pr_{13*})$ and let $\epsilon'$ be the counit of the adjoint pair $(\pr_{13}^{'*},{\pr_{13}'}_{*})$.  We claim that Figure 2,
\begin{sidewaysfigure}
$$
\xymatrix{
\mathbb{P}(\pr_{12}^{'*}p_{xy}^{*}\mathcal{E} \otimes \pr_{23}^{'*}p_{yz}^{*}\mathcal{F}) \ar[r]^{\hskip -.3in \mathbb{P}(\epsilon')} \ar[d] & \mathbb{P}(\pr_{13}^{'*}{\pr_{13}'}_{*}(\pr_{12}^{'*}p_{xy}^{*}\mathcal{E} \otimes \pr_{23}^{'*}p_{yz}^{*}\mathcal{F})) \ar[d] \ar[r] & \mathbb{P}({\pr_{13}'}_{*}\pr_{12}^{'*}p_{xy}^{*}\mathcal{E} \otimes \pr_{23}^{'*}p_{yz}^{*}\mathcal{F}) \ar[d] \\
\mathbb{P}(p_{xyz}^{*}\pr_{12}^{*}\mathcal{E} \otimes p_{xyz}^{*}\pr_{23}^{*}\mathcal{F}) \ar[d] &  \mathbb{P}({\pr_{13}^{'*}}{\pr_{13}'}_{*}(p_{xyz}^{*}\pr_{12}^{*}\mathcal{E} \otimes p_{xyz}^{*}\pr_{23}^{*}\mathcal{F})) \ar[d]  &  \mathbb{P}({\pr_{13}'}_{*}(p_{xyz}^{*}\pr_{12}^{*}\mathcal{E} \otimes p_{xyz}^{*}\pr_{23}^{*}\mathcal{F})) \ar[d] \\
\mathbb{P}(p_{xyz}^{*}(\pr_{12}^{*}\mathcal{E} \otimes \pr_{23}^{*}\mathcal{F})) \ar[r]^{\hskip -.3in \mathbb{P}(\epsilon')} \ar[dd] & \mathbb{P}(\pr_{13}^{'*}{\pr_{13}'}_{*}p_{xyz}^{*}(\pr_{12}^{*}\mathcal{E} \otimes \pr_{23}^{*}\mathcal{F})) \ar[r] & \mathbb{P}({\pr_{13}'}_{*}p_{xyz}^{*}(\pr_{12}^{*}\mathcal{E} \otimes \pr_{23}^{*}\mathcal{F})) \ar[d]^{\cong} \\
& & \mathbb{P}(p_{xz}^{*}(\mathcal{E} \otimes_{{\mathcal{O}}_{Y}} \mathcal{F})) \ar[d] \\
\mathbb{P}(\pr_{12}^{*}\mathcal{E} \otimes \pr_{23}^{*}\mathcal{F}) \ar[r]_{\hskip -.3in \mathbb{P}(\epsilon)} & \mathbb{P}(\pr_{13}^{*}\pr_{13*}(\pr_{12}^{*}\mathcal{E} \otimes \pr_{23}^{*}\mathcal{F})) \ar[r] & \mathbb{P}(\mathcal{E} \otimes_{{\mathcal{O}}_{Y}} \mathcal{F})
}
$$
\caption{}
\end{sidewaysfigure}
whose third right vertical is induced by the isomorphism of Lemma \ref{lem.genisommonoid}, and whose other unlabeled arrows are either induced by isomorphism of the various projectivized modules or are from Lemma \ref{lem.proj}, commutes.  For, the bottom rectangle commutes by Theorem \ref{theorem.nat}, the upper right square commutes by Lemma \ref{lem.proj}, and the upper left square commutes by the functoriality of $\mathbb{P}(-)$.  The theorem follows by noticing that if we put Figure 1 to the left of Figure 2, the outer circuit of this combined diagram is just (\ref{eqn.basechanger}).
\end{proof}\index{bimodule Segre embedding!is compatible with base change|)}

\section{$s$ is associative}\index{bimodule Segre embedding!is associative|(}
We show $s$ is associative.  Our proof employs the fact that closed immersions induce well behaved adjoints (Corollary \ref{cor.groth}).

\subsection{The proof that $s$ is associative}

\begin{lemma} \label{lem.diamond}
Let $A$, $B$, $C$, and $D$ be commutative rings such that there is a commutative diagram
\begin{equation} \label{eqn.diamond}
\xymatrix{
& A \ar[dl] \ar[d] \ar[dr] & \\
B \ar[dr] \ar[r] & C \ar[d] & D \ar[l] \ar[dl] \\
& A &
}
\end{equation}
whose verticals are identities.  Let
$$
\phi:(A \otimes_{B} C) \otimes_{C} (A \otimes_{D} C) \rightarrow A \otimes_{C}A
$$
be the $C$-module map defined by $\phi((a \otimes c) \otimes (a' \otimes c')) = c \cdot a \otimes c' \cdot a'$.  If every simple tensor of $(A \otimes_{B} C) \otimes_{C} (A \otimes_{D} C)$ may be written as $(a \otimes 1) \otimes (a' \otimes 1)$ for $a,a' \in A$ then $\phi$ is an isomorphism.
\end{lemma}

\begin{proof}
Since $\phi$ is obviously surjective, we need only show $\phi$ is injective.  Suppose, for $a_{i}, a_{i}'' \in A$ such that $\{a_{i}''\}$ generates $A$ as a $C$-module, that $\Sigma_{i}a_{i} \otimes a_{i}'' \in A \otimes_{C}A$ equals $0$.  By hypothesis, we need only show that $\Sigma_{i}(a_{i} \otimes 1) \otimes (a_{i}'' \otimes 1)=0$.  By Lemma \ref{lem.eisenbud}, there exist $r_{ij} \in C$ and $a_{j}' \in A$ such that
$$
\Sigma_{j}r_{ij} \cdot a_{j}'=a_{i}
$$
and
$$
\Sigma_{i}r_{ij} \cdot a_{i}''=0.
$$
We find
\begin{align*}
\Sigma_{i}(a_{i} \otimes 1) \otimes (a_{i}'' \otimes 1) & =\Sigma_{i}(\Sigma_{j}r_{ij} \cdot a_{j}' \otimes 1) \otimes (a_{i}'' \otimes 1) \\
& = \Sigma_{j}(\Sigma_{i}((r_{ij} \cdot a_{j}' \otimes 1) \otimes (a_{i}'' \otimes 1)) \\
& = \Sigma_{j}(a_{j}' \otimes \Sigma_{i}r_{ij} \cdot 1) \otimes (a_{i}'' \otimes 1) \\
& = \Sigma_{j}(a_{j}' \otimes 1) \otimes (\Sigma_{i}r_{ij} \cdot a_{i}'' \otimes 1) \\
& = 0.
\end{align*}
where we have used the commutativity of (\ref{eqn.diamond}) as well as the assumption that all vertical routes are the identity.
\end{proof}
Let $W$ be a scheme and let $\pr_{ij}^{W}:(W \times X)\times_{W}(W \times Y)\times_{W}(W \times Z) \rightarrow W \times X,W \times Y, W \times Z$ be the standard projections.
\begin{proposition} \label{prop.cruciso}
Let ${\mathcal{L}}_{1}$ and ${\mathcal{L}}_{2}$ be invertible ${\mathcal{O}}_{W}$-modules, and let $\Omega_{1}$ and $\Omega_{2}$ be relative versions of the maps given in (\ref{eqn.alphaone}) and (\ref{eqn.alphatwo}) respectively.  Then the map
\begin{equation} \label{eqn.map}
{\Omega_{1}}_{{\mathcal{L}}_{1}} \otimes {\Omega_{2}}_{{\mathcal{L}}_{2}}:{\pr_{12}}^{W*}{q_{12}^{W}}_{*}{\mathcal{L}}_{1} \otimes {\pr_{23}}^{W*}{q_{23}^{W}}_{*}{\mathcal{L}}_{2} \rightarrow
\end{equation}
$$
{q_{123}^{W}}_{*}{\mathcal{L}}_{1} \otimes  {q_{123}^{W}}_{*}{\mathcal{L}}_{2}
$$
is an isomorphism.
\end{proposition}

\begin{proof}
We check the assertion locally.  Let $p$ be a closed point in the support of the domain of (\ref{eqn.map}).  Then, since $q_{123}^{W}$ is a closed immersion, there exists a unique $w \in W$ such that $q_{123}^{W}(w)=p$.  Let $p_{1}=q_{12}^{W}(w)$ and let $p_{2}=q_{23}^{W}(w)$.  Then $({\Omega_{1}}_{{\mathcal{L}}_{1}} \otimes {\Omega_{2}}_{{\mathcal{L}}_{2}})_{p}$ is a map
$$
({\mathcal{L}}_{1w} \otimes_{{\mathcal{O}}_{p_{1}}} {\mathcal{O}}_{p}) \otimes_{{\mathcal{O}}_{p}} ({\mathcal{L}}_{2w} \otimes_{{\mathcal{O}}_{p_{2}}} {\mathcal{O}}_{p}) \rightarrow {\mathcal{L}}_{1w} \otimes_{{\mathcal{O}}_{p}} {\mathcal{L}}_{2w}
$$
sending the simple tensor $(a \otimes c) \otimes (a' \otimes c')$ to $c \cdot a \otimes c' \cdot a'$ for $a \in {\mathcal{L}}_{1w}$, $a' \in {\mathcal{L}}_{2w}$ and $c,c' \in {\mathcal{O}}_{p}$.  Since ${\mathcal{L}}_{1}$ and ${\mathcal{L}}_{2}$, are invertible ${\mathcal{L}}_{1w} \cong {\mathcal{L}}_{2w} \cong {\mathcal{O}}_{w}$.  In addition,
$$
{\mathcal{O}}_{p}={\mathcal{O}}_{q_{1}^{W}(w)} \otimes_{{\mathcal{O}}_{w}} {\mathcal{O}}_{q_{2}^{W}(w)} \otimes_{{\mathcal{O}}_{w}} {\mathcal{O}}_{q_{3}^{W}(w)},
$$
$$
{\mathcal{O}}_{p_{1}}={\mathcal{O}}_{q_{1}^{W}(w)} \otimes_{{\mathcal{O}}_{w}} {\mathcal{O}}_{q_{2}^{W}(w)},
$$
and
$$
{\mathcal{O}}_{p_{2}}={\mathcal{O}}_{q_{2}^{W}(w)} \otimes_{{\mathcal{O}}_{w}} {\mathcal{O}}_{q_{3}^{W}(w)}.
$$
The assertion now follows from Lemma \ref{lem.diamond} with $A = {\mathcal{O}}_{w}$, $B = {\mathcal{O}}_{p_{1}}$, $C = {\mathcal{O}}_{p}$ and $D = {\mathcal{O}}_{p_{2}}$.
\end{proof}

\begin{proposition} \label{prop.segretensor}
Let $W$ be a scheme, and let
$$
p_{xy}, p_{zy}:(W \times X) \times_{W} (W \times Y) \times (W \times Z) \rightarrow
$$
$$
X \times Y,Y \times Z
$$

be projections.  Suppose
$$
r^{W}:W \rightarrow {\mathbb{P}}_{(W \times X) \times_{W} (W \times Y)}(p_{xy}^{*}\mathcal{E}) \otimes_{W \times Y} {\mathbb{P}}_{(W \times Y) \times_{W} (W \times Z)}(p_{yz}^{*}\mathcal{F})
$$
has projection to $(W \times X) \times_{W} (W \times Y)$, $q_{12}^{W}$ and projection to $(W \times Y) \times_{W} (W \times Z)$, $q_{23}^{W}$.  If $r^{W}$ corresponds to the pair of epimorphisms
$$
\psi_{1}:p_{xy}^{*}\mathcal{E} \rightarrow {q_{12}^{W}}_{*}{\mathcal{L}}_{1}
$$
and
$$
\psi_{2}:p_{yz}^{*}\mathcal{F} \rightarrow {q_{23}^{W}}_{*}{\mathcal{L}}_{2},
$$
the composition $s \circ r^{W}$ corresponds to the epimorphism
\begin{equation} \label{eqn.shorter}
\psi_{1} \otimes_{{\mathcal{O}}_{Y}} \psi_{2}.
\end{equation}
\end{proposition}

\begin{proof}
By Proposition \ref{prop.cruciso}, $\Omega_{1{\mathcal{L}}_{1}} \otimes \Omega_{2{\mathcal{L}}_{2}}$ is an isomorphism.  In addition, since $q_{123}^{W}$ is a closed immersion, $\Upsilon$ ((\ref{eqn.betaone}) on page \pageref{eqn.betaone}) is also an isomorphism.  By Corollary \ref{cor.groth}, the maps (\ref{eqn.shorter}) and (\ref{eqn.bigsegre}), defined on page \pageref{eqn.bigsegre}, correspond to the same map $r^{W}$.
\end{proof}

\begin{theorem}
$s$ is associative:  if $T$ is a scheme, $\mathcal{E}$ is an ${\mathcal{O}}_{X \times Y}$-bimodule, $\mathcal{F}$ is an ${\mathcal{O}}_{Y \times Z}$-bimodule and $\mathcal{G}$ is an ${\mathcal{O}}_{Z \times T}$-bimodule, then the diagram
\begin{equation} \label{eqn.ass}
\xymatrix{
({\mathbb{P}}_{X \times Y}(\mathcal{E}) \otimes_{Y} {\mathbb{P}}_{Y \times Z}(\mathcal{F})) \otimes_{Z} {\mathbb{P}}_{Z \times T}(\mathcal{G}) \ar[r]^{\cong} \ar[d]_{s \times \id} & {\mathbb{P}}_{X \times Y}(\mathcal{E}) \otimes_{Y} ({\mathbb{P}}_{Y \times Z}(\mathcal{F}) \otimes_{Z} {\mathbb{P}}_{Z \times T}(\mathcal{G})) \ar[d]^{\id \times s} \\
{\mathbb{P}}_{X \times Z}(\mathcal{E}\otimes_{{\mathcal{O}}_{Y}} \mathcal{F}) \otimes_{Z} {\mathbb{P}}_{Z \times T}(\mathcal{G})  \ar[d]_{s} & {\mathbb{P}}_{X \times Y}(\mathcal{E}) \otimes_{Y} {\mathbb{P}}_{Y \times Z}(\mathcal{F} \otimes_{{\mathcal{O}}_{Z}} \mathcal{G})  \ar[d]^{s} \\
{\mathbb{P}}_{X \times T}((\mathcal{E}\otimes_{{\mathcal{O}}_{Y}} \mathcal{F})\otimes_{{\mathcal{O}}_{Z}} \mathcal{G}) \ar[r]_{\cong} & {\mathbb{P}}_{X \times T}(\mathcal{E} \otimes_{{\mathcal{O}}_{Y}}(\mathcal{F} \otimes_{{\mathcal{O}}_{Z}} \mathcal{G}))
}
\end{equation}
whose bottom row is induced by the associativity isomorphism (\ref{prop.tensor}), commutes.
\end{theorem}

\begin{proof}
To clarify exposition, we will sometimes drop subscripts from projectivizations.  Let $p_{xy}$ and $p_{yz}$ be defined as in Proposition \ref{prop.segretensor}.  Let $p_{zt}$ be defined similarly.  Suppose $W$ is an affine scheme and define $X' = W \times X$, $Y' = W \times Y$ and $Z' = W \times Z$.  Let $T' = W \times T$ and suppose
$$
r:W \rightarrow ({\mathbb{P}}_{X \times Y}(\mathcal{E}) \otimes_{Y} {\mathbb{P}}_{Y \times Z}(\mathcal{F})) \otimes_{Z} {\mathbb{P}}_{Z \times T}(\mathcal{G})
$$
is a morphism with projection $q_{12}$ to $X \times Y$, projection $q_{23}$ to $Y \times Z$ and projection $q_{34}$ to $Z \times T$.  To prove the Theorem, it suffices to show that the two routes of (\ref{eqn.ass}) precomposed with $r$ commute, since $W$ is an arbitrary affine scheme.  Suppose $r$ corresponds to epis
$$
\phi_{1}:q_{12}^{*}\mathcal{E} \rightarrow {\mathcal{L}}_{1}
$$
$$
\phi_{2}:q_{12}^{*}\mathcal{E} \rightarrow {\mathcal{L}}_{2}
$$
and
$$
\phi_{3}:q_{12}^{*}\mathcal{E} \rightarrow {\mathcal{L}}_{3}.
$$
Then there are epis
$$
\xymatrix{
\gamma_{1}:q_{12}^{W*}p_{xy}^{*}\mathcal{E} \ar[r] & q_{12}^{*}\mathcal{E} \ar[r]^{\phi_{1}} & {\mathcal{L}}_{1}
}
$$
$$
\xymatrix{
\gamma_{2}:q_{23}^{W*}p_{yz}^{*}\mathcal{F} \ar[r] & q_{23}^{*}\mathcal{F} \ar[r]^{\phi_{2}} & {\mathcal{L}}_{2}
}
$$
$$
\xymatrix{
\gamma_{3}:q_{34}^{W*}p_{zw}^{*}\mathcal{G} \ar[r] & q_{34}^{*}\mathcal{G} \ar[r]^{\phi_{3}} & {\mathcal{L}}_{3}.
}
$$
These epis give a map
$$
r^{W}:W \rightarrow ({\mathbb{P}}_{X' \times_{W} Y'}(p_{xy}^{*}\mathcal{E}) \otimes_{Y'} {\mathbb{P}}_{Y' \times_{W} Z'}(p_{yz}^{*}\mathcal{F})) \otimes_{Z'} {\mathbb{P}}_{Z' \times_{W} T'}(p_{zt}^{*}\mathcal{G})
$$
such that if the vertical map in the diagram
$$
\xymatrix{
W \ar[r] \ar[dr]_{r} & ({\mathbb{P}}_{X' \times_{W} Y'}(p_{xy}^{*}\mathcal{E}) \otimes_{Y'} {\mathbb{P}}_{Y' \times_{W} Z'}(p_{yz}^{*}\mathcal{F})) \otimes_{Z'} {\mathbb{P}}_{Z' \times_{W} T'}(p_{zt}^{*}\mathcal{G}) \ar[d] \\
& ({\mathbb{P}}_{X \times Y}(\mathcal{E}) \otimes_{Y} {\mathbb{P}}_{Y \times Z}(\mathcal{F})) \otimes_{Z} {\mathbb{P}}_{Z \times T}(\mathcal{G})
}
$$
is that constructed in Lemma \ref{lem.proj}, and the horizontal is $r^{W}$, then this diagram commutes (Lemma \ref{lem.algproj}).

We claim the outer circuit of the diagram
\begin{equation} \label{eqn.almostass}
\xymatrix{
W \ar[d]_{r^{W}} \ar[r]^{=} & W \ar[d] \\
(\mathbb{P}(p_{xy}^{*}\mathcal{E}) \otimes_{Y'} \mathbb{P}(p_{yz}^{*}\mathcal{F})) \otimes_{Z'} \mathbb{P}(p_{zt}^{*}\mathcal{G}) \ar[r]^{\cong} \ar[d]_{s \times \id} & \mathbb{P}(p_{xy}^{*}\mathcal{E}) \otimes_{Y'} (\mathbb{P}(p_{yz}^{*}\mathcal{F}) \otimes_{Z'} \mathbb{P}(p_{zt}^{*}\mathcal{G})) \ar[d]^{\id \times s} \\
\mathbb{P}(p_{xy}^{*}\mathcal{E}\otimes_{{\mathcal{O}}_{Y'}} p_{yz}^{*}\mathcal{F}) \otimes_{Z'} \mathbb{P}(p_{zt}^{*}\mathcal{G})  \ar[d]_{s} & \mathbb{P}(p_{xy}^{*}\mathcal{E}) \otimes_{Y'} \mathbb{P}(p_{yz}^{*}\mathcal{F} \otimes_{{\mathcal{O}}_{Z'}} p_{zt}^{*}\mathcal{G})  \ar[d]^{s} \\
\mathbb{P}((p_{xy}^{*}\mathcal{E}\otimes_{{\mathcal{O}}_{Y'}} p_{yz}^{*}\mathcal{F})\otimes_{{\mathcal{O}}_{Z'}} p_{zt}^{*}\mathcal{G}) \ar[r]_{\cong} & \mathbb{P}(p_{xy}^{*}\mathcal{E} \otimes_{{\mathcal{O}}_{Y'}}(p_{yz}^{*}\mathcal{F} \otimes_{{\mathcal{O}}_{Z'}} p_{zt}^{*}\mathcal{G}))
}
\end{equation}
whose upper right vertical is induced by $\gamma_{1}$, $\gamma_{2}$ and $\gamma_{3}$, commutes.  Let $\psi_{i}$ be the right adjunct of $\gamma_{i}$.  By Proposition \ref{prop.segretensor}, the left hand circuit corresponds to the left hand circuit of the commutative diagram
\begin{equation} \label{eqn.smallass}
\xymatrix{
(p_{xy}^{*}\mathcal{E} \otimes_{{\mathcal{O}}_{Y'}} p_{yz}^{*}\mathcal{F}) \otimes_{{\mathcal{O}}_{Z'}}p_{zw}^{*}\mathcal{G} \ar[d]_{(\psi_{1}' \otimes_{{\mathcal{O}}_{Y'}} \psi_{2}') \otimes_{{\mathcal{O}}_{Z'}} \psi_{3}'}  & p_{xy}^{*}\mathcal{E} \otimes_{{\mathcal{O}}_{Y'}} (p_{yz}^{*}\mathcal{F} \otimes_{{\mathcal{O}}_{Z'}}p_{zw}^{*}\mathcal{G}) \ar[d]^{\psi_{1}' \otimes_{{\mathcal{O}}_{Y'}} (\psi_{2}' \otimes_{{\mathcal{O}}_{Z'}} \psi_{3}')} \ar[l]_{\cong} \\
(q_{12*}^{W}{\mathcal{L}}_{1} \otimes_{{\mathcal{O}}_{Y'}} q_{23*}^{W}{\mathcal{L}}_{2}) \otimes_{{\mathcal{O}}_{Z'}}q_{34*}^{W}{\mathcal{L}}_{3}
& q_{12*}^{W}{\mathcal{L}}_{1} \otimes_{{\mathcal{O}}_{Y'}} (q_{23*}^{W}{\mathcal{L}}_{2} \otimes_{{\mathcal{O}}_{Z'}}q_{34*}^{W}{\mathcal{L}}_{3}) \ar[l]^{\cong}
}
\end{equation}
while the right hand vertical of (\ref{eqn.almostass}) corresponds to the right hand vertical of (\ref{eqn.smallass}).  By the commutativity of (\ref{eqn.smallass}), both the left hand circuit and the right vertical of (\ref{eqn.smallass}) have the same kernel.  Thus, by Corollary \ref{cor.kernel}, (\ref{eqn.almostass}) commutes.

To show (\ref{eqn.ass}) commutes, it suffices to construct a cube of schemes whose top is (\ref{eqn.almostass}), whose bottom is (\ref{eqn.ass}), whose edge
$$
\xymatrix{
({\mathbb{P}}_{X' \times_{W} Y'}(p_{xy}^{*}\mathcal{E}) \otimes_{Y'} {\mathbb{P}}_{Y' \times_{W} Z'}(p_{yz}^{*}\mathcal{F})) \otimes_{Z'} {\mathbb{P}}_{Z' \times_{W} T'}(p_{zt}^{*}\mathcal{G}) \ar[d] \\
({\mathbb{P}}_{X \times Y}(\mathcal{E}) \otimes_{Y} {\mathbb{P}}_{Y \times Z}(\mathcal{F})) \otimes_{Z} {\mathbb{P}}_{Z \times T}(\mathcal{G})
}
$$
is the map constructed in Lemma \ref{lem.proj} and whose sides commute.  With the orientation fixed so that (\ref{eqn.ass}) forms the bottom of a cube, let the upper left edge into the page be the map above and construct the upper right edge into the page similarly.  Let the lower left edge into the page be induced by the composition of isomorphisms
$$
\xymatrix{
p_{xt}^{*}((\mathcal{E}\otimes_{{\mathcal{O}}_{Y}} \mathcal{F})\otimes_{{\mathcal{O}}_{Z}} \mathcal{G}) \ar[d] \\
p_{xz}^{*}(\mathcal{E}\otimes_{{\mathcal{O}}_{Y}} \mathcal{F})\otimes_{{\mathcal{O}}_{Z'}} p_{zt}^{*}\mathcal{G} \ar[d] \\
(p_{xy}^{*}\mathcal{E}\otimes_{{\mathcal{O}}_{Y'}} p_{yz}^{*}\mathcal{F})\otimes_{{\mathcal{O}}_{Z'}} p_{zt}^{*}\mathcal{G}
}
$$
(Lemma \ref{lem.isommonoid}), followed by the map
$$
{\mathbb{P}}(p_{xt}^{*}((\mathcal{E}\otimes_{{\mathcal{O}}_{Y}} \mathcal{F})\otimes_{{\mathcal{O}}_{Z}} \mathcal{G})) \rightarrow {\mathbb{P}}((\mathcal{E}\otimes_{{\mathcal{O}}_{Y}} \mathcal{F})\otimes_{{\mathcal{O}}_{Z}} \mathcal{G})
$$
from Lemma \ref{lem.proj}, and let the lower right edge be constructed similarly.  The left and right sides of this cube commute since $s$ is compatible with base change (Theorem \ref{theorem.basechanger}) and since  $s$ is functorial (Theorem \ref{theorem.functorial}).  It is easy to see the back commutes by examining the algebraic description of the back edges.  To complete the proof of the assertion, we need only show that the front face of the cube, the diagram
$$
\xymatrix{
\mathbb{P}((p_{xy}^{*}\mathcal{E}\otimes_{{\mathcal{O}}_{Y'}} p_{yz}^{*}\mathcal{F})\otimes_{{\mathcal{O}}_{Z'}} p_{zt}^{*}\mathcal{G}) \ar[r]^{\cong} \ar[d] & \mathbb{P}(p_{xy}^{*}\mathcal{E} \otimes_{{\mathcal{O}}_{Y'}}(p_{yz}^{*}\mathcal{F} \otimes_{{\mathcal{O}}_{Z'}} p_{zt}^{*}\mathcal{G})) \ar[d] \\
\mathbb{P}((p_{xz}^{*}(\mathcal{E}\otimes_{{\mathcal{O}}_{Y}} \mathcal{F})\otimes_{{\mathcal{O}}_{Z'}} p_{zt}^{*}\mathcal{G}))  \ar[d] & \mathbb{P}(p_{xy}^{*}\mathcal{E} \otimes_{{\mathcal{O}}_{Y'}}(p_{yt}^{*}(\mathcal{F} \otimes_{{\mathcal{O}}_{Z}}\mathcal{G}))) \ar[d] \\
\mathbb{P}(p_{xt}^{*}((\mathcal{E}\otimes_{{\mathcal{O}}_{Y}} \mathcal{F})\otimes_{{\mathcal{O}}_{Z}} \mathcal{G}))  \ar[d] \ar[r]^{\cong} & \mathbb{P}(p_{xt}^{*}(\mathcal{E} \otimes_{{\mathcal{O}}_{Y}}(\mathcal{F} \otimes_{{\mathcal{O}}_{Z}}\mathcal{G}))) \ar[d] \\
\mathbb{P}((\mathcal{E}\otimes_{{\mathcal{O}}_{Y}} \mathcal{F})\otimes_{{\mathcal{O}}_{Z}} \mathcal{G}) \ar[r]_{\cong} & \mathbb{P}(\mathcal{E} \otimes_{{\mathcal{O}}_{Y}}(\mathcal{F} \otimes_{{\mathcal{O}}_{Z}} \mathcal{G}))
}
$$
whose third verticals are from Lemma \ref{lem.proj}, commutes.  By Lemma \ref{lem.proj}, the bottom square commutes, so to complete the proof of the theorem, we need only show that the diagram
$$
\xymatrix{
(p_{xy}^{*}\mathcal{E}\otimes_{{\mathcal{O}}_{Y'}} p_{yz}^{*}\mathcal{F})\otimes_{{\mathcal{O}}_{Z'}} p_{zt}^{*}\mathcal{G} \ar[r]^{\cong} \ar[d] & p_{xy}^{*}\mathcal{E} \otimes_{{\mathcal{O}}_{Y'}}(p_{yz}^{*}\mathcal{F} \otimes_{{\mathcal{O}}_{Z'}} p_{zt}^{*}\mathcal{G}) \ar[d] \\
p_{xz}^{*}(\mathcal{E}\otimes_{{\mathcal{O}}_{Y}} \mathcal{F})\otimes_{{\mathcal{O}}_{Z'}} p_{zt}^{*}\mathcal{G} \ar[d] & p_{xy}^{*}\mathcal{E} \otimes_{{\mathcal{O}}_{Y'}}p_{yt}^{*}(\mathcal{F} \otimes_{{\mathcal{O}}_{Z}}\mathcal{G}) \ar[d] \\
p_{xt}^{*}((\mathcal{E}\otimes_{{\mathcal{O}}_{Y}} \mathcal{F})\otimes_{{\mathcal{O}}_{Z}} \mathcal{G}) \ar[r]_{\cong} & p_{xt}^{*}(\mathcal{E} \otimes_{{\mathcal{O}}_{Y}}(\mathcal{F} \otimes_{{\mathcal{O}}_{Z}}\mathcal{G}))  \\
}
$$
whose verticals are the isomorphisms constructed in Lemma \ref{lem.isommonoid}, commutes.  But this just follows from the fact that the associativity of bimodules is an indexed natural transformation (by Proposition \ref{prop.assindex}).
\end{proof}\index{bimodule Segre embedding|)}\index{bimodule Segre embedding!is associative|)}

\chapter{The Representation of $\Gamma_{n}$ for High
$n$}\index{family!free, of truncated point modules|(} We utilize
the following notation throughout this chapter: suppose $U$ is an
affine, noetherian scheme with structure map $f:U \rightarrow S$,
$\tilde{f}=f \times \id_{X}:U \times X \rightarrow S \times X$,
$X$ is a separated, noetherian scheme, $\mathcal{E}$ is a coherent
${\mathcal{O}}_{X}$-bimodule, $\mathcal{I}$ is a graded ideal of
$T(\mathcal{E})$ and $\mathcal{B}=T(\mathcal{E})/\mathcal{I}$.  If
$\mathcal{G}$ is an ${\mathcal{O}}_{X^{2}}$-module, then we let
${\mathcal{G}}^{U}=\tilde{f}^{2*}\mathcal{G}$.  Assume all
unadorned tensor products are bimodule tensor products.  If
$\mathcal{M}$ is a free truncated $U$-family of length $n+1$
(Definition \ref{def.freefamily}), it is also a free truncated
$U$-family of length $n+1$ over $T(\mathcal{E})$.  We denote
$\mathcal{M}$'s multiplication over $T(\mathcal{E})$ by
$\mu_{i,j}$ and let its components be denoted by
${\mathcal{M}}_{i} \cong (\id_{U} \times
q_{i})_{*}{\mathcal{O}}_{U}$.  As before, we let $d:U \rightarrow
U \times U$ denote the diagonal morphism and we let $q_{ij} =
((\id_{U} \times q_{i}) \times (\id_{U} \times q_{j}))\circ (d
\times d)\circ d:U \rightarrow (U \times X)_{U}^{2}$.  For $0 \leq
k \leq n$, let
$$
i_{k}:\mbox{SSupp }\pr_{1}^{*}{\mathcal{M}}_{k} \rightarrow (U \times X)_{U}^{2}
$$
be inclusion and let

$$
F_{k}= (\pr_{2}i_{k})_{*} i_{k}^{*} (\pr_{1}^{*}{\mathcal{M}}_{k} \otimes_{{\mathcal{O}}_{(U \times X)_{U}^{2}}} - )
$$
and
$$
G_{k}= {\mathcal{H}}{\it om}_{{\mathcal{O}}_{(U \times X)_{U}^{2}}}(\pr_{1}^{*}{\mathcal{M}}_{k},-) i_{k*} (\pr_{2}i_{k})^{!}.
$$
Then $(F_{k},G_{k})$ is an adjoint pair by Lemma \ref{lem.unit}.  Finally, for $i \geq 0$, $j>0$, let the right adjunct of $\mu_{i,j} \in \operatorname{Hom }_{(U \times X)_{U}^{2}}(F_{i}{\mathcal{E}}^{U \otimes j},{\mathcal{M}}_{i+j})$ be $\operatorname{Ad }\mu_{i,j} \in \operatorname{Hom }_{(U \times X)_{U}^{2}}({\mathcal{E}}^{U \otimes j},G_{i}{\mathcal{M}}_{i+j})$.

We define the bimodule Segre embedding
$s:{\mathbb{P}}_{X^{2}}(\mathcal{E})^{\otimes 1} \rightarrow
{\mathbb{P}}_{X^{2}}({\mathcal{E}}^{\otimes 1})$ as the identity
map.  Our goal in this chapter is to prove the following theorem.
\begin{theorem} \label{theorem.bigone}
For $n \geq 1$, $\Gamma_{n}$ is represented by the pullback of the diagram
\begin{equation} \label{eqn.thisisthelastone}
\xymatrix{
& {\mathbb{P}}_{X^{2}}(\mathcal{E})^{\otimes n} \ar[d]^{s} \\
{\mathbb{P}}_{X^{2}}({\mathcal{E}}^{\otimes n}/{\mathcal{I}}_{n}) \ar[r] & {\mathbb{P}}_{X^{2}}({\mathcal{E}}^{\otimes n}).
}
\end{equation}
\end{theorem}
To prove this theorem, we begin with a technical lemma (Lemma \ref{prop.main}).  Then, in Proposition \ref{prop.injy} and Corollary \ref{cor.newinjy}, we construct a monomorphism
$$
\Upsilon:\Gamma_{n} \Longrightarrow \operatorname{Hom}_{S}(-,P)
$$
where $P$ is the pullback of the diagram (\ref{eqn.thisisthelastone}).  Finally, in Proposition \ref{prop.lastfinal}, we prove $\Upsilon$ is an epimorphism, and hence, an equivalence.
\begin{lemma} \label{prop.main}
With the above notation,
$$
\operatorname{ker }(\operatorname{Ad }\mu_{i,j})=\operatorname{ker }(\operatorname{Ad }\mu_{i,1} \otimes \operatorname{Ad }\mu_{i+1,j-1}).
$$
\end{lemma}

\begin{proof}
If the adjoint pair $(F_{k},G_{k})$ has counit $\eta_{k}$, then, to prove the Lemma, we must show that for $i,j \geq 0$, the two routes of

\begin{equation} \label{eqn.main}
\xymatrix{
{\mathcal{E}}^{U} \otimes {\mathcal{E}}^{U \otimes j-1} \ar[rr]^{=} \ar[d]_{\eta_{i{\mathcal{E}}^{U} \otimes {\mathcal{E}}^{U \otimes j-1}}} & & {\mathcal{E}}^{U} \otimes {\mathcal{E}}^{U \otimes j-1} \ar[d]^{\eta_{i{\mathcal{E}}^{U}} \otimes \eta_{i+1{\mathcal{E}}^{U \otimes j-1}}} \\
G_{i}F_{i}({\mathcal{E}}^{U} \otimes {\mathcal{E}}^{U \otimes j-1}) \ar[d]_{G_{i}(\mu_{i,1}\otimes{\mathcal{E}}^{U \otimes j-1})}  & & G_{i}F_{i}{\mathcal{E}}^{U} \otimes G_{i+1}F_{i+1}{\mathcal{E}}^{U \otimes j-1} \ar[d]^{G_{i}\mu_{i,1} \otimes G_{i+1}\mu_{i+1,j-1}} \\
G_{i}F_{i+1}{\mathcal{E}}^{U \otimes j-1} \ar[d]_{G_{i}\mu_{i+1,j-1}} & & G_{i}{\mathcal{M}}_{i+1} \otimes G_{i+1}{\mathcal{M}}_{i+j} \\
G_{i}{\mathcal{M}}_{i+j} & &
}
\end{equation}
have equal kernel.  For, the left hand map in (\ref{eqn.main}) equals $\operatorname{Ad }(\mu_{i+1,j-1} \circ \mu_{i,1} \otimes {\mathcal{E}}^{U \otimes j-1})$.  But since $\mathcal{B}$-module multiplication is associative,
$$
\operatorname{Ad }(\mu_{i+1,j-1} \circ \mu_{i,1} \otimes {\mathcal{E}}^{U \otimes j-1})=\operatorname{Ad }\mu_{i,j}.
$$
Now, by Corollary \ref{cor.bimodtensor} the kernel of the right hand map is
$$
\operatorname{ker }(\operatorname{Ad }\mu_{i,1}) \otimes {\mathcal{E}}^{U \otimes j-1}+ {\mathcal{E}}^{U} \otimes \operatorname{ker }\operatorname{Ad }\mu_{i+1,j-1}.
$$
We show that $\operatorname{ker }(\operatorname{Ad }\mu_{i,1}) \otimes {\mathcal{E}}^{U \otimes j-1} \subset \operatorname{ker }\operatorname{Ad }\mu_{i,j}$, the kernel of the left map.  To prove this fact, we claim it suffices to show
\begin{equation} \label{eqn.conditionly}
\mu_{i,1}(F_{i}\operatorname{ker }\operatorname{Ad }\mu_{i,1})=0.
\end{equation}
For, suppose (\ref{eqn.conditionly}) holds.  By naturality of $\eta_{i}$, the diagram
\begin{equation} \label{eqn.provey}
\xymatrix{
\operatorname{ker }\operatorname{Ad }\mu_{i,1} \otimes {\mathcal{E}}^{U \otimes j-1} \ar[d] \ar[r]^{\hskip -.2in \eta_{i}} & G_{i}F_{i}(\operatorname{ker }\operatorname{Ad }\mu_{i,1} \otimes {\mathcal{E}}^{U \otimes j-1}) \ar[d] \\
{\mathcal{E}}^{U} \otimes {\mathcal{E}}^{U \otimes j-1} \ar[r]_{\hskip -.2in \eta_{i}} & G_{i}F_{i}({\mathcal{E}}^{U} \otimes {\mathcal{E}}^{U \otimes j-1}) \ar[d]^{G_{i}(\mu_{i,1}\otimes {\mathcal{E}}^{U \otimes j-1})} \\
& G_{i}F_{i+1}{\mathcal{E}}^{U \otimes j-1}
}
\end{equation}
commutes.  Since $\mu_{i,1}(F_{i}\operatorname{ker }\operatorname{Ad }\mu_{i,1})=0$, the right vertical equals $0$.  Thus, the right route equals $0$ so that, by commutativity of (\ref{eqn.provey}), the left route equals $0$.  But the bottom row composed with the bottom right vertical is a composite of the left hand side of (\ref{eqn.main}) so that
$$
\operatorname{ker }\operatorname{Ad }\mu_{i,1} \otimes {\mathcal{E}}^{U \otimes j-1} \subset \operatorname{ker }\operatorname{Ad }\mu_{i,j}
$$
as desired.

We now show $\mu_{i,1}(F_{i}\operatorname{ker }\operatorname{Ad }\mu_{i,1})=0$.  By adjointness of $(F_{i},G_{i})$, the diagram
$$
\xymatrix{
\operatorname{Hom}_{U \times X}(F_{i}{\mathcal{E}}^{U},{\mathcal{M}}_{i+1}) \ar[r] \ar[d] & \operatorname{Hom}_{(U \times X)_{U}^{2}}({\mathcal{E}}^{U},G_{i}{\mathcal{M}}_{i+1}) \ar[d] \\
\operatorname{Hom}_{U \times X}(F_{i} \operatorname{ker }\operatorname{Ad }\mu_{i,1}, {\mathcal{M}}_{i+1}) \ar[r] & \operatorname{Hom}_{(U \times X)_{U}^{2}}(\operatorname{ker }\operatorname{Ad }\mu_{i,1}, G_{i}{\mathcal{M}}_{i+1})
}
$$
commutes.  Now, $\mu_{i,1}$ in the upper left goes to $\operatorname{Ad }\mu_{i,1}$ in the upper right.  This clearly goes to zero on the lower right.  Thus, $\mu_{i,1}$ goes to zero on the lower left, as desired.  This completes the proof that
$$
\operatorname{ker }(\operatorname{Ad }\mu_{i,1}) \otimes {\mathcal{E}}^{U \otimes j-1} \subset \operatorname{ker }\operatorname{Ad }\mu_{i,j}.
$$
In a similar fashion, one shows that
$$
{\mathcal{E}}^{U} \otimes \operatorname{ker }\operatorname{Ad }\mu_{i+1,j-1} \subset \operatorname{ker }\operatorname{Ad }\mu_{i,j}.
$$
Thus,
$$
{\mathcal{E}}^{U} \otimes \operatorname{ker }\operatorname{Ad }\mu_{i+1,j-1} + \operatorname{ker }\operatorname{Ad }\mu_{i,1} \otimes {\mathcal{E}}^{U \otimes j-1} \subset \operatorname{ker }\operatorname{Ad }\mu_{i,j}.
$$
To complete the proof of the Proposition, we must show that the inclusion is equality.  By Proposition \ref{prop.bij}, both routes in (\ref{eqn.main}) are epimorphisms.  Thus,
\begin{equation} \label{eqn.super1}
{\mathcal{E}}^{U \otimes j}/({\mathcal{E}}^{U} \otimes \operatorname{ker }\operatorname{Ad }\mu_{i+1,j-1}+\operatorname{ker }\operatorname{Ad }\mu_{i,1} \otimes {\mathcal{E}}^{U \otimes j-1}) \cong G_{i}{\mathcal{M}}_{i} \otimes G_{i+1}{\mathcal{M}}_{i+j}
\end{equation}
and
\begin{equation} \label{eqn.super2}
{\mathcal{E}}^{U \otimes j}/\operatorname{ker }\operatorname{Ad }\mu_{i,j} \cong G_{i}{\mathcal{M}}_{i+j}.
\end{equation}
Since $\mathcal{M}$ is free then, by Proposition \ref{prop.simpleiso}, $G_{i}{\mathcal{M}}_{i+j} \cong q_{i,i+j*}{\mathcal{O}_{U}}$, $G_{i}{\mathcal{M}}_{i+1} \cong q_{i,i+1*}{\mathcal{O}}_{U}$ and $G_{i+1}{\mathcal{M}}_{i+j} \cong q_{i+1,i+j*}{\mathcal{O}}_{U}$.  Since $q_{k,l}$ is a closed immersion for all $k$ and $l$,
$$
q_{i,i+1*}{\mathcal{O}}_{U} \otimes q_{i+1,i+j*}{\mathcal{O}}_{U} \cong q_{i,i+j*}{\mathcal{O}_{U}}
$$
by Proposition \ref{prop.cruciso}, so that
$$
G_{i}{\mathcal{M}}_{i+1} \otimes G_{i+1}{\mathcal{M}}_{i+j} \cong G_{i}{\mathcal{M}}_{i+j}.
$$
This, along with (\ref{eqn.super1}) and (\ref{eqn.super2}), implies there is an isomorphism
$$
{\mathcal{E}}^{U \otimes j}/({\mathcal{E}}^{U} \otimes \operatorname{ker }\operatorname{Ad }\mu_{i+1,j-1}+\operatorname{ker }\operatorname{Ad }\mu_{i,1} \otimes {\mathcal{E}}^{U \otimes j-1}) \cong {\mathcal{E}}^{U \otimes j}/\operatorname{ker }\operatorname{Ad }\mu_{i,j}.
$$
Thus, the projection
\begin{equation} \label{eqn.epinow}
{\mathcal{E}}^{U \otimes j}/({\mathcal{E}}^{U} \otimes \operatorname{ker }\operatorname{Ad }\mu_{i+1,j-1}+\operatorname{ker }\operatorname{Ad }\mu_{i,1} \otimes {\mathcal{E}}^{U \otimes j-1}) \rightarrow {\mathcal{E}}^{U \otimes j}/\operatorname{ker }\operatorname{Ad }\mu_{i,j}
\end{equation}
with kernel
\begin{equation} \label{eqn.kernelya}
\operatorname{ker }\operatorname{Ad }\mu_{i,j}/({\mathcal{E}}^{U} \otimes \operatorname{ker }\operatorname{Ad }\mu_{i+1,j-1}+\operatorname{ker }\operatorname{Ad }\mu_{i,1} \otimes {\mathcal{E}}^{U \otimes j-1})
\end{equation}
induces an epimorphism from ${\mathcal{E}}^{U \otimes j}/({\mathcal{E}}^{U} \otimes \operatorname{ker }\operatorname{Ad }\mu_{i+1,j-1}+\operatorname{ker }\operatorname{Ad }\mu_{i,1} \otimes {\mathcal{E}}^{U \otimes j-1})$ to itself.  Since $q_{i,i+j*}^{U} {\mathcal{O}}_{U}$ is coherent (\cite[II, ex. 5.5, p.124]{alggeo}) and $U$ and $X$ are noetherian, ${\mathcal{E}}^{U \otimes j}/({\mathcal{E}}^{U} \otimes \operatorname{ker }\operatorname{Ad }\mu_{i+1,j-1}+\operatorname{ker }\operatorname{Ad }\mu_{i,1} \otimes {\mathcal{E}}^{U \otimes j-1})$ is locally finitely generated.  Thus, (\ref{eqn.epinow}) is an isomorphism (\cite[Corollary 4.4a, p.120]{comalg}) so that (\ref{eqn.kernelya}) equals zero.  The proof follows.
\end{proof}
Repeated application of this result yields
\begin{corollary} \label{cor.segrekernel}
With the notation as in Lemma \ref{prop.main},
$$
\operatorname{ker }\operatorname{Ad }\mu_{i,j} = \Sigma_{l=0}^{j-1}{\mathcal{E}}^{U \otimes l} \otimes \operatorname{ker }\operatorname{Ad }\mu_{i+l,1} \otimes {\mathcal{E}}^{U \otimes j-l-1}.
$$
\end{corollary}

\begin{proposition} \label{prop.injy}
Let $\Phi$ be the natural transformation constructed in Theorem \ref{theorem.rep1}.  There exists a natural transformation
\begin{equation} \label{eqn.sigma}
\Sigma: \Gamma_{n} \Longrightarrow \operatorname{Hom}_{S}(-,{\mathbb{P}}_{X^{2}}({\mathcal{E}}^{\otimes n}/{\mathcal{I}}_{n}))
\end{equation}
making
$$
\xymatrix{
\Gamma_{n} \ar@{=>}[d]_{\Sigma} \ar@{=>}[r]^{\Phi} & \operatorname{Hom}_{S}(-,{\mathbb{P}}_{X^{2}}(\mathcal{E})^{\otimes n}) \ar@{=>}[d]^{s \circ -} \\
\operatorname{Hom}_{S}(-,{\mathbb{P}}_{X^{2}}({\mathcal{E}}^{\otimes n}/{\mathcal{I}}_{n})) \ar@{=>}[r] & \operatorname{Hom}_{S}(-,{\mathbb{P}}_{X^{2}}({\mathcal{E}}^{\otimes n}))
}
$$
commute.
\end{proposition}

\begin{proof}
By Theorem \ref{theorem.rep1}, there is a monomorphic natural transformation
$$
\Phi:\Gamma_{n} \Longrightarrow \operatorname{Hom}_{S}(-,{\mathbb{P}}_{X^{2}}(\mathcal{E})^{\otimes n}).
$$
Let $U$ be an affine scheme.  We construct an injection $\operatorname{Hom}_{S}(U,{\mathbb{P}}_{X^{2}}(\mathcal{E})^{\otimes n}) \rightarrow \operatorname{Hom}_{S}(U,{\mathbb{P}}_{(U \times X)_{U}^{2}}({\mathcal{E}}^{U})^{\otimes n})$.  We first note that
\begin{equation} \label{eqn.referencer}
{\mathbb{P}}_{(U \times X)_{U}^{2}}({\mathcal{E}}^{U})^{\otimes n} \cong {\mathbb{P}}_{X^{2}}(\mathcal{E})^{\otimes n} \times_{S} U.
\end{equation}
Since we may define an injection
\begin{equation} \label{eqn.newmap}
\operatorname{Hom}_{S}(U,{\mathbb{P}}_{X^{2}}(\mathcal{E})^{\otimes n}) \rightarrow \operatorname{Hom}_{S}(U,{\mathbb{P}}_{X^{2}}(\mathcal{E})^{\otimes n} \times_{S} U)
\end{equation}
by sending $f$ to $f \times \id_{U}$, (\ref{eqn.referencer}) allows us to define an injection
\begin{equation} \label{eqn.jj}
j: \operatorname{Hom}_{S}(U,{\mathbb{P}}_{X^{2}}(\mathcal{E})^{\otimes n}) \rightarrow \operatorname{Hom}_{S}(U,{\mathbb{P}}_{(U \times X)_{U}^{2}}({\mathcal{E}}^{U})^{\otimes n})
\end{equation}
as desired.  Furthermore, if $\mathcal{M}$ is a free truncated $U$-family of length $n+1$, then $j \circ \Phi_{U}([\mathcal{M}])$ corresponds, by Proposition \ref{prop.trans}, to a collection of $n$ epimorphisms $\psi_{i}$, which are the compositions
$$
\xymatrix{
{\mathcal{E}}^{U} \ar[r]^{\hskip -.16in \operatorname{Ad }\mu_{i,1}} & G_{i}{\mathcal{M}}_{i+1} \ar[r]^{\cong} &  q_{i,i+1*}{\mathcal{O}}_{U}
}
$$
whose right map is given by Proposition \ref{prop.simpleiso}.  By Proposition \ref{prop.segretensor}, $s \circ j \circ \Phi_{U}([\mathcal{M}])$ corresponds to the epimorphism $\otimes_{l=0}^{n-1}\psi_{l}$.  By Corollary \ref{cor.bimodtensor},
$$
\operatorname{ker }\otimes_{l=0}^{n-1}\psi_{l}=\Sigma_{l=0}^{n-1}{\mathcal{E}}^{U \otimes l} \otimes \operatorname{ker }\psi_{l} \otimes {\mathcal{E}}^{U \otimes n-l-1}=
$$
$$
\Sigma_{l=0}^{n-1}{\mathcal{E}}^{U \otimes l} \otimes \operatorname{ker }\operatorname{Ad }\mu_{l,1} \otimes {\mathcal{E}}^{U \otimes n-l-1}.
$$
By Corollary \ref{cor.segrekernel}, this module equals $\operatorname{ker }\operatorname{Ad }\mu_{0,n}$.  We claim ${\mathcal{I}}_{n}^{U} \subset \operatorname{ker }\operatorname{Ad }\mu_{0,n}$.  By adjointness of $(F_{i},G_{i})$, the diagram
$$
\xymatrix{
\operatorname{Hom}_{U \times X}(F_{i}{\mathcal{E}}^{U \otimes n-i}, {\mathcal{M}}_{n}) \ar[r] \ar[d] & \operatorname{Hom}_{(U \times X)_{U}^{2}}({\mathcal{E}}^{U \otimes n-i}, G_{i}{\mathcal{M}}_{n}) \ar[d] \\
\operatorname{Hom}_{U \times X}(F_{i}{\mathcal{I}}^{U}_{n-i}, {\mathcal{M}}_{n}) \ar[r]  & \operatorname{Hom}_{(U \times X)_{U}^{2}}({\mathcal{I}}^{U}_{n-i}, G_{i}{\mathcal{M}}_{n})
}
$$
commutes.  Now, $\mu_{i,n-i}$ in the upper left goes to zero in the lower left since $\mathcal{M}$ is a $T({\mathcal{E}}^{U})/{\mathcal{I}}^{U}$-module.  Thus, $\operatorname{Ad }\mu_{i,n-i}$ goes to zero in the lower right, as claimed.

We conclude ${\mathcal{I}}_{n}^{U} \subset \operatorname{ker } \otimes_{i=0}^{n-1}\psi_{i}.$  Thus, $\otimes_{i=0}^{n-1}\psi_{i}$ must factor as
$$
\xymatrix{
{\mathcal{E}}^{U \otimes n} \ar[rr]^{\otimes_{i=0}^{n-1}\psi_{i}} \ar[d] & & \otimes_{i=0}^{n-1}q_{i,i+1*}{\mathcal{O}}_{U} \\
{\mathcal{E}}^{U \otimes n}/{\mathcal{I}}_{n}^{U} \ar[rr]_{=} & & {\mathcal{E}}^{U \otimes n}/{\mathcal{I}}_{n}^{U} \ar[u]
}
$$
so there exists a map $\sigma_{U}([\mathcal{M}]):U \rightarrow  {\mathbb{P}}_{(U \times X)_{U}^{2}}({\mathcal{E}}^{U \otimes n}/{\mathcal{I}}_{n}^{U})$ such that $s \circ j \circ \Phi_{U}([\mathcal{M}]) \in \operatorname{Hom}_{S}(U,{\mathbb{P}}_{X^{2}}({\mathcal{E}}^{\otimes n}))$ factors as
$$
\xymatrix{
U \ar[d]_{\sigma_{U}([\mathcal{M}])} \ar[rr]^{s \circ j \circ \Phi_{U}([\mathcal{M}])} & & {\mathbb{P}}_{(U \times X)_{U}^{2}}({\mathcal{E}}^{U \otimes n})  \\
{\mathbb{P}}_{(U \times X)_{U}^{2}}({\mathcal{E}}^{U \otimes n}/{\mathcal{I}}_{n}^{U}) \ar[rr]_{=} & & {\mathbb{P}}_{(U \times X)_{U}^{2}}({\mathcal{E}}^{U \otimes n}/{\mathcal{I}}_{n}^{U}) \ar[u].
}
$$
Thus, there exists a map $\sigma_{U}: \Gamma_{n}^{\Fr}(U) \rightarrow \operatorname{Hom}_{S}^{\Fr}(U,{\mathbb{P}}_{(U \times X)_{U}^{2}}({\mathcal{E}}^{U\otimes n}/{\mathcal{I}}_{n}^{U}))$ such that the diagram
\begin{equation} \label{eqn.almost}
\xymatrix{
\Gamma_{n}^{\Fr}(U) \ar[d]_{\sigma_{U}} \ar[r]^{j \circ \Phi_{U}(-)} & \operatorname{Hom}_{S}^{\Fr}(U,{\mathbb{P}}_{(U \times X)_{U}^{2}}({\mathcal{E}}^{U})^{\otimes n}) \ar[d]^{s \circ -} \\
\operatorname{Hom}_{S}^{\Fr}(U,{\mathbb{P}}_{(U \times X)_{U}^{2}}({\mathcal{E}}^{U\otimes n}/{\mathcal{I}}_{n}^{U})) \ar[r] & \operatorname{Hom}_{S}^{\Fr}(U,{\mathbb{P}}_{(U \times X)_{U}^{2}}({\mathcal{E}}^{U\otimes n}))
}
\end{equation}
commutes as well.  Let $\Sigma_{U}^{\Fr}$ be the composition
$$
\xymatrix{
\Gamma_{n}^{\Fr}(U) \ar[r]^{\hskip -.7in \sigma_{U}} & \operatorname{Hom}_{S}^{\Fr}(U,{\mathbb{P}}_{(U \times X)_{U}^{2}}({\mathcal{E}}^{U\otimes n}/{\mathcal{I}}_{n}^{U})) \ar[r] & \operatorname{Hom }_{S}^{\Fr}(U,{\mathbb{P}}_{X^{2}}({\mathcal{E}}^{\otimes n}/{\mathcal{I}}_{n}))
}
$$
whose right composite is induced by the composition
\begin{equation} \label{eqn.ipromise}
{\mathbb{P}}_{(U \times X)_{U}^{2}}({\mathcal{E}}^{U\otimes n}/{\mathcal{I}}_{n}^{U}) \rightarrow {\mathbb{P}}_{(U \times X)_{U}^{2}}(({\mathcal{E}}^{\otimes n}/{\mathcal{I}}_{n})^{U}) \rightarrow {\mathbb{P}}_{X^{2}}({\mathcal{E}}^{\otimes n}/{\mathcal{I}}_{n}).
\end{equation}
If we let $(-,-)=\operatorname{Hom }_{S}^{\Fr}(-,-)$, the half-cube
\begin{equation} \label{eqn.helpthereader}
\xymatrix{
& \Gamma_{n}^{\Fr}(U) \ar[dl]_{\sigma_{U}} \ar[rr]^{j \circ \Phi_{U}(-)}  & & (U,{\mathbb{P}}_{(U \times X)_{U}^{2}}({\mathcal{E}}^{U})^{\otimes n}) \ar[dl]^{s \circ -} \ar[dd] \\
(U,{\mathbb{P}}_{(U \times X)_{U}^{2}}({\mathcal{E}}^{U\otimes n}/{\mathcal{I}}_{n}^{U}))    \ar[dd]  \ar[rr] & & (U,{\mathbb{P}}_{(U \times X)_{U}^{2}}({\mathcal{E}}^{U\otimes n})) \ar[dd] & \\
&  & &  (U,{\mathbb{P}}_{X^{2}}(\mathcal{E})^{\otimes n}) \ar[dl]^{s \circ -} \\
(U,{\mathbb{P}}_{X^{2}}({\mathcal{E}}^{\otimes n}/{\mathcal{I}}_{n})) \ar[rr]  & & (U,{\mathbb{P}}_{X^{2}}({\mathcal{E}}^{\otimes n}))  &
}
\end{equation}
commutes since the diagram
\begin{equation} \label{eqn.penulitimate}
\xymatrix{
{\mathbb{P}}_{(U \times X)_{U}^{2}}({\mathcal{E}}^{U \otimes n}/{\mathcal{I}}_{n}^{U}) \ar[r] \ar[d] & {\mathbb{P}}_{(U \times X)_{U}^{2}}({\mathcal{E}}^{U \otimes n})  \ar[d] & {\mathbb{P}}_{(U \times X)_{U}^{2}}({\mathcal{E}}^{U})^{ \otimes n}  \ar[l]_{s} \ar[dd] \\
{\mathbb{P}}_{(U \times X)_{U}^{2}}(({\mathcal{E}}^{\otimes n}/{\mathcal{I}}_{n})^{U}) \ar[r] \ar[d] & {\mathbb{P}}_{(U \times X)_{U}^{2}}(({\mathcal{E}}^{\otimes n})^{U}) \ar[d] & \\
{\mathbb{P}}_{X^{2}}({\mathcal{E}}^{\otimes n}/{\mathcal{I}}_{n}) \ar[r] & {\mathbb{P}}_{X^{2}}({\mathcal{E}}^{\otimes n}) &  {\mathbb{P}}_{X^{2}}({\mathcal{E}})^{\otimes n} \ar[l]^{s}
}
\end{equation}
whose left vertical equals (\ref{eqn.ipromise}), and whose bottom verticals and right-most vertical are maps from Lemma \ref{lem.proj}, commutes.  The right square commutes since $s$ is compatible with base change, the upper left square commutes by the naturality of the indexing of the tensor product, and the lower left square commutes by Lemma \ref{lem.proj}.  Thus, we have a commutative diagram
\begin{equation} \label{eqn.almost2}
\xymatrix{
\Gamma_{n}^{\Fr}(U) \ar[d]_{\Sigma_{U}^{\Fr}} \ar[r] & \operatorname{Hom}_{S}^{\Fr}(U,{\mathbb{P}}_{X^{2}}(\mathcal{E})^{\otimes n}) \ar[d]^{s \circ -} \\
\operatorname{Hom}_{S}^{\Fr}(U,{\mathbb{P}}_{X^{2}}({\mathcal{E}}^{\otimes n}/{\mathcal{I}}_{n})) \ar[r] & \operatorname{Hom}_{S}^{\Fr}(U,{\mathbb{P}}_{X^{2}}({\mathcal{E}}^{\otimes n}))
}
\end{equation}
for each $U \in \sf{S}$.

To complete the proof, we need only show that $\Sigma^{\Fr}$ is natural.  For, in that case, $\Sigma^{\Fr}$ extends to a natural transformation
$$
\Sigma:\Gamma_{n} \Longrightarrow \operatorname{Hom}_{S}(-,{\mathbb{P}}_{X^{2}}({\mathcal{E}}^{\otimes n}/{\mathcal{I}}_{n}))
$$
by Proposition \ref{prop.subfunctor}.  To show that $\Sigma^{\Fr}$ is natural, we must show that if $f:V \rightarrow U$ is a morphism of affine schemes, then the diagram
$$
\xymatrix{
\Gamma_{n}^{\Fr}(V) \ar[rr]^{\hskip -.5in \Sigma_{V}^{\Fr}} \ar[d]_{\Gamma_{n}^{\Fr}(f)} & & \operatorname{Hom}_{S}^{\Fr}(V,{\mathbb{P}}_{X^{2}}({\mathcal{E}}^{\otimes n}/{\mathcal{I}}_{n})) \ar[d]^{-\circ f} \\
\Gamma_{n}^{\Fr}(U) \ar[rr]_{\hskip -.5in \Sigma_{U}^{\Fr}} & & \operatorname{Hom}_{S}^{\Fr}(U,{\mathbb{P}}_{X^{2}}({\mathcal{E}}^{\otimes n}/{\mathcal{I}}_{n}))
}
$$
commutes.  If we let $(-,-)=\operatorname{Hom }_{S}^{\Fr}(-,-)$, this follows easily from the fact that all faces of the cube
$$
\xymatrix{
& \Gamma_{n}^{\Fr}(V) \ar[dl] \ar[rr] \ar[dd] & & (V,{\mathbb{P}}_{X^{2}}(\mathcal{E})^{\otimes n}) \ar[dl] \ar[dd] \\
(V,{\mathbb{P}}_{X^{2}}({\mathcal{E}}^{\otimes n}/{\mathcal{I}}_{n}))    \ar[dd]  \ar[rr] & & (V,{\mathbb{P}}_{X^{2}}({\mathcal{E}}^{\otimes n})) \ar[dd] & \\
& \Gamma_{n}^{\Fr}(U) \ar[rr] \ar[dl] & &  (U,{\mathbb{P}}_{X^{2}}(\mathcal{E})^{\otimes n}) \ar[dl] \\
(U,{\mathbb{P}}_{X^{2}}({\mathcal{E}}^{\otimes n}/{\mathcal{I}}_{n})) \ar[rr]  & & (U,{\mathbb{P}}_{X^{2}}({\mathcal{E}}^{\otimes n}))  &
}
$$
save possibly the left face, commute.  For, the top and bottom commute (\ref{eqn.almost2}), while the far square commutes by naturality of $\Phi^{\Fr}$.  Since the horizontal arrows are injections, it must be true that the left face commutes also.
\end{proof}

\begin{corollary} \label{cor.newinjy}
Let $P$ be the pullback in the diagram
$$
\xymatrix{
& {\mathbb{P}}_{X^{2}}(\mathcal{E})^{\otimes n} \ar[d]^{s} \\
{\mathbb{P}}_{X^{2}}({\mathcal{E}}^{\otimes n}/{\mathcal{I}}_{n}) \ar[r] & {\mathbb{P}}_{X^{2}}({\mathcal{E}}^{\otimes n})
}
$$
whose bottom vertical is induced by the epimorphism ${\mathcal{E}}^{\otimes n} \rightarrow {\mathcal{E}}^{\otimes n}/{\mathcal{I}}_{n}$.  Then the natural transformation $\Sigma$ ((\ref{eqn.sigma}) on page \pageref{eqn.sigma}) induces a monomorphism
\begin{equation} \label{eqn.upsilon}
\Upsilon:\Gamma_{n} \Longrightarrow \operatorname{Hom}_{S}(-,P).
\end{equation}
\end{corollary}

\begin{proof}
By Proposition \ref{prop.injy}, $\Sigma$ makes
$$
\xymatrix{
\Gamma_{n} \ar@{=>}[d]_{\Sigma} \ar@{=>}[r] & \operatorname{Hom}_{S}(-,{\mathbb{P}}_{X^{2}}(\mathcal{E})^{\otimes n}) \ar@{=>}[d]^{s \circ -} \\
\operatorname{Hom}_{S}(-,{\mathbb{P}}_{X^{2}}({\mathcal{E}}^{\otimes n}/{\mathcal{I}}_{n})) \ar@{=>}[r] & \operatorname{Hom}_{S}(-,{\mathbb{P}}_{X^{2}}({\mathcal{E}}^{\otimes n}))
}
$$
commute.  Since the right route is a monomorphism, the induced map to the pullback must also be a monomorphism
$$
\Gamma_{n} \Longrightarrow \operatorname{Hom}_{S}(-,{\mathbb{P}}_{X^{2}}(\mathcal{E})^{\otimes n}) \times_{ \operatorname{Hom}_{S}(-,{\mathbb{P}}_{X^{2}}({\mathcal{E}}^{\otimes n}))} \operatorname{Hom}_{S}(-,{\mathbb{P}}_{X^{2}}({\mathcal{E}}^{\otimes n}/{\mathcal{I}}_{n})).
$$
By \cite[p.20]{curves}, this is the same as a monomorphism
$$
\Upsilon:\Gamma_{n} \Longrightarrow \operatorname{Hom}_{S}(-,P)
$$
as desired.
\end{proof}

\begin{lemma} \label{lem.factorproj}
Suppose $Y$ is a scheme, $\mathcal{F}$ is an ${\mathcal{O}}_{Y}$-module, $\mathcal{I} \subset \mathcal{F}$ is a submodule, $p:{\mathbb{P}}_{Y}(\mathcal{F}) \rightarrow Y$ is the structure map,  and there is a commutative diagram
\begin{equation} \label{eqn.routenow}
\xymatrix{
U \ar[r]^{g} \ar[d]_{g'} & {\mathbb{P}}_{Y}(\mathcal{F})  \\
{\mathbb{P}}_{Y}(\mathcal{F}/\mathcal{I}) \ar[ur] &
}
\end{equation}
where $g$ and $g'$ are $Y$-morphisms.  If $g$ corresponds, by Proposition \ref{prop.groth}, to a morphism
$$
\psi:\mathcal{F} \rightarrow q_{*}\mathcal{L},
$$
then $\psi$ factors as
$$
\xymatrix{
\mathcal{F} \ar[d] \ar[r]^{\psi} & q_{*}\mathcal{L}  \\
\mathcal{F}/\mathcal{I} \ar[ur] &
}
$$
\end{lemma}

\begin{proof}
If $g'$ corresponds to $\psi':\mathcal{F}/\mathcal{I} \rightarrow q_{*}\mathcal{L}'$, then
$$
\xymatrix{
\mathcal{F} \ar[r] & \mathcal{F}/\mathcal{I} \ar[r]^{\psi'} & q_{*}\mathcal{L}'
}
$$
corresponds to the bottom route of (\ref{eqn.routenow}).  Since (\ref{eqn.routenow}) commutes, the top map must also correspond to the morphism.  Thus, by Proposition \ref{prop.groth}, there exists an isomorphism $\alpha:q_{*}\mathcal{L}' \rightarrow q_{*}\mathcal{L}$ such that
$$
\xymatrix{
\mathcal{F} \ar[d] \ar[r]^{\psi} & q_{*}\mathcal{L}  \\
\mathcal{F}/\mathcal{I} \ar[r]_{\psi'} & q_{*}\mathcal{L}' \ar[u]_{\alpha}
}
$$
commutes, as desired.
\end{proof}

\begin{proposition} \label{prop.lastfinal}
The natural transformation $\Upsilon:\Gamma_{n} \Longrightarrow \operatorname{Hom}_{S}(-,P)$ (\ref{eqn.upsilon}) is an equivalence.
\end{proposition}

\begin{proof}
We show $\Upsilon^{\Fr}:\Gamma_{n}^{\Fr} \Longrightarrow \operatorname{Hom}_{S}^{\Fr}(-,P)$, the restriction of $\Upsilon$ to free families, is an equivalence.  This suffices to prove the proposition in light of Proposition \ref{prop.subfunctor}.  Let $U$ be an affine scheme and suppose $f \in \operatorname{Hom}_{S}^{\Fr}(U,P)$.  Then $f$ is defined by a pair of maps $(f_{1},f_{2})$ making the diagram
\begin{equation} \label{eqn.lastya}
\xymatrix{
U \ar[r]^{f_{1}} \ar[d]_{f_{2}} & {\mathbb{P}}_{X^{2}}(\mathcal{E})^{\otimes n} \ar[d]^{s} \\
{\mathbb{P}}_{X^{2}}({\mathcal{E}}^{\otimes n}/{\mathcal{I}}_{n}) \ar[r] & {\mathbb{P}}_{X^{2}}({\mathcal{E}}^{\otimes n})
}
\end{equation}
commute.  If $\Delta_{n}$ is the functor of flat families of $T(\mathcal{E})$-modules of length $n+1$, then there exists a morphism $\Gamma_{n} \Longrightarrow \Delta_{n}$ induced by the epimorphism $T(\mathcal{E}) \rightarrow T(\mathcal{E})/\mathcal{I}$.  By Theorem \ref{theorem.rep1}, $f_{1}$ corresponds to an element of $\Delta_{n}^{\Fr}(U)$, $[\mathcal{M}]$, where $\mathcal{M}$ is a free $T({\mathcal{E}}^{U})$-module of length $n+1$.  By Lemma \ref{lem.inheritmod}, to show that $\mathcal{M}$ is actually a $T({\mathcal{E}}^{U})/{\mathcal{I}}^{U}$-module, we need to show, for all $i,j \geq 0$, that the composition
$$
\xymatrix{
{\mathcal{M}}_{i} \otimes {\mathcal{I}}_{j}^{U} \ar[r] & {\mathcal{M}}_{i} \otimes {\mathcal{E}}^{U\otimes j} \ar[r]^{\hskip .2in \mu_{i,j}} & {\mathcal{M}}_{i+j}
}
$$
equals zero.  It thus suffices to show that, for all $i, j \geq 0$,
$$
{\mathcal{I}}_{j}^{U} \subset \operatorname{ker }\operatorname{Ad }\mu_{i,j}.
$$
Since (\ref{eqn.lastya}) commutes, we claim the diagram
$$
\xymatrix{
U \ar[r]^{f_{1} \times \id_{U}} \ar[d]_{f_{2} \times \id_{U}} & {\mathbb{P}}_{X^{2}}(\mathcal{E})^{\otimes n} \times_{S}U \ar[d]^{s \times \id_{U}} \ar[dr]^{\cong} & \\
{\mathbb{P}}_{X^{2}}({\mathcal{E}}^{\otimes n}/{\mathcal{I}}_{n})\times_{S} U \ar[r] \ar[dr]_{\cong} & {\mathbb{P}}_{X^{2}}({\mathcal{E}}^{\otimes n}) \times_{S}U \ar[dr]^{\cong} &  {\mathbb{P}}_{(U \times X)_{U}^{2}}({\mathcal{E}}^{U})^{\otimes n} \ar[d]^{s} \\
& {\mathbb{P}}_{(U \times X)_{U}^{2}}({\mathcal{E}}^{U\otimes n}/{\mathcal{I}}_{n}^{U}) \ar[r] & {\mathbb{P}}_{(U \times X)_{U}^{2}}({\mathcal{E}}^{U\otimes n})
}
$$
whose diagonals are from Lemma \ref{lem.proj}, commutes.  The upper left commutes since (\ref{eqn.lastya}) commutes.  The lower left square commutes by Lemma \ref{lem.proj}, while the lower right square commutes since $s$ is a $U$-morphism which is compatible with base change.  Denote the top and bottom route of the diagram by $f_{1}^{U}$ and $f_{2}^{U}$ respectively.  By Proposition \ref{prop.segretensor}, $f_{1}^{U}$ corresponds to the epimorphism
$$
\xymatrix{
\upsilon: {\mathcal{E}}^{U \otimes n} \ar[rr]^{\bigotimes_{i=1}^{n}\psi_{i}} & & \bigotimes_{i=1}^{n}q_{i.i+1*}{\mathcal{O}}_{U}
}
$$
where $\psi_{i}$ is the composition
$$
\xymatrix{
{\mathcal{E}}^{U} \ar[rr]^{\operatorname{Ad }\mu_{i,1}} & & G_{i}{\mathcal{M}}_{i+1} \ar[r]^{\cong} &  q_{i,i+1*}{\mathcal{O}}_{U}
}
$$
whose right composite is the map in Proposition \ref{prop.simpleiso}.  Since $f_{1}^{U}$ factors through ${\mathbb{P}}_{(U \times X)_{U}^{2}}({\mathcal{E}}^{U\otimes n}/{\mathcal{I}}_{n}^{U})$, $\upsilon$ factors through ${\mathcal{E}}^{U\otimes n}/{\mathcal{I}}_{n}^{U}$ by Lemma \ref{lem.factorproj}.  Thus, ${\mathcal{I}}_{n}^{U} \subset \operatorname{ker }\bigotimes_{i=1}^{n}\psi_{i}$.  By Corollary \ref{cor.segrekernel},
$$
\operatorname{ker }\otimes_{i=1}^{n}\psi_{i} = \operatorname{ker }\operatorname{Ad }\mu_{i,j}
$$
as desired.
\end{proof}
Thus, we have established Theorem \ref{theorem.bigone}.\index{family!free, of truncated point modules|)}

\backmatter

\printindex

\end{document}